\documentclass[11pt]{amsart} 

\usepackage{pstricks}\psset{xunit=2.2ex,yunit=2.2ex}
\usepackage{latexsym}
\usepackage{amssymb}
\usepackage{youngtab}\Yboxdim{12pt}\Yinterspace{1pt}

\newcommand{\excise}[1]{}

\newtheorem{thm}{Theorem}
\newtheorem{lemma}[thm]{Lemma}

\newtheorem{cor}[thm]{Corollary}
\newtheorem{prop}[thm]{Proposition}
\newtheorem{conj}[thm]{Conjecture}

\theoremstyle{definition}
\newtheorem{example}[thm]{Example}
\newtheorem{remark}[thm]{Remark}

\newtheorem{defn}[thm]{Definition}

\newtheorem{conv}[thm]{Convention}

\newenvironment{example*}{\begin{trivlist}\item {\bf Example.\,}}
        {\end{trivlist}}

\newenvironment{movedThm}[1]{\begin{trivlist}\item {\bf
        {#1}\,}\it}{\end{trivlist}}

\newenvironment{noqed}{\begin{trivlist}\item {\it Proof.\,}}%
        {\end{trivlist}}

\newenvironment{numbered}%
        {\begin{list}
                {\noindent\makebox[0mm][r]{\arabic{enumi}.}}
                {\leftmargin=5.5ex \usecounter{enumi}}
        }
        {\end{list}}

        {\begin{list}
                {\noindent\makebox[0mm][r]{(\roman{enumi})}}
                {\leftmargin=5.5ex \usecounter{enumi}}
        }
        {\end{list}}

\newcounter{separated}

\newenvironment{stareqn}%
        {
         \setcounter{separated}{\value{equation}}
         \setcounter{equation}{0}
         \begin{equation}%
        }%
        {\end{equation}%
         \setcounter{equation}{\value{separated}}%
        }%

\newenvironment{rcgraph}{\begin{trivlist}\item\centering\footnotesize$}
                        {$\vspace{1ex}\end{trivlist}}

\def\tinyrc#1{\hbox{\tiny${#1}$}}

\newenvironment{vcenteredyoung}{\begin{array}{@{}c@{}}\\[-2.2ex]\young}
                                {\end{array}}
\def\cyoung#1{\begin{vcenteredyoung}(#1)\end{vcenteredyoung}}
\def\lyoung#1{\raisebox{-.5ex}{\young(#1)}}

\def\bem#1{\textbf{#1}}

\def\dis{\displaystyle}

\def\<{\langle}
\def\>{\rangle}
\def\0{{\mathbf 0}}
\def\1{{\mathbf 1}}
\def\AA{{\mathbb A}}

\def\CC{{\mathbb C}}
\def\cC{{\mathcal C}}

\def\da{{\dot a}}
\def\db{{\dot b}}
\def\dc{{\dot c}}
\def\dt{{\dot d}}

\def\EE{{\mathcal E}}
\def\FF{{\mathcal F}}
\def\GG{{\mathcal G}}

\def\II{{\mathcal I}}

\def\KK{{\mathcal K}}

\def\NN{{\mathbb N}}
\def\OO{{\mathcal O}}
\def\PP{{\mathcal P}}
\def\QQ{{\mathcal Q}}

\def\SS{{\mathfrak S}}
\def\VV{{\mathbf V}}

\def\ZZ{{\mathbb Z}}
\def\cZ{{\mathcal Z}}
\def\aa{{\mathbf a}}
\def\bb{{\mathbf b}}
\def\cc{{\mathbf c}}
\def\dd{{\mathbf d}}

\def\ff{{\mathbf f}}
\def\bg{{\mathbf g}}
\def\ii{{\mathbf i}}
\def\kk{{\Bbbk}}

\def\rr{{\mathbf r}}
\renewcommand\ss{{\mathbf{s}}}
\renewcommand\tt{{\mathbf{t}}}
\def\uu{{\mathbf u}}

\def\ww{{\mathbf w}}
\def\xx{{\mathbf x}}
\def\yy{{\mathbf y}}
\def\zz{{\mathbf z}}

\def\Gr{{\rm Gr}}

\def\st{{\rm st}}
\def\th{{\rm th}}
\def\Hom{{\rm Hom}}
\def\Top{{\rm top}}
\def\col{{\rm col}}
\def\row{{\rm row}}
\def\diag{{\rm diag}}
\def\rank{{\rm rank}}
\def\spec{{\rm Spec}}
\def\codim{{\rm codim}}
\def\cycle{{\rm cycle}}

\def\weight{{\rm wt}}

\def\K{$K$}
\def\BtB{{B_+ \times B_-}}

\def\ID{I_{\!\Delta}}
\def\IN{\mathsf{in}}
\def\SL{{\mathit{S\!L}}}
\def\UB{{U\hspace{-.4ex}B}}

\def\gl{{\mathit{G\!L}}}
\def\rp{\mathcal{RP}}
\def\sq{\square}
\def\ti{\times}
\def\tv{{T}}
\def\IDi{\til{I}_{\!\Delta_i}}

\def\dom{\backslash}
\def\gld{{G_{\!}L_d}}
\def\glv{{\mathit{G\!L}}}
\def\too{\longrightarrow}
\def\oxx{\text{$\stackrel{\raisebox{-.2ex}{$\scriptscriptstyle\circ$}}%
     {\hbox{$\xx$}}\hspace{-.7ex}^{}$}}
\def\oyy{\text{$\stackrel{\raisebox{-.2ex}{$\scriptscriptstyle\circ$}}%
     {\hbox{$\yy$}}\hspace{-.7ex}^{}$}}
\def\btbv{{B_+ \times_T B_-}}
\def\endv{{M_d}}
\def\from{\leftarrow}
\def\glvv{{\mathit{G\!L}^2}}
\def\homv{{\hspace{-.2ex}\text{\sl{Hom}}\hspace{.1ex}}}
\def\into{\hookrightarrow}
\def\spot{{\hbox{\raisebox{.33ex}{\large\bf .}}\hspace{-.10ex}}}

\def\pgld{P\dom\hspace{-.2ex}\gld}
\def\adots{{.\hspace{.2ex}\raisebox{.5ex}{.}\hspace{.2ex}\raisebox{1ex}{.}}}
\def\cross{\textrm{`+'} }

\def\hhomv{{\hspace{-.2ex}\text{\bf\textsl{Hom}}\hspace{.1ex}}}
\def\homvp{{\text{\hspace{-.2ex}\textsl{Hom}}'}}
\def\minus{\smallsetminus}
\def\OOmega{{\mathbf{\Omega}}}

\def\goesto{\rightsquigarrow}
\def\implies{\Rightarrow}
\def\nothing{\varnothing}

\def\flatddots{\raisebox{.6ex}{.}\hspace{.4ex}\raisebox{.3ex}{.}\hspace{.4ex}.}
\def\toplinedots{\hfill\raisebox{-5.1pt}[0pt][0pt]{\ $\ldots$\ }\hfill}

\def\ol#1{{\overline {#1}}}
\def\ub#1{\text{\underbar{$#1$\hspace{-.25ex}}\hspace{.25ex}}}

\def\wt#1{{\widetilde{#1}}}
\def\til#1{\rlap{\hspace{.4ex}\raisebox{-.25ex}[0pt][0pt]{%
        $\widetilde{\phantom{#1}\hspace{-.6ex}}$}}{#1}}
\def\ggg#1#2#3{{g_{#1}}^{\!\!#2#3}}
\def\fff#1#2#3{{f^{\makebox[0ex][l]{$\scriptstyle #1$}}%
         {}_{\!\scriptstyle #2#3}}}%
\def\sub#1#2{_{#1 \times #2}}
\def\xxx#1#2{{x^{#1}_{#2}}}
\def\ffff#1#2#3#4{{f_{#3#4}}^{\hspace{-1.85ex}#1#2}}
\def\lvec#1{\overleftarrow{#1}\!}
\def\rvec#1{\overrightarrow{#1}\!}
\def\lgam#1{\lvec{\gamma}_{\!{#1}}}
\def\rgam#1{\rvec{\gamma}_{\!{#1}}}
\def\lggg#1#2#3{{\lvec {g_{#1\!}}^{\,#2#3}}}
\def\rggg#1#2#3{{\rvec {g_{#1\!}}^{\,#2#3}}}
\def\twoline#1#2{\aoverb{\scriptstyle {#1}}{\scriptstyle {#2}}}

\newcommand{\aoverb}[2]{{\genfrac{}{}{0pt}{1}{#1}{#2}}}

\newcommand{\Xb}{\overline{X}}
\newcommand{\cd}{\cdot}

\newcommand{\Schub}{\mathfrak{S}}

\newcommand{\la}{\lambda}

\newcommand{\FFF}[3]{F_{#1}(\xx^{#2}\!-\!\xx^{#3})}
\newcommand{\SSS}[3]{\Schub_{#1}(\xx^{#2}\!-\!\xx^{#3})}

\newcommand{\Peel}{\mathrm{Peel}}
\newcommand{\vn}{\varnothing}
\newcommand{\kp}{\kappa}
\newcommand{\shape}{\mathrm{shape}}
\newcommand{\bla}{\ub\lambda}

\newcommand{\rc}{\mathrm{rc}}
\newcommand{\lc}{\mathrm{lc}}
\newcommand{\key}{\mathrm{key}}
\newcommand{\wgt}{\mathrm{wt}}
\newcommand{\north}{\!\uparrow}
\newcommand{\west}{\!\leftarrow}
\newcommand{\ldom}{\trianglelefteq}
\newcommand{\rdom}{\trianglerighteq}
\newcommand{\ptofs}{\Psi}
\newcommand{\hgt}{\mathrm{height}}
\newcommand{\wid}{\mathrm{width}}
\newcommand{\rrh}{{\widehat{\rr}}}
\newcommand{\Rh}{\widehat{\mathbf{R}}}

\newcommand{\Image}{\mathrm{Im}}
\newcommand{\code}{\mathrm{code}}

\newcommand{\leftkey}{K_-}
\newcommand{\rightkey}{K_+}
\newcommand{\Kh}{\widehat{K}}
\newcommand{\Xh}{\widehat{X}}
\newcommand{\Zh}{\widehat{Z}}

\font\co=lcircle10
\def\petit#1{{\scriptstyle #1}}

\def\jr{\smash{\raise2pt\hbox{\co \rlap{\rlap{\char'005} \char'007}}
               \raise6pt\hbox{\rlap{\vrule height5pt}}
               \raise2pt\hbox{\rlap{\hskip4pt \vrule height0.4pt depth0pt
                width5.7pt}}
               \raise2pt\hbox{\rlap{\hskip-9.5pt \vrule height.4pt depth0pt
                width6.2pt}}
               \lower6pt\hbox{\rlap{\vrule height4.5pt}}}}
\def\je{\smash{\raise2pt\hbox{\co \rlap{\rlap{\char'005}
                \phantom{\char'007}}}\raise6pt\hbox{\rlap{\vrule height5pt}}
               \raise2pt\hbox{\rlap{\hskip-9.5pt \vrule height.4pt depth0pt
                width6.2pt}}}}
\def\er{\smash{\raise2pt\hbox{\co \rlap{\rlap{\phantom{\char'005}} \char'007}}
               \raise2pt\hbox{\rlap{\hskip4pt \vrule height0.4pt depth0pt
                width5.7pt}}
               \lower6pt\hbox{\rlap{\vrule height4.5pt}}}}
\def\+{\smash{\lower6pt\hbox{\rlap{\vrule height17pt}}
                \raise2pt\hbox{\rlap{\hskip-9pt \vrule height.4pt depth0pt
                width18.7pt}}}}
\def\hor{\smash{\raise2pt\hbox{\rlap{\hskip-9.5pt \vrule height.4pt depth0pt
                width19.2pt}}}}
\def\ver{\smash{\lower6pt\hbox{\rlap{\vrule height17pt}}}}

\def\ph{\phantom{+}}

\def\sr{\makebox[0ex]{$*$}}

\def\mcc{\multicolumn{1}{@{}c@{}}}

\def\textcross{\ \smash{\lower4pt\hbox{\rlap{\hskip4.15pt\vrule height14pt}}
                \raise2.8pt\hbox{\rlap{\hskip-3pt \vrule height.4pt depth0pt
                width14.7pt}}}\hskip12.7pt}
\def\textelbow{\ \hskip.1pt\smash{\raise2.8pt%
                \hbox{\co \hskip 4.15pt\rlap{\rlap{\char'005} \char'007}
                \lower6.8pt\rlap{\vrule height3.5pt}
                \raise3.6pt\rlap{\vrule height3.5pt}}
                \raise2.8pt\hbox{%
                  \rlap{\hskip-7.15pt \vrule height.4pt depth0pt width3.5pt}%
                  \rlap{\hskip4.05pt \vrule height.4pt depth0pt width3.5pt}}}
                \hskip8.7pt}

\begin{document}

\title{Four positive formulae for type~\(A\) quiver polynomials}
\author{Allen Knutson}
\thanks{AK was partly supported by the Sloan Foundation and NSF}
\address{Mathematics Department\\ UC Berkeley\\ Berkeley, California}
\email{allenk@math.berkeley.edu}
\author{Ezra Miller}
\thanks{EM and MS were partly supported by the NSF}
\address{Mathematical Sciences Research Institute\\ Berkeley, California}
\email{ezra@math.umn.edu}
\author{Mark Shimozono}
\address{Mathematics Department\\ Virginia Tech \\ Blacksburg, Virginia}
\email{mshimo@math.vt.edu}
\date{31 August 2003}

\begin{abstract}

\noindent
We give four positive formulae for the (equioriented type~$A$) quiver
polynomials of Buch and Fulton \cite{BF}.  All four formulae are
combinatorial, in the sense that they are expressed in terms of
combinatorial objects of certain types: Zelevinsky permutations,
lacing diagrams, Young tableaux, and pipe dreams (also known as
rc-graphs).  Three of our formulae are multiplicity-free and
geometric, meaning that their summands have coefficient~$1$, and
correspond bijectively to components of a torus-invariant scheme.  The
remaining (presently non-geometric) formula was conjectured for by
Buch and Fulton in terms of factor sequences of Young
tableaux~\cite{BF}; our proof of it proceeds by way of a new
characterization of the tableaux counted by quiver constants.  All
four formulae come naturally in ``doubled'' versions, two for {\em
double quiver polynomials}, and the other two for their stable
versions, the {\em double quiver functions}, where setting half the
variables equal to the other half specializes to the ordinary case.

Our method begins by identifying quiver polynomials as multidegrees
\cite{joseph,rossmann} via equivariant Chow groups \cite{EG98}.  Then
we make use of Zelevinsky's map from quiver loci to open subvarieties
of Schubert varieties in partial flag manifolds \cite{Zel85}.
Interpreted in equivariant cohomology, this lets us write double
quiver polynomials as ratios of double Schubert polynomials
\cite{LS82} associated to Zelevinsky permutations; this is our first
formula.  In~the process, we provide a simple argument that Zelevinsky
maps are scheme-theoretic isomorphisms (originally proved in~\cite{LM}
using standard monomial theory).  Writing double Schubert polynomials
in terms of pipe dreams \cite{FK96} then provides another geometric
formula for double quiver polynomials, via \cite{grobGeom}.  The
combinatorics of pipe dreams for Zelevinsky permutations implies an
expression for double quiver functions in terms of products of Stanley
symmetric functions \cite{St}.  A~degeneration of quiver loci (orbit
closures of $\hspace{-.2ex}\glv$ on quiver representations) to unions
of products of matrix Schubert varieties \cite{Ful92,grobGeom}
identifies the summands in our Stanley function formula
combinatorially, as lacing diagrams that we construct based on the
strands of Abeasis and Del\thinspace{}\thinspace{}Fra in the
representation theory of quivers \cite{AD}.  Finally, we apply the
combinatorial theory of key polynomials to pass from our lacing
diagram formula to a double Schur function formula in terms of
peelable tableaux \cite{RSplact95,RSpeelable98}, and from there to a
`double stable' generalization of the Buch--Fulton conjecture.
\end{abstract}

\maketitle

{}

\tableofcontents
\raggedbottom

\part*{Introduction}

\subsection*{Overview}

Universal formulae for cohomology classes appear in a number of
different guises.  In topology, classes with attached names of
Pontrjagin, Chern, and Stiefel--Whitney arise as obstructions to
vector bundles having linearly independent sections.  Universality of
formulae for these and other more general classes can be traced, for
our purposes, to the fact that they---the formulae as well as the
classes---live canonically on classifying spaces, from which they are
pulled back to arbitrary spaces by classifying~maps.

In algebra, cohomology classes on certain kinds of varieties
(projective space, for example) are called `degrees', and many
universal formulae show up in degree calculations for classical
ideals, such as those generated by minors of fixed size in generic
matrices \cite{giambelli} (or see \cite[Chapter~14]{FulIT}).

In geometry, pioneering work of R.\thinspace{}Thom \cite{thom55}
associated cohomology classes to sets of {\em critical points}\/ of
generic maps between manifolds, where the differential drops rank.
Using slightly different language, the set of critical points can be
characterized as the {\em degeneracy locus}\/ for the associated
morphism of tangent bundles.  Subsequently there have appeared
numerous extensions of Thom's notion of degeneracy locus, such as to
maps between pairs of arbitrary vector bundles.

The cohomology classes Poincar\'e dual to degeneracy loci for (certain
collections of) morphisms of complex vector bundles are expressible as
polynomials in the Chern roots of the given vector bundles.  It seems
to be a general principle that coefficients in universal such
expressions as sums of simpler polynomials always seem to be governed
by combinatorial rules.  This occurs, for example, in
\cite{Ful92,Bu01a,BKTY02}.
However, even when it is possible to prove an explicit combinatorial
formula, it has usually been unclear how the geometry of degeneracy
loci reflects the combinatorics~directly.

Our original motivation for this work was to close this gap, by
bringing combinatorial universal formulae for cohomology classes of
degeneracy loci into the realm of geometry, continuing the point of
view set forth in \cite{grobGeom}.  To do so, we reduce to the
algebraic perspective mentioned above in terms of degrees of
determinantal ideals, by way of the topological perspective mentioned
above in terms of classifying spaces, as in
\cite{Kaz97,FRthomPoly,grobGeom}.  This reduction then allows us to
manufacture formulae that are simultaneously geometric as well as
combinatorial, by applying flat degenerations of subvarieties inside
finite-dimensional torus representations.  The cohomology class
remains unchanged in the degeneration, and then it is given by summing
the classes of all components of the special~fiber.

The degenerations we employ result from one-parameter linear torus
actions, and hence have Gr\"obner bases as their natural language.
Our theorems on orbit closures and their combinatorics explicitly
generalize results from the theory of ideals generated by minors of
fixed size in generic matrices, where the prototypical result is
Giambelli's degree formula \cite{giambelli}.  The ideals that interest
us here are generated by minors in {\em products}\/ of generic
matrices.  It is intriguing that our combinatorial analysis of these
general determinantal ideals proceeds (using the `Zelevinsky map') via
Schubert determinantal ideals \cite{Ful92,grobGeom}, which are
generated by minors of varying sizes in a {\em single}\/ matrix of
variables.  Arbitrary Schubert conditions may seem, consequently, to
be in some sense more general; but while this may be true, the special
form taken by the Schubert rank conditions on single matrices here
give rise to a substantially richer combinatorial structure.  This
richness is fully borne out only by combining methods based on
Schubert determinantal ideals with an approach in terms of minors in
products of generic matrices, which is rooted more directly in
representation theory of quivers.

Our technique of orbit degeneration automatically produces positive
geometric formulae that are universal, since they essentially live on
classifying spaces.  However, there still remains the elucidation of
what combinatorics most naturally describes these formulae, or indexes
its components.  Well over half of our exposition in this paper is
dedicated to unearthing rich and often surprising interconnections
between a number of combinatorial objects, some known and some new.
Combinatorics---of objects including Zelevinsky permutations, lacing
diagrams, Young tableaux, and pipe dreams (also known as
rc-graphs)---forms the bridge between positive geometric formulae and
the explicit algebra of several families of polynomials, including
Schur functions, Schubert polynomials, Stanley symmetric functions,
and now quiver polynomials.  Working in the broader context of double
quiver polynomials and their stable versions, the double quiver
functions, allows the extra flexibility required for our proof of the
Buch--Fulton conjecture~\cite{BF}.  Special cases of this conjecture
appear in \cite{Bu01a,BKTY02}.

\subsection*{Acknowledgments}
We would like to thank Anders Buch and Bill Fulton for suggesting the
problem that concerns us here, and for patiently listening, with
Sergey Fomin, Richard Stanley, and Alex Yong, to endless hours of us
presenting potential and actual meth\-ods toward its solution.  Frank
Sottile and Chris Woodward invited
us
to their superb AMS Special Session on Modern Schubert Calculus at
Northeastern University in October 2002; we are indebted to them for
the opportunity at that meeting both to present our positivity proof
for the quiver constants and our conjectural component formula,
as well as to learn of
\cite{BKTY02} for the first time.  We are grateful to the Banff
International Research Station (BIRS), as well as the organizers of
the May 2003 meeting on Algebraic Combinatorics there, for bringing
the three of us together into the same room for the first time,
thereby allowing us to complete the first version of this paper.
Comments that improved the exposition were provided by Anders Buch,
Bill Fulton, Peter Magyar, \mbox{Rich{\'a}rd Rim{\'a}nyi, and
Frank~Sottile}.

\subsection*{The four formulae}

The formulae we consider in this paper are for cohomology classes of
degeneracy loci for type~$A$ equioriented quivers of vector bundles.
This simply means we start with sequences \mbox{$E_0 \to E_1 \to
\cdots \to E_n$} of vector bundle morphisms over a fixed base.  Since
we immediately reduce in Theorem~\ref{t:BFtoKMS} to considering the
equivariant classes of universal degeneracy loci, or quiver loci (to
be defined shortly), we present our exposition in that context,
refering the reader to \cite{BF} for an introduction in the language
of vector bundles.  We work over an arbitrary field~$\kk$.

Consider the vector space $\homv$ of sequences
\mbox{$V_0\stackrel{\phi_1}\too V_1\stackrel{\phi_2}\too\cdots\
\stackrel{\phi_{n-1}}\too V_{n-1}\stackrel{\phi_n}\too V_n$}
of linear transformations between vector spaces of dimensions
$r_0,\ldots,r_n$.  We think of {\em quiver representations} $\phi \in
\homv$ as sequences of matrices.  Each matrix list~$\phi$ determines
its {\em rank array}\/ $\rr(\phi) = (r_{ij}(\phi))_{i \leq j}$, where
$r_{ij}(\phi)$ for $i < j$ equals the rank of the composite map $V_i
\to V_j$, and $r_{ii} = \dim(V_i)$.  The data of a rank array
determines a {\em quiver locus}\/ $\Omega_\rr$, defined as the
subscheme of matrix lists $\phi \in \homv$ with rank array dominated
by~$\rr$ entrywise: $r_{ij}(\phi) \leq r_{ij}$ for all $i \leq j$.
Thus~$\Omega_\rr$ is the zero scheme of the ideal~$I_\rr$ in the
coordinate ring of~$\homv$ generated by all rank $(1+r_{ij})$ minors
in the product of the appropriate $j - i$ generic matrices, for all $i
< j$ (Definition~\ref{d:quivlocus}).

We assume throughout that $\rr = \rr(\phi)$ for some matrix list $\phi
\in \Omega_\rr$.  This is a nontrivial condition equivalent to the
irreducibility of~$\Omega_\rr$.  Alternatively, it means that
$\Omega_\rr$ is an orbit closure for the group $\glv = \prod_{i=0}^n
\gl(V_i)$ that acts on~$\homv$ by change of basis in each vector
space~$V_i$.  In particular, $\Omega_\rr$ is stable under the action
of any torus in~$\glv$.  Picking a maximal torus~$\tv$, we define the
{\em quiver polynomial}\/ as the $T$-multidegree of~$\Omega_\rr$.  We show
in Proposition~\ref{p:chow} that this polynomial coincides with the
$T$-equivariant class of~$\Omega_\rr$ in the Chow ring $A^*_T(\homv)$,
or alternatively in cohomology when $\kk = \CC$.  This Chow ring is a
polynomial ring~$\ZZ[\xx_\rr]$ over the integers in an
alphabet~$\xx_\rr$ of size $r_0 + \cdots + r_n$, the union of bases
$\xx^0,\ldots,\xx^n$ for the weight lattices of the maximal tori
in~$\glv(V_i)$.  Since~$\Omega_\rr$ is stable under~$\glv$ and not
just the torus~$T$, its equivariant class is symmetric in each of the
alphabets~$\xx^i$.  More naturally, this is precisely the statement
that the quiver polynomial lies inside $A^*_\glv(\homv) \subset
A^*_T(\homv)$, and in this sense does not depend on our choice of~$T$.
We express the quiver polynomial determined by the rank array~$\rr$ as
$\QQ_\rr(\xx-\oxx)$ for reasons that should become clear in the ratio
formula and double~versions.

Next we present our four positive combinatorial formulae for quiver
polynomials, in the same order that we will prove them in the main
body of the text.  Besides making the overview of their proofs more
coherent, this choice will serve to emphasize an important point that
is worth bearing in mind before seeing the statements: while each
formula stands well enough on its own, this paper is not merely a
catalog of four different perspectives on quiver loci.  The
principles underlying each formula seep throughout the others,
infusing all of them with just enough additional insight to make the
proofs work.  The interdependencies will become increasingly clear as
we sketch the proofs later in this Introduction.

For each formula, we shall give here a precise statement, although
many details of definitions will be left for later.  Each formula will
be presented along with a pointer to its location in the text, and an
example using the following $n=3$ rank array~$\rr = (r_{ij})$.
\begin{stareqn} \label{rr}
\rr \ \:=\:\
  \begin{array}{cccc|c}
  3 & 2 & 1 & 0 & i \diagup j \\ \hline
    &   &   & 1 & 0 \\
    &   & 3 & 1 & 1 \\
    & 3 & 2 & 1 & 2 \\
  1 & 1 & 1 & 0 & 3
  \end{array}
\end{stareqn}%
We use alphabets $\aa = \xx^0$, $\bb = \xx^1$, $\cc = \xx^2$, and $\dd
= \xx^3$ in the examples, keeping in mind that for the above
ranks~$\rr$, this means that $\aa = \{a_1\}$, $\bb = \{b_1,b_2,b_3\}$,
$\cc = \{c_1,c_2,c_3\}$, and $\dd = \{d_1\}$ are alphabets of sizes
$1$, $3$, $3$, and~$1$.  With these conventions, the quiver polynomial
$\QQ_\rr(\xx-\oxx)$ can be calculated to be
\begin{eqnarray*}
  \QQ_\rr &=& (b_1+b_2+b_3-c_1-c_2-c_3)(a_1-d_1).
\end{eqnarray*}

In the general setting, we distinguish between the symbol~$\xx$, in
which every alphabet in the sequence $\xx^0,\ldots,\xx^n$ is {\em
infinite}, and~$\xx_\rr$, where $\xx^i$ has cardinality~$r_i$.  In the
notation $\QQ_\rr(\xx-\oxx)$ for quiver polynomials, only the finitely
many variables $\xx_\rr \subset \xx$ actually appear; that is, we can
view $\QQ_\rr$ as taking infinitely many variables as input, even
though only finitely many contribute.

\subsubsection*{Ratio formula (Zelevinsky permutation)}

Associated to each rank array~$\rr$ is a permutation $v(\rr)$ in the
symmetric group~$S_d$ for $d = r_0 + \cdots + r_n$.  This {\em
Zelevinsky permutation}\/~$v(\rr)$ has a block structure determined
by~$\rr$, with block rows of heights $r_0,\ldots,r_n$, and block
columns of widths $r_0,\ldots,r_n$ proceeding from {\em right to
left}.  Zelevinsky permutations enjoy a number of pleasant properties,
such as having no descents within any block row, and all entries zero
above the block superantidiagonal.  Proposition~\ref{p:zel} gives a
complete characterization.

Given the dimensions $r_0,\ldots,r_n$, there is a unique Zelevinsky
permutation $v(\homv)$ of minimal length.  It is associated to the
maximal rank array, whose $i\hspace{-.2ex}j$ entry is defined to be
$r_{ij}(\homv) = \min_{i\leq k \leq j}\{r_k\}$.  This is the maximum
possible rank for a map $V_i \to V_j$ that factors through $V_k$ for
$i \leq k \leq j$.  The Zelevinsky permutation~$v(\homv)$ is best
characterized by its {\em diagram}.  Generally, the diagram of a
permutation~$v$ consists of all cells neither due south of nor due
east of a nonzero entry in the permutation matrix of~$v$.  For
$v(\homv)$, the diagram consists precisely of all cells strictly above
the block superantidiagonal.  These cells also lie inside the diagram
of~$v(\rr)$ for every rank array~$\rr$.  The particularly simple form
of this diagram says that $v(\homv)$ is a dominant permutation, whose
matrix Schubert variety (Section~\ref{sec:schub}) is a linear
subspace.

The last ingredients for our first formula are the {\em double
Schubert polynomials} of Lascoux and Sch\"utzenberger \cite{LS82}
(Section~\ref{sec:schub}).  These represent torus-equivariant
cohomology classes of Schubert varieties in flag manifolds, and are
characterized geometrically as being the multidegrees of matrix
Schubert varieties \cite{grobGeom}.  Double Schubert polynomials
\mbox{$\SS_w(\xx-\oyy)$} take as input two alphabets $\xx$ and~$\yy$,
of sizes at least $r$ and~$r'$, along with an $r \times r'$ {\em
partial permutation}~$w$, by which we mean an $r \times r'$ matrix
with no more than one~$1$ in each row and column, and zeros elsewhere.
If~$w$ happens to lie in~$S_d$, and both alphabets have size~$d$, then
we write~$\oxx_\rr$ for the block reversal of~$\xx_\rr$.  In general,
$\oxx$ denotes the sequence $\xx^n,\ldots,\xx^0$ of alphabets.  The
denominator in what follows is a product of linear forms, because
$v(\homv)$ has such a simple diagram.

\begin{movedThm}{Theorem (Ratio formula).}
$\displaystyle
\QQ_\rr(\xx-\oxx)\ =\ \frac{\SS_{v(\rr)}(\xx_\rr-\oxx_\rr)}
{\SS_{v(\homv)}(\xx_\rr-\oxx_\rr)}$.
\end{movedThm}

This result is proved in Theorem~\ref{c:ratio}.  It is geometric in
the sense that it interprets in equivariant cohomology the {\em
Zelevinsky map}\/ (Section~\ref{sec:zel}), which takes quiver
loci~$\Omega_\rr$ isomorphically to open subvarieties of Schubert
varieties in partial flag manifolds.  The geometry of the Zelevinsky
map explains why Schubert polynomials appear.

\begin{example*}
The rank array~$\rr$ in~\eqref{rr} has the associated
Zelevinsky permutation
\begin{eqnarray*}
v(\rr) &=&
  \begin{array}{@{}l|%
        @{}  c@{\:}|%
        @{\:}c@{\:}@{\:}c@{\:}@{\:}c@{\:}|%
        @{\:}c@{\:}@{\:}c@{\:}@{\:}c@{\:}|%
        @{\:}c@{\:}|@{}}
  \cline{2-9}
   5 &\: *  &  * & * & *  & \ti&\cd&\cd & \cd
  \\\cline{2-9}
   2 &\: *  & \ti&\cd&\cd & \cd&\cd&\cd & \cd
  \\
   3 &\: *  & \cd&\ti&\cd & \cd&\cd&\cd & \cd
  \\
   6 &\: *  & \cd&\cd&\sq & \cd&\ti&\cd & \cd
  \\\cline{2-9}
   1 &\:\ti & \cd&\cd&\cd & \cd&\cd&\cd & \cd
  \\
   4 &\:\cd & \cd&\cd&\ti & \cd&\cd&\cd & \cd
  \\
   8 &\:\cd & \cd&\cd&\cd & \cd&\cd&\sq & \ti
  \\\cline{2-9}
   7 &\:\cd & \cd&\cd&\cd & \cd&\cd&\ti & \cd
  \\\cline{2-9}
  \end{array}
\end{eqnarray*}
The numbers running down the left side of~$v(\rr)$ constitute one-line
notation for~$v(\rr)$, so $1 \mapsto 5$, $2 \mapsto 2$, $3 \mapsto 3$,
$4 \mapsto 6$, and so on.  We have replaced all $1$~entries
in~$v(\rr)$ with $\ti$~entries.  The asterisks~$*$ and boxes~$\sq$ all
denote cells in the diagram of~$v(\rr)$, whereas only the $*$~entries
denote cells in the diagram of the Zelevinsky permutation $v(\homv) =
52341678$ of the maximal rank array~$\rr(\homv)$.  That $v(\homv)$ has
such a simple diagram means that in fact the denominator in the
formula
\begin{eqnarray*}
  \QQ_\rr &=& \frac{\SS_{52361487}(\aa,\bb,\cc,\dd-\dd,\cc,\bb,\aa)}
                {\SS_{52341678}(\aa,\bb,\cc,\dd-\dd,\cc,\bb,\aa)}
\end{eqnarray*}
equals the product
$
  (a_1-c_3)(a_1-c_2)(a_1-c_1)(a_1-d_1)(b_1-d_1)(b_2-d_1)(b_3-d_1)
$
of linear factors corresponding to the locations of $*$~entries.
\end{example*}
\begin{excise}{%
  Double quiver polynomial, via ratio formula:
  $(b_1-\dc_3+b_2-\dc_2+b_3-\dc_1)
  (a_1-\db_3+b_1-\db_2+b_2-b_1+b_3-\dc_3+c_1-\dc_2+c_2-\dc_1+c_3-\dt_1)$
}\end{excise}%

\subsubsection*{Pipe formula}

One way to see that the denominator in the ratio formula is always the
corresponding product of linear factors, and divides the numerator, is
through our next combinatorial formula.  It is in fact little more
than an application of the Fomin--Kirillov double version \cite{FK96}
of the Billey--Jockusch--Stanley formula for Schubert polynomials
\cite{BJS,FS}.

Briefly, for a permutation $v \in S_d$, denote by~$\rp(v)$ its set of
{\em reduced pipe dreams}\/ (also known as {\em planar histories}\/ or
{\em rc-graphs}\/).  These are fillings the $d \times d$ grid with
square tiles of the form $\textcross$ or~$\textelbow$, in which two
pipes either cross or avoid each other
(Section~\ref{sec:pipe2quiver}).  The adjective `reduced' indicates
that no pair of pipes crosses more than once, and the permutation~$v$
dictates that the pipe entering row~$i$ exits out of column~$v(i)$.
Label the rows of the grid from top to bottom with the ordered
alphabet~$\xx_\rr$, and label the columns from left to right with the
ordered alphabet~$\oxx_\rr$.  Given a pipe dream~$D$, define the
monomial $(\xx-\oxx)^D$ as the product of linear binomials
$(x_{\row(+)} - x_{\col(+)})$, one for each \mbox{crossing
tile~$\textcross$ in~$D$}.

For example, the permutation $v(\homv)$ has only one reduced pipe
dream~$D_\homv$, with crossing tiles~$\textcross$ in every square
strictly above the block superantidiagonal.  Hence its monomial
$(\xx-\oxx)^{D_\homv}$ equals the product of linear binomials
corresponding to $*$~entries in the diagram of {\em any}\/ Zelevinsky
permutation~$v(\rr)$, as long as $\rr$ gives rank conditions on quiver
representations with dimension vector $(r_0,\ldots,r_n)$.  We identify
pipe dreams with their sets of crossing tiles, as subsets of the $d
\times d$ grid.

\begin{movedThm}{Theorem (Pipe formula).}
$\displaystyle \QQ_\rr(\xx-\oxx)\ =\sum_{D \in \rp(v(\rr))}
(\xx-\oxx)^{D \minus D_\homv}$.
\end{movedThm}

This follows from Theorem~\ref{t:QQrr}.  It is geometric in the sense
that summands on the right hand side are equivariant classes of
coordinate subspaces in a flat (Gr\"obner) degeneration
\cite[Theorem~B]{grobGeom} of the matrix Schubert variety for~$v(\rr)$
(Section~\ref{sec:schub}).  This is the largest matrix Schubert
variety whose quotient modulo the appropriate parabolic subgroup
contains the Zelevinsky image of~$\Omega_\rr$ as an open dense
subvariety.

\begin{example*}
Here is a typical pipe dream for the Zelevinsky permutation $v(\rr)$
from the previous Example, drawn in two ways.
$$
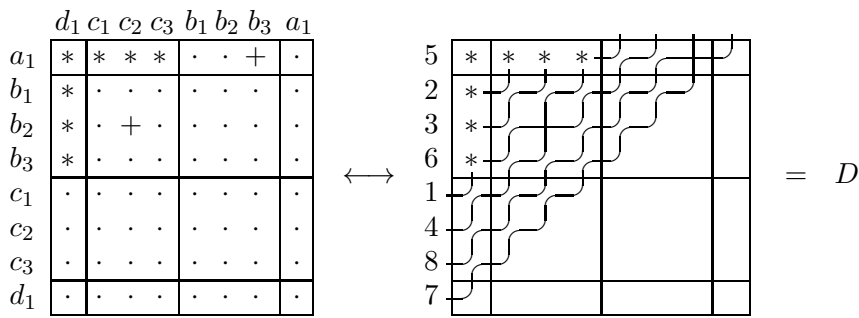
\def\kj{\hspace{.1ex}}
\begin{array}{@{}l|c@{\ \,}|@{}c@{\:}c@{\:}c@{\ }%
                |@{}c@{\kj}c@{\,}c|@{\ }c@{\ }|@{}}
\multicolumn{9}{@{}c@{}}{}\\[-3ex]
\mcc{}&\mcc{\ph}&\ph&\ph&\mcc{\ph}&\ph&\ph&\mcc{\ph}&\mcc{\ph}
\\[-2ex]
\mcc{}&\mcc{d_1}&c_1&c_2&\mcc{c_3\,}&\,b_1&\kj\kj\kj b_2&\mcc{\!\!\!b_3}
                                                            &\mcc{\!a_1}
\\\cline{2-9}
a_1&
  \sr&\sr&\sr&\sr&\cd&\cd& + &\cd
\\\cline{2-9}
b_1&
  \sr&\cd&\cd&\cd&\cd&\cd&\cd&\cd
\\
b_2&
  \sr&\cd& + &\cd&\cd&\cd&\cd&\cd
\\
b_3&
  \sr&\cd&\cd&\cd&\cd&\cd&\cd&\cd
\\\cline{2-9}
c_1&
  \cd&\cd&\cd&\cd&\cd&\cd&\cd&\cd
\\
c_2&
  \cd&\cd&\cd&\cd&\cd&\cd&\cd&\cd
\\
c_3&
  \cd&\cd&\cd&\cd&\cd&\cd&\cd&\cd
\\\cline{2-9}
d_1&
  \cd&\cd&\cd&\cd&\cd&\cd&\cd&\cd
\\\cline{2-9}
\end{array}
\quad\raisebox{-1.5ex}{$\longleftrightarrow$}\quad
\begin{array}{@{}l|c|ccc|ccc|c|@{}}
\multicolumn{9}{@{}c@{}}{}\\[-2.2ex]
\mcc{}&\mcc{\ \,\ \,\ }&\ &\ &\mcc{\ \,\ \,\ }&\ &\ &\mcc{\ \,\ \,\ }&
                                                     \mcc{\ \,\ \,\ }
\\\cline{2-9}
 5 &\sr&\sr&\sr&\sr&\jr&\jr&\+ &\je
\\\cline{2-9}
 2 &\sr&\jr&\jr&\jr&\jr&\jr&\je&
\\
 3 &\sr&\jr&\+ &\jr&\jr&\je&   &
\\
 6 &\sr&\jr&\jr&\jr&\je&   &   &
\\\cline{2-9}
 1 &\jr&\jr&\jr&\je&   &   &   &
\\
 4 &\jr&\jr&\je&   &   &   &   &
\\
 8 &\jr&\je&   &   &   &   &   &
\\\cline{2-9}
 7 &\je&   &   &   &   &   &   &
\\\cline{2-9}
\end{array}
\quad\, \raisebox{-1.5ex}{$= \ \ D$}
$$
The right hand diagram depicts all the tiles except those in the
diagram of~$v(\homv)$, whose $*$'s are all really $\textcross$ tiles,
and those parts of pipes below the main antidiagonal (since the
``sea'' of~$\textelbow$ tiles there can be confusing to look at).  The
elbows have been entirely left out of the left hand diagram, but its
row and column labels are present.  The monomial $(\xx-\oxx)^{D \minus
D_\homv}$ for the above pipe dream~$D$ is simply $(a_1-b_3)(b_2-c_2)$.

For this permutation $v(\rr)$, moving the~\cross
in the $\bb\cc$ block to either of the two remaining available cells
on the antidiagonal of the $\bb\cc$ block produces another reduced
pipe dream for~$v(\rr)$.  Independently, the other~\cross can move
freely along the antidiagonal on which it sits.  None of the
$*$~entries can move.  Therefore the right hand side of the pipe
formula becomes a product of two linear forms, namely the sums
\begin{eqnarray*}
& \big((a_1\!-\!b_3) + (b_1\!-\!b_2) + (b_2\!-\!b_1) + (b_3\!-\!c_3) +
  (c_1\!-\!c_2) + (c_2\!-\!c_1) + (c_3\!-\!d_1)\big)
\\
& \llap{and}\qquad \big((b_1\!-\!c_3) + (b_2\!-\!c_2) +
  (b_3\!-\!c_1)\big)
\end{eqnarray*}
of the binomials associated to cells on the corresponding
antidiagonals.  Cancelation occurs in the longer form above, to give
$(a_1-d_1)$, so the product of these two forms returns $\QQ_\rr =
(b_1+b_2+b_3-c_1-c_2-c_3)(a_1-d_1)$ again.
\end{example*}

\subsubsection*{Component formula (lacing diagrams)}

Double Schubert polynomials are not symmetric,
in general.  However, Stanley introduced certain symmetrized versions
\cite{St}, now called {\em double Stanley symmetric functions}\/ or
{\em stable double Schubert polynomials}\/ and denoted
by~$F_w(\xx-\oyy)$.  They are indexed again by partial
permutations~$w$, and take as arguments a pair of infinite alphabets
$\xx$ and~$\yy$.  When we evaluate such a symmetric function on a
finite alphabet, we mean to set all remaining variables to zero.
Stanley functions are produced by an algebraic limiting procedure
(Proposition~\ref{p:doublestanley}) applied to Schubert polynomials;
hence the term `stable polynomial'.  In general, when the double
Schubert polynomial~$\SS_w$ and double Stanley function~$F_w$ are
evaluated on the same pair of finite alphabets, the polynomial~$F_w$,
which is symmetric separately in each of the two alphabets, will tend
to have many more terms.

Partial permutations as we defined them are matrices with entries in
the field~$\kk$, so it makes sense to say that a list $\ww =
(w_1,\ldots,w_n)$ of partial permutations lies in~$\homv$, if each
$w_i$ has size $r_{i-1} \times r_i$.  Such lists of partial
permutation matrices can be identified with nonembedded graphs drawn
in the plane, called {\em lacing diagrams}, that we define in
Section~\ref{part:lacing}.  The vertex set of~$\ww$ consists of $n+1$
columns of dots, where column~$i$ has $r_i$~dots.  An edge of~$w$
connects the dot at height~$\alpha$ in column~\mbox{$i-1$} to the dot
at height~$\beta$ in column~$i$ if the $\alpha\beta$ entry of~$w_i$
is~$1$.  Lacing diagrams come with a natural notion of length derived
from the Bruhat order, and the minimum possible length for a lacing
diagram with rank array~$\rr$ is the codimension of the quiver
locus~$\Omega_\rr$, which equals the degree of~$\QQ_\rr(\xx-\oxx)$.
Denote by~$W(\rr)$ the set of minimal length lacing diagrams with rank
array~$\rr$ (characterized combinatorially in
Theorem~\ref{t:mindiag}).

Since Schubert polynomials can be indexed by partial permutations, the
notation
\begin{eqnarray*}
  \SS_\ww(\xx-\oxx) &=& \SSS{w_1}01 \cdots \SSS{w_n}{n-1}n
\end{eqnarray*}
for products of double Schubert polynomials indexed by the partial
permutations in~$\ww$ makes sense.  (For consistency with
$\QQ_\rr(\xx-\oxx)$, we write $\SS_\ww(\xx-\oxx)$ even though
$\SS_{w_i}$ takes $\xx^{i-1}-\xx^i$ for input, not
$\xx^{i-1}-\xx^{n-i}$.)  Similarly, we have the product notation
\begin{eqnarray*}
  F_\ww(\xx-\oxx) &=& \FFF{w_1}01 \cdots \FFF{w_n}{n-1}n
\end{eqnarray*}
for double Stanley symmetric functions.  Again, the $n+1$ alphabets in
the latter equation are infinite, but we are allowed to write
$F_\ww(\xx_\rr-\oxx_\rr)$ if we want to evaluate $F_\ww$ on sequences
of finite alphabets.

\begin{movedThm}{Theorem (Component formula).}
$\begin{array}[t]{@{}rcl@{}}
\QQ_\rr(\xx-\oxx)
  &=&\displaystyle \sum_{\ww \in W(\rr)} \SS_\ww(\xx-\oxx)
\\[3.5ex]\mbox{}
  &=&\displaystyle \sum_{\ww \in W(\rr)} F_\ww(\xx_\rr-\oxx_\rr)
\end{array}$
\end{movedThm}

The two versions are Corollary~\ref{c:component} and
Corollary~\ref{c:specialize}.
The theorem implies, in particular, that all the extra terms
in the Stanley version cancel.

The first sum in this theorem is combinatorially positive in a manner
that most directly (of all our formulae) reflects the geometry of
quiver loci.  The basic idea is to flatly degenerate the group
action of~$\glv$ under which the quiver locus~$\Omega_\rr$ is an orbit
closure (Section~\ref{part:degen}).  As $\glv$ degenerates, so do its
orbits, and the flat limits of the orbits are stable under the action
of the limiting group.  (The general version of this statement is
Proposition~\ref{p:flat}; the specific case of interest to us is
Proposition~\ref{p:action}.)  However, the limit need not be
irreducible, and its components are precisely the closures of orbits
through minimal length lacing diagrams (recall that lacing diagrams
are matrix lists in~$\homv$).  This, together with the statement that
the components in the degenerate limit of~$\Omega_\rr$ are generically
reduced, is precisely the content of Theorem~\ref{t:components}.  It
is our main geometric theorem concerning quiver loci, and it
immediately implies the Schubert version of the component formula.

\begin{example*} \label{ex:blah}
Continuing with $\rr$ from~\eqref{rr}, the set~$W(\rr)$ consists of
the following three minimal length lacing diagrams, with their partial
permutation lists underneath:
$$
\begin{array}{@{}c@{\qquad}c@{\qquad}c@{}}
  \pspicture[.1](0,1)(3,3)
  \psdots(0,1)(1,1)(1,2)(1,3)(2,1)(2,2)(2,3)(3,1)
  \psline(0,1)(1,1)(2,1)
  \psline(1,2)(2,2)(3,1)
  \endpspicture
&
  \pspicture[.1](0,1)(3,3)
  \psdots(0,1)(1,1)(1,2)(1,3)(2,1)(2,2)(2,3)(3,1)
  \psline(0,1)(1,1)(2,2)
  \psline(1,2)(2,1)(3,1)
  \endpspicture
&
  \pspicture[.1](0,1)(3,3)
  \psdots(0,1)(1,1)(1,2)(1,3)(2,1)(2,2)(2,3)(3,1)
  \psline(0,1)(1,2)(2,2)
  \psline(1,1)(2,1)(3,1)
  \endpspicture
\\[2ex]
\left(
        \left[\tinyrc{\begin{array}{@{}c@{\ }c@{\ }c@{}}
                        1&0&0\end{array}}\right],
        \left[\tinyrc{\begin{array}{@{}c@{\ }c@{\ }c@{}}
                        1&0&0\\
                        0&1&0\\
                        0&0&0\end{array}}\right]\!,
        \left[\tinyrc{\begin{array}{@{}c@{}}
                        0    \\
                        1    \\
                        0    \end{array}}\right]
\right)
&
\left(
        \left[\tinyrc{\begin{array}{@{}c@{\ }c@{\ }c@{}}
                        1&0&0\end{array}}\right],
        \left[\tinyrc{\begin{array}{@{}c@{\ }c@{\ }c@{}}
                        0&1&0\\
                        1&0&0\\
                        0&0&0\end{array}}\right]\!,
        \left[\tinyrc{\begin{array}{@{}c@{}}
                        1    \\
                        0    \\
                        0    \end{array}}\right]
\right)
&
\left(
        \left[\tinyrc{\begin{array}{@{}c@{\ }c@{\ }c@{}}
                        0&1&0\end{array}}\right],
        \left[\tinyrc{\begin{array}{@{}c@{\ }c@{\ }c@{}}
                        1&0&0\\
                        0&1&0\\
                        0&0&0\end{array}}\right]\!,
        \left[\tinyrc{\begin{array}{@{}c@{}}
                        1    \\
                        0    \\
                        0    \end{array}}\right]
\right)
\end{array}
$$
%
%
%
%
%
The two versions (Schubert and Stanley) of the component formula read
$$
\begin{array}{@{}rcc@{\ }c@{\ }c@{\ }c@{\ }c@{}}
\QQ_\rr
  &=& \SS_{1243}(\bb-\cc)\SS_{213}(\cc-\dd) &+& \SS_{2143}(\bb-\cc)
  &+& \SS_{213}(\aa-\bb)\SS_{1243}(\bb-\cc)
\\[.5ex]
  &=& F_{1243}(\bb-\cc)F_{213}(\cc-\dd) &+& F_{2143}(\bb-\cc)
  &+& F_{213}(\aa-\bb)F_{1243}(\bb-\cc)
\end{array}
$$
Note that the Stanley symmetric functions are evaluated on finite
alphabets $\aa,\bb,\cc,\dd$.  We have indexed the Schubert polynomials
and Stanley functions by permutations instead of partial permutations
by canonically completing each partial permutation $w$ to a
permutation~$\wt w$; this does not alter the polynomials in any way
(Section~\ref{sec:schub}).
\end{example*}

\subsubsection*{Tableau formula (peelable tableaux and factor sequences)}

Our previous formulae write quiver polynomials in terms of binomials,
double Schubert polynomials (in two completely different ways), and
double Stanley symmetric functions.  Now we turn to Schur functions.
For a list $\bla = (\lambda_1,\ldots,\lambda_n)$ of partitions, write
\begin{eqnarray*}
  s_\bla(\xx-\oxx) &=& s_{\la_1}(\xx^0-\xx^1) \cdots
  s_{\la_n}(\xx^{n-1}-\xx^n)
\end{eqnarray*}
for the corresponding product of double Schur functions
(Section~\ref{sec:quiverpolys}).  Using our Theorem~\ref{t:BFtoKMS} to
identify our quiver polynomials with those in \cite{BF}, the Main
Theorem of Buch and Fulton states that there exist unique
integers~$c_\bla(\rr)$ \mbox{satisfying}
\begin{eqnarray*}
  Q_\rr(\xx-\oxx) &=& \sum_\bla c_\bla(\rr) s_\bla(\xx_\rr-\oxx_\rr),
\end{eqnarray*}
and exhibits an explicit way to generate them.  Although this
procedure involves negative integers, Buch and Fulton conjectured the
positivity of all the {\em quiver constants}\/~$c_\bla(\rr)$.  In
addition, starting with a rank array~$\rr$, they produced a concrete
recursive method for gener\-ating lists of semistandard Young tableaux
(Section~\ref{part:coeffs}), called {\em factor sequences}\/
(Section~\ref{sec:factor}).  Defining $\Phi(\rr)$ to be the set of
factor sequences coming from the rank array~$\rr$, every list $W \in
\Phi(\rr)$ of tableaux has an associated list~$\bla(W)$ of partitions.
The combinatorial conjec\-ture of \cite{BF} says that $c_\bla(\rr)$
counts the number of factor sequences $W \in \Phi(\rr)$ of
shape~$\bla$.

The direct connection between factor sequences and the combinatorial
geometry of quiver polynomials was from the beginning---and still
remains now---a mystery to us.  It seems that the question should come
down to finding an appropriate geometric explanation for
Schur-positivity of Stanley symmetric functions, which is known both
combinatorially and algebraically from various points of view.
Although we lack a geometric framework, one of these other points of
view, namely that of {\em Demazure characters}\/
(Section~\ref{sec:demazure}), still allows us to deduce a
combinatorial formula for the quiver constants~$c_\bla(\rr)$.

The argument is based on the fact that every Stanley function~$F_w$
expands as a sum $\sum_\la \alpha^\la_w s_\la$ of Schur
functions~$s_\la$ with nonnegative integer {\em Stanley
coefficients}~$\alpha^\la_w$ (Section~\ref{sec:stanley}).  The main
point, stated precisely in Theorem~\ref{t:qc=ss}, is that every quiver
constant is a Stanley coefficient for the Zelevinsky permutation:
$c_\bla(\rr) = \alpha^\la_{v(\rr)}$ for a certain partition~$\la$.

It is known that every Stanley coefficient counts a set of {\em
peelable tableaux} (Section~\ref{sec:peelable}).
Our case focuses on the set~$\Peel(D_\rr)$ of peelable tableaux for
the diagram~$D_\rr$
of the Zelevinsky permutation~$v(\rr)$.  A semistandard tableau~$P$
lies in~$\Peel(D_\rr)$ if it satisfies a readily checked, mildly
recursive condition.  The shape of every such tableau contains the
partition whose shape is the diagram~$D_\homv$ of~$v(\homv)$, so
removing~$D_\homv$ leaves a skew tableau $P - D_\homv$.  The connected
components of this skew tableau
form a list $\ptofs_\rr(P)$ of tableaux, read northeast to southwest
(Definition~\ref{d:ptofs}), whose list of shapes
we denote by~$\bla(P)$.

\begin{movedThm}{Theorem (Tableau formula; Buch--Fulton factor
sequence conjecture).}
\begin{eqnarray*}
  \QQ_\rr(\xx-\oxx)
  &=& \sum_{P \in\Peel(\rr)} s_{\bla(P)}(\xx_\rr-\oxx_\rr)
\\
  &=& \sum_{W \in \Phi(\rr)} s_{\bla(W)}(\xx_\rr-\oxx_\rr)
\end{eqnarray*}
\end{movedThm}

The two formulae in this theorem come from Theorem~\ref{t:PeelQuiver}
and Corollary~\ref{c:bla} (see Theorem~\ref{t:bla}, as well, which
says that our $c_\bla(\rr)$ agree with those in \cite{BF}).  In
contrast to the Schubert and Stanley versions of the component
formula, the peelable and factor sequence versions of the tableau
formula are equal summand by summand; this is immediate from
Theorem~\ref{t:ptofs}, which says simply that $P \mapsto
\ptofs_\rr(P)$ is a bijection from peelable tableaux to factor
sequences.  In fact, this is how our proof of the factor sequence
conjecture proceeds.  We shall mention more about the proof of the
tableau formula later in the Introduction.  For now, we reiterate that
we know of no direct geometric interpretation for it, in
either~version.

\begin{example*}
The set $\Peel(D_\rr)$ of peelable tableaux for
$D_\rr$ consists of four tableaux:
$$
\begin{array}{@{}c@{}}
\mbox{}\hspace{\textwidth}\mbox{}\\[-3ex]
\Yboxdim{12pt} \Yinterspace{1pt}
 \hfill
 \begin{array}{@{}c@{}}\\[-13.7ex]\lyoung{11117,24,3,4}\end{array}\hfill
 \begin{array}{@{}c@{}}\\[-13.7ex]\lyoung{1111,247,3,4}\end{array}\hfill
 \begin{array}{@{}c@{}}\\[-13.7ex]\lyoung{1111,24,37,4}\end{array}\hfill
 \lyoung{1111,24,3,4,7}
 \hfill
\end{array}
$$
Hence, under the bijection~$\ptofs_\rr$
that acts by removing entries whose locations correspond to $*$
entries in the diagram and reduced pipe dreams for~$v(\rr)$, there are
four factor sequences:
$$
\begin{array}{@{}c@{}}
\mbox{}\hspace{\textwidth}\mbox{}\\[-3ex]
\Yboxdim{12pt} \Yinterspace{1pt}
  \hfill\quad
  \big(\,\cyoung{7}\,,\cyoung{4}\,,\vn\big)\hfill
  \big(\vn\,,\cyoung{47}\,,\vn\big)\hfill
  \big(\vn\,,\cyoung{4,7}\,,\vn\big)\hfill
  \big(\vn\,,\cyoung{4}\,,\cyoung{7}\,\big)
  \hfill\ \
\end{array}
$$
The resulting tableau formula for $\QQ_\rr$ in terms of Schur
polynomials is
\begin{eqnarray*}
  Q_\rr
&=&\Yboxdim{3.5pt} \Yinterspace{1pt}
  s_{\yng(1)}(\aa-\bb)s_{\yng(1)}(\bb-\cc) +
  s_{\yng(2)}(\bb-\cc) +
  s_{\yng(1,1)}(\cc-\dd) +
  s_{\yng(1)}(\bb-\cc) s_{\yng(1)}(\cc-\dd).
\end{eqnarray*}
\end{example*}

\subsection*{Proofs via double and stable generalizations}

The formulae we presented above are all for {\em ordinary}\/ quiver
polynomials.  One of the main innovations in this paper is the idea of
working in the context of two
kinds of {\em double}\/ generalizations of quiver polynomials.  One is
obtained by applying limits to the other in the manner that Stanley
symmetric functions are obtained from Schubert polynomials.  Our need
for these limits, as well as our need for double versions at all,
will
surface later, when we sketch the proof of the Buch--Fulton factor
sequence conjecture.  In the meantime, let us begin by explaining the
double versions of the ratio and pipe formulae.

\subsubsection*{Double quiver polynomials: ratio and pipe formulae}

The definitions of $\SS_v(\xx-\oxx)$ for permutations $v \in S_d$, and
$s_\bla(\xx-\oxx)$ for partition lists~$\bla$, as well as
$\SS_\ww(\xx-\oxx)$~and $F_\ww(\xx-\oxx)$ for lacing diagrams~$\ww \in
\homv$ all required two sequences of $n+1$ alphabets.  Instead of
taking the second sequence~$\oxx$ to be the reverse of the first, use
a new list $\yy = \yy^n,\ldots,\yy^0$ for the second arguments.  All
of the definitions go through verbatim after replacing each
alphabet~$\oxx$ or~$\xx^j$ that has a minus sign in front of it by the
corresponding $\yy$ alphabet.  We define the {\em double quiver
polynomial}\/ $\QQ_\rr(\xx-\oyy)$ by making this replacement in the
ratio formula (Definition~\ref{d:QQrr}).

Given the way we define double quiver polynomials, the content of the
ratio formula for ordinary quiver polynomials is that specializing
$\yy = \xx$ in the ratio $\QQ_\rr(\xx-\oyy)$ of double Schubert
polynomials yields the multidegree of the quiver locus~$\Omega_\rr$
inside~$\homv$.  This is accomplished in Section~\ref{sec:ratio} using
the {\em Zelevinsky map}\/ from Section~\ref{sec:zel}.  The
equivariant class of the matrix Schubert variety~$\ol X_{v(\rr)}$ for
the Zelevinsky permutation equals the double Schubert
polynomial~$\SS_{v(\rr)}(\xx-\oyy)$ by \cite{grobGeom}.  The
Zelevinsky map takes~$\Omega_\rr$ isomorphically to the intersection
of the $\ol X_{v(\rr)}$ with the `opposite big cell'~$Y_0$ for partial
the flag manifold corresponding to the block decomposition of~$\endv$.
Restriction of equivariant cohomology classes on~$\endv$ to the affine
subspace~$Y_0$ of~$\endv$ results in the specialization $\yy = \xx$.
The division by~$\SS_{v(\homv)}$ comes from comparison of classes
on~$\homv$ with their images on~$Y_0$ pushed forward under the
Zelevinsky map.

The Zelevinsky map was introduced in a slightly different form
in~\cite{Zel85}.  The fact that it is a scheme-theoretic isomorphism
is equivalent to the main result in~\cite{LM}.  However, since we need
additionally the explicit connection to Zelevinsky permutations
(Theorem~\ref{t:zel}), we derive a new proof, which happens to be
quite elementary.  It does not even rely on primality of the
ideal~$I_\rr$ defining $\Omega_\rr$ as a subscheme of~$\homv$.

The double version of the pipe formula, Theorem~\ref{t:QQrr}, is
automatic from the formula for double Schubert polynomials
in~\cite{FK96}, as we explain in Section~\ref{sec:pipe2quiver}.

\subsubsection*{A combinatorial hierarchy}

Before explaining how stabilization enters the picture, it will be
helpful to have in mind a natural ordering of our four formulae,
different from the order in which we presented them above (which was
determined by our proofs).  The combinatorial objects indexing
summands in our four formulae can all be described in terms of the set
$\rp(v(\rr))$ of reduced pipe dreams for the Zelevinsky
permutation~$v(\rr)$.  More precisely, there are four partitions of
the set~$\rp(v(\rr))$ such that each type of combinatorial object
corresponds to an equivalence class in one of these partitions.  The
hierarchy is as follows.
\begin{numbered}
\item
Pipe formula: equivalence classes are singletons.

\item
Tableau formula: equivalence classes are {\em key classes}.  There is
a bijection from pipe dreams in $\rp(v(\rr))$ to `compatible pairs'
$(\aa,\ii)$ in which $\aa$ is a reduced word for the Zelevinsky
permutation and $\ii$ is a `compatible sequence' of positive integers
\cite{BB}.  The set of compatible pairs can be partitioned according
to the Coxeter--Knuth classes of the reduced words~$\aa$.  The
generating function of each equivalence class with respect to the
weight of the compatible sequences~$\ii$ is known as a `key
polynomial' or `Demazure character of type~$A$' \cite{RS95K}.  These
polynomials are crucial to our derivation of the tableau formula
(Section~\ref{part:coeffs}).

\item
Component formula:
equivalence classes are strip classes, indexed by lacing diagrams.
One of our favorite combinatorial results in this paper is the
combination of Theorem~\ref{t:lacing} and Corollary~\ref{l:W}.  These
say that out of every reduced pipe dream for~$v(\rr)$ comes a minimal
length lacing diagram with rank array~$\rr$.  The example we used
before is not quite complicated enough to get the feel, so the reader
should instead see the short Section~\ref{sec:pipe2lace} for the
simple construction.  Briefly, we divide pipe dreams into horizontal
strips (block rows) and say that two pipe dreams are equivalent if the
pipes ``do the same thing''---that is, accomplish the same partial
permutation---in each strip.  The resulting map from pipe dreams to
lacing diagrams happens to be fundamental to the logic, methods, and
overall perspective of the paper, as we shall see shortly.

\item
Ratio formula:
only one equivalence class.
\end{numbered}

Roughly speaking, summing monomials $(\xx-\oxx)^{D \minus D_\homv}$
over an equivalence class in one of these four partitions yields a
summand in one of the four formulae.  This is content-free for the
ratio and pipe formulae.  Moreover, we do not explicitly use this idea
for the tableau formula, although it is faintly evident in the proof
of Proposition~\ref{p:QQrr}.  That proof is also suggestive of the
essentially true statement that key equivalence refines lacing
equivalence, which is the fundamental (though not explicit) principle
behind our derivation of the peelable tableau formula from the
component formula.  For the component formula itself, summation over a
lacing equivalence classes is an essential operation, but only {\em
roughly}\/ yields a lacing diagram summand.  This is the first of many
situations where stable limits come to the rescue.

\subsubsection*{Limits and stabilization: component formula}

There are two primary aspects to our proof of the component formula:
we must show that
\begin{itemize}
\item
components in the degeneration of~$\Omega_\rr$ have multiplicity~$1$,
and are exactly the closures of orbits through minimal length lacing
diagrams with ranks~$\rr$; and

\item
the summands can be taken to be products of Stanley symmetric
functions rather than products of double Schubert polynomials.
\end{itemize}
As we mentioned during our presentation of the component formula, the
Schubert version for ordinary quiver polynomials follows immediately
from the first two statements above.  But keep in mind that we need to
prove a component formula for {\em double}\/ quiver polynomials, for
reasons to be revealed a little later.  The bad news
(Remark~\ref{rk:SSww}) is the subtle fact that replacing $\oxx$
with~$\oyy$ in the Schubert version {\em always}\/ yields a false
statement whenever $W(\rr)$ contains more than one lacing diagram!

The two aspects itemized above split the proof into two essentially
disjoint parts: a geometric argument producing a lower bound, and a
combinatorial argument producing an upper bound.  The conclusion comes
by noticing that the lower bound is visibly at least as high as the
upper bound, so that both must be equal.  Let us make this vague
description more precise, starting with the geometric lower bound.

The mere existence of the flat degeneration we introduce in
Section~\ref{sec:quivdegen} automatically implies that the ordinary
quiver polynomial $\QQ_\rr(\xx-\oxx)$ is a sum of double Schubert
products $\SS_\ww(\xx-\oxx)$, with positive integer coefficients
possibly greater than~$1$ (Corollary~\ref{c:quivSchub}).  This is
because the components in the degeneration are products of matrix
Schubert varieties, possibly nonreduced a~priori
(Theorem~\ref{t:degen}).  Moreover, it is easy to verify that the
orbit through each lacing diagram in~$W(\rr)$ appears with
multiplicity at least~$1$ (Proposition~\ref{p:lace}).

One of the fundamental observations in \cite{BF} gained by viewing
quiver polynomials as ordinary classes of degeneracy loci for vector
bundle maps on arbitary schemes is that nothing vital changes when the
ranks of all the bundles are increased by~$1$, or by any fixed
integer~$m$.  The key observation for us along these lines is that
this stability under uniformly increasing ranks can be proved directly
for the total families of quiver locus degenerations---and hence for
components in the special fiber, including their scheme
structures---by a direct geometric argument
(Proposition~\ref{p:stability}).  It does not require knowing anything
in advance about the special fibers of such degenerations.  As a
consequence, the previous paragraph, which a~priori only gives
existence of a Schubert polynomial component formula along with lower
bounds on its coefficients, in fact proves the {\em same}\/ existence
and lower bounds for an ordinary Stanley function component formula,
simply by taking limits.

The combinatorial argument producing an upper bound is based on the
partition of the pipe dreams in~$\rp(v(\rr))$ into strip classes, and
it automatically works in the `double' setting, with $\xx$ and~$\yy$
alphabets.  Summing $(\xx-\oyy)^{D \minus D_\homv}$ over reduced pipe
dreams~$D$ with strip equivalence class~$\ww$ yields a polynomial
that, while not actually equal to~$\SS_\ww(\xx-\oyy)$ because of the
bad news above, agrees with $\SS_\ww(\xx-\oyy)$ in all coefficients on
monomials all of whose variables appear ``not too late'' in the
alphabets $\xx^0,\ldots,\xx^n$ and $\yy^0,\ldots,\yy^n$ (this follows
from Corollary~\ref{c:bottom}).  This observation depends subtly on
the combinatorics of reduced pipe dreams and the symmetry of
$\QQ_\rr(\xx-\oyy)$ in each of its $2n+2$ alphabets
(Proposition~\ref{p:sym}).

Just as in the geometric argument, we have a stability statement;
here, it says that strip classes are stable under uniform increase in
rank (Corollary~\ref{c:W}).  Consequently, summing over all strip
classes and taking the limit as ranks uniformly become infinite yields
Theorem~\ref{t:stanley}, the ``existence and upper bound theorem'' for
{\em double quiver functions}~\mbox{$\FF_\rr(\xx-\oyy)$}, and the
double Stanley symmetric function component formula for them.  The
formula says that~$\FF_\rr(\xx-\oyy) = \lim_{m \to \infty}
\QQ_{m+\rr}(\xx-\oyy)$ is a multiplicity-free sum of double Stanley
symmetric function products $F_\ww(\xx-\oyy)$, where $m+\rr$ is the
rank array obtained from~$\rr$ by adding~$m$ to every entry.  That the
sum is over strip classes, whose lacing diagrams automatically lie
inside~$W(\rr)$ by Theorem~\ref{t:lacing}, implies the upper bound:
the coefficients on~$F_\ww(\xx-\oyy)$ is~$1$ if $\ww$ indexes a strip
class, and zero otherwise.

After specializing the strip class component formula for double quiver
functions to $\yy = \xx$, the lower and upper bounds are forced to
match.  Therefore, the preceding arguments culminate in the following
master component formula (Theorem~\ref{t:lace}), which is arguably the
single most important result in the paper.

\begin{movedThm}{Theorem (Stable double component formula for double
quiver functions).}
$$
\begin{array}{rcccl}
  \displaystyle\FF_\rr(\xx-\oyy)
&:=&
  \displaystyle\lim_{m \to \infty} \QQ_{m+\rr}(\xx-\oyy)
&=&
  \displaystyle\sum_{\ww \in W(\rr)} F_\ww(\xx-\oyy)
\end{array}
$$
\end{movedThm}
Part of this statement is the important combinatorial observation that
every
minimal lacing diagram indexes a strip class (Corollary~\ref{l:W}).
For another corollary, the somewhat mysteriously constructed
combinatorial upper bound magically implies the geometric
multiplicity~$1$ statement we were after: components in the
degenerated quiver locus are generically reduced
(Theorem~\ref{t:components}).  This finally yields the Schubert
version of the component formula for ordinary quiver polynomials; and
we of course also recover the ordinary Stanley version.

To convey a better feel for what limits have accomplished here, let us
summarize the fundamental principle that drives most of our logic.
There are a number of geometric, algebraic, and combinatorial objects
that we are unable to get a precise handle on directly.  These objects
include the components in quiver degenerations, as well as their
multiplicities; the sum of all monomials associated to pipe dreams in
a fixed strip class; and the set of lacing diagrams indexing strip
classes for a fixed Zelevinsky permutation.  The key realization is
that, even though we fail to identify these objects precisely, we can
prove their stability under uniformly increasing ranks, and in the
limit it becomes possible to identify
the desired objects exactly.  Hence our failures occurred only because
we fixed the ranks to begin with.

Let us note that our failures in fixed finite rank are not always
artifacts of our proof technique.  For example, the sum of monomials
for pipe dreams in a fixed strip class~$\ww$ really {\em isn't} equal to the
Stanley product~$F_\ww$ or the Schubert product~$\SS_\ww$; it {\em
only}\/ becomes so in the limit.  We take this opportunity to remark
that a stability statement holds also for peelable tableaux
(Corollary~\ref{c:PeelLim}); though it turned out to be unnecessary
for the proofs in this final version, it was essential in an earlier
draft.

\subsubsection*{Why double and stable versions: the factor sequence
conjecture}

The substitution $\oxx \goesto \oyy$ has no real geometric motivation
for us, even though we knew from the start that
\mbox{$\QQ_\rr(\xx-\oyy)$} can be interpreted directly in terms of
equivariant classes of Zelevinsky matrix Schubert varieties, and in
hindsight it is likely possible to make a functorial construction on
vector bundles that produces double quiver polynomials as degeneracy
locus formulae.  Our reason for seriously considering double quiver
polynomials is different: we need to be able to set the second set of
alphabets to zero.  Observe that when $\yy = \xx$, such a procedure
succeeds in killing almost all terms in the ordinary versions.

At a crucial point in our derivation of the peelable tableau formula
from a component formula (see Remark~\ref{rk:y=0}), we need certain
key polynomials
to equal the products of Schur polynomials that they a~priori only
approximate.  We check this by applying the theory of Demazure
characters via isobaric divided differences (Section~\ref{sec:sums}).
But these have only been developed for ordinary (not double) Schubert
polynomials, so we must work with ratios
$\SS_{v(\rr)}(\xx)/\SS_{v(\homv)}(\xx) = \QQ_\rr(\xx)$ of ordinary
Schubert polynomials obtained from double quiver polynomials by
setting~\mbox{$\yy = \0$}.

To be a little more precise, consider the two sides of the equation
\begin{eqnarray*}
  \frac{\SS_{v(\rr)}(\xx)}
  {\SS_{v(\homv)}(\xx)}
&=&
  \sum_{\ww \in W(\rr)} F_\ww(\xx_\rr),
\end{eqnarray*}
which we prove in Proposition~\ref{p:QQrr} from the stable double
component formula.  Both sides expand into sums of products
$s_\bla(\xx_\rr)$ of ordinary Schur polynomials, the left side by
Demazure characters, and the right side by using the fact that each
summand is a product of symmetric polynomials in the
sequence~$\xx_\rr$.  This formula explains why we need double
versions, and simultaneously provides yet another reason why we need
stable versions: Schur-expansion of the left hand side demands double
versions because it requires $\yy=\0$, while Schur-expansion of the
right hand side demands stability because it requires the summands to
be symmetric functions.  Note that if the products of Stanley
functions in the right hand sum are replaced by products of Schubert
polynomials, then the resulting statement is always false whenever
there are at least two summands (Remark~\ref{rk:SSww}).

We show in Theorem~\ref{t:bla} that the coefficients in the
Schur-expansion of the right hand side are the quiver constants from
\cite{BF}, while the above Demazure character argument implies that
the coefficients on the left hand side are Stanley coefficients
(Theorem~\ref{t:qc=ss}).
These count peelable tableaux, proving the peelable tableau formula.
The factor sequence conjecture follows from the peelable tableau
formula via the bijection~$\ptofs_\rr$.  The proof of this bijection
in Section~\ref{part:fs} is a complex but essentially elementary
computation.  It works by breaking the diagram~$D_\rr$ into pieces
small enough so that putting them back together in different ways
yields peelable tableaux on the one hand, and \mbox{factor sequences
on the~other}.

\subsection*{Related notions and extensions}

Quiver loci have been studied by a number of authors in their roles
as universal degeneracy loci for vector bundle morphisms, as
generalizations of Schubert varieties, and in representation theory of
quivers; see \cite{AD,ADK,LM,FP,Ful99,BF,FRthomPoly}, for a sample.
The particular rank conditions we consider here appeared first in
\cite{BF}, as the culmination of an increasingly general progression
of degeneracy loci.  Snapshots from this progression include the
following; for more on the history,
see \cite{FP}.
\begin{itemize}
\item
The case of two vector bundles and a fixed map required to have rank
bounded by a given integer is known as the Giambelli--Thom--Porteous
formula; see \cite[Chapter~14.4]{FulIT}.  (In this case, both the
ratio formula and the component formula reduce directly to the
Giambelli--Thom--Porteous formula: the denominator in the ratio
formula equals~$1$, while its numerator is the appropriate Schur
polynomial; and there is just one lacing diagram, with no crossing
laces.)

\item
Given two filtered vector bundles, there are degeneracy loci defined
by Schubert conditions bounding the ranks of induced morphisms from
subbundles in the filtration of the source to quotients by those in
the target.  The resulting universal formulae are double Schubert
polynomials~\cite{Ful92}.

\item
The same Schubert conditions can be placed on a quiver like those
considered here, with $2n$ bundles of ranks increasing from~$1$ up
to~$n$ and then decreasing from~$n$ down to~$1$.  The degeneracy locus
formula produces polynomials that Fulton called `universal Schubert
polynomials' \cite{Ful99}, because they specialize to quantum and
double Schubert polynomials.  We propose to call them {\em Fulton
polynomials}.  Results similar to the ones we prove here have been
obtained independently for the special case of Fulton polynomials
in~\cite{BKTY02}.
\end{itemize}

One can consider even more general quiver loci and degeneracy loci,
for arbitrary quivers and arbitrary rank conditions.  We expect that
the resulting quiver polynomials should have interesting combinatorial
descriptions, at least for finite type quivers.  Some indications of
this come from \cite{FRthomPoly}, where the same general idea of
reducing to the equivariant study of quiver loci also appears (but
the specifics differ in
what kinds of statements are made concerning quiver polynomials),
along with Rim{\'a}nyi's Thom polynomial proof of the component
formula \cite{BFR03}, which he produced
in response to seeing the formula.

For nilpotent cyclic quivers the orbit degeneration works just as for
the type $A_{n+1}$ quiver considered in this paper, and the Zelevinsky
map can be replaced by Lusztig's embedding of a nilpotent cyclic
quiver into a partial flag variety \cite{Lu90}.  The technique of
orbit degeneration might also extend to arbitrary finite type quivers,
but Zelevinsky maps do not extend in any straightforward way.
%
%
One of the advantages of orbit degeneration is that it not only
implies existence of a combinatorial formula, but provides strong
hints (Proposition~\ref{p:lace}) as to the format of such a formula.
Thus, in contrast to \cite[p.~668]{BF}, where much of the work lay in
``discovering the shape of the formula'', the geometry of degeneration
provides a blueprint~automatically.

Among the cohomological statements in this paper lurks a single result
in \mbox{\K-theory}: Theorem~\ref{t:ratio} gives a combinatorial
formula for quiver \K-classes as ratios of double Grothendieck
polynomials.  This \K-theoretic analogue of the ratio formula (which
is cohomological) immediately implies the \K-theoretic analogue of our
pipe formula by \cite{FK94} (see \cite[Section~5]{subword}).  There is
also a \K-theoretic generalization of our component formula (including
the stable double version),
proved independently in two papers subsequent to this one
\cite{BuchAltSign,altSign}.  It implies Buch's conjecture
\cite{Buch02} that the \K-theoretic analogues of the quiver
constants~$c_\bla(\rr)$
exhibit a certain sign alternation,
generalizing the positivity~of~$c_\bla(\rr)$.

It is an interesting problem to find a \K-theoretic analogue of the
factor sequence tableau formula.  Substantial work in this direction
has already been done by Buch \cite{Buch02}.  In analogy with our
proof of the cohomological Buch--Fulton factor sequence conjecture, a
complete combinatorial description of its \K-theoretic analogue will
involve a better understanding of how stable Grothendieck polynomials
indexed by permutations expand into stable Grothendieck polynomials
indexed by partitions.  Lascoux has given an algorithm for this
expansion based on a transition formula for Grothendieck polynomials
\cite{Las01}.

\subsection*{Organization}

The logical structure of the paper---that is, which formulae imply
others and how---is summarized by the following diagram:
$$
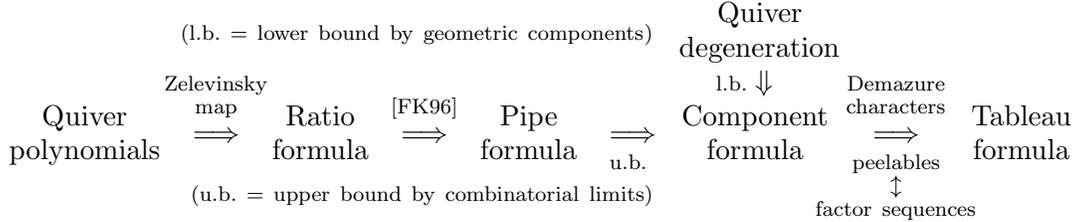
\begin{array}{@{}c@{\ \ }c@{\ \ }c@{\ \ }c@{\ \ }c@{\ \ }c@{\ \ }
                 c@{\ \ }c@{\ \ }c@{}}
\\[-3ex]
\multicolumn{6}{r}{\scriptstyle%
        ({\rm l.b.\ =\ lower\ bound\ by\ geometric\ components})}
 &
  \begin{array}{c}\hbox{Quiver}\\\hbox{degeneration}\end{array}
&&
\\
&\raisebox{-1ex}[0pt][0pt]{\makebox[0pt][c]{$\twoline{\rm Zelevinsky}
                           {\rm map}$}}&
&\raisebox{-1ex}[0pt][0pt]{\makebox[0pt][c]{$\twoline{}
                           {\rm [FK96]}$}}&
&&
  \llap{$\scriptstyle {\rm l.b.}$}\Downarrow
&\raisebox{-1ex}[0pt][0pt]{\makebox[0pt][c]{$\twoline{\rm Demazure}
                           {\rm characters}$}}&
\\
  \begin{array}{@{}c}\hbox{Quiver}\\\hbox{polynomials}\end{array}
&\Longrightarrow&
  \begin{array}{c}\hbox{Ratio}\\\hbox{formula}\end{array}
&\Longrightarrow&
  \begin{array}{c}\hbox{Pipe}\\\hbox{formula}\end{array}
&\Longrightarrow\!\!&
  \begin{array}{c}\hbox{Component}\\\hbox{formula}\end{array}\
&\Longrightarrow&
 \ \ \begin{array}{c@{}}\hbox{Tableau}\\\hbox{formula}\end{array}
\\
&&
&&
&\raisebox{2ex}[0pt][0pt]{$\scriptstyle {\rm u.b.}$}&
&\raisebox{-.5ex}[0pt][0pt]{\makebox[0pt][c]{$\twoline%
                {\twoline{\rm peelables}
                         {\updownarrow}}
                {\rm factor\ sequences}$}}&
\\[-2.2ex]
\multicolumn{6}{r}{\scriptstyle%
        ({\rm u.b.\ =\ upper\ bound\ by\ combinatorial\ limits})}
&&
\\[1ex]
\end{array}
$$

Although to simplify the exposition we use language at times as if the
field~$\kk$ were algebraically closed, all of our results hold for an
arbitrary field~$\kk$.  In fact the schemes we consider are all
defined over the integers~$\ZZ$, and our geometric statements
concerning flat (Gr\"obner) degenerations work in that context.


\part{Geometry of quiver loci}

\label{part:quiver}
\setcounter{section}{0}
\setcounter{thm}{0}
\setcounter{equation}{0}

\section{Quiver loci and ideals}

\label{sec:quiver}

Fix a sequence $V = (V_0,V_1,\ldots,V_n)$ of vector spaces with
$\dim(V_i)=r_i$, and denote by $\homv$ the variety $\Hom(V_0,V_1)
\times \cdots \times \Hom(V_{n-1},V_n)$ of type~$A_{n+1}$ \bem{quiver
representations} on~$V$.  That is, $\homv$ equals the vector space of
sequences
\begin{eqnarray*}
  \phi: && V_0\ \stackrel{\phi_1}\too\ V_1\ \stackrel{\phi_2}\too\
  \cdots\ \stackrel{\phi_{n-1}}\too\ V_{n-1}\ \stackrel{\phi_n}\too\ V_n
\end{eqnarray*}
of linear transformations.  By convention, set $V_{-1} = 0 = V_{n+1}$,
and $\phi_0 = 0 = \phi_{n+1}$.  Once and for all fix a basis%
    \footnote{This is not so bad, since we shall in due course require
    a maximal torus in $\gl(V_0) \times \cdots \times \gl(V_n)$.}
for each vector space~$V_i$, and express elements of~$V_i$ as row
vectors of length~$r_i$.  Doing so identifies each map $\phi_i$ in a
quiver representation~$\phi$ with a matrix of size $r_{i-1} \times r_i$.
Thus the coordinate ring~$\kk[\homv]$ becomes a polynomial ring in
variables $(\fff 1\alpha\beta), \ldots, (\fff n\alpha\beta)$, where the
$i^\th$ index~$\beta$ and the $(i+1)^\st$ index~$\alpha$ run from
$1$~to~$r_i$.  Let $\Phi$ be the \bem{generic} quiver representation, in
which the entries in the matrices $\Phi_i: V_{i-1} \to V_i$ are the
variables $\fff i\alpha\beta$.

\begin{defn} \label{d:quivlocus}
Given an array $\rr=(r_{ij})_{0 \leq i \leq j \leq n}$ of nonnegative
integers, the \bem{quiver locus} $\Omega_\rr$ is the zero scheme of the
\bem{quiver ideal} $I_\rr \subset \kk[\homv]$ generated by the union
over all $i < j$ of the minors of size $(1+r_{ij})$ in the product\/
$\Phi_{i+1} \cdots \Phi_j$ of matrices:
\begin{eqnarray*}
  I_\rr &=& \<\text{minors of size } (1+r_{ij}) \text{ in } \Phi_{i+1}
  \cdots \Phi_j \hbox{ for } i < j\>.
\end{eqnarray*}
\end{defn}

Thus the quiver locus $\Omega_\rr$ gives a natural scheme structure to
the set of quiver representations whose composite maps $V_i \to V_j$
have rank at most~$r_{ij}$ for all $i < j$.

Two quiver representations $(V,\phi)$ and $(W,\psi)$ on sequences $V$
and~$W$ of $n$ vector spaces are isomorphic if there are linear
isomorphisms $\eta_i : V_i \to W_i$ for $i = 0,\ldots,n$ commuting with
$\phi$ and~$\psi$.  Also, we can take the direct sum $V \oplus W$ in the
obvious manner.  Quiver representations that cannot be expressed
nontrivially as direct sums are called \bem{indecomposable}.  Whenever
$0 \leq p \leq q \leq n$, there is an (obviously) indecomposable
representation
\begin{eqnarray*}
  I_{p,q}: && 0 \to \cdots \to 0 \to
  \begin{array}[t]{@{}c@{}}\kk\\[-1.2ex]\scriptstyle p\end{array} =
  \cdots = \begin{array}[t]{@{}c@{}}\kk\\[-1.2ex]\scriptstyle
  q\end{array} \to 0 \to \cdots \to 0
\end{eqnarray*}
having copies of the field~$\kk$ in spots between $p$ and~$q$, with
identity maps between them and zeros elsewhere.  Here is a standard
result, c.f.\ \cite[Section~1.1]{LM}, generalizing the rank--nullity
theorem from linear algebra.

\begin{prop} \label{p:indec}
Every indecomposable quiver representation is isomorphic to some
$I_{p,q}$, and every quiver representation $\phi \in \homv$ is
isomorphic to a direct sum of these.
\end{prop}%
\begin{excise}{%
  \begin{proof}
  Given~$\phi$, we may as well assume $V_0 \neq 0$, and let $j$ be the
  largest index for which the composite $V_0 \to V_j$ is nonzero.
  Choose $y \in V_0$ mapping to a nonzero element of~$V_j$.  Set $V_i' =
  V_i/y$ for $i \leq j$, set $V_i' = V_i$ for $i > j$, and let $\phi'$
  be the induced quiver representation on~$V' = V_0',\ldots,V_n'$.  By
  induction on~$r_0 + \cdots + r_n$, the desired result holds
  for~$\phi'$.  Taking the direct sum of~$\phi'$ with the quiver
  $I_{0,q}$ determined by restricting~$\phi$ to the subspaces spanned by
  the images of~$y$ in $V_0,\cdots,V_j$, we conclude that the desired
  result holds for~$\phi$, as well.%
  \end{proof}
}\end{excise}%

The space $\homv$ of quiver representations carries an action of the
group
\begin{eqnarray} \label{glvv}
  \glvv &=& \gl(V_0)^2 \times \gl(V_1)^2 \times \cdots \times
  \gl(V_{n-1})^2 \times \gl(V_n)^2
\end{eqnarray}
of linear transformations.  Specifically, if we think of elements in each $V_i$
as row vectors and writing $\gl(V_i)^2 = \lvec{\gl}\,(V_i) \times
\rvec{\gl}\,(V_i)$ for $i = 0,\ldots,n$, the factor $\lvec\gl\,(V_i)$
acts by inverse multiplication on the right of $\Hom(V_{i-1},V_i)$,
while the factor $\rvec\gl\,(V_i)$ acts by multiplication on the left of
$\Hom(V_i,V_{i+1})$.  The group
\begin{eqnarray*}
  \glv &=& \gl(V_0) \times \cdots \times \gl(V_n)
\end{eqnarray*}
embeds diagonally inside $\glvv$, so that $\gamma =
(\gamma_0,\ldots,\gamma_n) \in \glv$ acts by
\begin{eqnarray*}
  \gamma \cdot \phi &=& (\,\ldots,\gamma_{i-1} \phi_i
  \gamma_i^{-1},\gamma_i \phi_i \gamma_{i+1}^{-1},\ldots\,).
\end{eqnarray*}
This action of $\glv$ preserves ranks, in the sense that $\glv \cdot
\Omega_\rr = \Omega_\rr$ for any ranks~$\rr$, because $\gamma_i^{-1}$
cancels with $\gamma_i$ when the maps in $\gamma \cdot \phi$ are
composed.

\begin{lemma} \label{l:orbits}
The group $\glv$ has finitely many orbits on~$\homv$, and every quiver
locus $\Omega_\rr$ is supported on a union of closures of such orbits.
\end{lemma}
\begin{proof}
Since $\glv \cdot \Omega_\rr = \Omega_\rr$, quiver loci are unions of
orbit closures for~$\glv$, so it is enough to prove there are finitely
many orbits.  Proposition~\ref{p:indec} implies that given any quiver
representation $\phi \in \homv$, we can choose new bases for
$V_0,\ldots,V_n$ in which $\phi$ is expressed as a list of partial
permutation matrices (all zero entries, except for at most one~$1$ in
each row and column).  Hence $\phi$ lies in the $\glv$ orbit of one of
the finitely many lists of partial permutation matrices in~$\homv$.%
\end{proof}

Definition~\ref{d:quivlocus} makes no assumptions about the array~$\rr$
of nonnegative integers, but only certain arrays can actually occur as
ranks of quivers.  More precisely, associated to each quiver
representation $\phi\in\homv$ is its \bem{rank array}~$\rr(\phi)$, whose
nonnegative integer entries are given by the ranks of composite maps
$V_i\rightarrow V_j$:
\begin{eqnarray*}
  r_{ij}(\phi) &=& \rank(\phi_{i+1}\phi_{i+2}\dotsm\phi_j)
  \quad\text{for}\quad i < j,
\end{eqnarray*}
and $r_{ii}(\phi) = r_i$ for $i = 0,\ldots,n$.  It is a consequence of
the discussion above that the support of a quiver locus $\Omega_\rr$
is irreducible if and only if $\rr = \rr(\phi)$ occurs as the rank
array of some quiver representation~$\phi \in \homv$, and this happens
if and only if $\Omega_\rr$ equals the closure of the orbit of $\glv$
through~$\phi$.
\begin{conv}
We make the convention here, once and for all, that in this paper we
consider only rank arrays~$\rr$ that can occur, so $\rr = \rr(\phi)$
for some $\phi \in \Omega_\rr$.
\end{conv}
That being said, we remark that this assumption is unnecessary in the
part of Section~\ref{sec:quiverpolys} preceding
Theorem~\ref{t:BFtoKMS}, as well as in Section~\ref{part:degen}.

\section{Rank, rectangle, and lace arrays}

\label{sec:zelperm}

Given a rank array $\rr = \rr(\phi)$ that can occur,
Proposition~\ref{p:indec} implies that there exists a unique array of
nonnegative integers $\ss = (s_{ij})_{0 \leq i \leq j \leq n}$ such
that $\phi = \bigoplus_{p \leq q} I_{p,q}^{\oplus s_{pq}}$.  In other
words, the ranks~$\rr$ can occur if and only if there is an associated
array~$\ss$.  We call this data the \bem{lace array}; the nomenclature
is explained by Lemma~\ref{l:ranklace}, and the symbol `$\ss$' recalls
that our construction is based on the `strands' appearing
in~\cite{AD}.  Taking ranks, we find that
\begin{eqnarray} \label{eq:laces}
  r_{ij} &=& \sum_\twoline{k \leq i}{j \leq m} s_{km}
  \quad\text{for}\quad i \leq j.
\end{eqnarray}
The (transposed) Buch--Fulton \bem{rectangle array} is the array
of rectangles $(R_{ij})$ for $0 \leq i < j \leq n$, such that
$R_{ij}$ is the rectangle of height $r_{i,j-1}-r_{ij}$ and width
$r_{i+1,j}-r_{ij}$.

\begin{example} \label{ex:rank}
Consider the rank array $\rr=(r_{ij})$, its lace array $\ss=(s_{ij})$,
and its rectangle array $\mathbf{R}=(R_{ij})$, which we depict as
follows.
\begin{equation*}
\Yboxdim{5pt} \Yinterspace{1pt}
\rr \ \:=\:\ \begin{array}{cccc|c}
3 & 2 & 1 & 0 & i \diagup j \\ \hline
  &   &   & 2 & 0 \\
  &   & 3 & 2 & 1 \\
  & 4 & 2 & 1 & 2 \\
3 & 2 & 1 & 0 & 3
\end{array} \qquad
\ss \ \:=\:\ \begin{array}{cccc|c}
3 & 2 & 1 & 0 & i \diagup j \\ \hline
  &   &   & 0 & 0 \\
  &   & 0 & 1 & 1 \\
  & 1 & 0 & 1 & 2 \\
1 & 1 & 1 & 0 & 3
\end{array} \qquad
\mathbf{R} \ \:=\:\
\begin{array}{ccc|c}
     2     &     1   &    0    & i \diagup j \\ \hline
           &         & &     1       \\
           & \yng(2) & \yng(1) &     2       \\
 \yng(1,1) & \yng(1) & \yng(1) &     3
\end{array}
\end{equation*}
The relation \eqref{eq:laces} says that an entry of~$\rr$ is the
sum of the entries in~$\ss$ that are weakly southeast of the
corresponding location.  The height of $R_{ij}$ is obtained by
subtracting the entry $r_{ij}$ from the one above it, while the
width of $R_{ij}$ is obtained by subtracting the entry $r_{ij}$
from the one to its left.\qed
\end{example}

It follows from the definition of $R_{ij}$ that
\begin{alignat}{2}
\label{eq:heightsums} \sum_{k\ge j} \hgt(R_{ik}) &= r_{i,j-1}-r_{i,n}
  \le r_{i,j-1} & \qquad &\text{for all $i$} \\
\label{eq:widthsums} \sum_{\ell\le i} \wid(R_{\ell j})
  &=r_{i+1,j}-r_{0,j} \le r_{i+1,j} & &\text{for all $j$.}
\end{alignat}
(This will be applied in Proposition~\ref{pp:ZPeelShape}.)  The
relation~\eqref{eq:laces} can be inverted to obtain
\begin{eqnarray} \label{eq:sr}
  s_{ij} &=& r_{ij}-r_{i-1,j}-r_{i,j+1}+r_{i-1,j+1}
\end{eqnarray}
for $i \leq j$, where $r_{ij} = 0$ if $i$ and~$j$ do not both lie
between $0$ and~$n$.  Therefore the rank array $\rr$ can occur for~$V$
if and only if $r_{ii}=r_i$ for $i = 0,\ldots,n$ and
\begin{eqnarray} \label{eq:rankineq}
  r_{ij}-r_{i-1,j}-r_{i,j+1}+r_{i-1,j+1} &\ge& 0
\end{eqnarray}
for $i \leq j$, the left hand side being simply~$s_{ij}$.  We shall
interchangeably use a lace array~$\ss$ or its corresponding rank
array~$\rr$ to describe a given irreducible quiver locus.

It follows from \eqref{eq:rankineq} that
\begin{align}
\label{eq:heightincr} &\text{$\hgt(R_{ij})$ decreases as $i$ decreases.}
\\
\label{eq:widthincr} &\text{$\wid(R_{ij})$ decreases as $j$ increases.}
\end{align}

Given a rank array~$\rr$ or equivalently a lace array~$\ss$, we shall
construct a permutation $v(\rr)\in S_d$, where $d = r_0 + \cdots +
r_n$. In~general, any matrix in the space~$\endv$ of \mbox{$d \times
d$} matrices comes with a decomposition into block rows of heights
$r_0,\ldots,r_n$ (from top to bottom) and block columns of lengths
$r_n,\ldots,r_0$ (from left to right).  Note that our indexing
convention may be unexpected, with the square blocks lying along the
main block {\em anti\/}diagonal rather than on the diagonal as
usual. With these conventions, the $i^\th$~block column refers to the
block column of length~$r_i$, which sits $i$~blocks from the {\em
right}.

We draw the matrix for the permutation $w$ with a symbol $\times$
(instead of a $1$) at each position $(q,w(q))$ and zeros elsewhere.

\begin{prop} \label{p:zel}
Given a rank array $\rr$ for $V$, there exists a unique element
$v(\Omega_\rr)=v(\rr)\in S_d$ satisfying the following conditions.
Consider the block in the $i^\th$ block {\em column} and $j^\th$
block~{\em row}.
\begin{numbered}
\item
If $i \leq j$ (that is, the block sits on or below the main block
antidiagonal) then the number of~$\times$ entries in that block
equals~$s_{ij}$.

\item
If $i = j+1$ (that is, the block sits on the main block
superantidiagonal) then the number of~$\times$ entries in that block
equals~$r_{j,j+1}$.

\item
If $i \geq j+2$ (that is, the block lies strictly above the main block
superantidiagonal) then there are no~$\times$ entries in that block.

\item
Within every block row or block column, the $\times$ entries proceed
from northwest to southeast, that is, no~$\times$~entry is northeast
of another~$\times$~entry.
\end{numbered}
\end{prop}
\begin{defn} \label{d:zel}
$v(\rr)$ is the \bem{Zelevinsky permutation} for the rank array~$\rr$.
\end{defn}
\begin{proof}
We must show that the number of~$\times$~entries in any block row, as
dictated by conditions~\mbox{1--3}, equals the height of that block row
(and transposed for columns), since condition~4 then stipulates uniquely
how to arrange the~$\times$~entries within each block.  In other words
the height $r_{\!j}$ of the $j^\th$ block row must equal the number
$r_{\!j,j+1}$ of $\times$ entries in the superantidiagonal block in that
block row, plus the sum $\sum_i s_{ij}$ of the number of $\times$
entries in the rest of the blocks in that block row (and a similar
statement must hold for block columns).  These statements follow
from~(\ref{eq:laces}).%
\end{proof}

The \bem{diagram} $D(v)$ of a permutation $v\in S_d$, which by
definition consists of those cells in the \mbox{$d \times d$} grid
having no~$\times$ in~$v$ due north or due west of~it, refines the
data contained in the Buch--Fulton rectangle array \cite{BF}.  The
following Lemma is a straightforward consequence of the definition of
the Zelevinsky permutation $v(\rr)$ and equation~\eqref{eq:sr}; to
simplify notation, we often write $D_\rr = D(v(\rr))$.

\begin{lemma} \label{l:diagram}
In any block of the diagram $D_\rr = D(v(\rr))$, the cells form a
rectangle justified in the southeast corner of the block. If the block
is above the superantidiagonal, this rectangle consists of all the
cells in the block. If the block is on or below the superantidiagonal,
say in the $i^\th$~block column and $j^\th$~block row, this rectangle
is the Buch--Fulton rectangle $R_{i-1,j+1}$.
\end{lemma}

\begin{example} \label{ex:zel}
Let $\rr,\ss,\mathbf{R}$ be as in Example \ref{ex:rank}. The Zelevinsky
permutation is
\setcounter{MaxMatrixCols}{20}
\begin{eqnarray*}
v(\rr) &=&  \begin{pmatrix}
        1&2&3&4&5&6&7&8&9&10&11&12\\
        8&9&4&5&11&1&2&6&12&3&7&10
        \end{pmatrix},
\end{eqnarray*}
whose permutation matrix is indicated by $\times$ entries in
Fig.~\ref{f:zel},
\begin{figure}
$$
\begin{array}{@{}l|%
        @{}  c@{\:}@{\:}c@{\:}@{\:}c@{\:}|%
        @{\:}c@{\:}@{\:}c@{\:}@{\:}c@{\:}@{\:}c@{\:}|%
        @{\:}c@{\:}@{\:}c@{\:}@{\:}c@{\:}|%
        @{\:}c@{\:}@{\:}c@{\:}|@{}}
\cline{2-13}
 8 &\: * & * & * & * & * & * & * &\ti&\cd&\cd&\cd&\cd
\\
 9 &\: * & * & * & * & * & * & * &\cd&\ti&\cd&\cd&\cd
\\\cline{2-13}
 4 &\: * & * & * &\ti&\cd&\cd&\cd&\cd&\cd&\cd&\cd&\cd
\\
 5 &\: * & * & * &\cd&\ti&\cd&\cd&\cd&\cd&\cd&\cd&\cd
\\
 11&\: * & * & * &\cd&\cd&\sq&\sq&\cd&\cd&\sq&\ti&\cd
\\\cline{2-13}
 1 &\:\ti&\cd&\cd&\cd&\cd&\cd&\cd&\cd&\cd&\cd&\cd&\cd
\\
 2 &\:\cd&\ti&\cd&\cd&\cd&\cd&\cd&\cd&\cd&\cd&\cd&\cd
\\
 6 &\:\cd&\cd&\sq&\cd&\cd&\ti&\cd&\cd&\cd&\cd&\cd&\cd
\\
 12&\:\cd&\cd&\sq&\cd&\cd&\cd&\sq&\cd&\cd&\sq&\cd&\ti
\\\cline{2-13}
 3 &\:\cd&\cd&\ti&\cd&\cd&\cd&\cd&\cd&\cd&\cd&\cd&\cd
\\
 7 &\:\cd&\cd&\cd&\cd&\cd&\cd&\ti&\cd&\cd&\cd&\cd&\cd
\\
 10&\:\cd&\cd&\cd&\cd&\cd&\cd&\cd&\cd&\cd&\ti&\cd&\cd
\\\cline{2-13}
\end{array}
$$
\caption{Zelevinsky permutation and its diagram}
\label{f:zel}
\end{figure}
and whose diagram $D_\rr = D(v(\rr))$ is indicated by the union of
the~$*$ and~$\square$ entries there.\qed%
\end{example}

\begin{defn} \label{d:diagram}
Consider a fixed dimension vector $(r_0,r_1,\dotsc,r_n)$ and a varying
rank array $\rr$ with that dimension vector.
\begin{numbered}
\item
The diagram $D_\homv := D(v(\homv))$ for the quiver locus $\homv$
consists of all cells above the block superantidiagonal.  Thus $D_\rr
\supset D_\homv$ for all occurring rank arrays~$\rr$.  Define the
\bem{(skew) diagram} $D_\rr^*$ for $\rr$ as the difference
\begin{eqnarray*}
  D_\rr^* &=& D_\rr\setminus D_\homv.
\end{eqnarray*}
In our depictions of~$D_\rr$, each cell of~$D_\rr^*$ is indicated by a
box~$\sq$, while each cell of~$D_\homv$ is indicated by an
asterisk~$*$.
\item
For the zero quiver~$\Omega_0$, whose rank array $\rr_0$ is zero
except for the $r_{ii}$ entries, the diagram $D(\Omega_0)$ consists of
all locations strictly above the antidiagonal.
\end{numbered}
\end{defn}

\section{The Zelevinsky map}

\label{sec:zel}

Let $\gld$ be the invertible $d \times d$ matrices.  Denote
by~$P\subset\gld$ the parabolic subgroup of block {\em lower}\/
triangular matrices, where the diagonal blocks have sizes
$r_0,\ldots,r_n$, and let $B_+$ be the group of {\em upper\/}
triangular matrices.

The group~$P$ acts by multiplication on the left of~$\gld$, and the
quotient~$\pgld$ is the manifold of partial flags with dimension jumps
$r_0, \ldots, r_n$.  Schubert varieties $X_v$ in~$\pgld$ are orbit
closures for the (left) action of~$B_+$ by inverse right
multiplication.  In~writing $X_v \subseteq \pgld$ we~shall always
assume that the permutation~$v \in S_d$ has minimal length in its
right coset $(S_{r_0} \times \cdots \times S_{r_n}) v$, which means
that $v$~has no descents within a block.  Graphically, this condition
is reflected in the permutation matrix for~$v$ by saying that
no~$\times$~entry is northeast of another in the same block row.
Observe that this condition holds by definition for Zelevinsky
permutations~$v(\rr)$.

The preimage in~$\gld$ of a Schubert variety $X_v \subseteq \pgld$ is
the closure in~$\gld$ of the double coset $P v B_+$.  The closure
$\ol{P v B_+}$ of that inside~$\endv$ is the `matrix Schubert
variety'~$\Xb_v$, whose definition we now recall.  Since we shall need
matrix Schubert varieties in the slightly more general setting of
partial permutations later on, we use language here compatible with
that level of generality.

\begin{defn}[\cite{Ful92}] \label{d:matSchub}
Let $w$ be a \bem{partial permutation} of size $k \times \ell$, meaning
that $w$ is a matrix with $k$~rows and $\ell$~columns having all entries
equal to~$0$ except for at most one entry equal to~$1$ in each row and
column.  Let $M_{k\ell}$ be the vector space of $k \times \ell$
matrices.  The \bem{matrix Schubert variety} $\Xb_w$ is the subvariety
\begin{eqnarray*}
   \Xb_w &=& \{Z \in M_{k\ell} \mid \rank(Z\sub qp) \leq
   \rank(w\sub qp) \text{ for all } q \text{ and } p\}
\end{eqnarray*}
inside $M_{k\ell}$, where $Z\sub qp$ consists of the top $q$~rows and
left $p$~columns of~$Z$.
\end{defn}

If $v \in S_d$ is a permutation, we do not require~$v$ to be a minimal
element in its coset of \mbox{$S_{r_0} \times \cdots \times S_{r_n}$}
when we write~$\Xb_v$ (as opposed to~$X_v$); however, the matrix
Schubert variety will fail to be a $P \times B_+$ orbit closure
in~$\endv$ unless $v$ is minimal.

Let $Y_0$ be the variety of all matrices in~$\gld$ whose antidiagonal
blocks are all identity matrices, and whose other blocks below the
antidiagonal are zero.  Thus $Y_0$ is obtained from the unipotent
radical $U(P_+)$ of the block upper-triangular subgroup~$P_+$ by a
global left-to-right reflection followed by block left-to-right
reflection.  This has the net effect of reversing the order of the
block columns of~$U(P_+)$.  If~${w_0}$ is the \bem{long element}
in~$S_d$, represented as the antidiagonal permutation matrix, and
$\ww_0$ is the \bem{block long element}, with antidiagonal permutation
matrices in each diagonal block, then $Y_0 = U(P_+){w_0}\ww_0$.  The
variety~$Y_0$ maps isomorphically (scheme-theoretically) to the
\bem{opposite big cell} $U_0$ in~$\pgld$ under projection modulo~$P$.
In other words, $U_0 \to Y_0$ is a section of the projection $\gld \to
\pgld$ over the open set~$U_0$.

Using the isomorphism $U_0 \cong Y_0$, the intersection $X_v \cap U_0$
of a Schubert variety in~$\pgld$ with the opposite big cell is a closed
subvariety~$Y_v$ of~$Y_0$ (we assume $v$ is minimal in its coset when we
write~$Y_v$).  It will be important for Theorem~\ref{t:zel} to make the
(standard) comparison between the equations defining~$Y_v$ inside~$Y_0$
and the equations defining the corresponding matrix Schubert
variety~$\ol X_v$ inside~$\endv$.

\begin{lemma} \label{l:zero}
The ideal of\/ $Y_v$ as a subscheme of\/~$Y_0$ is obtained from the
ideal of\/~$\ol X_v$ as a subscheme of~$\endv$ by setting
all variables
on the diagonal of each antidiagonal block equal to\/~$1$ and all other
variables in or below the block antidiagonal to zero.
\end{lemma}
\begin{proof}
The statement is equivalent to its geometric version, which says that
$Y_v$ equals the scheme-theoretic intersection of $\ol X_v$ with~$Y_0$.
To prove this version, note that $\ol X_v \cap Y_0 = (\ol X_v \cap \gld)
\cap Y_0$, and that $\ol X_v \cap \gld$ projects to~$X_v$ as a fiber
bundle over $X_v$ with fiber~$P$.  The intersection of~$\ol X_v \cap
\gld$ with the image~$Y_0$ of the section $U_0 \to \gld$ is
scheme-theoretically the image~$Y_v$ of the section restricted
to~\mbox{$X_v \cap U_0$}.%
\end{proof}

\begin{defn}
The \bem{Zelevinsky map}~$\cZ$ takes $\phi \in \homv$ to the block
matrix
\begin{eqnarray} \label{zel}
  (\phi_1,\phi_2,\ldots,\phi_n)
  &\stackrel {\textstyle \cZ}{\displaystyle\longmapsto}&
  \left[
  \begin{array}{ccccc}
       0   &        &    0   & \phi_1 & \1 \\
       0   &        & \phi_2 &   \1   &  0 \\
       0   & \adots &   \1   &    0   &  0 \\
    \phi_n & \adots &    0   &    0   &  0 \\
      \1   &        &    0   &    0   &  0
  \end{array}
  \right].
\end{eqnarray}
Set $Y_\homv \subset \gld$ equal to the image in~$\gld$ of the quiver
space~$\homv$ under the Zelevinsky map~$\cZ$, and denote by $Y_\rr =
\cZ(\Omega_\rr)$ the image of the quiver locus.
\end{defn}

Zelevinsky's original map \cite{Zel85} sent each $\phi \in \homv$ to
the matrix $(\cZ(\phi){w_0}\ww_0)^{-1}$ essentially inverse
to~$\cZ(\phi)$, and Zelevinsky proved that it was bijective.  The
following theorem is therefore equivalent
to the main theorem in~\cite{LM}, which proved using standard monomial
theory that Zelevinsky's original map is an isomorphism of
schemes---that is, the quiver ideal is prime.
Using the simpler form~$\cZ$ of Zelevinsky's map clarifies the proof,
and connects quiver loci to the explicit combinatorics of matrix
Schubert varieties for Zelevinsky permutations, as we shall require
later.

\begin{thm} \label{t:zel}
The Zelevinsky map $\cZ$ induces a scheme isomorphism from each quiver
locus~$\Omega_\rr$ to the closed subvariety\/~$Y_{v(\rr)}$ inside the
Schubert subvariety~$X_{v(\rr)}$ of the partial flag manifold~$\pgld$.
In other words, $Y_\rr = Y_{v(\rr)}$ as subschemes of\/~$Y_0$.
\end{thm}
\begin{proof}
Both of the schemes $Y_\rr$ and~$Y_{v(\rr)}$ are contained in~$Y_\homv$:
for $Y_\rr$ this is by definition, and for $Y_{v(\rr)}$ this is because
the diagram of~$v(\rr)$ contains the union of all blocks strictly above
the block superantidiagonal, whence every coordinate above the block
superantidiagonal is zero on~$Y_{v(\rr)}$.  Since the Zelevinsky map
$\Omega_\rr \to Y_\rr$ is obviously an isomorphism of schemes, we must
show that $Y_\rr$ is defined by the same equations in~$\kk[Y_\homv]$
defining~$Y_{v(\rr)}$.

The set-theoretic description of~$\Xb_v$ in Definition~\ref{d:matSchub}
implies that the Schubert determinantal ideal $I(\Xb_v) \subseteq
\kk[\endv]$ contains the union (over $q,p = 1,\ldots,d$) of minors with
size \mbox{$1 + \rank(v\sub qp)$} in the northwest \mbox{$q \times p$}
submatrix of the \mbox{$d \times d$} matrix~$\Phi$ of variables.  In
fact, using Fulton's `essential set' \cite[Section~3]{Ful92}, the ideal
$I(\ol X_v)$ is generated by those minors arising from cells $(q,p)$ at
the southeast corner of some block, along with all the variables
strictly above the block superantidiagonal.

Consider a box $(q,p)$ at the southeast corner of $B_{i+1,j-1}$, the
intersection of block column~\mbox{$i+1$} and block row~\mbox{$j-1$}, so
that by definition of Zelevinsky permutation,
\begin{eqnarray} \label{rank}
\rank(v(\rr)\sub qp) &=&
\sum_\twoline{\alpha>i}{\beta<j}s_{\alpha\beta}+\sum_{k=i+1}^jr_{k-1,k}.
\end{eqnarray}

\begin{lemma} \label{l:rank}
The number\/ $\rank(v(\rr)\sub qp)$ in \eqref{rank} equals $r_{ij} +
\sum_{k=i+1}^{j-1} r_k$.
\end{lemma}
\begin{proof}
The coefficient on~$s_{\alpha\beta}$ in $r_{ij} + \sum_{k=i+1}^{j-1}
r_k$ is the number of elements in
$\{r_{ij}\}\cup\{r_{i+1,i+1},\ldots,r_{j-1,j-1}\}$ that are weakly
northwest of $r_{\alpha\beta}$ in the rank array~$\rr$ (when the
array~$\rr$ is oriented so that its southeast corner is~$r_{0n}$).
This number equals the number of elements
in~$\{r_{i,i+1},\ldots,r_{j-1,j}\}$ that are weakly northwest
of~$r_{\alpha\beta}$, unless $r_{\alpha\beta}$ happens to lie strictly
north and strictly west of~$r_{ij}$, in which case we get one fewer.
This one fewer is exactly made up by the sum of entries from~$\ss$
in~\eqref{rank}.%
\end{proof}

Resuming the proof of Theorem~\ref{t:zel}, we set the appropriate
variables in the generic matrix~$\Phi$ to $0$ or~$1$ by
Lemma~\ref{l:zero}, and consider the equations coming from the northwest
\mbox{$q \times p$} submatrix for~$(q,p)$ in the southeast corner
of~$B_{i+1,j-1}$.  Since $Y_{v(\rr)} \subseteq Y_\homv$, these
equations are minors of~\eqref{zel}.  In particular, using
Lemma~\ref{l:rank}, and assuming that $i,j \in \{0,\ldots,n\}$,
we~find that these equations in~$\kk[Y_\homv]$ are the minors of size
$1 + u + r_{ij}$ in the generic \mbox{$(u+r_i) \times (u+r_i)$} block
matrix
\begin{eqnarray} \label{rij}
  \left[
  \begin{array}{ccccc}
     0  &   0       &      &   0      &\phi_{i+1}\\
     0  &   0       &      &\phi_{i+2}&  \1      \\
     0  &   0       &\adots&  \1      &   0      \\
     0  &\phi_{j-1} &\adots&   0      &   0      \\
  \phi_j&  \1       &      &   0      &   0
  \end{array}
  \right],
\end{eqnarray}
where $u = \sum_{k=i+1}^{j-1} r_k$ is the sum of the ranks of the
subantidiagonal $\1$~blocks.  The ideal generated by these minors of
size \mbox{$1+u+r_{ij}$} is preserved under multiplication
of~\eqref{rij} by any~matrix in $\SL_{u+r_i}(\kk[Y_\homv])$.  In
particular, multiply~\eqref{rij} on~the~left~by
\begin{eqnarray*}
  \left[
  \begin{array}{c@{\quad\ }c@{\quad\ }c@{\quad\ }c@{\quad\ }c@{\quad\ }c}
    \ \ \1\quad
  & \quad\ -\phi_{i+1}\quad\
  & \phi_{i+1}\phi_{i+2}
  & \cdots
  & \pm\phi_{i+1,j-2}
  & \mp\phi_{i+1,j-1}
  \\
       & \1 &    &          &    &    \\
       &    & \1 &          &    &    \\[-1ex]
       &    &    &\flatddots&    &    \\[-.5ex]
       &    &    &          & \1 &    \\
       &    &    &          &    & \1
  \end{array}
  \right],
\end{eqnarray*}
where $\phi_{i+1,k} = \phi_{i+1} \cdots \phi_k$ for $i+1 \leq k$.  The
result agrees with~\eqref{rij} except in its first block row, which has
left block $(-1)^{j-1-i}\phi_{i+1} \cdots \phi_j$ and all other blocks
zero.  Therefore the minors coming from the southeast corner $(q,p)$ of
the block~$B_{i+1,j-1}$ generate the same ideal in~$\kk[Y_\homv]$ as the
size \mbox{$1+r_{ij}$} minors of $\phi_{i+1} \cdots \phi_j$, for
all~\mbox{$i \leq j$}.  We conclude that the ideals of~$Y_\rr$ and
$Y_{v(\rr)}$ inside $\kk[Y_\homv]$ coincide.%
\end{proof}

The proof of Theorem~\ref{t:zel} never really uses the fact that the
minors vanishing on matrix Schubert varieties generate a prime ideal.
Nonetheless, they do \cite{Ful92,grobGeom}, and in fact one can
conclude much more by citing (as in \cite{LM}) other properties of
Schubert varieties from \cite{Ram85,RR85}, and using
Lemma~\ref{l:diagram}.

\begin{cor}[\cite{LM,AD}] \label{c:codim}
Irreducible quiver loci $\Omega_\rr$ are reduced, normal, and
Cohen--Macaulay with rational singularities.  The total area
$$
\begin{array}{@{}rcccl@{}}
  \dis d(\rr) &=& \dis \ell(v(\rr)) - \ell(v(\homv)) &=& \dis \sum_{0
  \leq k \leq i < j \leq m \leq n} s_{k,j-1} s_{i+1,m}
\end{array}
$$
of the rectangles~$R_{ij}$ equals the codimension of\/~$\Omega_\rr$
in~$\homv$.
\end{cor}

\begin{remark}
One can prove even more simply that the Zelevinsky map~$\cZ$ induces
an isomorphism of the reduced variety underlying the
scheme~$\Omega_\rr$ with the opposite `cell' in {\em some}\/ Schubert
variety~$X_v$ inside~$\pgld$, without identifying $v = v(\rr)$.
Indeed, $\cZ$ is obviously an isomorphism (in fact, a linear map) from
$\homv$ to~$Y_\homv$, which is easily seen to equal $Y_{v(\homv)}$.
Then one uses equivariance of~$\cZ$
under appropriate actions of $\glv$ to conclude that it takes orbit
closures (the reduced varieties underlying quiver loci) to orbit closures
(opposite Schubert `cells').
\end{remark}

\section{Quiver polynomials as multidegrees}

\label{sec:quiverpolys}

In this section we define quiver polynomials, our title characters, as
multidegrees of quiver loci~$\Omega_\rr$, and prove that these agree
with the polynomials appearing in~\cite{BF}.  For the reader's
convenience, we present the definition of multidegree
\cite{joseph,rossmann} and introduce related notions here; those
wishing more background and explicit examples can consult
\cite[Sections~1.2 and~1.7]{grobGeom} or~\cite[Chapter~8]{CCA}.

Let $\kk[\ff]$ be a polynomial ring graded by~$\ZZ^d$.  The gradings
we consider here are all%
    \footnote{Actually, a nonpositive grading appears in the proof
    of Proposition~\ref{p:grothBigCell}, but we apply no theorems to
    it---the nonpositive grading there is only a transition between
    two positive gradings.}
\bem{positive} in the sense of \cite[Section~1.7]{grobGeom}, which implies
that the graded pieces of $\kk[\ff]$ have finite dimension as vector
spaces.  Writing $\xx = x_1,\ldots,x_d$ for a basis of~$\ZZ^d$, each
variable $f \in \kk[\ff]$ has an \bem{ordinary weight} $\deg(f) = \uu(f)
= \sum_{i=1}^d u_i(f)\cdot x_i \in \ZZ^d$, and the corresponding
\bem{exponential weight} $\weight(f) = \xx^{\uu(f)} = \prod_{i=1}^d
x_{\makebox[0ex][l]{$\scriptstyle i$}}{}^{u_i(f)}$, which is a Laurent
monomial in the ring $\ZZ[x_1^{\pm 1}, \ldots, x_d^{\pm 1}]$.  Every
finitely generated $\ZZ^d$-graded module $\Gamma = \bigoplus_{\uu \in
\ZZ^d} \Gamma_{\!\uu}$ over $\kk[\ff]$ has a finite $\ZZ^d$-graded free
resolution
$$
  \EE_\spot:\ 0 \from \EE_0 \from \EE_1 \from \EE_2 \from \cdots
  \qquad\hbox{where}\qquad
  \EE_i = \bigoplus_{j=1}^{\beta_i} \kk[\ff](-\tt_{ij})
$$
is $\ZZ^d$-graded, with the $j^\th$ summand of $\EE_i$ generated in
degree $\tt_{ij} \in \ZZ^d$.  In this case, the \bem{\K-polynomial} of
$\Gamma$ is
\begin{eqnarray*}
  \KK(\Gamma;\xx) &=& \sum_i (-1)^i \sum_j \xx^{\tt_{ij}}.
\end{eqnarray*}

Given any Laurent monomial $\xx^\uu = x_1^{u_1} \cdots x_d^{u_d}$, the
rational function $\prod_{j=1}^d(1-x_j)^{u_j}$ can be expanded as a
well-defined (that is, convergent in the $\xx$-adic topology) formal
power series \mbox{$\prod_{j=1}^d(1 - u_jx_j + \cdots)$} in~$\xx$.
Doing the same for each monomial in an arbitrary Laurent polynomial
$\KK(\xx)$ results in a power series denoted by $\KK(\1-\xx)$.  The
\bem{multidegree} of a $\ZZ^d$-graded $\kk[\ff]$-module $\Gamma$ is
\begin{eqnarray*}
  \cC(\Gamma;\xx) &=& \text{the sum of terms of degree codim$(\Gamma)$
  in } \KK(\Gamma;\1-\xx).
\end{eqnarray*}
If $\Gamma = \kk[\ff]/I$ is the coordinate ring of a subscheme $X
\subseteq \spec(\kk[\ff])$, then we may also write $[X]_{\ZZ^d}$ or
$\cC(X;\xx)$ to mean $\cC(\Gamma;\xx)$.

In the context of quiver loci, $d = r_0 + \cdots + r_n$, and the
variables $\xx$ are split into a sequence of $n+1$ alphabets
$\xx^0,\ldots,\xx^n$ of sizes $r_0,\ldots,r_n$, so that the $i^\th$
alphabet is $\xx^i = x^i_1,\ldots,x^i_{r_i}$.  The coordinate ring
$\kk[\homv]$ is graded by~$\ZZ^d$, with the variable $\fff i\alpha\beta
\in \kk[\homv]$ having
\begin{eqnarray} \label{eq:weight}
\begin{array}{r@{\ }r@{\ }c@{\ }l}
  \text{ordinary weight}
  &\deg(\fff i\alpha\beta) &=& \xxx{i-1}\alpha-\xxx i\beta
\\[.5ex]
  \hbox{and exponential weight}
  &\weight(\fff i\alpha\beta) &=& \xxx {i-1}\alpha\!/\xxx i\beta
\end{array}
\end{eqnarray}
for each $i = 1,\ldots,n$.  Under this grading by $\ZZ^d$, the quiver
ideal~$I_\rr$ is homogeneous.

\begin{defn} \label{d:quivpoly}
The \bem{(ordinary) quiver polynomial} is the multidegree
\begin{eqnarray*}
  \QQ_\rr(\xx-\oxx) &=& \cC(\Omega_\rr;\xx)
\end{eqnarray*}
of the subvariety $\Omega_\rr$ inside $\homv$, under the $\ZZ^d$-grading
in which $\deg(\fff i\alpha\beta) = \xxx{i-1}\alpha-\xxx i\beta$.
\end{defn}

For the moment, the argument $\xx - \oxx$ of\/~$\QQ_\rr$ can be
regarded as a formal symbol, denoting that $n+1$ alphabets $\xx =
\xx^0,\ldots,\xx^n$ are required as input.  Later, in
Section~\ref{part:zel}, we shall define `double quiver polynomials',
with arguments $\xx-\oyy$, indicating two sequences of alphabets as input,
and then the symbol $\oxx$ will take on additional meaning as the
sequence $\xx^n,\ldots,\xx^0$ of alphabets constructed from~$\xx$.

Our discussion of quiver polynomials will require supersymmetric Schur
functions.  Given (finite or infinite) alphabets $X$ and~$Y$, define
for $r \in \NN$ the \bem{homogeneous symmetric function} $h_r(X-Y)$ in
the difference $X-Y$ of alphabets by the generating series
\begin{eqnarray*}
  \dfrac{\prod_{y\in Y} (1-zy)}{\prod_{x\in X} (1-zx)} &=&
  \sum_{r\ge0} z^r h_r(X-Y),
\end{eqnarray*}
in which $h_r(X-Y)=0$ for $r<0$.  Then the \bem{Schur function}
$s_\la(X-Y)$ in the difference $X-Y$ of alphabets is defined by the
Jacobi--Trudi determinant
\begin{eqnarray*}
  s_\la(X-Y) &=& \det\big(h_{\la_i-i+j}(X-Y)\big).
\end{eqnarray*}
The right hand side is the determinant of a $|\lambda|\times|\lambda|$
matrix, where $|\lambda|$ is the number of parts in the
partition~$\la$.

Buch and Fulton defined polynomials in \cite{BF} that also deserve
rightly to be called `quiver polynomials' (though they did not use
this appellation).  Given a sequence of maps between vector bundles
$E_0,\ldots,E_n$ on a scheme, the polynomials appearing in their Main
Theorem take as arguments the Chern classes of the bundles.  Here, we
shall express these polynomials formally as symmetric functions of the
Chern roots
\begin{eqnarray} \label{xr}
  \xx_\rr &=& \{x^i_j \mid i = 0,\ldots,n \text{ and } j =
  1,\ldots,r_i\}
\end{eqnarray}
of the dual bundles $E_0^\vee,\ldots,E_n^\vee$.  The Main Theorem
of~\cite{BF} is expressed in terms of integers $c_\bla(\rr)$ called
\bem{quiver constants}, each depending on a rank array~$\rr$ and a
sequence $\bla=(\la_1,\la_2,\dotsc,\la_n)$ of $n$ partitions.  Let
\mbox{$\sum c_\bla(\rr)s_\bla(\xx-\oxx)$} be the \mbox{corresponding
sum of products}
\begin{eqnarray} \label{eq:slax}
  s_\bla(\xx-\oxx) &=& \prod_{i=1}^n s_{\la_j}(\xx^{i-1}-\xx^i)
\end{eqnarray}
of Schur functions in differences of consecutive alphabets from
the sequence $\xx = \xx^0,\dotsc,\xx^n$.  In this notation, each
of the alphabets $\xx^i$ is taken to be {\em infinite}, so that
$s_\bla(\xx-\oxx)$ is a power series symmetric in each of $n+1$
sets of infinitely many variables.

This allows us to evaluate $s_\bla$ on finite alphabets of any given
size, by using $\xx_\rr$ as in~\eqref{xr}.  It results in our notation
\mbox{$\sum c_\bla(\rr)s_\bla(\xx_\rr-\oxx_\rr)$} for the polynomials
appearing in the Main Theorem of~\cite{BF}.

Note that quiver polynomials $\QQ_\rr(\xx-\oxx)$, as we defined them
in Definition~\ref{d:quivpoly}, automatically use only the finitely
many variables in~\eqref{xr}, even if evaluated on sequences of
a~priori infinite alphabets.

Before getting to Theorem~\ref{t:BFtoKMS}, we must first say more about
the geometry underlying multidegrees, which are algebraic substitutes
for equivariant cohomology classes of sheaves on torus representations.
Here, the torus is a maximal torus~$\tv$ inside~$\glv$ acting
on~$\homv$.  We choose~$\tv$ so the bases for $V_0,\ldots,V_n$ that were
fixed in Section~\ref{sec:quiver} consist of eigenvectors.  Thus $\tv$
consists of sequences of diagonal matrices, of sizes $r_0,\ldots,r_n$.
The group~$\ZZ^d$ in this context is the weight lattice of the
torus~$\tv$, and $x^i_j$ is the weight corresponding to the action of
the $j^\th$ diagonal entry in the $i^\th$~matrix.  The action of~$\tv$
on~$\homv$ via the inclusion $T \subset \glv$ induces the
grading~\eqref{eq:weight} on~$\kk[\homv]$.

The geometric data of a torus action on~$\homv$ determines an
equivariant Chow ring $A^*_T(\homv)$.  Let us recall briefly the
relevant definitions from \cite{EG98} in the context of a linear
algebraic group~$G$ acting on a smooth scheme~$M$.  Let $W$ be a
representation of~$G$ containing an open subset $U \subset W$ on which
$G$ acts freely.  For any integer~$k$, we can choose $W$ and~$U$ so
that $W \minus U$ has codimension at least $k+1$ inside~$W$
\cite[Lemma~9 in Section~6]{EG98}, and the construction is independent
of our particular choices by \cite[Definition--Proposition~1 and
Proposition~4]{EG98}.  The \bem{$G$-equivariant Chow ring} $A^*_G(M)$
has degree~$k$ piece
\begin{eqnarray*}
  A^k_G(M) &=& A^k((M \times U)/G)
\end{eqnarray*}
equal to the degree~$k$ piece of the ordinary Chow ring of the
quotient of~$M \times U$ modulo the diagonal action of~$G$.  We have
used smoothness of~$M$ via \cite[Corollary~2]{EG98}.

The equivariant Chow rings $A_T^*(M)$ for torus actions on vector
spaces~$M$ all equal the integral symmetric algebra ${\rm
Sym}^\spot_\ZZ({\mathfrak t}^*_\ZZ)$ of the weight lattice ${\mathfrak
t}^*_\ZZ$, independently of~$M$.  Thus, once we pick a basis $\xx =
x_1,\ldots,x_d$ for ${\mathfrak t}^*_\ZZ$, the class of any $T$-stable
subvariety is uniquely expressed as a polynomial in~$\ZZ[\xx]$.  More
generally, every $\ZZ^d$-graded module~$\Gamma$ (equivalently,
$T$-equivariant sheaf~$\tilde\Gamma$) determines a polynomial
$[\cycle(\Gamma)]_T$, where the \bem{cycle} of~$\Gamma$ is the sum of
its maximal-dimensional support components with coefficients given by
the generic multiplicity of~$\Gamma$ along each.

\begin{prop} \label{p:chow}
The multidegree\/ $\cC(\Gamma;\xx)$ of a $\ZZ^d$-graded
module\/~$\Gamma$ over $\kk[\ff]$ is the class\/ $[\cycle(\Gamma)]_T$
in the equivariant Chow ring $A^*_T(M)$ of $M = \spec(\kk[\ff])$.
\end{prop}
\begin{proof}
The function sending each graded module to its multidegree is
characterized uniquely as having the `Additivity', `Degeneration', and
`Normalization' properties in \cite[Theorem~1.7.1]{grobGeom}.
Thus it suffices to show that the function sending each module to the
equivariant Chow class of its cycle also has these properties.
Additivity for Chow classes is by definition.  Degeneration is the
invariance of Chow classes under rational equivalence of cycles.
Normalization, which says that a coordinate subspace~$L \subseteq M$
has class equal to the product of the weights of the variables in its
defining ideal, follows from \cite[Theorem~2.1]{Bri97} and
\cite[Proposition~4]{EG98} by \mbox{downward induction on the
dimension of~$L$}.%
\end{proof}

\begin{thm} \label{t:BFtoKMS}
The quiver polynomial $\QQ_\rr(\xx-\oxx)$ in
Definition~\ref{d:quivpoly} coincides with the polynomial $\sum
c_\bla(\rr)s_\bla(\xx_\rr-\oxx_\rr)$ appearing in the Main Theorem
of\/~\cite{BF}.
\end{thm}
\begin{proof}
Let $\Gr(\rr,\ell)$ be the product $\Gr(r_0,\ell) \times \cdots \times
\Gr(r_n,\ell)$ of Grassmannians of subspaces of dimensions
$r_0,\ldots,r_n$ inside $\kk^\ell$ for some large~$\ell$.  This
variety $\Gr(\rr,\ell)$ comes endowed with $n+1$ universal vector
bundles $\VV_0,\ldots,\VV_n$ of ranks $r_0,\ldots,r_n$.
Denote by $\hhomv$ the vector bundle $\prod_{j=1}^n
\Hom(\VV_{i-1},\VV_i)$ over $\Gr(\rr,\ell)$.  Pulling back $\VV_i$ to
$\hhomv$ for all~$i$ yields vector bundles with a universal (or
``tautological'') sequence
\begin{eqnarray*}
  \wt\Phi: && \wt\VV_0\ \stackrel{\wt\Phi_1}\too\ \wt\VV_1\
  \stackrel{\wt\Phi_2}\too\ \cdots\ \stackrel{\wt\Phi_{n-1}}\too\
  \wt\VV_{n-1}\ \stackrel{\wt\Phi_n}\too\ \wt\VV_n.
\end{eqnarray*}
of maps between them.  Let $\OOmega_\rr \subseteq \hhomv$ denote the
degeneracy locus of~$\wt\Phi$ for ranks~$\rr$.  If~$k = d(\rr)$ is the
codimension of $\OOmega_\rr$ from Corollary~\ref{c:codim}, then we can
assume $\ell$\/ has been chosen so large that the graded component
$A^k(\hhomv)$ of the Chow ring of~$\hhomv$ consists of all symmetric
polynomials of degree~$k$ in~$\ZZ[\xx_\rr]$, for $\xx_\rr$ as
in~\eqref{xr}.

The class $[\OOmega_\rr] \in A^k(\hhomv)$ is by definition the
polynomial in the Main Theorem of~\cite{BF}.  On the other hand,
$A^k(\hhomv)$ is by \cite[Corollary~2]{EG98} the degree~$k$ piece
$A^k_\glv(\homv)$ of the $\glv$-equivariant Chow ring of~$\homv$, and
the ordinary class $[\OOmega_\rr] \in A^k(\hhomv)$ coincides with the
equivariant class $[\Omega_\rr]_\glv \in A^k_\glv(\homv)$.  The result
now follows from Proposition~\ref{p:chow}, given that
$A^*_\glv(\homv)$ consists of the symmetric group invariants
inside $A^*_T(\homv)$ by \cite[Proposition~6]{EG98}.%
\end{proof}

\begin{remark}
Thinking of quiver polynomials as equivariant Chow or equivariant
cohomology classes makes sense from the topological viewpoint taken
originally by Thom \cite{thom55} in the context of what we now know as
the Giambelli--Thom--Porteous formula.  In the present context, a
quiver of vector bundles on a space~$Y$ is induced by a unique
homotopy class of maps from~$Y$ to the universal bundle
$\homv_{\!\glv}$ with fiber~$\homv$ on the classifying space
of~$\glv$.  The induced map $H^*_\glv(\text{point}) = H^*(\glv) \to
H^*(Y)$ on cohomology simply evaluates the universal degeneracy locus
polynomial for~$\rr$, namely the quiver polynomial, on the Chern
classes of the vector bundles in the quiver.  Replacing the
classifying space for~$\glv$ by the corresponding Artin stack
accomplishes the same thing in the Chow ring context instead of
cohomology.
\end{remark}

\part{Double quiver polynomials}

\label{part:zel}
\setcounter{section}{0}
\setcounter{thm}{0}
\setcounter{equation}{0}

\section{Double Schubert polynomials}

\label{sec:schub}

Let $d$ be a positive integer, $\xx=(x_1,x_2,\dotsc,x_d)$ and
$\yy=(y_1,y_2,\dotsc,y_d)$ two sets of variables, and $w\in S_d$ a
permutation.  The double Schubert polynomial $\SS_w(\xx-\oyy)$ of
Lascoux and Sch\"utzenberger \cite{LS82} is defined by downward
induction on the length of~$w$ as follows.  For the long permutation
$w_0\in S_d$ that reverses the order of $1,\ldots,d$, set
\begin{eqnarray*}
  \SS_{w_0}(\xx-\oyy) &=& \prod_{i+j\le d} (x_i-y_j).
\end{eqnarray*}
If $w \neq w_0$, then choose $q$ so that $w(q) < w(q+1)$, and define
\begin{eqnarray*}
  \SS_w(\xx-\oyy) &=& \partial_q \SS_{ws_q}(\xx-\oyy).
\end{eqnarray*}
The operator $\partial_q$ is the $q^\th$ divided difference acting on
the $\xx$ variables, so
\begin{eqnarray}\label{eq:partial}
  \partial_q P(x_1, x_2, \ldots) &=& \frac{P(x_1, x_2, \ldots\,) -
  P(x_1, \ldots, x_{q-1}, x_{q+1}, x_q, x_{q+2}, \ldots\,)}{x_q-x_{q+1}}
\end{eqnarray}
for any polynomial $P$ taking the $\xx$ variables and possibly other
variables for input.

In Section~\ref{part:pos} we will require Schubert polynomials for
{\em partial}\/ permutations~$w$ (Definition~\ref{d:matSchub}).  By
convention, when we write $\SS_w$ for a partial permutation~$w$ of
size $k \times \ell$, we shall mean $\SS_\wt w$, where $\wt w$ is a
minimal-length completion of~$w$ to an actual permutation.  Thus $\wt
w$ can be be any permutation whose matrix agrees with~$w$ in the
northwest $k \times \ell$ rectangle, and such that in the columns
strictly to the right of~$\ell$, as well as in the rows strictly
below~$k$, the nonzero entries of~$\wt w$ proceed from northwest to
southeast (no nonzero entry is northeast of another).  This convention
is consistent with Definition~\ref{d:matSchub}, in view of the next
result, because the matrix Schubert varieties $\ol X_w$ and $\ol X_\wt
w$ have equal multidegrees, as they are defined by the {\em same}\/
equations (though in a~priori different polynomial rings).

\begin{thm}{\cite[Theorem~A]{grobGeom}} \label{t:dubSchub}
Let $w$ be a partial permutation matrix of size $k \times \ell$, so
$\Xb_w$ is a subvariety of $M_{k\ell} = \spec(\kk[\ff])$, where\/ $\ff =
(f_{\alpha\beta})$ for $\alpha = 1,\ldots,k$ and $\beta =
1,\ldots,\ell$.  If~$f_{\alpha\beta}$ has ordinary weight
$\deg(f_{\alpha\beta}) = x_\alpha-y_\beta$, then $\ol X_w$ has
multidegree
\begin{eqnarray*}
  \cC(\ol X_w;\xx,\yy) &=& \SS_w(\xx-\oyy)
\end{eqnarray*}
equal to the double Schubert polynomial for~$w$.  (We remind the reader
that the matrix for a partial permutation~$w$ has a $1$ in row~$q$ and
column~$p$ if and only if $w(q) = p$.)
\end{thm}

When the partial permutation~$w$ is an honest permutation $v \in S_d$,
we usually break the variables $\xx$ and~$\yy$ in the argument of
$\SS_v(\xx-\oyy)$ into two sequences of alphabets
\begin{eqnarray} \label{xx}
  \xx = \xx^0,\ldots,\xx^n &\text{and}& \yy = \yy^n,\ldots,\yy^0,
\\ \label{xxj}
  \makebox[0ex][r]{where\quad} \xx^j = x^j_1,\ldots,x^j_{r_{\!j}}
  &\text{and}& \yy^j = y^j_1,\ldots,y^j_{r_{\!j}}.
\end{eqnarray}
On the other hand, most partial permutations~$w$ that occur in the
sequel will have size $r_{j-1} \times r_{\!j}$ for some~$j \in
\{1,\ldots,n\}$; in that case we consider $\SS_w(\xx^{j-1}-\yy^j)$.

\section{Double quiver polynomials}

\label{sec:doublequiv}

The coordinate ring $\kk[\endv]$ of the $d \times d$ matrices is
multigraded by the group~$\ZZ^{2d} = (\ZZ^{r_0} \oplus \cdots \oplus
\ZZ^{r_n})^2$, which we take to have basis as in~\eqref{xx}
and~\eqref{xxj}.  In our context, it is most natural to index the
variables $\ffff ij\alpha\beta\,$ in the generic \mbox{$d \times d$}
matrix in a slightly unusual manner: $\ffff ij\alpha\beta \in
\kk[\endv]$ occupies the spot in row~$\alpha$ and column~$\beta$ within
the rectangle at the intersection of the $i^\th$~block row and the
\mbox{$j^\th$}~block column, where $i,j = 0,\ldots,n$ and we label block
columns starting from the right.  Set
\begin{eqnarray*}
  \deg(\ffff ij\alpha\beta\,) &=& x^i_\alpha - y^j_\beta
\end{eqnarray*}
as the ordinary weight of~$\ffff ij\alpha\beta\,$.

More pictorially, label the rows of the \mbox{$d \times d$} grid with
the~$\xx$ basis vectors in the order they are given, from top to bottom,
and similarly label the columns with~$\yy = \yy^n,\ldots,\yy^0$, from
left to right (go left to right within each alphabet).  Under our
multigrading, the ordinary weight of a variable is its row basis vector
minus its column basis vector.  The different alphabets on rows and
columns line up as follows.
\begin{itemize}
\item
Rows in the $j^\th$ block from the top are labeled $x^j_1,\ldots,x^j_r$
for $r = r_{\!jj}$, starting with $x^j_1$ in the highest row of that
block.

\item
Columns in the $j^\th$ block {\em from the right}\/ are labeled
$y^j_1,\ldots,y^j_r$ for $r = r_{\!jj}$, starting with $y^j_1$ above the
leftmost column of that block.
\end{itemize}

For notational clarity in examples, it is usually convenient to rename
the alphabets $\xx^0,\xx^1,\xx^2,\ldots\ $ using distinct Latin names,
such as
\begin{equation*}
  \xx^0 = \aa = a_1,a_2,a_3,\ldots\ \hbox{ and }\ \xx^1 = \bb =
  b_1,b_2,b_3,\ldots\ \hbox{ and }\ \xx^2 = \cc = c_1,c_2,c_3,\ldots
\end{equation*}
and so on, rather than upper indices.  Then, we rename the alphabets
$\yy^0,\yy^1,\yy^2,\ldots\ $~as
\begin{equation*}
  \yy^0 = \dot\aa = \da_1,\da_2,\da_3,\ldots\ \hbox{ and }\ \yy^1 =
  \dot\bb = \db_1,\db_2,\db_3,\ldots\ \hbox{ and }\ \yy^2 = \dot\cc =
  \dc_1,\dc_2,\dc_3,\ldots
\end{equation*}
and so on, the same as the $\xx$'s but with dots on top.

All of the notation above should be made clearer by the following
example.

\begin{example} \label{ex:weights}
Consider quivers on the vector spaces $V = (V_0,V_1,V_2)$ with
dimensions $(r_0,r_1,r_2) = (2,3,1)$.  The coordinate ring $\kk[\endv]$
has variables $\ffff ij\alpha\beta$ as they appear in the matrices below
(the $f$~variables are the same in both):
\begin{rcgraph}
\begin{array}{@{}l@{\:}|c@{\ }|@{\ }c@{\ }c@{\ }c@{\ }|@{\ }c@{\ }%
    c|l@{\hspace{-4ex}}}
\mcc{}&\mcc{\ph}&\ph&\ph&\mcc{\ph}&\ph&\mcc{\ph}
\\[-2ex]
\mcc{}&\mcc{\!\!y^2_1}&y^1_1&y^1_2&\mcc{\!\!y^1_3}&y^0_1&\mcc{\!\!y^0_2}
\\[.25ex]\cline{2-7}
x^0_1&\ffff 0211&\ffff 0111&\ffff 0112&\ffff 0113&\ffff 0011&\ffff 0012
&\phantom{\displaystyle\sum}\\
x^0_2&\ffff 0221&\ffff 0121&\ffff 0122&\ffff 0123&\ffff 0021&\ffff 0022
&\phantom{\displaystyle\sum}\\\cline{2-7}
x^1_1&\ffff 1211&\ffff 1111&\ffff 1112&\ffff 1113&\ffff 1011&\ffff 1012
&\phantom{\displaystyle\sum}\\
x^1_2&\ffff 1221&\ffff 1121&\ffff 1122&\ffff 1123&\ffff 1021&\ffff 1022
&\phantom{\displaystyle\sum}\\
x^1_3&\ffff 1231&\ffff 1131&\ffff 1132&\ffff 1133&\ffff 1031&\ffff 1032
&\phantom{\displaystyle\sum}\\\cline{2-7}
x^2_1&\ffff 2211&\ffff 2111&\ffff 2112&\ffff 2113&\ffff 2011&\ffff 2012
&\phantom{\displaystyle\sum}\\\cline{2-7}
\end{array}
\quad\begin{array}{c}\\[2ex]=\end{array}\quad
\begin{array}{@{}l@{\:}|c@{\ }|@{\ }c@{\ }c@{\ }c@{\ }|@{\ }c@{\ }c|l@{\hspace{-4ex}}}
\mcc{}&\mcc{\ph}&\ph&\ph&\mcc{\ph}&\ph&\mcc{\ph}
\\[-2ex]
\mcc{}&\mcc{\!\!\dc_1}&\db_1&\db_2&\mcc{\!\!\db_3}&\da_1&\mcc{\!\!\da_2}
\\[.25ex]\cline{2-7}
a_1&\ffff 0211&\ffff 0111&\ffff 0112&\ffff 0113&\ffff 0011&\ffff 0012
&\phantom{\displaystyle\sum}\\
a_2&\ffff 0221&\ffff 0121&\ffff 0122&\ffff 0123&\ffff 0021&\ffff 0022
&\phantom{\displaystyle\sum}\\\cline{2-7}
b_1&\ffff 1211&\ffff 1111&\ffff 1112&\ffff 1113&\ffff 1011&\ffff 1012
&\phantom{\displaystyle\sum}\\
b_2&\ffff 1221&\ffff 1121&\ffff 1122&\ffff 1123&\ffff 1021&\ffff 1022
&\phantom{\displaystyle\sum}\\
b_3&\ffff 1231&\ffff 1131&\ffff 1132&\ffff 1133&\ffff 1031&\ffff 1032
&\phantom{\displaystyle\sum}\\\cline{2-7}
c_1&\ffff 2211&\ffff 2111&\ffff 2112&\ffff 2113&\ffff 2011&\ffff 2012
&\phantom{\displaystyle\sum}\\\cline{2-7}
\end{array}
\end{rcgraph}
The ordinary weight of each $f$ variable equals its row label minus its
column label, and its exponential weight is the ratio instead of the
difference.  For example, the variable $\ffff {\,0}123$ has ordinary
weight~$x^0_2 - y^1_3 = a_2 - \db_3$ and exponential weight $x^0_2/y^1_3
= a_2/\db_3$.\qed%
\end{example}

The coordinate ring $\kk[Y_0]$ of the variety $Y_0$ from
Section~\ref{sec:zel} is not naturally multigraded by all of $\ZZ^{2d}$,
but only by~$\ZZ^d$, with the variable $\ffff ij\alpha\beta \in
\kk[Y_0]$ having ordinary weight $x^i_\alpha-x^j_\beta$.  This
convention is consistent with the multigrading of~$\kk[\homv]$
in~\eqref{eq:weight} under the homomorphism to~$\kk[Y_0]$ induced by the
Zelevinsky map, since the variable $\ffff {i-1,}i\alpha\beta \in
\kk[Y_0]$ maps to $\fff i\alpha\beta \in \kk[\homv]$, and their ordinary
weights $\xxx{i-1}\alpha-\xxx i\beta$ agree.  Our $\ZZ^d$-grading
of~$\kk[Y_0]$ is positive, so every quotient of~$\kk[Y_0]$ has a Hilbert
series.

\begin{example} \label{ex:U}
In the situation of Example~\ref{ex:weights}, the coordinate ring
$\kk[Y_0]$ has only the variables $\ffff ij\alpha\beta$ that appear
in the matrices below:
\begin{rcgraph}
\begin{array}{@{}l@{\:}|c@{\ }|@{\ }c@{\ }c@{\ }c@{\ }|@{\ }c@{\ }c|@{}}
\mcc{}&\mcc{\ph}&\ph&\ph&\mcc{\ph}&\ph&\mcc{\ph}
\\[-2ex]
\mcc{}&\mcc{\!\!x^2_1}&x^1_1&x^1_2&\mcc{\!\!x^1_3}&x^0_1&\mcc{\!\!x^0_2}
\\[.25ex]\cline{2-7}
x^0_1&\ffff 0211&\ffff 0111&\ffff 0112&\ffff 0113&1&\phantom{\displaystyle\sum}
\\
x^0_2&\ffff 0221&\ffff 0121&\ffff 0122&\ffff 0123&\phantom{\displaystyle\sum}&1
\\\cline{2-7}
x^1_1&\ffff 1211&     1    &          &          & &\phantom{\displaystyle\sum}
\\
x^1_2&\ffff 1221&          &     1    &          & &\phantom{\displaystyle\sum}
\\
x^1_3&\ffff 1231&          &          &     1    & &\phantom{\displaystyle\sum}
\\\cline{2-7}
x^2_1&1&\phantom{\ffff 2211}&         &          & &\phantom{\displaystyle\sum}
\\\cline{2-7}
\end{array}
\quad\begin{array}{c}\\[2ex]=\end{array}\quad
\begin{array}{@{}l@{\,}|c@{\ }|@{\ }c@{\ }c@{\ }c@{\ }|@{\ }c@{\ }c|@{}}
\mcc{}&\mcc{\ph}&\ph&\ph&\mcc{\ph}&\ph&\mcc{\ph}
\\[-2ex]
\mcc{}&\mcc{\!\!c_1}&b_1&b_2&\mcc{\!\!b_3}&a_1&\mcc{\!\!a_2}
\\[.25ex]\cline{2-7}
a_1&\ffff 0211&\ffff 0111&\ffff 0112&\ffff 0113&1&\phantom{\displaystyle\sum}
\\
a_2&\ffff 0221&\ffff 0121&\ffff 0122&\ffff 0123&\phantom{\displaystyle\sum}&1
\\\cline{2-7}
b_1&\ffff 1211&     1    &          &          & &\phantom{\displaystyle\sum}
\\
b_2&\ffff 1221&          &     1    &          & &\phantom{\displaystyle\sum}
\\
b_3&\ffff 1231&          &          &     1    & &\phantom{\displaystyle\sum}
\\\cline{2-7}
c_1&1&\phantom{\ffff 2211}&         &          & &\phantom{\displaystyle\sum}
\\\cline{2-7}
\end{array}
\end{rcgraph}
In this case, the variable $\ffff {\,0}123$ has ordinary weight~$x^0_2 -
x^1_3 = a_2 - b_3$ and exponential weight $x^0_2/x^1_3 = a_2/b_3$.\qed%
\end{example}

The multigrading giving rise to the ordinary quiver polynomials in
Section~\ref{part:quiver} is derived from the action of~$\glv$ on~$\homv$.
In~the current context of~$\gld$ acting on~$\endv$, we identify $\glv =
\prod_{i=0}^n \gl(V_i)$ as the Levi factor~$L$, consisting of block
diagonal matrices inside the parabolic subgroup $P \subset \gld$, by
sticking $\gl(V_i)$ in the $i^\th$ diagonal block.  This identifies the
torus~$T$ of diagonal matrices in~$L$ with the torus~$\tv$ from
Section~\ref{sec:quiverpolys}.

\begin{lemma} \label{l:PY}
The matrix Schubert variety $\ol X_{v(\rr)}$ is the closure
inside~$\endv$ of the variety $P Y_\rr$, by which we mean the image
of~$P \times Y_\rr$ under multiplication.
\end{lemma}
\begin{proof}
In general, if~$X_v$ is a Schubert subvariety of~$\pgld$, and~$\tilde
X_v$ is the preimage of~$X_v$ inside~$\gld$ (under projection
modulo~$P$), then $\ol X_v$ equals the closure of~$\tilde X_v$
inside~$\endv$.  Theorem~\ref{t:zel} implies that $PY_\rr$ is dense
inside the preimage of~$X_{v(\rr)}$.%
\end{proof}

The quiver locus~$\Omega_\rr$ has only one copy of the torus $T
\subset \glv$ acting on it, as does its Zelevinsky image~$Y_\rr$
(conjugation by $T \subset L$).  Multiplying by~$P$ on the left as in
Lemma~\ref{l:PY} frees up the left and right hand tori (and even the
left and right hand Levi factors), allowing them to act independently.
Consequently, we can define a {\em doubly}\/ equivariant cohomology
class on the smeared Zelevinsky image~$\ol X_{v(\rr)}$
of~$\Omega_\rr$.  In terms of multidegrees, we can now use the two
sequences $\xx$ and~$\yy$ of alphabets instead of only the $\xx$
alphabets in ordinary quiver polynomials.

\begin{defn} \label{d:QQrr}
The \bem{double quiver polynomial~$\QQ_\rr$} is the ratio
\begin{eqnarray*}
  \QQ_\rr(\xx-\oyy) &=&
  \frac{\SS_{v(\rr)}(\xx-\oyy)}{\SS_{v(\homv)}(\xx-\oyy)}
\end{eqnarray*}
of double Schubert polynomials in the concatenations of the two
sequences of {\em finite}\/ alphabets described in~\eqref{xx}
and~\eqref{xxj}.
\end{defn}

The denominator $\SS_{v(\homv)}(\xx-\oyy)$ should be regarded as a fudge
factor, being simply the product of all weights $(x_* - y_*)$ of
variables lying strictly above the block superantidiagonal.  These
variables lie in locations corresponding to~$*$~entries in the diagram
of every Zelevinsky permutation, so $\SS_{v(\homv)}$ obviously
divides~$\SS_{v(\rr)}$ (see Theorem~\ref{t:BJS}).

The simple relation between double and ordinary quiver polynomials, to
be presented in Theorem~\ref{c:ratio}, justifies the notation
$\QQ_\rr(\xx-\oxx)$ for the ordinary case.

\section{Ratio formula for quiver polynomials}

\label{sec:ratio}

Our goal in this section is to express quiver polynomials as the
specializations of double quiver polynomials obtained by setting $\yy^j
= \xx^j$ for all~$j$.

More generally, we shall express the quiver \K-polynomial as a ratio
of corresponding specialized \bem{double Grothendieck polynomials}.
These are usually defined by `isobaric divided difference operators'
(also known as `Demazure operators'), which we shall review in
Section~\ref{sec:demazure}.  However, since Grothendieck polynomials
are tangential to the focus here, we do not review the background,
which can be found in \cite{LS82,Las90} or \cite{grobGeom}, the latter
being well-suited to the current applications.  Readers unfamiliar
with Grothendieck polynomials will lose nothing in the present context
by taking the first sentence of Proposition~\ref{p:grothBigCell} as a
definition.

Double Grothendieck polynomials $\GG_v(\xx/\yy)$ for permutations $v \in
S_d$ take as arguments two alphabets $\xx$ and~$\yy$, each having
size~$d$.  For purposes relating to quiver loci, $\xx$ and~$\yy$ are
broken into sequences of alphabets as in~\eqref{xx} and~\eqref{xxj}.
Write
\begin{eqnarray*}
  \oxx &=& \xx^n,\ldots,\xx^0
\end{eqnarray*}
to mean the reverse of the finite list~$\xx$ of alphabets
from~\eqref{xx}, which will arise when calculating the $K$-polynomials
of opposite `cells' in Schubert subvarieties of partial flag manifolds
$\pgld$, as we shall see shortly.  Note that in the case of full flags,
where $P = B$ is the Borel subgroup of lower-triangular matrices
in~$\gld$, each alphabet in the list~$\xx$ consists of just one variable
(as opposed to there being only one alphabet in the list), so the
reversed list~$\oxx$ is really just a reversed alphabet in that case.

For notation, let $\KK_M(\ol X;\xx,\yy)$ denote the \K-polynomial of a
multigraded subvariety $\ol X \subseteq \endv$ under the
$\ZZ^{2d}$-grading from Section~\ref{sec:doublequiv}, and let
$\KK_{Y}(Z;\xx)$ denote the \K-polynomial of a multigraded subvariety $Z
\subseteq Y_0$ under the $\ZZ^d$-grading above.

\begin{prop} \label{p:grothBigCell}
The \K-polynomial $\KK_M(\ol X_v;\xx,\yy)$ of a matrix Schubert variety
$\ol X_v \subseteq \endv$ is the double Grothendieck
polynomial~$\GG_v(\xx/\yy)$.  Intersecting $\ol X_{v(\rr)}$ with~$Y_0$
yields the variety~$Y_\rr$, whose \K-polynomial $\KK_{Y}(Y_v;\xx)$ is the
specialization~$\GG_{v(\rr)}(\xx/\oxx)$.
\end{prop}
\begin{proof}
The first sentence comes from \cite[Theorem~A]{grobGeom}, one of the
main results there.  For the second sentence, observe that $Y_\rr$ has
the same codimension inside~$Y_0$ as does $\ol X_{v(\rr)}$
inside~$\endv$.  Hence the equations defining~$Y_0$ as a subvariety
of~$\endv$ (namely $f - 1$ for diagonal~$f$ in antidiagonal blocks and
$f - 0$ for the other variables in or below the main block antidiagonal)
form a regular sequence on~$\ol X_{v(\rr)}$.  Coarsening the
\mbox{$\ZZ^{2d}$-grading} on $\kk[\endv]$ to the grading by~$\ZZ^d$ in
which $\ffff ij\alpha\beta\,$ has ordinary weight $x^i_\alpha -
x^j_\beta$, by setting $y^j_\beta = x^j_\beta$, makes the equations
defining~$Y_0$ become homogeneous, because the variables set equal
to~$1$ have weight zero.  Therefore, if $\FF_\spot$ is any
$\ZZ^{2d}$-graded free resolution of~$\kk[\ol X_{v(\rr)}]$
over~$\kk[\endv]$, then $\FF_\spot \otimes_{\kk[\endv]} \kk[Y_0]$ is a
$\ZZ^d$-graded free resolution of~$\kk[Y_\rr]$ over~$\kk[Y_0]$.  It
follows that $\KK_M(\ol X_{v(\rr)};\xx,\oxx) = \KK_{Y}(Y_\rr;\xx)$,
because the relevant $\ZZ^d$-graded Betti numbers coincide.  Applying
the first sentence of the Proposition and substituting $y^j_\beta =
x^j_\beta$ now proves the second sentence.%
\end{proof}

\begin{thm} \label{t:ratio}
The \K-polynomial $\KK_\homv(\Omega_\rr;\xx)$ of\/ $\Omega_\rr$
inside\/~$\homv$ is the ratio
\begin{eqnarray*}
  \KK_\homv(\Omega_\rr;\xx) &=& \frac{\GG_{v(\rr)}(\xx/\oxx)}
  {\GG_{v(\homv)}(\xx/\oxx)}
\end{eqnarray*}
of specialized double Grothendieck polynomials for $v(\rr)$
and~$v(\homv)$.%
\end{thm}
\begin{proof}
The equality $H(\Omega_\rr;\xx) = \KK_\homv(\Omega_\rr;\xx)
H(\homv;\xx)$ of Hilbert series (which are well-defined by positivity of
the grading of $\kk[\homv]$ by~$\ZZ^d$) follows from the definition of
\K-polynomial.  For the same reason, we have
\begin{eqnarray*}
  H(\homv;\xx) = \KK_{Y}(Y_\homv;\xx) H(Y_0;\xx) &\hbox{and also}&
  H(\Omega_\rr;\xx) = \KK_{Y}(Y_\rr;\xx) H(Y_0;\xx).
\end{eqnarray*}
Therefore $\KK_{Y}(Y_\rr;\xx)H(Y_0;\xx) = \KK_\homv(\Omega_\rr;\xx)
\KK_{Y}(Y_\homv;\xx)H(Y_0;\xx)$, and hence
\begin{eqnarray*}
  \KK_\homv(\Omega_\rr;\xx) &=&
  \frac{\KK_{Y}(Y_\rr;\xx)}{\KK_{Y}(Y_\homv;\xx)}.
\end{eqnarray*}
Now apply Theorem~\ref{t:zel} and Proposition~\ref{p:grothBigCell} to
complete the proof.%
\end{proof}

\begin{remark} \label{rk:v0}
Here is the geometric content of Proposition~\ref{p:grothBigCell}.
The torus $T \times T$ of dimension~$2d$ acts on~$\endv$, the left
factor by multiplication on the left, and the right factor by inverse
multiplication on the right.  Viewing $\GG_v(\xx/\yy)$ as an
equivariant class for this action of $T \times T$, specializing the
$\yy$~variables to~$\oxx$ does not arise by restricting the $T \times
T$ action to some subtorus.  Rather, it uses a different
$d$-dimensional torus~$T$ entirely.  The main point is that the
opposite big cell is an orbit for an {\em opposite}\/ parabolic group,
whose torus has been conjugated by the identity-in-each-%
antidiagonal-block permutation~${w_0}\ww_0$ from the usual one.

In more detail, the subvariety~$Y_0 \subset \gld$ is not stable under
the usual action of $T \times T$ above, because we need to preserve
the~$1$s in antidiagonal blocks.  Instead, it is stable under the
action of the diagonal torus~$T$ inside the copy of $T \times T$ whose
left factor acts by usual multiplication on the left, but where a
diagonal matrix~$\tau$ in the right factor acts as right
multiplication by ${w_0}\ww_0 \tau^{-1} \ww_0{w_0}$.  Thus the
diagonal $T$ acts with the weight on the variable at position
$(\alpha,\beta)$ in block row~$i$ and block column~$j$ (from the {\em
right}) being $x^i_\alpha - x^j_\beta$

As a reality check, note that this action of the right hand~$T$ factor
is indeed the restriction of the ``opposite'' parabolic action
on~$\pgld$ to its torus.  Indeed, our ``usual'' parabolic acting
on~$\pgld$ is the block upper-triangular $P_+ \subset \gld$ acting by
inverse multiplication on the right (we want $P \times P_+$ to act via
a left group action), so the opposite parabolic is ${w_0}\ww_0 P_+
\ww_0{w_0}$ acting by inverse multiplication on the right.%

One might observe that the torus in ${w_0}\ww_0 P_+ \ww_0{w_0}$ acts
on $\pgld$ simply by inverse right multiplication, not by twisted
conjugation.  All this means is that multiplying by a diagonal matrix
${w_0}\ww_0 \tau^{-1} \ww_0{w_0}$ moves the cosets in $\pgld$
correctly.  However, the subvariety $Y_0 \subset \gld$ consists of
specific coset {\em representatives}\/ for points in the opposite big
cell of~$\pgld$, and we need to know which coset representative
corresponds to the image of a matrix after right multiplication by
${w_0}\ww_0 \tau^{-1} \ww_0{w_0}$.  For this, we have to left-multiply
by~$\tau$, bringing the antidiagonal blocks back to a list of identity
matrices without altering the coset.%
\end{remark}

The proofs of both Proposition~\ref{p:grothBigCell} and
Theorem~\ref{t:ratio} could as well have been presented entirely in
geometric language, referring to pullbacks and pushforwards in
equivariant \K-theory rather than ratios of numerators of Hilbert
series.

Next comes what we were after from the beginning, the consequence of
Theorem~\ref{t:ratio} on multidegrees (or equivariant cohomology or Chow
groups, by Proposition~\ref{p:chow}).

\begin{thm} \label{c:ratio}
The quiver polynomial $\QQ_\rr(\xx-\oxx)$ for a rank array\/~$\rr$
equals the ratio
\begin{eqnarray*}
  \QQ_\rr(\xx-\oxx) &=& \frac{\SS_{v(\rr)}(\xx-\oxx)}
  {\SS_{v(\homv)}(\xx-\oxx)}
\end{eqnarray*}
of specialized double Schubert polynomials.  In other words, the
ordinary quiver polynomial $\QQ_\rr(\xx-\oxx)$ is the $\yy = \xx$
specialization of the double quiver polynomial~$\QQ_\rr(\xx-\oyy)$.%
\end{thm}
\begin{proof}
Clear denominators in Theorem~\ref{t:ratio},
substitute $1-x$ for every occurrence of each variable~$x$, take lowest
degree terms, and divide the result by $\SS_{v(\homv)}(\xx-\oxx)$.  Note that
$\SS_{v(\homv)}(\xx-\oxx)$ is nonzero, being simply the product of the
$\ZZ^d$-graded weights of the variables in~$\kk[Y_0]$ strictly above the
block superantidiagonal.%
\end{proof}

In view of Theorem~\ref{c:ratio}, we consider Definition~\ref{d:QQrr}
to be our first positive formula for double quiver polynomials, since
it expresses them in terms of known polynomials (Schubert polynomials)
that depend on combinatorial data (the Zelevinsky permutation)
associated to given ranks~$\rr$.  This is not so bad: we could have
(without motivation) defined double quiver polynomials starting from
the Buch--Fulton formula \cite{BF} by using independent dual bundles
in the second half of each argument in their Schur functions, and then
Definition~\ref{d:QQrr} would be a theorem (but not one in this paper;
see Remark~\ref{rk:double}).  As we have chosen to present things, the
double version of the Buch--Fulton formula will come out as a theorem
instead, in Theorems~\ref{t:bla}, \ref{t:PeelQuiver},
and~\ref{t:ptofs}; see also Corollary~\ref{c:bla}.

\part{Lacing diagrams}

\label{part:lacing}
\setcounter{section}{0}
\setcounter{thm}{0}
\setcounter{equation}{0}

\section{Geometry of lacing diagrams}

\label{sec:lacing}

In this section we use diagrams derived from those in~\cite{AD} to
describe the orbits of certain groups related to~$\glv$.  Recall that we
have fixed bases for the vector spaces $V_0,\ldots,V_n$, and this
essentially amounted to fixing a maximal torus~$\tv$ inside the
group~$\gl$ acting on~$V$.  Suppose that, in this basis, the quiver
representation $\phi \in \homv$ is represented by a sequence $\ww =
(w_1,\ldots,w_n)$ of $n$ partial permutation matrices with respect to
the above bases.  We shall freely identify $\phi$ with~$\ww$.

Each list $\ww = (w_1,\ldots,w_n)\in\homv$ of partial permutation
matrices can be represented by a (nonembedded) graph in the plane
called a \bem{lacing diagram}.  The vertex set of the graph consists
of $r_i$ bottom-justified dots in column~$i$ for $i = 0,\ldots,n$,
with an edge connecting the dot at height~$\alpha$ (from the bottom)
in column~$i-1$ with the dot at height~$\beta$ in column~$i$ if and
only if the entry of~$w_i$ at $(\alpha,\beta)$ is~$1$.  A~\bem{lace}
is a connected component of a lacing diagram.  An $(i,j)$-lace is one
that starts in column~$i$ and ends in column~$j$.  We shall also
identify~$\ww\in\homv$ with \mbox{its associated lacing diagram}.

\begin{example} \label{ex:lace}
Here is a lacing diagram for the dimension vector $(2,3,4,3)$ and
its corresponding list of partial permutation matrices.
\begin{eqnarray*}
\pspicture[.1](0,0)(3,3)
\psdots(0,-1)(0,0)(1,-1)(1,0)(1,1)(2,-1)(2,0)(2,1)(2,2)(3,-1)(3,0)(3,1)
\psline(0,-1)(1,-1)(2,0)
\psline(0,0)(1,0)
\psline(1,1)(2,-1)(3,-1)
\psline(2,1)(3,0)
\endpspicture\
&\quad\longleftrightarrow\quad&
\left(\begin{bmatrix}
  1&0&0 \\
  0&1&0
\end{bmatrix},
\begin{bmatrix}
0&1&0&0\\
0&0&0&0\\
1&0&0&0
\end{bmatrix},
\begin{bmatrix}
1&0&0\\
0&0&0\\
0&1&0\\
0&0&0
\end{bmatrix}\right)
\end{eqnarray*}
This lacing diagram has rank array~$\rr$ and lace array~$\ss$ from
Example~\ref{ex:rank}.\qed%
\end{example}

\begin{lemma} \label{l:ranklace}
Fix a rank array $\rr$ that can occur for~$V$, and let $\ss$ be the
lace array related to~$\rr$ by~\eqref{eq:laces}.  For a lacing diagram
$\ww\in\homv$, the following are equivalent:
\begin{numbered}
\item
$\ww$ has rank array $\rr$.

\item
The quiver locus\/ $\Omega_\rr$ is the closure $\overline{\glv\cdot\ww}$
of the orbit of\/~$\glv$ through\/~$\ww$.

\item $s_{ij}$ is the number of $(i,j)$-laces in $\ww$ for all $0
\leq i \leq j \leq n$.
\end{numbered}
\end{lemma}
\begin{proof}
This is a consequence of Proposition~\ref{p:indec} plus
Lemma~\ref{l:orbits} and its proof.%
\end{proof}

Lemma~\ref{l:ranklace} identifies orbits indexed by a collection
of lacing diagrams.  It will be important for us to consider
another group action, for which each lacing diagram uniquely
determines its own orbit. As a starting point, we have the
following standard result; see \cite{Ful92} or \cite{grobGeom}.
As a matter of notation, we write $w' \leq w$ for partial permutation
matrices $w$ and~$w'$ of the same shape if their associated
permutations satisfy $\wt w' \leq \wt w$ in Bruhat order.

\begin{lemma} \label{l:matSchuborbit}
Let the product $B_-(k) \times B_+(\ell)$ of lower- and upper-triangular
Borel groups inside $\gl_k \times \gl_\ell$ act on matrices $Z \in
M_{k\ell}$ by $(b_-,b_+) \cdot Z = b_- Z b_+^{-1}$.  The matrix Schubert
variety $\ol X_w$ for a partial permutation matrix~$w$ is the closure
in~$M_{k\ell}$ of the orbit of $B_-(k) \times B_+(\ell)$ through~$w$.
Moreover, $\Xb_{w'} \subset \Xb_w$ if and only if $w' \leq w$.
\end{lemma}

The observation in Lemma~\ref{l:matSchuborbit} easily extends to
sequences of partial permutations---that is, lacing diagrams $\ww
= (w_1,\ldots,w_n)$---using the subgroup
\begin{eqnarray*}
  \BtB &=& \prod_{i=0}^n B_+(V_i) \times B_-(V_i)
\end{eqnarray*}
of the group $\glvv$ from~\eqref{glvv}, where $B_+(V_i)$ and $B_-(V_i)$
are the upper and lower triangular Borel groups in $\gl(V_i) =
\gl_{r_i}$, respectively.  The action of $\BtB$ on~$\homv$ is inherited
from the $\glvv$ action after~\eqref{glvv}.  We also need to
\mbox{consider the subgroup}
\begin{eqnarray} \label{btbv}
  \btbv &=& \prod_{i=0}^n B_+(V_i) \times_{T(V_i)} B_-(V_i) \
  \:\subset\:\ \BtB,
\end{eqnarray}
where $B_+(V_i) \times_{T(V_i)} B_-(V_i)$ consists of those pairs
\mbox{$(b_+,b_-) \in B_+(V_i) \times B_-(V_i)$} in which the diagonals
of $b_+$ and~$b_-$ are equal.

\begin{prop} \label{p:sameorbit}
The groups $\btbv$ and $\BtB$ have the same orbits on\/~$\homv$,
and each orbit contains a unique lacing diagram\/~$\ww$.  The
closure in\/~$\homv$ of the orbit through\/~$\ww$ is the product
\begin{eqnarray*}
  \OO(\ww) &=& \prod_{j=1}^n \Xb_{w_j}
\end{eqnarray*}
of matrix Schubert varieties, with\/ $\ol X_{w_j}$ lying inside the
factor\/ $\Hom(V_{j-1},V_j)$ of\/~$\homv$.
\end{prop}
\begin{proof}
By Lemma~\ref{l:matSchuborbit} the action of $\BtB$ on $\homv$ has
orbit representatives given by lacing diagrams.  Hence it suffices
to show that every orbit of the subgroup $\btbv\subset \BtB$
acting on $\homv$ contains a lacing diagram.

Note that $\btbv = T \cdot U$ for the maximal torus $T \subset \glv$
(diagonally embedded in~$\glvv$) and the unipotent group $U = \prod_i
U_+(V_i) \times U_-(V_i)$, where $U_+(V_i) \subset B_+(V_i)$ and
$U_-(V_i) \subset B_-(V_i)$ are the unipotent radicals (unitriangular
matrices).  Acting on the factor $\Hom(V_{i-1},V_i)$ by the group
$U_-(V_{i-1}) \times U_+(V_i)$ via the action of~$\glvv$
after~\eqref{glvv}, any matrix can be transformed into a monomial matrix
(one having at most one nonzero element in each row and each column) by
performing downward row elimination and rightward column elimination
with no scaling of the pivot elements.  Thus each orbit of~$U$
in~$\homv$ contains a list of monomial matrices
$\phi=(\phi_1,\dotsc,\phi_n)$.

Let us view a list of monomial matrices as a lacing diagram whose
edges have been labeled with nonzero elements of~$\kk$.  The
underlying lacing diagram $\ww$ is determined by the list of
partial permutation matrices obtained by changing each nonzero
entry of $\phi$ to a $1$, and if the coordinate $\fff
i\alpha\beta$ is nonzero on $\phi$, then its value $\fff
i\alpha\beta(\phi)$ is the label on the edge connecting the
height~$\alpha$ vertex in column~$i-1$ to the height~$\beta$
vertex in column~$i$.  It now suffices to find $\tau\in T$ with
$\tau\phi = \ww$.  We construct~$\tau$ by setting the $j^\th$
diagonal entry $\tau_i^j$ in its $i^\th$ factor~$\tau_i$ to be the
product of the edge labels that are to the left of column~$i$ on
the lace in~$\phi$ containing the
height~$j$ vertex in column~$i$.%
\end{proof}

The entry $(i,j)$ of a partial permutation matrix is called an
\bem{inversion} if there is no~$1$ due north or due west of it, i.e.,
there is no~$1$ in position $(i',j)$ for $i'<i$ or $(i,j')$ for $j'<j$.
Define the \bem{length}~$\ell(w)$ of a partial permutation matrix~$w$ to
be the number of inversions of~$w$.  Equivalently, this is the length of
any permutation of minimal length whose matrix contains~$w$ as its
northwest corner.  Given a list~$\ww=(w_1,\dotsc,w_n)$ of partial
permutation matrices, define the length of~$\ww$ to be $\ell(\ww) =
\sum_{j=1}^n \ell(w_j)$.

\begin{defn} \label{d:W}
Denote by~$W(\rr)$ the subset of lacing diagrams in~$\homv$ with
rank array $\rr$ and length $d(\rr) = \codim(\Omega_\rr)$ from
Corollary~\ref{c:codim}.
\end{defn}

\begin{example} \label{ex:lacediags}
For $\rr$ as in Example \ref{ex:rank}, the set $W(\rr)$ is given by
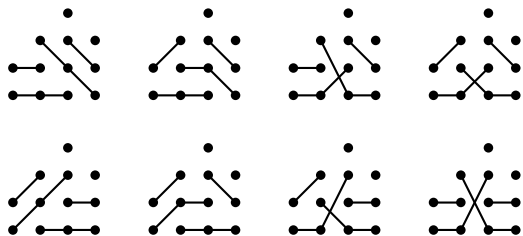
\begin{equation*}
\begin{split}
&
\pspicture[.1](0,0)(3,3)
\psdots(0,0)(0,1)(1,0)(1,1)(1,2)(2,0)(2,1)(2,2)(2,3)(3,0)(3,1)(3,2)
\psline(0,0)(1,0)(2,0)
\psline(0,1)(1,1)
\psline(1,2)(2,1)(3,0)
\psline(2,2)(3,1)
\endpspicture \qquad
\pspicture[.1](0,0)(3,3)
\psdots(0,0)(0,1)(1,0)(1,1)(1,2)(2,0)(2,1)(2,2)(2,3)(3,0)(3,1)(3,2)
\psline(0,0)(1,0)(2,0)
\psline(0,1)(1,2)
\psline(1,1)(2,1)(3,0)
\psline(2,2)(3,1)
\endpspicture \qquad
\pspicture[.1](0,0)(3,3)
\psdots(0,0)(0,1)(1,0)(1,1)(1,2)(2,0)(2,1)(2,2)(2,3)(3,0)(3,1)(3,2)
\psline(0,0)(1,0)(2,1)
\psline(0,1)(1,1)
\psline(1,2)(2,0)(3,0)
\psline(2,2)(3,1)
\endpspicture \qquad
\pspicture[.1](0,0)(3,3)
\psdots(0,0)(0,1)(1,0)(1,1)(1,2)(2,0)(2,1)(2,2)(2,3)(3,0)(3,1)(3,2)
\psline(0,0)(1,0)(2,1)
\psline(0,1)(1,2)
\psline(1,1)(2,0)(3,0)
\psline(2,2)(3,1)
\endpspicture
\\[5mm]
&
\pspicture[.1](0,0)(3,3)
\psdots(0,0)(0,1)(1,0)(1,1)(1,2)(2,0)(2,1)(2,2)(2,3)(3,0)(3,1)(3,2)
\psline(0,0)(1,1)(2,2)
\psline(0,1)(1,2)
\psline(1,0)(2,0)(3,0)
\psline(2,1)(3,1)
\endpspicture \qquad
\pspicture[.1](0,0)(3,3)
\psdots(0,0)(0,1)(1,0)(1,1)(1,2)(2,0)(2,1)(2,2)(2,3)(3,0)(3,1)(3,2)
\psline(0,0)(1,1)(2,1)
\psline(0,1)(1,2)
\psline(1,0)(2,0)(3,0)
\psline(2,2)(3,1)
\endpspicture \qquad
\pspicture[.1](0,0)(3,3)
\psdots(0,0)(0,1)(1,0)(1,1)(1,2)(2,0)(2,1)(2,2)(2,3)(3,0)(3,1)(3,2)
\psline(0,0)(1,0)(2,2)
\psline(0,1)(1,2)
\psline(1,1)(2,0)(3,0)
\psline(2,1)(3,1)
\endpspicture \qquad
\pspicture[.1](0,0)(3,3)
\psdots(0,0)(0,1)(1,0)(1,1)(1,2)(2,0)(2,1)(2,2)(2,3)(3,0)(3,1)(3,2)
\psline(0,0)(1,0)(2,2)
\psline(0,1)(1,1)
\psline(1,2)(2,0)(3,0)
\psline(2,1)(3,1)
\endpspicture
\end{split}
\end{equation*}
The codimension of this quiver locus~$\Omega_\rr$ is $7$, which will
be reflected in the number of `virtual crossings' in these diagrams,
as defined in this next section.\qed
\end{example}

\section{Minimum length lacing diagrams}

\label{sec:min}

In Theorem~\ref{t:mindiag} of this section we give a combinatorial
characterization of the lacing diagrams~$W(\rr)$ from
Definition~\ref{d:W}.

For partial permutation matrix lists (lacing diagrams) $\ww$
and~$\ww'$, write $\ww' \leq \ww$ if $w'_j \leq w_j$ in Bruhat order
for all $1 \leq j \leq n$.

To compute the length of a lacing diagram we introduce the notion of a
\bem{crossing}.  The vertex at height~$i$ (from the bottom) in
column~$k$ is written as the ordered pair~\mbox{$(i,k)$}.  For the
purposes of the following definition one may imagine that every
$(a,b)$-lace is extended so that its left endpoint is connected by an
invisible edge to an invisible vertex $(\infty,a-1)$ and its right
endpoint is connected by an invisible edge to an invisible vertex
$(\infty,b+1)$.  Two laces in a lacing diagram are said to \bem{cross}
if one contains an edge of the form $(i,k-1)\rightarrow(j,k)$ and the
other an edge of the form $(i',k-1)\rightarrow(j',k)$ such that $i>i'$
and $j<j'$.  Such a crossing shall be labeled by the triple
$(j,j',k)$.  Looking at these two edges, this crossing is classified
as follows.

Call the crossing \bem{ordinary} if neither edge is invisible;
\bem{left-virtual} if the edge touching $(j,k)$ is invisible but the
other edge is not; \bem{right-virtual} if the edge touching $(j',k)$
is invisible but the other is not; and \bem{doubly virtual} if both
edges are invisible.  Two laces have a doubly virtual crossing
(necessarily unique) if and only if one is an $(a,b)$-lace and the
other a $(b+1,c)$-lace for some $a\le b<c$.

For the next proposition and theorem, a \bem{crossing} refers to any
of the above four kinds of crossings.

\begin{prop} \label{pp:lencross}
The number of crossings in any lacing diagram\/ $\ww\in\homv$
is\/~$\ell(\ww)$.
\end{prop}
\begin{proof}
Fix $1\le k\le n$. The inversions $(i,j)$ of $w_k$ come in four flavors:
\begin{numbered}
\item
$w_k$ has a $1$ in positions $(i,j')$ and $(i',j)$ with $j'>j$ and
$i'>i$.
\item
$w_k$ has a $1$ in position $(i,j')$ with $j'>j$ but no $1$ in the
$j^\th$ column.
\item
$w_k$ has a $1$ in position $(i',j)$ with $i'>i$ but no $1$ in the
$i^\th$ row.
\item
$w_k$ has no $1$ in either the $i^\th$ row or $j^\th$ column.
\end{numbered}
Consider the two laces in $\ww$ containing the vertices $(i,k-1)$ and
$(j,k)$ respectively; they are distinct since $(i,j)$ is an inversion
of $w_k$. Then the four kinds of inversions correspond respectively to
crossings $(j,j',k)$ of these two laces that are ordinary,
left-virtual, right-virtual, and doubly virtual respectively.%
\end{proof}

It is a bit difficult to use this proposition to identify by eye
the lacing diagrams in $W(\rr)$, because it is hard to make an accurate
count of virtual crossings. The following theorem gives a more
easily checked criterion (though it still involves virtual crossings).

\begin{thm} \label{t:mindiag}
Suppose the lacing diagram\/ $\ww$ has rank array\/~$\rr$.  Then
$\ell(\ww) \geq d(\rr)$, with equality if and only if
\begin{enumerate}
\item[(C1)] no pair of laces in $\ww$ crosses more than once; and
\item[(C2)] two laces in~$\ww$ have no crossings if both start or both
end in the same column.
\end{enumerate}
\end{thm}

Theorem~\ref{t:mindiag} follows immediately from Proposition
\ref{pp:lencross} and the two following results, which may be of
independent interest.

\begin{prop} \label{pp:nonminlength}
Suppose the lacing diagram\/~$\ww$ has rank array\/~$\rr$ but fails to
satisfy (C1) or (C2).  Then there is a lacing diagram\/~$\ww'$ with
rank array\/~$\rr$ such that\/ $\ww' < \ww$.
\end{prop}
\begin{proof}
Suppose first that $\ww$ has a pair of laces~$L$ and~$L'$ that cross
at least twice.  Then undoing adjacent crossings in these two laces
produces the desired lacing diagram~$\ww'$.  More precisely, consider
a pair of crossings $(j,j',k)$ and $(j_1,j'_1,k_1)$ of~$L$ and~$L'$
with $k<k_1$ and $k_1-k$ minimal.  Let $\ww'$ be the lacing diagram
obtained from $\ww$ by changing the vertex sets of~$L$ and~$L'$ by
trading vertices in columns~$c$ such that $k\le c < k_1$ (and leaving
all other laces unchanged).  It is straightforward to verify that
$\ww'$ has rank array~$\rr$, that $w'_k<w_k$, that $w'_{k_1}<w_{k_1}$,
and that $w'_c=w_c$ for $c\not\in\{k,k_1\}$.

Suppose $\ww$ has an $(a,b)$-lace~$L$ and an $(a,b')$-lace~$L'$ that
cross (the case that two crossing laces end in the same column is
entirely similar).  Let $(j,j',k)$ be the crossing of~$L$ and~$L'$
with $k$ minimal.  As before we may change the vertex sets of~$L$
and~$L'$ by trading vertices in columns~$c$ such that $a\le c<k$,
resulting in the lacing diagram~$\ww'$.  Then $\ww'$ has rank
array~$\rr$, with $w'_k<w_k$, and $w'_c=w_c$ for $a\le c<k$.%
\end{proof}

\begin{prop} \label{pp:minlength}
If\/ $\ww$ has rank array\/ $\rr$ then\/ $\ww$ has at least $d(\rr)$
crossings.  Moreover if\/ $\ww$ satisfies (C1) and (C2) then\/ $\ww$
has exactly\/ $d(\rr)$ crossings.
\end{prop}
\begin{proof}
Consider any $(a,b)$-lace $L$ and $(a',b')$-lace~$L'$ in~$\ww$ such
that $a<a'$ and $b<b'$.  It is easy to check that the laces cross at
least once if and only if $b \geq a'-1$.  {}From this and
Corollary~\ref{c:codim} it follows that $\ww$ has at least
$d(\rr)$-crossings.  Suppose $\ww$ has rank array~$\rr$ and satisfies
(C1) and (C2).  Let $L$ be an $(a,b)$ lace and~$L'$ an $(a',b')$-lace
with $a\le a'$.  If $a=a'$ the laces cannot cross by (C2), so suppose
$a<a'$.  If $b>b'$ then the laces cannot cross, for if they did then
they would have to cross twice, violating (C1).  If $b=b'$ then again
by (C2) the laces cannot cross.  Finally, for $b<b'$, the laces cross
at least once and only if $b\ge a'-1$, and they cannot cross twice
by~(C1).  This gives precisely the desired $d(\rr)$ crossings.%
\end{proof}

\part{Degeneration of quiver loci}

\label{part:degen}
\setcounter{section}{0}
\setcounter{thm}{0}
\setcounter{equation}{0}

\section{Orbit degenerations}

\label{sec:groupdegen}

Tracing the degeneration of a group acting on a space can allow one to
degenerate the orbits of the group.  More precisely, the next result
says that under appropriate conditions, an action of a one-dimensional
flat family $\ol G$ of groups on the general fiber~$X$ of a flat
family~$\ol X$ extends to a $\ol G$-action on all of~$\ol X$.  This
implies, in particular, that the special fiber of $\ol G$ acts on the
special fiber of~$\ol X$.

\begin{prop} \label{p:flat}
Let the group scheme $\Gamma{}_{\!0}$ act on affine $n$-space $\AA^n$
(over the integers, say).  Fix a regular scheme $\ol Y$ of dimension~$1$
along with a closed point $y \in \ol Y$, and consider the family $\ol V
= \AA^n \times \ol Y$.  Suppose $\ol G$ is a flat closed subgroup scheme
of\/ $\ol \Gamma \cong \Gamma{}_{\!0} \times \ol Y$ over~$\ol Y$.  Let
$V$, $G$, and $\Gamma$ denote the restrictions of\/ $\ol V$, $\ol G$,
and $\ol \Gamma$ to schemes over $Y = \ol Y \minus y$.  If $X \subseteq
V$ is a $G$-subscheme flat over~$Y$, then there is a unique equivariant
completion of~$X$ to a closed $\ol G$-subscheme $\ol X \subseteq \ol V$
flat over~$\ol Y$.
\end{prop}
\begin{proof}
The flat family $X \subseteq V$ has a unique extension to a flat family
$\ol X$ whose total space is closed inside $\ol V$; this is standard,
taking $\ol X$ to be the closure of $\ol X$ inside~$\ol V$.  (The ideal
sheaf $\II_\ol X$ is thus $\II_X \cap \OO_\ol V$, the intersection
taking place inside~$\OO_V$.)  It remains to show that the map $\ol G
\times_\ol Y \ol X \to \ol V$ has image contained in $\ol X$ (more
formally, the map factors through the inclusion $\ol X \into \ol V$).

The product $G \times_Y X$ is a closed subscheme of $\Gamma \times_Y V$,
which is itself an open subscheme of $\ol \Gamma \times_\ol Y \ol V$.
As such we can take the closure $\ol{G \times_Y X}$ inside $\ol \Gamma
\times_\ol Y \ol V$.  The composite map
$$
  \ol{G \times_{Y\!}X}\ \too\ \ol\Gamma\times_{\ol Y} \ol V\ \too\ \ol V
$$
lands inside $\ol X$, because (i)~removing all the overbars in the above
composite map yields a map whose image is precisely the $G$-scheme~$X$,
and (ii)~$\ol X$ is closed in~$\ol V$.  It suffices to show now that
$\ol{G \times_Y X} = \ol G \times_\ol Y \ol X$ as subschemes of $\ol
\Gamma \times_\ol Y \ol V$.

Clearly $\ol{G \times_{Y\!} X} \subseteq \ol G \times_\ol Y \ol X$,
because $\ol G \times_\ol Y \ol X$ is closed in $\ol V$ and contains $G
\times_{Y\!} X$.  Moreover, $\ol G \times_\ol Y \ol X$ is flat over $\ol
Y$ (because both $\ol G$ and $\ol X$~are) and agrees with $G
\times_{Y\!}  X$ over~$Y$.  Thus $\ol G \times_\ol Y \ol X$ and $\ol{G
\times_{Y\!} X}$ both equal the unique extension of the family $G
\times_{Y\!} X \to Y$ to a flat family over $\ol Y$ that is closed
inside $\ol \Gamma \times_\ol Y \ol V$.%
\end{proof}

In our case, the group degeneration is a group scheme $\til\glv$ over
$\AA^1$ (parametrized by the coordinate~$t$ on~$\AA^1$).  $\til\glv$
is defined as follows via Gr\"obner degeneration of the diagonal
embedding
$$
  \Delta : \glv \into \glvv.
$$
Let $\bg_i = (\ggg i\alpha\beta)_{\alpha,\beta = 1}^{r_i}$ be
coordinates on $\gl(V_i)$ for each index~$i$, so the coordinate ring
\begin{eqnarray*}
  \kk[\gl(V_i)^2] &=&
  \kk[\lvec{\bg_i}\,,\rvec{\bg_i}\,][\lvec\det\,^{-1},\rvec\det\,^{-1}]
\end{eqnarray*}
is the polynomial ring with the determinants of the left and right
pointing variables inverted.  The ideal of the $i^\th$ component of the
diagonal embedding,
\begin{eqnarray*}
  I_{\!\Delta_i} &=& \<\lggg i\alpha\beta - \rggg i\alpha\beta \mid
  \alpha,\beta = 1,\ldots,r_i\>,
\end{eqnarray*}
simply sets the right and left pointing coordinates equal to each other
in the $i^\th$ com\-ponent.  Choose a weight $\omega$ on the variables
in the polynomial ring $\kk[\lvec{\aa_i}\,,\rvec{\aa_i}\,]$~so~that
$$
  \omega(\lggg i\alpha\beta) = \alpha-\beta \quad\hbox{ and }\quad
  \omega(\rggg i\alpha\beta) = \beta-\alpha.
$$
Using this weight we get a family
\begin{eqnarray*}
  \IDi &=& \<\lggg i\alpha\beta - \rggg i\alpha\beta
  t^{2\alpha-2\beta} \mid \alpha \geq \beta\>\ +\ \<\lggg i\alpha\beta
  t^{2\beta-2\alpha} - \rggg i\alpha\beta \mid \alpha \leq \beta\>
\end{eqnarray*}
of ideals over $\AA^1$, whose coordinate we call~$t$.

\begin{lemma} \label{l:flatG}
$\til{I}_{\!\Delta} = \sum_{i=0}^n \IDi$ defines a flat group
subscheme $\til\glv$ of $\glvv\times \AA^1$ whose special fiber at $t
= 0$ is the group $\glv(0) = \btbv$ from~\eqref{btbv}.
\end{lemma}
\begin{proof}
The special fiber $\glv(0)$ over $t = 0$ is defined by $\ID(0)$, which
coincides with
\begin{eqnarray*}
  \IN_\omega(\ID)
&=&
  \sum_{i=0}^n \Bigl(
    \<\lggg i\alpha\beta \mid \alpha > \beta\>
  + \<\rggg i\alpha\beta \mid \alpha < \beta\>
  + \<\lggg i\alpha\alpha-\rggg i\alpha\alpha \mid \alpha=1,\ldots,r_i\>
  \Bigr)
\end{eqnarray*}
by definition of the initial ideal $\IN_\omega$ for the weight~$\omega$
\cite[Chapter~15]{Eis}.  Thus $\til\glv$ is a flat family (specifically,
a Gr\"obner
degeneration) over~$\AA^1$.  Flatness still holds after inverting the
left and right pointing determinants.

It remains to show that $\til\glv$ is a subgroup scheme over $\AA^1$,
i.e.\ that fiberwise multiplication in $\glvv\times \AA^1$ preserves
$\til\glv \subset \glvv\times\AA^1$.  This can be done
set-theoretically because (for instance) all fibers of $\til\glv$ are
geometrically reduced.

The fiber $\glv(0)$ over $t = 0$ is the subgroup $\btbv \subset \glvv$.
On the other hand, the union of all fibers with $t \neq 0$ is isomorphic
to $\glv \times (\AA^1 \minus \{0\})$ as a group subscheme of $\glvv
\times (\AA^1 \minus \{0\})$, the isomorphism being
\begin{eqnarray} \label{eq:tau}
  (\lgam{\,}\,,\,\rgam{\,}) \times \tau &\mapsto&
  (\tau_\spot\lgam{\,}\tau_\spot^{-1}\,,
  \,\tau_\spot^{-1}\rgam{\,}\tau_\spot)
\end{eqnarray}
where $\tau_\spot \in \tv$ has $T(V_i)$ component $\tau_i$ with
diagonal entries $1,\tau,\tau^2,\ldots,\tau^{r_i-1}$.  Indeed, the
matrices $(\tau_\spot\lgam{\,}\tau_\spot^{-1}\,,\,
\tau_\spot^{-1}\rgam{\,}\tau_\spot)$ are precisely those satisfying
the relations defining the ideal $\til{I}_{\!\Delta}$, which can be
rewritten more symmetrically~as
\begin{eqnarray*}
  \ID(t \neq 0) &=& \<t^{-\alpha} \, \lggg i\alpha\beta \, t^\beta
  - t^\alpha \, \rggg i\alpha\beta \, t^{-\beta} \mid i = 0,\ldots,n\>
\end{eqnarray*}
away from $t=0$.
\end{proof}

\section{Quiver degenerations}

\label{sec:quivdegen}

The orbits we wish to degenerate are of course the quiver loci.
Recall the notation from Section~\ref{sec:quiver} regarding the
indexing on variables $\fff i\alpha\beta$ in the coordinate
ring~$\kk[\homv]$.

\begin{defn}
For $p(\ff) \in \kk[\homv]$ and a coordinate $t$ on $\AA^1 = \kk$, let
\begin{eqnarray*}
  p_t(\ff) &=& \frac{p(t^{\alpha+\beta}\fff i\alpha\beta)}{\hbox{as
  much $t$ as possible}}.
\end{eqnarray*}
Set $\til{I}_\rr = \<p_t(\ff) \mid p(\ff) \in I_\rr\> \subseteq
\kk[\homv\times \AA^1]$, and let $\til\Omega_\rr \subseteq
\homv\times \AA^1$ be the zero scheme of~$\til{I}_\rr$.  For $\tau\in
\kk$ let $\Omega_\rr(\tau)\subset \homv$ be the zero scheme of the
ideal $I_\rr(\tau)\subset\kk[\homv]$ obtained from $\til{I}_\rr$ by
setting $t=\tau$.  The special fiber $\Omega_\rr(0)\subset \homv$ is
called the \bem{quiver degeneration}, and we refer to $\til\Omega_\rr$
as the \bem{family} of $\Omega_\rr(0)$.
\end{defn}

\begin{lemma} \label{l:flatO}
The family\/ $\til\Omega_\rr \subseteq \homv \times \AA^1$ is flat over
$\AA^1$.
\end{lemma}
\begin{proof}
The polynomial $p_t(\ff)$ is obtained from $p(\ff)$ by homogenizing with
respect to the weight function on $\kk[\homv]$ that assigns weight
$-\alpha-\beta$ to the variable $\fff i\alpha\beta$.  The resulting
initial ideal is~$I_\rr(0)$, so \cite[Theorem~15.17]{Eis} applies.%
\end{proof}

Note that $I_\rr(0)$ is rarely a monomial ideal.  Thus $\Omega_\rr(1)
\goesto \Omega_\rr(0)$ is only a ``partial'' Gr\"obner degeneration of
$\Omega_\rr = \Omega_\rr(1)$.

\begin{prop} \label{p:action}
The flat family $\til\Omega_\rr$ over $\AA^1$ is acted upon
(fiberwise over $\AA^1$) by the group scheme $\til\glv$ over $\AA^1$.
\end{prop}
\begin{proof}
This will follow from Proposition~\ref{p:flat} applied to
$$
  \AA^n = V, \quad \Gamma{}_{\!0} = \glvv, \quad \ol Y = \AA^1,
  \quad y = 0, \quad \ol G = \til\glv, \quad \hbox{and} \quad \ol X =
  \til\Omega_\rr
$$
by Lemmas~\ref{l:flatO} and~\ref{l:flatG}, as soon as we show that
$\til\glv(\tau \neq 0)$ acts on $\til\Omega_\rr(\tau \neq 0)$, these
being the restrictions of $\til\glv$ and $\til\Omega_\rr$ to families
over $\AA^1 \minus \{0\}$.

Resume the notation from the last paragraph of the proof of
Lemma~\ref{l:flatG}.  Observe that $\Omega_\rr \times (\AA^1 \minus
\{0\}) \to \til\Omega_\rr(\tau \neq 0)$ is an isomorphism of varieties
under the map
\begin{eqnarray} \label{eq:phi}
  (\phi,\tau)\ =\ ((\phi_i)_{i=1}^n,\tau) &\mapsto& (\tau_{i-1}^{-1}
  \phi_i \tau_i^{-1})\ =:\ \tau_{\spot}^{-1} \phi \tau_\spot^{-1}.
\end{eqnarray}
Now combine Eq.~\eqref{eq:tau} and~\eqref{eq:phi},
\begin{eqnarray*}
  (\tau_\spot\lgam{\,}\tau_\spot^{-1}\,,\,
  \tau_\spot^{-1}\rgam{\,}\tau_\spot) \cdot (\tau_{\spot}^{-1}
  \phi\,\tau_\spot^{-1})
&=&
  (\tau_{i-1}^{-1} \rgam{i-1} \phi_i \lgam{i}^{-1} \tau_i^{-1})\\
&=&
  \tau_{\spot}^{-1}[(\lgam{\,},\rgam{\,}) \cdot \phi]
  \tau_\spot^{-1},
\end{eqnarray*}
and use that $(\lgam{\,},\rgam{\,}) \cdot \phi \in \Omega_\rr$ to get the
desired result.
\end{proof}

Theorem~\ref{t:degen}, the main result of Section~\ref{part:degen},
relates the quiver degeneration $\Omega_\rr(0)$ to matrix Schubert
varieties as in Proposition~\ref{p:sameorbit}.

\begin{thm} \label{t:degen}
Each irreducible component of\/ $\Omega_\rr(0)$ is a possibly
nonreduced product\/~$\OO(\ww)$ of matrix Schubert varieties for a
lacing diagram\/~$\ww\in\homv$.
\end{thm}
\begin{proof}
By Proposition~\ref{p:action}, the special fiber $\Omega_\rr(0)$ is
stable under the action of $\glv(0)=\btbv$ on~$\homv$, and hence is a
union of orbit closures for this group.  By
Proposition~\ref{p:sameorbit} these orbit closures are direct products
of matrix Schubert varieties.%
\end{proof}

\begin{remark} \label{r:cyclic}
The methods of the current and the previous sections can be applied with
no essential changes to the more general case of nilpotent cyclic
quivers.
\end{remark}

\section{Multidegrees of quiver degenerations}

\label{sec:multidegen}

Knowing multidegrees of quiver degenerations gives us multidegrees of
quiver loci:

\begin{lemma} \label{l:multideg}
The quiver degeneration $\Omega_\rr(0)$ has the same multidegree
as~$\Omega_\rr$.
\end{lemma}
\begin{proof}
By Lemma~\ref{l:flatO} we may apply `Degeneration' in
\cite[Theorem~1.7.1]{grobGeom}.%
\end{proof}

Now we can conclude from Theorem~\ref{t:degen} that the quiver
polynomial has a positive expression in terms of products of double
Schubert polynomials.  This will provide a ``lower bound'' on the
quiver polynomial $\QQ_\rr(\xx-\oxx)$, in Proposition~\ref{p:lace}.

\begin{cor} \label{c:quivSchub}
If $\SS_\ww(\xx-\oxx)$ for a lacing diagram~$\ww$ denotes the
corresponding product $\prod_{j=1}^n \SS_{w_j}(\xx^{j-1}-\xx^j)$
of double Schubert polynomials, then
\begin{eqnarray*}
  \QQ_\rr(\xx-\oxx) &=& \sum_{\ww\in\,\homv} c_\ww(\rr)\,
  \SS_\ww(\xx-\oxx),
\end{eqnarray*}
where $c_\ww(\rr)$ is the multiplicity of the component\/~$\OO(\ww)$
in\/~$\Omega_\rr(0)$.
In particular each $c_\ww(\rr)\geq 0$.
\end{cor}
\begin{proof}
In light of Lemma~\ref{l:multideg} and the additivity of
multidegrees under taking unions of components of equal dimension
\cite[Theorem~1.7.1]{grobGeom}, the quiver polynomial
$\QQ_\rr(\xx-\oxx)$, which is by definition the multidegree of the
quiver locus~$\Omega_\rr$, breaks up as a sum over lacing diagrams
by Theorem~\ref{t:degen}.  The summands are products of double
Schubert
polynomials by Theorem~\ref{t:dubSchub}.%
\end{proof}

To make Corollary~\ref{c:quivSchub} more explicit, it is necessary to
identify which lacing diagrams~$\ww$ give components $\OO(\ww)$, and
to find their multiplicities.  We shall accomplish this goal in
Theorem~\ref{t:components}, whose proof requires the less precise
statement in Proposition~\ref{p:lace}, below, as a stepping stone (in
the end, it turns out to have identified all the components, not just
some of them).  First, let us observe that the quiver
degeneration~$\Omega_\rr(0)$ contains every lacing diagram with rank
array~$\rr$.

\begin{lemma} \label{l:paths}
If the lacing diagram $\ww$ has rank array\/~$\rr$, then $\ww \in
\Omega_\rr(\tau)$ for all~$\tau\!\in\!\kk$.
\end{lemma}
\begin{proof}
Let $\ww\in\homv$ have ranks~$\rr$.  It suffices to show that $\ww\in
\Omega_\rr(\tau)$ for $\tau\in\AA^1\minus0$; indeed, it follows that it
holds for all~$\tau$ because the family $\til\Omega_\rr$ is closed
in~$\homv \times \AA^1$.

By Lemma~\ref{l:ranklace}, $\ww$ lies inside~$\Omega_\rr$.  Using
the notation of \eqref{eq:tau} and~\eqref{eq:phi}, for every
$\tau\in\AA^1\minus0$ the list $\tau_\spot \ww \tau_\spot$ of
monomial matrices also lies in~$\Omega_\rr$.  Now consider the
curve
$$
  (\tau_\spot \ww \tau_\spot) \times \tau
$$
inside $\Omega_\rr \times (\AA^1 \minus 0)$.  Under the map given in
\eqref{eq:phi}, this curve is sent to the constant curve $\ww\times
(\AA^1\minus0)$ inside $\til\Omega_\rr(\tau \neq 0)$, proving the
lemma.
\end{proof}

\begin{prop} \label{p:lace}
Among the irreducible components of the quiver degeneration
$\Omega_\rr(0)$ are the orbit closures $\OO(\ww)$ for\/ $\ww\in W(\rr)$.
In particular these $c_\ww(\rr)$ are at least $1$.
\end{prop}
\begin{proof} Let $\ww\in W(\rr)$. Then $\OO(\ww)\subset\Omega_\rr(0)$,
by Lemma~\ref{l:paths} along with Propositions~\ref{p:sameorbit}
and~\ref{p:action}.  Consider a component of $\Omega_\rr(0)$
containing $\OO(\ww)$; it must have the form $\OO(\ww')$ by
Theorem~\ref{t:degen}.  But $\OO(\ww)\subset\OO(\ww')$ implies that
$w'_j\le w_j$ in the strong Bruhat order for all~$j$, by
Lemma~\ref{l:matSchuborbit}.  Thus $\ell(\ww)\ge\ell(\ww')$.  But
$\ell(\ww)$ is the codimension of $\OO(\ww)$ in $\homv$ for all~$\ww$,
and $d(\rr)$ is the codimension of the flat degeneration
$\Omega_\rr(0)$ of $\Omega_\rr$ in~$\homv$.  It follows that
$\ww'=\ww$.
\end{proof}

\section{Rank stability of components}

\label{sec:stability}

Later, in Section~\ref{part:pos}, it will be crucial to understand how
the components of~$\Omega_\rr(0)$ and their multiplicities behave when
a constant~$m$ is added to all the ranks in~$\rr$.

For notation, given $m \in \NN$ and a partial permutation~$w$, let
$m+w$ be the new partial permutation obtained by letting~$w$ act on $m
+ \ZZ_{> 0} = \{m+1,m+2,\ldots\}$ instead of $\ZZ_{> 0} =
\{1,2,\ldots\}$ in the obvious manner.  In particular $m+w$ fixes
$1,2,\dotsc,m$.  For a list $\ww = (w_1,\ldots,w_n)$ of partial
permutations, set $m + \ww = (m+w_1,\ldots,m+w_n)$.  If $m$ is a
nonnegative integer and~$\rr = (r_{ij})$ is a rank array, let $m +
\rr$ be the new rank array obtained by adding~$m$ to every
rank~$r_{ij}$ in the array~$\rr$.  This transformation is accomplished
simply by adding~$m$ to the bottom entry~$s_{0n}$ in the lace
array~$\ss$ appearing in~\eqref{eq:laces}.

In dealing with stability of components of quiver
degenerations~$\Omega_\rr(0)$, we take our cue from the
set~$W(\rr)$, whose behavior in response to uniformly increasing
ranks we deduce immediately from Theorem~\ref{t:mindiag}.

\begin{cor} \label{c:lace+1}
There is a bijection $W(\rr)\rightarrow W(m+\rr)$ given by adding
$m$ $(0,n)$-laces along the bottom.
\end{cor}

The corresponding statement for the components of~$\Omega_\rr(0)$ is
also true; moreover, the multiplicities of these components are also
stable.

\begin{prop} \label{p:stability}
The orbit closure $\OO(\ww)$ is a component of the quiver
degeneration~$\Omega_\rr(0)$ if and only if~$\OO(m+\ww)$ is a
component of the quiver degeneration~$\Omega_{m+\rr}(0)$, and the
multiplicities at their generic points are equal.
\end{prop}
\begin{proof}
It is enough by induction to prove the case~$m=1$.  Let $V_j' = \kk
\oplus V_j$ for $j = 0,\ldots,n$, write $V' = V_0' \oplus \cdots
\oplus V_n'$, and set $\homvp = \prod_{j=1}^n \Hom(V_{j-1}',V_j')$.
Define $V'_\row$ to be the subset of~$\homvp$ consisting of matrix
lists~$\phi'$ in which the only nonzero entries in each
factor~$\phi_i'$ may lie along its top row.  Similarly, define the
group~$U'$ to consist of those transformations $(u_0',\ldots,u_n')$ in
$\prod_{i=0}^n \gl(V_i')$ each of whose factors~$u_i'$ has $1$'s on
the diagonal and zeros outside the left column and main diagonal.  To
avoid notational nightmares, we identify~$\homv$ with the subset
of~$\homvp$ consisting of matrix lists~$\phi'$ in which the left
column and top row of every $\phi_i'$ is~zero.
\begin{lemma} \label{l:dense}
The sum $U'V'_\row + \Omega_\rr$ taken inside~$\homvp$ is dense
inside\/~$\Omega_{1+\rr}$.%
\end{lemma}
\begin{proof}
There are two claims here: that $U'V'_\row+\Omega_\rr$ is contained
in~$\Omega_{1+\rr}$, and that their dimensions are equal.  For
containment, observe that $V'_\row+\Omega_\rr$ lies
inside~$\Omega_{1+\rr}$ whenever $\phi$ lies in~$\Omega_\rr$ (this is
elementary linear algebra).
%
Then use the fact that, given a fixed quiver representation $\phi' \in
\homvp$, the image in~$V_{i+1}'$ of each component~$\phi_i'$ remains
unchanged when $\phi_i'$ is multiplied on the left by the $i^\th$
factor $u_i'$ of an element in~$U'$.  The containment follows because
$U'(V'_\row+\phi) = U'V'_\row+\phi$.

For the dimension count, use that the codimension of~$\Omega_\rr$
inside~$\homv$ equals the codimension of~$\Omega_{1+\rr}$
inside~$\homvp$, together with the fact that $\dim(\homvp) -
\dim(\homv) = \dim(U') + \dim(V'_\row) = \dim(U'V'_\row)$.%
\end{proof}

Returning to the proof of Proposition~\ref{p:stability}, the idea is
to show how the degeneration of $U'V'_\row + \Omega_\rr$ breaks up as
the sum of $U'V'_\row$ with the family $\til\Omega_\rr$, so the
limiting sum $U'V'_\row + \Omega_\rr(0)$ is dense
in~$\Omega_{1+\rr}(0)$.

Let $\tau_\spot' = (\tau_0'\ldots,\tau_n')$ for each nonzero $\tau \in
\kk$ be the list of diagonal matrices with $\tau_i'$ having
$(\tau^{-1},1,\tau,\tau^2,\ldots,\tau^{r_i-1})$ down its diagonal.
By~\eqref{eq:phi}, the fiber $\Omega_{1+\rr}(\tau)$ of the family
$\til\Omega_{1+\rr}$ over $t = \tau$ is given by
$(\tau_\spot')^{-1}(\Omega_{1+\rr})(\tau_\spot')^{-1}$, at least when
$\tau \neq 0$.  We have used that fact that the family
$\til\Omega_{1+\rr}$ is stable under a global scaling by~$\tau^{-1}$
in order to set $\tau_i' =
\diag(\tau^{-1},1,\tau,\tau^2,\ldots,\tau^{r_i-1})$ rather than
$\tau_i' = \diag(1,\tau,\tau^2,\ldots,\tau^{r_i-1},\tau^{r_i})$.

Using Lemma~\ref{l:dense}, we conclude that
$(\tau_\spot')^{-1}(U'V'_\row + \Omega_\rr)(\tau_\spot')^{-1}$ is
dense inside $\Omega_{1+\rr}(\tau)$ for all nonzero~$\tau \in \kk$.
Now calculate
\begin{eqnarray*}
  (\tau_\spot')^{-1}(U'V'_\row + \Omega_\rr)(\tau_\spot')^{-1}
&=&
  (\tau_\spot')^{-1}U'V'_\row(\tau_\spot')^{-1} +
  (\tau_\spot')^{-1}\Omega_\rr(\tau_\spot')^{-1}
\\
&=&
  U'V'_\row + \Omega_\rr(\tau),
\end{eqnarray*}
where $\Omega_\rr(\tau)$ is identified as a subvariety of~$\homv$
inside $\homvp$.  To justify removing the $(\tau_\spot')^{-1}$ factors
from~$U'V'_\row$, we have used that
\begin{eqnarray}
  (\tau_\spot')^{-1}(U'V'_\row)(\tau_\spot')^{-1}
&=&
  \Big((\tau_\spot')^{-1}U'\tau_\spot'\Big) \,
  \Big((\tau_\spot')^{-1}V'_\row(\tau_\spot')^{-1}\Big) \nonumber
\\
&=&
  U'V'_\row, \label{eq:UV}
\end{eqnarray}
which holds because the $\tau$'s in each big parenthesized group only
result in a scaling of each nonzero (and in the case of~$U'$,
nondiagonal) entry.  The $t = 0$ limit of the dense subfamily
$U'V'_\row + \til\Omega_\rr$ in~$\til\Omega_{1+\rr}$ is still dense in
the limit fiber~$\Omega_{1+\rr}(0)$, and equals $U'V'_\row +
\Omega_\rr(0)$ by~\eqref{eq:UV}.

It remains to identify the components of~$U'V'_\row + \Omega_\rr(0)$
and check multiplicities.  Let~$W'_\row$ be the subvariety of quiver
representations~$\phi_i' \in V'_\row$ in which every northwest entry
of every factor $\phi_i'$ is nonzero.  Then $U'W'_\row$ is dense
in~$U'V'_\row$, and we make the further claim that the map $U'W'_\row
\times \homv \to U'W'_\row + \homv$ is an isomorphism
scheme-theoretically.  Indeed, if $\phi' = u'v'+\phi$ lies in
$U'W'_\row + \homv$, then~$u_i'$ can be recovered as the unique matrix
in the $i^\th$~component of~$U'$ that subtracts the appropriate
multiple of the top row of~$\phi_i'$ to cancel the left entry in each
row below the top row, and subsequently $v_i'$ is the top row.
Subtraction of~$u'v'$ is the scheme-theoretic inverse to $U'W'_\row
\times \homv \to U'W'_\row+\homv$.

The scheme isomorphism just constructed implies that the product
\mbox{$U'W'_\row \times \Omega_\rr(0)$} of an open subvariety
$U'W'_\row$ of a vector space with the quiver degeneration for rank
array~$\rr$ is isomorphic to a dense subscheme of~$\Omega_{1+\rr}(0)$.
This is enough to conclude that the components of maximal dimension in
$\Omega_{1+\rr}(0)$ are in bijection with the components of maximal
dimension in~$\Omega_\rr(0)$, and furthermore that the multiplicities
at their generic points in corresponding components match.  Finally,
it is easy to verify that the list $1+\ww$ of partial permutation
matrices lies in $U'W'_\row+\Omega_\rr(\tau)$ whenever~$\ww$ lies
in~$\Omega_\rr(\tau)$, so that $U'W'_\row+\OO(\ww)$ is dense
in~$\OO(1+\ww)$.
\end{proof}

\part{Pipe dreams for Zelevinsky permutations}

\label{part:pipe}
\setcounter{section}{0}
\setcounter{thm}{0}
\setcounter{equation}{0}

\section{Pipe dream formula for double quiver polynomials}

\label{sec:pipe2quiver}

Schubert polynomials possess a beautiful combinatorial description,
proved independently in \cite{BJS,FS}, in terms of diagrams defined
in~\cite{FK96} that we call `reduced pipe dreams', although they have
also been called `rc-graphs' in the literature \cite{BB,Len02}.  These
diagrams are the main characters in our next combinatorial formula,
Theorem~\ref{t:QQrr}, which---although also of independent
interest---is an essential ingredient for proving our remaining
formulae.  {}From a purely combinatorial perspective, we find the
direct connection between reduced pipe dreams for Zelevinsky
permutations and lacing diagrams particularly appealing
(Section~\ref{sec:pipe2lace}).

Consider a square grid of size \mbox{$\ell \times \ell$}, with the box
in row~$q$ and column~$p$ labeled $(q,p)$, as in \mbox{$\ell \times
\ell$} matrix.  If each box in the grid is covered with a square tile
containing either $\textcross$ or $\textelbow$, then one can think of
the tiled grid as a network of pipes.

\begin{defn} \label{d:pipe}
A \bem{pipe dream} is a subset of the \mbox{$\ell \times \ell$} grid,
identified as the set of crosses in a tiling by \bem{crosses}
$\textcross$ and \bem{elbow joints} $\textelbow$.  A pipe dream is
\bem{reduced} if each pair of pipes crosses at most once.  The set
$\rp(v)$ of reduced pipe dreams for a permutation $v \in S_\ell$ is
the set of reduced pipe dreams~$D$ such that the pipe entering row~$q$
flows northeast to exit from column~$v(q)$.%
\end{defn}

We usually draw crossing tiles as some sort of cross, either \cross or
`$\textcross$', although we shall frequently use an asterisk~$*$ to
denote a~$\textcross$ strictly above the block superantdiagonal.
We~often leave elbow tiles blank or denote them by dots, to make the
diagrams less cluttered.  We shall only be interested in pipe dreams
whose $\textcross$ entries lie in the triangular region strictly above
the main antidiagonal (in spots $(q,p)$ with $q+p \leq \ell$), with
elbow joints elsewhere in the square grid of size~$\ell$.

\begin{example} \label{ex:pipe}
Using $\ell = 12$, a typical reduced pipe dream for the Zelevinsky
permutation~$v$ in Example~\ref{ex:zel} ($v$ is written down the left
side) looks like
$$
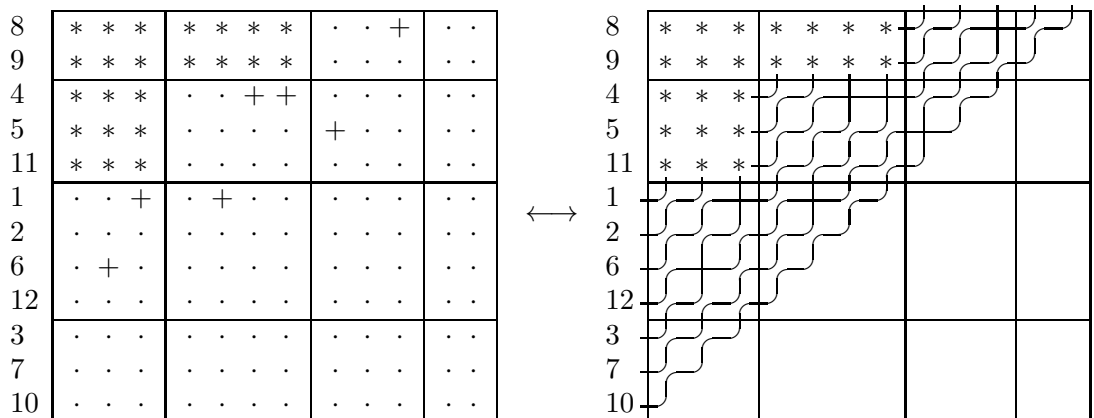
\begin{array}{@{}l|c@{\ }c@{\ }c|c@{\ }c@{\ }c@{\ }c|c@{\ }c@{\ }c|c@{\ }c|@{}}
\multicolumn{13}{@{}c@{}}{}\\[-5ex]
\mcc{}&\ph&\ph&\mcc{\ph}&\ph&\ph&\ph&\mcc{\ph}&\ph&\ph&\mcc{\ph}&\ph&\mcc{\ph}
\\\cline{2-13}
 8 &\sr&\sr&\sr&\sr&\sr&\sr&\sr&\cd&\cd& + &\cd&\cd\:
\\
 9 &\sr&\sr&\sr&\sr&\sr&\sr&\sr&\cd&\cd&\cd&\cd&\cd\:
\\\cline{2-13}
 4 &\sr&\sr&\sr&\cd&\cd& + & + &\cd&\cd&\cd&\cd&\cd\:
\\
 5 &\sr&\sr&\sr&\cd&\cd&\cd&\cd& + &\cd&\cd&\cd&\cd\:
\\
 11&\sr&\sr&\sr&\cd&\cd&\cd&\cd&\cd&\cd&\cd&\cd&\cd\:
\\\cline{2-13}
 1 &\cd&\cd& + &\cd& + &\cd&\cd&\cd&\cd&\cd&\cd&\cd\:
\\
 2 &\cd&\cd&\cd&\cd&\cd&\cd&\cd&\cd&\cd&\cd&\cd&\cd\:
\\
 6 &\cd& + &\cd&\cd&\cd&\cd&\cd&\cd&\cd&\cd&\cd&\cd\:
\\
 12&\cd&\cd&\cd&\cd&\cd&\cd&\cd&\cd&\cd&\cd&\cd&\cd\:
\\\cline{2-13}
 3 &\cd&\cd&\cd&\cd&\cd&\cd&\cd&\cd&\cd&\cd&\cd&\cd\:
\\
 7 &\cd&\cd&\cd&\cd&\cd&\cd&\cd&\cd&\cd&\cd&\cd&\cd\:
\\
 10&\cd&\cd&\cd&\cd&\cd&\cd&\cd&\cd&\cd&\cd&\cd&\cd\:
\\\cline{2-13}
\end{array}
\ \ \longleftrightarrow\ \
\begin{array}{@{}l|ccc|cccc|ccc|cc|@{}}
\multicolumn{13}{@{}c@{}}{}\\[-5ex]
\mcc{}&\ &\ &\mcc{\ \,\ \,\ }&\ &\ &\ &\mcc{\ \,\ \,\ }&\ &\ &\mcc{\ \,\ \,\ }&\ &\mcc{\ \,\ \,\ }
\\\cline{2-13}
 8 &\sr&\sr&\sr&\sr&\sr&\sr&\sr&\jr&\jr&\+ &\jr&\je
\\
 9 &\sr&\sr&\sr&\sr&\sr&\sr&\sr&\jr&\jr&\jr&\je&
\\\cline{2-13}
 4 &\sr&\sr&\sr&\jr&\jr&\+ &\+ &\jr&\jr&\je&   &
\\
 5 &\sr&\sr&\sr&\jr&\jr&\jr&\jr&\+ &\je&   &   &
\\
 11&\sr&\sr&\sr&\jr&\jr&\jr&\jr&\je&   &   &   &
\\\cline{2-13}
 1 &\jr&\jr&\+ &\jr&\+ &\jr&\je&   &   &   &   &
\\
 2 &\jr&\jr&\jr&\jr&\jr&\je&   &   &   &   &   &
\\
 6 &\jr&\+ &\jr&\jr&\je&   &   &   &   &   &   &
\\
 12&\jr&\jr&\jr&\je&   &   &   &   &   &   &   &
\\\cline{2-13}
 3 &\jr&\jr&\je&   &   &   &   &   &   &   &   &
\\
 7 &\jr&\je&   &   &   &   &   &   &   &   &   &
\\
 10&\je&   &   &   &   &   &   &   &   &   &   &
\\\cline{2-13}
\end{array}
$$
when the $*$'s remain as in the diagram.  Although each $*$~represents
a $\textcross$ in every pipe dream for~$v$, the~$*$'s will be just as
irrelevant here as they were for the diagram.\qed%
\end{example}

Each pipe dream~$D$ yields a monomial $(\xx-\oyy)^D$, defined as the
product over all $\textcross$ entries in~$D$ of $(x_+ - y_+)$, where
$x_+$ sits at the left end of the row containing~$\textcross$, and
$y_+$ sits atop the column containing~$\textcross$.  At this level of
generality, treat $\xx$ and~$\yy$ simply as alphabets (the sequence
structure is not important yet).  Here is the fundamental result
relating pipe dreams to Schubert polynomials; its `single' form was
proved independently in \cite{BJS,FS}, but it is the `double' form
that we need.

\begin{thm}\cite{FK96} \label{t:BJS}
The double Schubert polynomial for~$v$ equals the sum
\begin{eqnarray*}
  \SS_v(\xx-\oyy) &=& \sum_{D \in \rp(v)} (\xx-\oyy)^D
\end{eqnarray*}
of all monomials associated to reduced pipe dreams for~$v$.
\end{thm}

In~pictures of pipe dreams~$D$ for Zelevinsky permutations, we retain
the notation from Section~\ref{part:zel} regarding block matrices of
total size~$d$.  Also, don't count the~`$*$'~entries above the block
superantidiagonal of~$D$ as crosses.  That is, for the purposes of
seeing where crosses lie, we shall use the pipe dream $D \minus
D_\homv$ instead of~$D$, where~$D_\homv$ is identified here with the
unique reduced pipe dream consisting entirely of~$*$~locations (the
reduced pipe dream for the dominant Zelevinsky permutation
of~$\homv$).

\begin{example} \label{ex:monomial}
The pipe dream from Example~\ref{ex:pipe} has row and column labels
\begin{rcgraph}
\begin{array}{@{}l|c@{\ }c@{\ }c|c@{\ }c@{\ }c@{\ }c|c@{\ }c@{\ }c|c@{\ }c|@{}}
\multicolumn{13}{@{}c@{}}{}\\[-3ex]
\mcc{}&\ph&\ph&\mcc{\ph}&\ph&\ph&\ph&\mcc{\ph}&\ph&\ph&\mcc{\ph}&\ph&\mcc{\ph}
\\[-2ex]
\mcc{}&\!y^3_1&y^3_2&\mcc{\!y^3_3}&y^2_1&y^2_2&y^2_3&\mcc{\!y^2_4}&y^1_1&y^1_2&\mcc{\!y^1_3}&y^0_1&\mcc{\!y^0_2}
\\\cline{2-13}
x^0_1&
  \sr&\sr&\sr&\sr&\sr&\sr&\sr&\cd&\cd& + &\cd&\cd\:
\\
x^0_2&
  \sr&\sr&\sr&\sr&\sr&\sr&\sr&\cd&\cd&\cd&\cd&\cd\:
\\\cline{2-13}
x^1_1&
  \sr&\sr&\sr&\cd&\cd& + & + &\cd&\cd&\cd&\cd&\cd\:
\\
x^1_2&
  \sr&\sr&\sr&\cd&\cd&\cd&\cd& + &\cd&\cd&\cd&\cd\:
\\
x^1_3&
  \sr&\sr&\sr&\cd&\cd&\cd&\cd&\cd&\cd&\cd&\cd&\cd\:
\\\cline{2-13}
x^2_1&
  \cd&\cd& + &\cd& + &\cd&\cd&\cd&\cd&\cd&\cd&\cd\:
\\
x^2_2&
  \cd&\cd&\cd&\cd&\cd&\cd&\cd&\cd&\cd&\cd&\cd&\cd\:
\\
x^2_3&
  \cd& + &\cd&\cd&\cd&\cd&\cd&\cd&\cd&\cd&\cd&\cd\:
\\
x^2_4&
  \cd&\cd&\cd&\cd&\cd&\cd&\cd&\cd&\cd&\cd&\cd&\cd\:
\\\cline{2-13}
x^3_1&
  \cd&\cd&\cd&\cd&\cd&\cd&\cd&\cd&\cd&\cd&\cd&\cd\:
\\
x^3_2&
  \cd&\cd&\cd&\cd&\cd&\cd&\cd&\cd&\cd&\cd&\cd&\cd\:
\\
x^3_3&
  \cd&\cd&\cd&\cd&\cd&\cd&\cd&\cd&\cd&\cd&\cd&\cd\:
\\\cline{2-13}
\end{array}
\quad\begin{array}{c}\\[1ex]=\end{array}\quad
\begin{array}{@{}l|c@{\ }c@{\ }c|c@{\ }c@{\ }c@{\ }c|c@{\ }c@{\ }c|c@{\ }c|@{}}
\multicolumn{13}{@{}c@{}}{}\\[-3ex]
\mcc{}&\ph&\ph&\mcc{\ph}&\ph&\ph&\ph&\mcc{\ph}&\ph&\ph&\mcc{\ph}&\ph&\mcc{\ph}
\\[-2ex]
\mcc{}&\dt_1&\dt_2&\mcc{\!\dt_3}&\dc_1&\dc_2&\dc_3&\mcc{\!\dc_4}&\db_1&\db_2&\mcc{\!\db_3}&\da_1&\mcc{\!\da_2}
\\\cline{2-13}
a_1&
  \sr&\sr&\sr&\sr&\sr&\sr&\sr&\cd&\cd& + &\cd&\cd\:
\\
a_2&
  \sr&\sr&\sr&\sr&\sr&\sr&\sr&\cd&\cd&\cd&\cd&\cd\:
\\\cline{2-13}
b_1&
  \sr&\sr&\sr&\cd&\cd& + & + &\cd&\cd&\cd&\cd&\cd\:
\\
b_2&
  \sr&\sr&\sr&\cd&\cd&\cd&\cd& + &\cd&\cd&\cd&\cd\:
\\
b_3&
  \sr&\sr&\sr&\cd&\cd&\cd&\cd&\cd&\cd&\cd&\cd&\cd\:
\\\cline{2-13}
c_1&
  \cd&\cd& + &\cd& + &\cd&\cd&\cd&\cd&\cd&\cd&\cd\:
\\
c_2&
  \cd&\cd&\cd&\cd&\cd&\cd&\cd&\cd&\cd&\cd&\cd&\cd\:
\\
c_3&
  \cd& + &\cd&\cd&\cd&\cd&\cd&\cd&\cd&\cd&\cd&\cd\:
\\
c_4&
  \cd&\cd&\cd&\cd&\cd&\cd&\cd&\cd&\cd&\cd&\cd&\cd\:
\\\cline{2-13}
d_1&
  \cd&\cd&\cd&\cd&\cd&\cd&\cd&\cd&\cd&\cd&\cd&\cd\:
\\
d_2&
  \cd&\cd&\cd&\cd&\cd&\cd&\cd&\cd&\cd&\cd&\cd&\cd\:
\\
d_3&
  \cd&\cd&\cd&\cd&\cd&\cd&\cd&\cd&\cd&\cd&\cd&\cd\:
\\\cline{2-13}
\end{array}
\end{rcgraph}
where we use alphabets $\aa,\bb,\cc,\dd$ interchangeably with
$\xx^0,\xx^1,\xx^2,\xx^3$, and also $\dot\dd,\dot\cc,\dot\bb,\dot\aa$
interchangeably with $\yy^3,\yy^2,\yy^1,\yy^0$.  The monomial for the
above pipe dream is
$$
(a_1-\db_3)(b_1-\dc_3)(b_1-\dc_4)(b_2-\db_1)(c_1-\dt_3)(c_1-\dc_2)(c_3-\dt_2),
$$
ignoring all $*$~entries as required.  Removing the dots yields this
pipe dream's contribution to the ordinary quiver polynomial.\qed%
\end{example}

Now we have enough ingredients for our second positive combinatorial
formula for double quiver polynomials, the first being
Definition~\ref{d:QQrr}.

\begin{thm} \label{t:QQrr}
The double quiver polynomial for ranks\/~$\rr$ is the sum
\begin{eqnarray*}
  \QQ_\rr(\xx-\oyy) &=& \sum_{D \in \rp(v(\rr))}
  (\xx-\oyy)^{D \minus D_\homv}
\end{eqnarray*}
over all reduced pipe dreams~$D$ for the Zelevinsky permutation
$v(\rr)$ of the monomials $(\xx-\oyy)^D$ made from the crosses in~$D$
occupying the block antidiagonal and the block superantidiagonal (that
is, ignoring the~$*$~entries occupying every box strictly above the
block superantidiagonal).
\end{thm}
\begin{proof}
This follows from Definition~\ref{d:QQrr} and Theorem~\ref{t:BJS},
using the fact that every pipe dream $D \in \rp(v(\rr))$ contains the
subdiagram~$D_\homv$, and that $\rp(v(\homv))$ consists of the single
pipe dream~$D_\homv$.%
\end{proof}

In other words, the double quiver polynomial $\QQ_\rr(\xx-\oyy)$ is the
sum over monomials for all ``skew reduced pipe dreams'' $D \minus
D_\homv$ with $D \in \rp(v(\rr))$.

Since it will be important in Section~\ref{sec:pipestab} (and adds
even more concreteness to Theorem~\ref{t:BJS}), let us provide one way
of generating all reduced pipe dreams for~$v$.

\begin{defn}[\cite{BB}] \label{d:chute}
A \bem{chutable rectangle} is a connected $2 \times k$ rectangle $C$
inside a pipe dream $D$ such that $k \geq 2$ and every location in~$C$
is a~$\textcross$ except the following three: the northwest,
southwest, and southeast corners.  Applying a \bem{chute move} to~$D$
is accomplished by placing a $\textcross$ in the southwest corner of a
chutable rectangle~$C$ and removing the $\textcross$ from the
northeast corner of the same~$C$.%
\end{defn}

Heuristically, a chute move therefore looks like:
\begin{rcgraph}
\begin{array}{@{}c@{}}\\[-5ex]
\begin{array}{@{}r|c|c|c|@{}c@{}|c|c|c|l@{}}
  \multicolumn{7}{c}{}&\multicolumn{1}{c}{
                                          \phantom{\!+\!}}&
  \multicolumn{1}{c}{\begin{array}{@{}c@{}}\\\adots\end{array}}
  \\\cline{2-8}
           &\cdot&\!+\!&\!+\!& \toplinedots &\!+\!&\!+\!&\!+\!
  \\\cline{2-4}\cline{6-8}
           &\cdot&\!+\!&\!+\!&              &\!+\!&\!+\!&\cdot
  \\\cline{2-8}
  \multicolumn{1}{c}{\begin{array}{@{}c@{}}\adots\\ \\ \end{array}}&
  \multicolumn{1}{c}{\phantom{\!+\!}}
\end{array}
\quad\stackrel{\rm chute}\rightsquigarrow\quad
%
\begin{array}{@{}r|c|c|c|@{}c@{}|c|c|c|l@{}}
  \multicolumn{7}{c}{}&\multicolumn{1}{c}{
                                          \phantom{\!+\!}}&
  \multicolumn{1}{c}{\begin{array}{@{}c@{}}\\\adots\end{array}}
  \\\cline{2-8}
           &\cdot&\!+\!&\!+\!& \toplinedots &\!+\!&\!+\!&\cdot
  \\\cline{2-4}\cline{6-8}
           &\!+\!&\!+\!&\!+\!&              &\!+\!&\!+\!&\cdot
  \\\cline{2-8}
  \multicolumn{1}{c}{\begin{array}{@{}c@{}}\adots\\ \\ \end{array}}&
  \multicolumn{1}{c}{\phantom{\!+\!}}
\end{array}
\\[-3ex]
\end{array}
\end{rcgraph}

\begin{prop}\cite{BB} \label{p:chute}
There is a unique \bem{top reduced pipe dream} $D_\Top$ for~$v$, in
which no~$\textcross$ has any elbows~$\textelbow$ due north of it.
Every reduced pipe dream in~$\rp(v)$ can be obtained by starting
with~$D_\Top$ and applying some number of chute moves.%
\end{prop}

\section{Pipes to laces}

\label{sec:pipe2lace}

Pipe dreams for Zelevinsky permutations give rise to lacing
diagrams.

\begin{defn} \label{d:w(D)}
The $j^\th$ antidiagonal block is the block of size \mbox{$r_{\!j}
\times r_{\!j}$} along the main antidiagonal in the $j^\th$ block row.
%
%
Given a reduced pipe dream~$D$ for the Zelevinsky
permutation~$v(\rr)$, define the partial permutation $w_j = w_j(D)$
sending~$p$ to~$q$ if the pipe entering the $p^\th$ column from the
{\em right}\/ of the $(j-1)^\st$ antidiagonal block enters the
$j^\th$~antidiagonal block in its $q^\th$~column from the {\em right}.
Set $\ww(D) = (w_1,\ldots,w_n)$.  Equivalently, $\ww(D)$ is the lacing
diagram determined by~$D$.%
\end{defn}

\begin{example} \label{ex:strip}
The partial permutations arising from the reduced pipe dream in
Examples~\ref{ex:pipe} and~\ref{ex:monomial} come from the following
partial reduced pipe dreams,
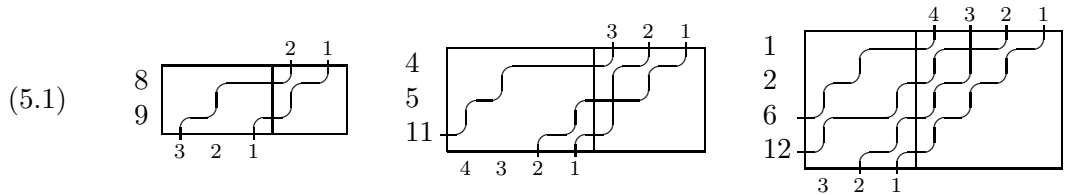
\begin{eqnarray} \label{eq:partrc}
\quad\qquad
\begin{array}{@{}l|ccc|cc|@{}}
\mcc{}&\ &\ &\mcc{\ \,\ \ \,}&\ &\mcc{\ \,\ \,\ }
\\[-3ex]\mcc{}&\multicolumn{3}{c}{}&\mcc{\petit 2}&\mcc{\petit 1}
\\\cline{2-6}
  8&&\er&\hor&\jr&\je
\\
  9&\er&\je&\er&\je&\
\\\cline{2-6}
\mcc{}&\mcc{\petit 3}&\mcc{\petit 2}&\mcc{\petit 1}&\multicolumn{2}{c}{}
\end{array}
\qquad
\begin{array}{@{}l|cccc|ccc|@{}}
\mcc{}&\ &\ &\ &\mcc{\ \,\ \ \,}&\ &\ &\mcc{\ \,\ \,\ }
\\[-3ex]\mcc{}&\multicolumn{4}{c}{}&\mcc{\petit 3}&\mcc{\petit 2}&\mcc{\petit 1}
\\\cline{2-8}
  4&   &\er&\hor&\hor&\jr&\jr&\je
\\
  5&\er&\je&    &\er &\+ &\je&
\\
 11&\je&   &\er &\jr &\je&   &
\\\cline{2-8}
\mcc{}&\mcc{\petit 4}&\mcc{\petit 3}&\mcc{\petit 2}&\mcc{\petit 1}&\multicolumn{3}{c}{}
\end{array}
\qquad
\begin{array}{@{}l|ccc|cccc|@{}}
\mcc{}&\ &\ &\mcc{\ \,\ \ \,}&\ &\ &\ &\mcc{\ \,\ \,\ }
\\[-3ex]\mcc{}&\multicolumn{3}{c}{}&\mcc{\petit 4}&\mcc{\petit 3}&\mcc{\petit 2}&\mcc{\petit 1}
\\\cline{2-8}
  1&   &\er &\hor&\jr&\+ &\jr&\je
\\
  2&\er&\je &\er &\jr&\jr&\je&
\\
  6&\jr&\hor&\jr &\jr&\je&   &
\\
 12&\je&\er &\jr &\je&   &   &
\\\cline{2-8}
\mcc{}&\mcc{\petit 3}&\mcc{\petit 2}&\mcc{\petit 1}&\multicolumn{4}{c}{}
\end{array}
\end{eqnarray}
where they send each number along the top to the number along the
bottom connected to it by a pipe (if such a pipe exists), or to
nowhere.  It is easy to see the lacing diagram from these
pictures.  Indeed, removing all segments of all pipes not
contributing to one of the partial permutations leaves some pipes
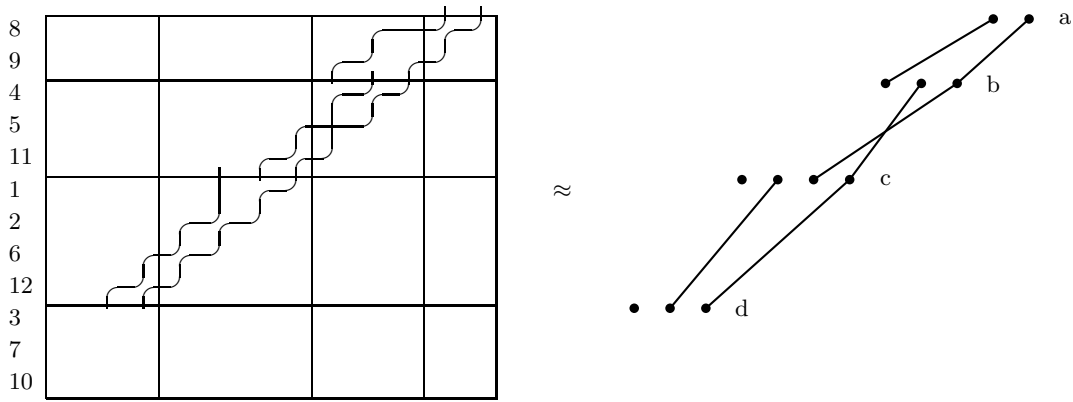
\psset{xunit=2.88ex,yunit=2.58ex}
\begin{rcgraph}
\begin{array}{@{}l|ccc|cccc|ccc|cc|@{}}
\mcc{}&\ \:&\ \:&\mcc{\ \,\ \,\ }&\ \:&\ \:&\ \:&\mcc{\ \,\ \,\ }&\ \:&\ \:&\mcc{\ \,\ \,\ }&\ \:&\mcc{\ \,\ \,\ }
\\[.3ex]\cline{2-13}
  8&   &   &   &   &   &   &   &   &\er&\hor&\jr&\je
\\[.3ex]
  9&   &   &   &   &   &   &   &\er&\je&\er &\je&
\\[.3ex]\cline{2-13}
  4&   &   &   &   &   &   &   &\er&\jr&\je &   &
\\[.3ex]
  5&   &   &   &   &   &   &\er&\+ &\je&    &   &
\\[.3ex]
 11&   &   &   &   &   &\er&\jr&\je&   &    &   &
\\[.3ex]\cline{2-13}
  1&   &   &   &  &\ver&\er&\je&   &   &    &   &
\\[.3ex]
  2&   &   &   &\er&\jr&\je&   &   &   &    &   &
\\[.3ex]
  6&   &   &\er&\jr&\je&   &   &   &   &    &   &
\\[.3ex]
 12&   &\er&\jr&\je&   &   &   &   &   &    &   &
\\[.3ex]\cline{2-13}
  3&   &   &   &   &   &   &   &   &   &    &   &
\\[.3ex]
  7&   &   &   &   &   &   &   &   &   &    &   &
\\[.3ex]
 10&   &   &   &   &   &   &   &   &   &    &   &
\\\cline{2-13}
\end{array}
\qquad\approx\quad
\begin{array}{c}
\pspicture[.1](-.5,-.5)(12,7.8)
\psdots(0,0)(1,0)(2,0)(3,4)(4,4)(5,4)(6,4)(7,7)(8,7)(9,7)(10,9)(11,9)
\psline(10,9)(7,7)
\psline(11,9)(9,7)(5,4)
\psline(8,7)(6,4)(2,0)
\psline(4,4)(1,0)
\rput(12,9)a
\rput(10,7)b
\rput(7,4)c
\rput(3,0)d
\endpspicture
\end{array}
\end{rcgraph}
\psset{xunit=2.2ex,yunit=2.2ex} that can be interpreted directly
as the desired lacing diagram
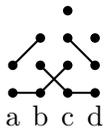
\begin{eqnarray*}
\begin{array}{c}
\pspicture[.1](-.5,-2)(3.5,3.5)
\psdots(0,0)(0,1)(1,0)(1,1)(1,2)(2,0)(2,1)(2,2)(2,3)(3,0)(3,1)(3,2)
\psline(0,0)(1,0)(2,1)
\psline(0,1)(1,2)
\psline(1,1)(2,0)(3,0)
\psline(2,2)(3,1)
\rput[b](0,-1.2)a
\rput[b](1,-1.2)b
\rput[b](2,-1.2)c
\rput[b](3,-1.2)d
\endpspicture
\end{array}
\end{eqnarray*}
by shearing to make the rightmost dots in each row line up vertically,
and then reflecting through the diagonal line $\diagdown$ of
slope~$-1$.\qed%
\end{example}

Warning: the $d \times d$ `pipe' representation of the lacing diagram
(at left above) does {\em not}\/ determine a reduced pipe dream
uniquely, because although the ordinary and singly virtual crossings
have been placed in the $d \times d$ grid, the doubly virtual
crossings still need to be placed (in the lower two superantidiagonal
blocks).

\begin{thm} \label{t:lacing}
Every reduced pipe dream $D \in \rp(v(\rr))$ gives rise to a minimally
crossing lacing diagram $\ww(D)$ with rank array~$\rr$; that is,
$\ww(D)\in W(\rr)$.
\end{thm}
\begin{proof}
Each $\times$~entry in the permutation matrix for~$v(\rr)$ corresponds
to a pipe in~$D$ entering due north of it and exiting due west of it.
The permutation~$v(\rr)$ was specifically constructed to have
exactly~$s_{ij}$ entries~$\times$ (for $i,j = 0,\ldots,n$) in the
intersection of the $i^\th$~block row and the $j^\th$~block column from
the right.%
\end{proof}

%
%
%
%
%
We shall see in Corollary~\ref{l:W} that in fact every minimally
crossing lacing diagram arises in this manner from a reduced pipe
dream.

\section{Rank stability of pipe dreams}

\label{sec:pipestab}

Just as the components in the quiver degeneration are stable for
uniformly increasing ranks~$\rr$ (Proposition~\ref{p:stability}), so are
the sets of pipe dreams for~$v(\rr)$.  Before stating the precise
result, we must analyze the changes to~$v(\rr)$ and its diagram upon
adding~$1$ to all ranks in~$\rr$.

\begin{lemma} \label{l:enlarge}
The passage from $\rr$ to~\mbox{$1+\rr$} enlarges every block
submatrix in the permutation matrix for $v(\rr)$ by one row and one
column.  Similarly, the diagram of~\mbox{$v(1+\rr)$} has the same
rectangular collections of cells in the southeast corners of the same
blocks as does the diagram of~$v(\rr)$.  (For an example, see
Fig.~\ref{fig:bigger}.)
\begin{figure}[ht]
$$
\begin{array}{@{}l|%
        @{}  c@{\:}@{\:}c@{\:}@{\:}c@{\:}@{\:}c@{\:}|%
        @{\:}c@{\:}@{\:}c@{\:}@{\:}c@{\:}@{\:}c@{\:}@{\:}c@{\:}|%
        @{\:}c@{\:}@{\:}c@{\:}@{\:}c@{\:}@{\:}c@{\:}|%
        @{\:}c@{\:}@{\:}c@{\:}@{\:}c@{\:}|@{}}
\cline{2-17}
 10&\: * & * & * & * & * & * & * & * & * &\ti&\cd&\cd&\cd&\cd&\cd&\cd
\\
 11&\: * & * & * & * & * & * & * & * & * &\cd&\ti&\cd&\cd&\cd&\cd&\cd
\\
 12&\: * & * & * & * & * & * & * & * & * &\cd&\cd&\ti&\cd&\cd&\cd&\cd
\\\cline{2-17}
 5 &\: * & * & * & * &\ti&\cd&\cd&\cd&\cd&\cd&\cd&\cd&\cd&\cd&\cd&\cd
\\
 6 &\: * & * & * & * &\cd&\ti&\cd&\cd&\cd&\cd&\cd&\cd&\cd&\cd&\cd&\cd
\\
 7 &\: * & * & * & * &\cd&\cd&\ti&\cd&\cd&\cd&\cd&\cd&\cd&\cd&\cd&\cd
\\
 14&\: * & * & * & * &\cd&\cd&\cd&\sq&\sq&\cd&\cd&\cd&\sq&\ti&\cd&\cd
\\\cline{2-17}
 1 &\:\ti&\cd&\cd&\cd&\cd&\cd&\cd&\cd&\cd&\cd&\cd&\cd&\cd&\cd&\cd&\cd
\\
 2 &\:\cd&\ti&\cd&\cd&\cd&\cd&\cd&\cd&\cd&\cd&\cd&\cd&\cd&\cd&\cd&\cd
\\
 3 &\:\cd&\cd&\ti&\cd&\cd&\cd&\cd&\cd&\cd&\cd&\cd&\cd&\cd&\cd&\cd&\cd
\\
 8 &\:\cd&\cd&\cd&\sq&\cd&\cd&\cd&\ti&\cd&\cd&\cd&\cd&\cd&\cd&\cd&\cd
\\
 15&\:\cd&\cd&\cd&\sq&\cd&\cd&\cd&\cd&\sq&\cd&\cd&\cd&\sq&\cd&\ti&\cd
\\\cline{2-17}
 4 &\:\cd&\cd&\cd&\ti&\cd&\cd&\cd&\cd&\cd&\cd&\cd&\cd&\cd&\cd&\cd&\cd
\\
 9 &\:\cd&\cd&\cd&\cd&\cd&\cd&\cd&\cd&\ti&\cd&\cd&\cd&\cd&\cd&\cd&\cd
\\
 13&\:\cd&\cd&\cd&\cd&\cd&\cd&\cd&\cd&\cd&\cd&\cd&\cd&\ti&\cd&\cd&\cd
\\
 16&\:\cd&\cd&\cd&\cd&\cd&\cd&\cd&\cd&\cd&\cd&\cd&\cd&\cd&\cd&\cd&\ti
\\\cline{2-17}
\end{array}
$$
\caption{The permutation \(v(1+\rr)\) and its diagram, for \(v(\rr)\)
as in Example~\ref{ex:zel}}\label{fig:bigger}
\end{figure}
\end{lemma}
\begin{proof}
The lace array giving rise to $1+\rr$ is the same as the array~$\ss$
giving rise to~$\rr$ except that $s_{0n}$ must be replaced by
$1+s_{0n}$.  The rest is straightforward.%
\end{proof}

Suppose $D \in \rp(v(\rr))$, and split $D$ into horizontal strips
(block rows).  In each strip, chute every~$\textcross$ upward (that
is, apply inverse chute moves) as far as possible.
By Proposition~\ref{p:chute}, this yields a top
reduced pipe dream in each horizontal strip.
Since every relevant~$\textcross$ in a
pipe dream $D \in \rp(v(\rr))$ is confined to the antidiagonal blocks
and those immediately to their left, we consider the $j^\th$
horizontal strip of a pipe dream to consist only of those~two~blocks,
with a single vertical dividing line.

\begin{example} \label{ex:top}
For example, in our running example of a reduced pipe dream
(Examples~\ref{ex:pipe}, \ref{ex:monomial}, and~\ref{ex:strip}), the
top pipe dream in the second (that is, $j=2$) horizontal strip is
depicted on the left side of Fig.~\ref{fig:bottom}.
\begin{figure}[ht]
\begin{rcgraph}
\begin{array}{@{}l|c@{\ }c@{\ }c|c@{\ }c@{\ }c@{\ }c|@{}}
\mcc{}&\ph&\ph&\mcc{\ph\ \:}&\ph&\ph&\ph&\mcc{\ph\ \:}
\\[-2ex]
\mcc{}&\dt_1&\dt_2&\mcc{\!\dt_3}&\dc_1&\dc_2&\dc_3&\mcc{\!\dc_4}
\\\cline{2-8}
c_1&\cd&\cd& + &\cd& + &\cd&\cd
\\
c_2&\cd&\cd& + &\cd&\cd&\cd&\cd
\\
c_3&\cd&\cd&\cd&\cd&\cd&\cd&\cd
\\
c_4&\cd&\cd&\cd&\cd&\cd&\cd&\cd
\\\cline{2-8}
\end{array}
\qquad\qquad
\begin{array}{@{}l|c@{\:\,}c@{\:\,}c@{\ }c|c@{\ }c@{\ }c@{\ }c@{\ }c|@{}}
\mcc{}&\ph&\ph&\ph&\mcc{\ph\ \:}&\ph&\ph&\ph&\ph&\mcc{\ph\ \:}
\\[-2ex]
\mcc{}&\dt_1&\dt_2&\dt_3&\mcc{\!\dt_4}&\dc_1&\dc_2&\dc_3&\dc_4&\mcc{\!\dc_5}
\\\cline{2-10}
c_1&\cd&\cd&\cd& + &\cd& + &\cd&\cd&\cd
\\
c_2&\cd&\cd&\cd& + &\cd&\cd&\cd&\cd&\cd
\\
c_3&\cd&\cd&\cd&\cd&\cd&\cd&\cd&\cd&\cd
\\
c_4&\cd&\cd&\cd&\cd&\cd&\cd&\cd&\cd&\cd
\\
c_5&\cd&\cd&\cd&\cd&\cd&\cd&\cd&\cd&\cd
\\\cline{2-10}
\end{array}
\end{rcgraph}
\caption{Top pipe dreams in \(D \in \rp(v)\), for \(v = v(\rr)\) or
\(v(1+\rr)\)}\label{fig:bottom}
\end{figure}
The corresponding top pipe dream (see Proposition~\ref{p:top}) in that
horizontal strip for $v(1+\rr)$ is on the right.\qed%
\end{example}

\begin{lemma} \label{l:blank}
Suppose $D \in \rp(v(\rr))$ consists of miniature top reduced pipe
dreams in its horizontal strips.  Inside the restriction of~$D$ to its
superantidiagonal block in block row~$j$, the\/ $2 \times 2$
configurations\/
$\begin{tinyrc}{
  \begin{array}{@{}|@{\,}c@{\,}|@{\,}c@{\,}|@{}}
      \hline   +   &  +
    \\\hline   +   &\cdot
    \\\hline
  \end{array}
  }\end{tinyrc}$
and\/~%
$\begin{tinyrc}{
  \begin{array}{@{}|@{\,}c@{\,}|@{\,}c@{\,}|@{}}
      \hline   *  &  *
    \\\hline   +  &\:\cdot\;
    \\\hline
  \end{array}
  }\end{tinyrc}$
in consecutive rows and columns never occur.
\end{lemma}
\begin{proof}
Suppose there is such a $2 \times 2$ configuration.  Since the
restriction of~$D$ to the superantidiagonal block in block row~$j$ is
a top pipe dream, all the cells in~$D$ due north of the $2 \times 2$
configuration are either $\textcross$ or~$*$, and all of these
represent crossings in~$D$.  Therefore the southwest~$\textcross$ in
the $2 \times 2$ configuration represents the crossing of two pipes
exiting through adjacent columns in block column $j+1$ from the right.
This is impossible because the transpose of~$D$ is a reduced pipe
dream for the permutation $v(\rr)^{-1}$, which is the Zelevinsky
permutation~$v(\rr')$ for the transpose ranks~$\rr'$ and hence can
have no descents within a block~row.%
\end{proof}

The easy proof of the next lemma is omitted.

\begin{lemma} \label{l:top}
Moving any reduced pipe dream for a permutation~$w$ one column due east,
or one row due south, yields a reduced pipe dream for\/~\mbox{$1+w$}.\qed
\end{lemma}

Here is the rank stability statement for pipe dreams, to be applied
in~Section~\ref{part:pos}.

\begin{prop} \label{p:top}
Let $\Top(\rr) \subseteq \rp(v(\rr))$ denote the set of reduced pipe
dreams~$D$ for~$v(\rr)$ consisting of a miniature top pipe dream $D_j$
in the $j^\th$ horizontal strip, for $j=0,\ldots,n-1$.  There is a
canonical bijection $\Top(\rr) \cong \Top(m+\rr)$ for all $m \geq 0$.
Moreover, corresponding miniature top pipe dreams $D_j$ and~$D'_{\!j}$
are translates of each other, in the same relative position to the
\mbox{vertical dividing line in horizontal strip~$j$}.%
\end{prop}
\begin{proof}
It suffices by induction to prove the case~$m=1$.  Starting with $D
\in \Top(\rr)$, inflate the blocks as in Lemma~\ref{l:enlarge}, and
construct a new pipe dream~$D'$ by placing each miniature top pipe
dream~$D_j$ in the same location relative to the vertical dividing
line in its inflated strip.  First we claim that $D' \in
\rp(v(1+\rr))$, for then automatically $D' \in \Top(1+\rr)$.  Then we
claim that every reduced pipe dream in~$\Top(1+\rr)$ comes from some
$D \in \Top(\rr)$ in this manner.

What we show precisely for the first claim is that for every~$\times$ in
the permutation matrix for~$v(1+\rr)$, the pipe (in $D'$) entering into
its row exits through its column.  Given a~$\times$ entry of~$v(1+\rr)$,
exactly one of the following three statements must hold, by construction
of Zelevinsky permutations: (i)~the~$\times$ lies in the northwest
corner of some superantidiagonal block; (ii)~the~$\times$ lies in the
southeast corner of the whole matrix; or (iii)~there is a
corresponding~$\times$ in~$v(\rr)$.

For an~$\times$~entry of type~(i), the pipe entering its row encounters
precisely one elbow, located in the same position as~$\times$, that
turns the pipe north to exit through its column.  For the~$\times$ of
type~(ii), observe that every entry along the superantidiagonal (not the
block superantidiagonal) is an elbow, because of the extra row added
into each block row (Lemma~\ref{l:enlarge}).  Thus its pipe proceeds
from the southwest corner to the northeast with no obstructions.

For an~$\times$ entry of type~(iii), we claim that the pipe~$\PP'$
starting in the row of~$\times$ enters the first antidiagonal block
(through which it passes) in the same location relative to the block's
northwest corner as does the pipe~$\PP$ for the corresponding~$\times$
in~$v(\rr)$.  If the~$\times$ lies in the bottom block row, then this is
trivial, and if~$\PP'$ never enters an antidiagonal block, then this is
vacuous.  Otherwise, the~$\times$ itself lies one row farther south from
the top of its block row, and we must show why this difference is
recovered by~$\PP'$.  Think of each $D'_{\!j}$ as having a blank column
(no~$\textcross$ entries) along the western edge of its strip, and a
blank row along the southern edge of its strip.  The pipe~$\PP'$
recovers one unit of height as it passes through each such region, and
it passes through exactly one on the way to its first antidiagonal
block.

Suppose that both~$\PP$ and~$\PP'$ enter an antidiagonal block at the
same vertical distance from the block's northwest corner.  Since $D'$
agrees with~$D$ in each antidiagonal block (as measured from their
northwest corners), the pipe~$\PP'$ exits the antidiagonal block in
the same position relative to its entry point as does~$\PP$.
Therefore we can use induction on the number of antidiagonal blocks
traversed by~$\PP$ to show that the pipe~$\PP'$ exits {\em every}\/
antidiagonal block in the same location relative to its northwest
corner as does~$\PP$.  If~the~$\times$ lies in the rightmost block
column, then $\PP'$ has already finished in the correct column.
Otherwise, $\PP'$ must recover one more horizontal unit, which it does
by traversing the blank row and column in the top horizontal strip.

For the second claim, we need only show that for any given $D' \in
\Top(1+\rr)$, each horizontal strip in~$D'$ has a column of elbows
along its western edge, a row of elbows along its southern edge, and
an antidiagonal's worth of elbows on the superantidiagonal, for then
each miniature top pipe dream $D'_{\!j}$ fits inside the deflated
block matrix for~$v(\rr)$ to make~$D \in \Top(\rr)$.  The antidiagonal
elbows come from the pipe added to make $s_{0n} \geq 1$: by
Theorem~\ref{t:lacing} there is a pipe going from the southwest block
to the northeast block, and by Theorem~\ref{t:mindiag} this pipe
crosses no others.  In particular, the southeast corner of each
superantidiagonal block in~$D'$ is an elbow.  Lemma~\ref{l:blank}
implies the southern elbow row, and its transpose (applied to
$v(1+\rr)^{-1} = v(1+\rr')$ for the transpose~$\rr'$ of~$\rr$) implies
the western elbow column.
\end{proof}

Interpreting Proposition~\ref{p:top} on lacing diagrams for pipe
dreams yields another useful stability statement that we shall
apply in Section~\ref{part:pos}; compare Corollary~\ref{c:lace+1}.

\begin{cor} \label{c:W}
Let $W_\rp(\rr) = \{\ww(D) \mid D \in \rp(v(\rr))\}$ be the set of
lacing diagrams obtained by decomposition of reduced pipe dreams
for~$v(\rr)$ into horizontal strips as in Definiton~\ref{d:w(D)} and
Theorem~\ref{t:lacing}.  Then $W_\rp(m+\rr) = \{m+\ww \mid \ww \in
W_\rp(\rr)\}$ is obtained by shifting each partial permutation list
in~$W_\rp(\rr)$ up by~$m$.
\end{cor}
\begin{proof}
Let $D \in \Top(\rr)$ correspond to $D' \in \Top(1+\rr)$.  The
miniature top pipe dream $D'_{\!j}$ in the $j^\th$ horizontal strip
of~$D'$ is one row farther from the bottom row of that strip
than~$D_j$ is from the bottom row in its strip, by
Lemma~\ref{l:enlarge} and Proposition~\ref{p:top}.  Reading (partial)
permutations in each horizontal strip with the southeast corner as the
origin, the result follows from Lemma~\ref{l:top}.%
\end{proof}

\part{Double quiver functions}

\label{part:pos}
\setcounter{section}{0}
\setcounter{thm}{0}
\setcounter{equation}{0}

\section{Double Stanley symmetric functions}

\label{sec:stanley}

Symmetric functions are formal power series in infinite alphabets,
which may be specialized to polynomials in alphabets of finite size.
We prepare conventions for such specializations, since they will be
crucial to our line of reasoning.

\begin{conv}
If a polynomial requiring $r$~or more variables for input is evaluated
on an alphabet with $r' < r$ letters, then the remaining $r - r'$
variables are to be set to zero.  Conversely, if a polynomial
requiring no more than $r$~variables for input is evaluated on an
alphabet of size $r' > r$, then the last $r'-r$ variables are
ignored.%
\end{conv}

\begin{example}
Consider the four $\yy\zz$ blocks obtained from the big pipe dream in
Example~\ref{ex:monomial}:
\begin{rcgraph}
\begin{array}{@{}l|c@{\ }c@{\ }c@{\ }c|c@{\ }c@{\ }c|@{}}
\multicolumn{8}{@{}c@{}}{}\\[-5ex]
\mcc{}&\ph&\ph&\ph&\mcc{\ph}&\ph&\ph&\mcc{\ph}
\\[-2ex]
\mcc{}&y^2_1&y^2_2&y^2_3&\mcc{\!y^2_4}&y^1_1&y^1_2&\mcc{\!y^1_3}
\\\cline{2-8}
x^1_1&
  \cd&\cd& + & + &\cd&\cd&\:\cd\:
\\[.2ex]
x^1_2&
  \cd&\cd&\cd&\cd& + &\cd&\:\cd\:
\\[.2ex]
x^1_3&
  \cd&\cd&\cd&\cd&\cd&\cd&\:\cd\:
\\[.2ex]\cline{2-8}
x^2_1&
  \cd& + &\cd&\cd&\cd&\cd&\:\cd\:
\\[.2ex]
x^2_2&
  \cd&\cd&\cd&\cd&\cd&\cd&\:\cd\:
\\[.2ex]
x^2_3&
  \cd&\cd&\cd&\cd&\cd&\cd&\:\cd\:
\\[.2ex]
x^2_4&
  \cd&\cd&\cd&\cd&\cd&\cd&\:\cd\:
\\\cline{2-8}
\end{array}
\end{rcgraph}
This is a reduced pipe dream for the permutation $w = 1253746$.  The
double Schubert polynomial $\SS_{1253746}(\xx^1,\xx^2-\yy^2,\yy^1)$ is
a sum of `monomials' of the form \mbox{$(x-y)$}, one for each reduced
pipe dream $D \in \rp(1253746)$.  Setting $\xx^2 = \yy^1 = 0$ in
\mbox{$\SS_{1253746}(\xx^1,\xx^2-\yy^2,\yy^1)$} yields
\mbox{$\SS_{1253746}(\xx^1-\yy^2)$}.  In the expression of this
polynomial as a sum of `monomials', reduced pipe dreams with crosses
outside the upper-left block contribute differently
to~$\SS_w(\xx^1-\yy^2)$ than they did
to~\mbox{$\SS_w(\xx^1,\xx^2-\yy^2,\yy^1)$}.  The pipe dream depicted
above, for example, contributes
\mbox{$(x^1_1-y^2_3)(x^1_1-y^2_4)(x^1_2-y^1_1)(x^2_1-y^2_2)$} to
\mbox{$\SS_w(\xx^1,\xx^2-\yy^1,\yy^2)$}, but only
\mbox{$(x^1_1-y^2_3)(x^1_1-y^2_4)(x^1_2)(-y^2_2)$} to
\mbox{$\SS_w(\xx^1-\yy^2)$}.\qed%
\end{example}

For each natural number $r \in \NN$ and any alphabet $\zz =
z_1,z_2,\ldots$ of infinite size, let $\zz_r = z_1, \ldots, z_r$ be an
alphabet of size~$r$.  By convention, define $\zz_r = \nothing$ to be
the empty alphabet for integers $r < 0$.  Suppose that $\xx =
\xx^0,\xx^1,\ldots,\xx^n$ and $\yy = \yy^n,\ldots,\yy^1,\yy^0$ are two
ordered finite lists of infinite alphabets, and that for each $m \in
\NN$ we are given truncations $\xx_{\rr(m)} =
\xx^0_{r_0(m)},\ldots,\xx^n_{r_n(m)}$ and $\yy_{\rr(m)} =
\yy^n_{r_n(m)},\ldots,\yy^0_{r_0(m)}$.  A~sequence
$p_m(\xx_{\rr(m)},\yy_{\rr(m)})$ of polynomials in these lists of
truncated alphabets is said to \mbox{\bem{converge}} to a power
series~$p(\xx,\yy)$ if the coefficient on any fixed monomial is
eventually constant as a function of~$m$.  Equivalently, we~say the
\bem{limit exists},~and~write
\begin{eqnarray*}
  \lim_{m\to\infty} p_m(\xx_{\rr(m)},\yy_{\rr(m)}) &=& p(\xx,\yy).
\end{eqnarray*}
Here, $p(\xx,\yy)$ is allowed to be an arbitrary formal sum of monomials
with integer coefficients in variables from the union of alphabets in
the lists $\xx$ and~$\yy$.

Recall from Section~\ref{sec:stability} the notation $m + w$ for
nonnegative integers~$m$ and partial permutations~$w$.

\begin{prop} \label{p:stanley}
Given a partial permutation~$w$, an infinite set of variables~$\zz$,
and a function $r:\NN\to\ZZ$ approaching~$\infty$, the following limit
exists and is symmetric~in~$\zz$:
\begin{eqnarray*}
  F_w(\zz) &=& \lim_{m\rightarrow\infty} \SS_{m+w}(\zz_{r(m)}).
\end{eqnarray*}
\end{prop}
\begin{defn}
The limit $F_w(\zz)$ is called the \bem{Stanley symmetric function} or
\bem{stable Schubert polynomial} for~$w$. (N.B. In the notation of
\cite{St} the permutation $w^{-1}$ is used for indexing instead of
$w$).
\end{defn}
\begin{proof}
After harmlessly extending $w$ to an honest permutation in the minimal
way (which leaves the Schubert polynomials unchanged), see \cite[Section
VII]{Mac}.%
\end{proof}

Stanley's original definition of~$F_w$ was combinatorial, and
motivated by the fact that the coefficient of any squarefree monomial
in~$F_w$ is the number of reduced decompositions of the
permutation~$w$.  The fact that Stanley's formulation of~$F_w$ is a
stable Schubert polynomial in the above sense follows immediately from
the formula for the Schubert polynomial given in
\cite[Theorem~1.1]{BJS}.

Being a symmetric function, $F_w$ can be written as a sum
\begin{eqnarray} \label{eq:FtoS}
  F_w(\zz) &=& \sum_\la \alpha^\la_w s_\la(\zz)
\end{eqnarray}
of Schur functions~$s_\la$ with coefficients~$\alpha^\la_w$, which we
call \bem{Stanley coefficients}.
Stanley
credits Edelman and Greene \cite{EGn87} with the first proof that
these coefficients~$\alpha^\la_w$ are {\em nonnegative}.  There are
various descriptions for $\alpha^\la_w$, all combinatorial: reduced
word tableaux \cite{LS89,EGn87}, the leaves of shape~$\la$ in the
transition tree for~$w$ \cite{LS85}, certain promotion sequences
\cite{Ha92}, and peelable tableaux \cite{RSplact95}.  The
coefficients~$\alpha^\la_w$ are also special cases of the quiver
constants~$c_\bla(\rr)$ from Theorem~\ref{t:BFtoKMS}~\cite{Bu01}.  We
shall prove the converse in Theorem~\ref{t:qc=ss}, below: the quiver
constants~$c_\bla(\rr)$ are special cases of Stanley
coefficients~$\alpha^\la_w$.

\begin{prop} \label{p:doublestanley}
Given functions $r,\dot r : \NN \to \ZZ$ approaching~$\infty$, a
partial permutation~$w$, and infinite alphabets $\zz,\dot\zz$, the
limit of double Schubert polynomials
\begin{eqnarray*}
  F_w(\zz-\dot\zz) &=&
  \lim_{m\to\infty}\SS_{m+w}\big(\zz_{r(m)}-\dot\zz_{\dot r(m)}\big)
\end{eqnarray*}
exists, does not depend on $r$ or~$\dot r$, and is symmetric
separately in $\zz$ and~$\dot\zz$.
\end{prop}
\begin{defn}
The limit $F_w(\zz-\dot\zz)$ is the \bem{double Stanley symmetric
function} or \bem{stable double Schubert polynomial}.
\end{defn}
\begin{proof}
Again, we may as well assume $w$ is an honest permutation.  One may
reduce to the ``single'' stable Schubert case by applying the
coproduct formula \cite[(6.3)]{Mac}
\begin{eqnarray*}
  \SS_w(\xx-\oyy) &=& \sum_{u,v} (-1)^{\ell(v)} \SS_u(\xx)
  \SS_{v^{-1}}(\yy)
\end{eqnarray*}
where the sum runs over the pairs of permutations $u,v$ such that $vu=w$
and the sum $\ell(u)+\ell(v)$ of the lengths of~$v$ and~$u$ equals the
length~$\ell(w)$ of~$w$.%
\end{proof}

For a sequence $\ww = w_1,\ldots,w_n$ of partial permutations of size
\mbox{$r_0\!\times\!r_1,\ldots,r_{n-1}\!\times\!r_n$}, and two
sequences of {\em finite}\/ alphabets as in~\eqref{xx}
and~\eqref{xxj}, denote by
\begin{eqnarray*}
  \SS_\ww(\xx-\oyy) &=& \prod_{j=1}^n \SS_{w_j}(\xx^{j-1}\!-\yy^j)
\end{eqnarray*}
the product of double Schubert polynomials in consecutive $\xx$
and~$\yy$ alphabets.  Write
\begin{eqnarray} \label{F}
  F_\ww(\xx-\oyy) &=& \prod_{j=1}^n F_{w_j}(\xx^{j-1}\!-\yy^j)
\end{eqnarray}
for the corresponding product of double Stanley symmetric functions,
where now the alphabets are infinite.  If we need to specify finite
alphabets~\eqref{xxj}, when ranks~$\rr$ have been set so that $\xx^j$
and~$\yy^j$ have size~$r_{\!j}$, then we shall write explicitly
$F_\ww(\xx_\rr-\yy_{\!\rr})$.  Keep in mind that the element
in~\eqref{F} remains unchanged when $\ww = (w_1,\ldots,w_n)$ is
replaced by $m+\ww$ for~$m \in \NN$.

\section{Double quiver functions}

\label{sec:stablequiv}

Stanley symmetric functions arise in our context by taking limits of
double quiver polynomials via the pipe formula, Theorem~\ref{t:QQrr}.
Each reduced pipe dream~$D$ for the Zelevinsky permutation~$v(\rr)$
determines a polynomial~$\SS_{\ww(D)}(\xx-\oyy)$ by
Definition~\ref{d:w(D)}, to which $D$ contributes a ``monomial'', as
in Section~\ref{sec:pipe2quiver}.  Here, in Proposition~\ref{p:sym},
we shall that each reduced pipe dream for~$v(\rr)$ can contribute a
different monomial by way of a crucial symmetry, as follows.

\begin{defn} \label{d:reverse}
The \bem{reverse monomial} associated to a pipe dream~$D$ is the product
$(\tilde\xx-\tilde\yy)^D$ over all $\textcross$ entries in~$D$ of
$(\tilde x_+ - \tilde y_+)$, where the variable~$\tilde x_+$ sits at the
left end of the row containing~$\textcross$ after reversing each of the
$\xx$~alphabets, and the variable~$\tilde y_+$ sits atop the column
containing~$\textcross$ after reversing each of the $\yy$~alphabets.
\end{defn}

\begin{example}
Reversing each of the row and column alphabets in
Example~\ref{ex:monomial} gives
\begin{rcgraph}
\begin{array}{@{}l|c@{\ }c@{\ }c|c@{\ }c@{\ }c@{\ }c|c@{\ }c@{\ }c|c@{\ }c|@{}}
\multicolumn{13}{@{}c@{}}{}\\[-3ex]
\mcc{}&\ph&\ph&\mcc{\ph}&\ph&\ph&\ph&\mcc{\ph}&\ph&\ph&\mcc{\ph}&\ph&\mcc{\ph}
\\[-1ex]
\mcc{}&y^3_3&y^3_2&\mcc{\!y^3_1}&y^2_4&y^2_3&y^2_2&\mcc{\!y^2_1}&y^1_3&y^1_2&\mcc{\!y^1_1}&
y^0_2&\mcc{\!y^0_1}
\\\cline{2-13}
x^0_2&
  \sr&\sr&\sr&\sr&\sr&\sr&\sr&\cd&\cd& + &\cd&\cd\:
\\
x^0_1&
  \sr&\sr&\sr&\sr&\sr&\sr&\sr&\cd&\cd&\cd&\cd&\cd\:
\\\cline{2-13}
x^1_3&
  \sr&\sr&\sr&\cd&\cd& + & + &\cd&\cd&\cd&\cd&\cd\:
\\
x^1_2&
  \sr&\sr&\sr&\cd&\cd&\cd&\cd& + &\cd&\cd&\cd&\cd\:
\\
x^1_1&
  \sr&\sr&\sr&\cd&\cd&\cd&\cd&\cd&\cd&\cd&\cd&\cd\:
\\\cline{2-13}
x^2_4&
  \cd&\cd& + &\cd& + &\cd&\cd&\cd&\cd&\cd&\cd&\cd\:
\\
x^2_3&
  \cd&\cd&\cd&\cd&\cd&\cd&\cd&\cd&\cd&\cd&\cd&\cd\:
\\
x^2_2&
  \cd& + &\cd&\cd&\cd&\cd&\cd&\cd&\cd&\cd&\cd&\cd\:
\\
x^2_1&
  \cd&\cd&\cd&\cd&\cd&\cd&\cd&\cd&\cd&\cd&\cd&\cd\:
\\\cline{2-13}
x^3_3&
  \cd&\cd&\cd&\cd&\cd&\cd&\cd&\cd&\cd&\cd&\cd&\cd\:
\\
x^3_2&
  \cd&\cd&\cd&\cd&\cd&\cd&\cd&\cd&\cd&\cd&\cd&\cd\:
\\
x^3_1&
  \cd&\cd&\cd&\cd&\cd&\cd&\cd&\cd&\cd&\cd&\cd&\cd\:
\\\cline{2-13}
\end{array}
\quad\begin{array}{c}\\[1ex]=\end{array}\quad
\begin{array}{@{}l|c@{\ }c@{\ }c|c@{\ }c@{\ }c@{\ }c|c@{\ }c@{\ }c|c@{\ }c|@{}}
\multicolumn{13}{@{}c@{}}{}\\[-3ex]
\mcc{}&\ph&\ph&\mcc{\ph}&\ph&\ph&\ph&\mcc{\ph}&\ph&\ph&\mcc{\ph}&\ph&\mcc{\ph}
\\[-1ex]
\mcc{}&\dt_3&\dt_2&\mcc{\!\dt_1}&\dc_4&\dc_3&\dc_2&\mcc{\!\dc_1}&\db_3&\db_2&\mcc{\!\db_1}&\da_2&\mcc{\!\da_1}
\\\cline{2-13}
a_2&
  \sr&\sr&\sr&\sr&\sr&\sr&\sr&\cd&\cd& + &\cd&\cd\:
\\
a_1&
  \sr&\sr&\sr&\sr&\sr&\sr&\sr&\cd&\cd&\cd&\cd&\cd\:
\\\cline{2-13}
b_3&
  \sr&\sr&\sr&\cd&\cd& + & + &\cd&\cd&\cd&\cd&\cd\:
\\
b_2&
  \sr&\sr&\sr&\cd&\cd&\cd&\cd& + &\cd&\cd&\cd&\cd\:
\\
b_1&
  \sr&\sr&\sr&\cd&\cd&\cd&\cd&\cd&\cd&\cd&\cd&\cd\:
\\\cline{2-13}
c_4&
  \cd&\cd& + &\cd& + &\cd&\cd&\cd&\cd&\cd&\cd&\cd\:
\\
c_3&
  \cd&\cd&\cd&\cd&\cd&\cd&\cd&\cd&\cd&\cd&\cd&\cd\:
\\
c_2&
  \cd& + &\cd&\cd&\cd&\cd&\cd&\cd&\cd&\cd&\cd&\cd\:
\\
c_1&
  \cd&\cd&\cd&\cd&\cd&\cd&\cd&\cd&\cd&\cd&\cd&\cd\:
\\\cline{2-13}
d_3&
  \cd&\cd&\cd&\cd&\cd&\cd&\cd&\cd&\cd&\cd&\cd&\cd\:
\\
d_2&
  \cd&\cd&\cd&\cd&\cd&\cd&\cd&\cd&\cd&\cd&\cd&\cd\:
\\
d_1&
  \cd&\cd&\cd&\cd&\cd&\cd&\cd&\cd&\cd&\cd&\cd&\cd\:
\\\cline{2-13}
\end{array}
\end{rcgraph}
This reduced pipe dream contributes the reverse monomial
$$
(a_2-\db_1)(b_3-\dc_2)(b_3-\dc_1)(b_2-\db_3)(c_4-\dt_1)(c_4-\dc_3)(c_2-\dt_2)
$$
to the double quiver polynomial $\QQ_\rr(\aa,\bb,\cc,\dd -
\dot\dd,\dot\cc,\dot\bb,\dot\aa)$.  Reversing the alphabets in the
argument of this polynomial leaves it invariant, because
by~Proposition~\ref{p:sym} this polynomial already equals the sum of
all reverse monomials from reduced pipe dreams for~$v(\rr)$.\qed%
\end{example}

In conjunction with Corollary~\ref{c:bottom}, the next observation
will allow us to take limits of double quiver polynomials effectively,
via Proposition~\ref{p:doublestanley}.

\begin{prop} \label{p:sym}
The double quiver polynomial\/~$\QQ_\rr(\xx-\oyy)$ is symmetric
separately in each of the alphabets\/~$\xx^i$ as well as in each of
the alphabets\/~$\yy^j$.  In particular,
\begin{eqnarray*}
  \QQ_\rr(\xx-\oyy) &=& \sum_{D \in \rp(v(\rr))} (\tilde\xx-\tilde\yy)^D,
\end{eqnarray*}
where $\tilde\xx$ and $\tilde\yy$ mean to reverse each alphabet in
each list (but not to reverse the sequential order of the alphabets
themselves).
\end{prop}
\begin{proof}
Matrix Schubert varieties for Zelevinsky permutations are preserved by
the action of block diagonal matrices---including block diagonal
permutation matrices---on the right~and~left.  Now apply the symmetry
to Theorem~\ref{t:QQrr}.%
\end{proof}

Next we record an important consequence of Proposition~\ref{p:top}.

\begin{cor} \label{c:bottom}
There is a fixed integer~$\ell$, independent of~$m$, such that setting
the last~$\ell$ variables to zero in every finite alphabet from the
lists $\xx_{m+\rr}$ and\/~$\yy_{\!m+\rr}$ kills the reverse monomial
for every reduced pipe dream $D \in \rp(v(m+\rr))$ containing at least
one cross~$\textcross$ in an antidiagonal block.%
\end{cor}
\begin{proof}
Suppose~$D$ is a reduced pipe dream for $v(m+\rr)$.  Divide~$D$ into
horizontal strips (block rows) as in Section~\ref{sec:pipestab}, and
consider the reduced pipe dream~$D'$ in the same chute class as~$D$
that consists of a top pipe dream in each strip.  The miniature top
pipe dream in horizontal strip~$j$ has an eastern half~$D'_{\!j}$ in
the $j^\th$ antidiagonal block.  Inside the $j^\th$~antidiagonal
block, consider the highest antidiagonal~$A_j$ above which every
$\textcross$ in~$D'_{\!j}$~lies.  Thus $A_j$ bounds an isosceles right
triangular region $\triangledown_{\!j}(D)$ in the northwest corner of
the $j^\th$ antidiagonal block that contains every $\textcross$
in~$D'_{\!j}$.

We claim that $\triangledown_{\!j}(D)$ in fact contains every
$\textcross$ of~$D$ itself that lies in the $j^\th$ antidiagonal
block.  Indeed, suppose that~$D$ contains a $\textcross$ in the
$j^\th$ antidiagonal block.  Using (any sequence of) upward chute
moves to get from~$D$ to~$D'$, the antidiagonal $\textcross$ we are
considering must end up inside $\triangledown_{\!j}(D)$, because
chuting any $\textcross$ up pushes it north and east.  Our claim
follows because upward chutes push each $\textcross$ at least as far
east as north.

Now define $\ell(D)$ to be the maximum of the leg lengths of the
isosceles right triangles~$\triangledown_{\!j}(D)$.  Chuting up as in
the previous paragraph, both the row variable~$\tilde x_+$ and the
column variable~$\tilde y_+$ of each antidiagonal $\textcross$ in~$D$
must lie within~$\ell(D)$ of the end of its finite alphabet.
Therefore, pick~$\ell$ to be the maximum of the numbers $\ell(D)$ for
reduced pipe dreams $D \in \rp(v(m+\rr))$.  The key (and final) point
is that~$\ell$ does not depend on~$m$ for $m \geq 0$, by
Proposition~\ref{p:top}.%
\end{proof}

The main result of this section concerning double quiver polynomials
$\QQ_\rr(\xx-\oyy)$ says that limits exist as the ranks~$\rr$ get
replaced by~\mbox{$m+\rr$} for $m \to \infty$.  We consider this an
algebraic rather than a combinatorial result (the latter being
Theorem~\ref{t:lace}), since we are unable at this stage to identify
the set~$W_\rp(\rr)$ from Corollary~\ref{c:W} in a satisfactory way.
Note that $W_\rp(\rr)$ is a set of distinct partial permutation
lists---not a multiset---contained in the set~$W(\rr)$ of minimum
length lacing diagrams with ranks~$\rr$, by Theorem~\ref{t:lacing}.
We shall see in Corollary~\ref{l:W} that in fact $W_\rp(\rr) =
W(\rr)$; but until then it is important for the proofs to distinguish
$W_\rp(\rr)$ from~$W(\rr)$.

\begin{thm} \label{t:stanley}
The limit of double quiver polynomials $\QQ_{m+\rr}(\xx-\oyy)$ for\/
$m$ approaching\/~$\infty$ exists and equals the sum
\begin{eqnarray*}
  \FF_\rr(\xx-\oyy) \ \::=\ \:
  \lim_{m\to\infty}\QQ_{m+\rr}(\xx-\oyy) &=& \!\!\!\sum_{\ww \in
  W_\rp(\rr)} F_\ww(\xx-\oyy)
\end{eqnarray*}
of products of double Stanley symmetric functions.  The limit power
series~$\FF_\rr(\xx-\oyy)$ is symmetric separately in each of the
$2n+2$ infinite alphabets $\xx^0,\ldots,\xx^n$
and~$\yy^n,\ldots,\yy^0$.%
\end{thm}
\begin{defn} \label{d:stanley}
$\FF_\rr(\xx-\oyy)$ is called the \bem{double quiver function} for
the rank array~$\rr$.%
\end{defn}
\begin{proof}
For all $m \geq 0$, the set\/~$W_\rp(m+\rr)$ is obtained
from\/~\mbox{$W_\rp(\rr)$} by replacing each permutation list\/~$\ww$
with\/ $m+\ww$, by Corollary~\ref{c:W}.  For each partial permutation
list $\ww \in W_\rp(\rr)$, consider the set $\rp_\ww(v(m+\rr))$ of
reduced pipe dreams~$D$ for~$v(m+\rr)$ whose lacing diagrams~$\ww(D)$
equal~$m+\ww$.  The sum of all reverse monomials for pipe dreams $D
\in \rp_\ww(v(m+\rr))$ is a product of polynomials, each factor coming
from the miniature pipe dreams in a single horizontal strip.  The
double quiver polynomial $\QQ_{m+\rr}(\xx-\oyy)$ is the sum of these
products, the sum being over lacing diagrams $\ww \in W_\rp(\rr)$.
A~similar statement can be made after setting the last~$\ell$
variables in each alphabet to zero as in Corollary~\ref{c:bottom},
although now $\QQ_\rr(\xx-\oyy)$ and the polynomials for each lacing
diagram $\ww \in W_\rp(\rr)$ have changed a little (but only in those
terms involving a variable near the end of some alphabet).

We claim that now, after setting the last few variables in each
alphabet to zero, the $j^\th$~polynomial in the product for any fixed
diagram~$\ww$ equals the result of setting the last few variables (at
most~$\ell$) to zero in the honest double Schubert polynomial
$\SS_{m+w_j}(\xx^{j-1} - \yy^j)$.  For this, note that every reduced
pipe dream $D$ for~$v(\rr)$ contributing a nonzero reverse monomial
now has $\textcross$ locations only in superantidiagonal blocks.
Hence we can read the restriction of~$D$ to the
$j^\th$~superantidiagonal block as a reduced pipe dream with its
southeast corner as the origin.  Since the variables in this one-block
rectangle that we have set to zero are far from the southeast corner,
our claim is proved.

Taking the limit as $m \to \infty$ produces the desired product of
double Stanley symmetric functions by
Proposition~\ref{p:doublestanley}.%
\end{proof}

The limit in Theorem~\ref{t:stanley} stabilizes in the following
strong sense.  The analogous stabilization for the limit defining
Stanley symmetric functions from Schubert polynomials
(Proposition~\ref{p:stanley}) almost never occurs, since the finite
variable specialization is symmetric, while Schubert polynomials
rarely are, even considering only variables that actually~occur.

\begin{cor} \label{c:stanley}
For all $m \gg 0$, the double quiver polynomial $\QQ_{m+\rr}(\xx-\oyy)$
actually equals the specialization $\FF_\rr(\xx_{m+\rr}-\yy_{\!m+\rr})$
of the double quiver function.%
\end{cor}
\begin{proof}
The double quiver polynomial $\QQ_{m+\rr}(\xx-\oyy)$ has fixed degree
and is symmetric in each of its finite alphabets for all~$m$
(Proposition~\ref{p:sym}).  Hence we need only show that for any fixed
monomial of degree at most $\deg(\QQ_{m+\rr}(\xx-\oyy))$, applying a
symmetry yields a monomial whose coefficient equals the coefficient on
the corresponding monomial in \mbox{$\FF_\rr(\xx-\oyy)$}.  This
follows immediately from Theorem~\ref{t:stanley}, because for $m \gg
0$, all monomials involving variables whose lower indices are at most
$\deg(\QQ_{m+\rr}(\xx-\oyy))$ have the correct coefficients.%
\end{proof}

The previous result
led us to suspect that the `$\gg$' sign may be replaced by~`$\geq$'.

\begin{conj} \label{m>0}
$\QQ_{m+\rr}(\xx-\oyy) = \FF_\rr(\xx_{m+\rr}-\oyy_{\!m+\rr})$ for all
$m \geq 0$.
\end{conj}

After we made Conjecture~\ref{m>0} in the first draft of this paper,
A.\thinspace{}Buch responded by proving it in \cite{BuchAltSign}.  We
have no need for the full strength of Conjecture~\ref{m>0} in this
paper, since Corollary~\ref{c:stanley} suffices; see
Proposition~\ref{p:QQrr}, below, and the last line of its proof.  We
shall comment on related issues in Remarks~\ref{rk:stable},
\ref{rk:double}, and~\ref{rk:SSww}.

\section{All components are lacing diagram orbit closures}

\label{sec:components}

Recall from Theorem~\ref{c:ratio} that the quiver polynomial
$\QQ_\rr(\xx-\oxx)$ is obtained from the double quiver polynomial
$\QQ_\rr(\xx-\oyy)$ by setting $\yy^j = \xx^j$ for all~$j$, so $\yy =
\oxx$.  Here, we shall need to specialize products $F_\ww(\xx-\oyy)$
and $\SS_\ww(\xx-\oyy)$ of double Schubert and Stanley symmetric
functions by setting $\yy = \xx$, to get $F_\ww(\xx-\oxx)$ and
$\SS_\ww(\xx-\oxx)$. Note that the alphabet~$\xx^j$ obtained by
specializing~$\yy^j$ never interferes in a catastrophic cancelative
manner with the original alphabet~$\xx^j$, because for any lacing
diagram~$\ww$, the $j^\th$~factors of~$\SS_\ww(\xx-\oxx)$
and~$F_\ww(\xx-\oxx)$ involve the distinct alphabets $\xx^{j-1}$
and~$\xx^j$.

\begin{lemma} \label{l:nonzero}
A~nonempty positive sum of products $F_\ww(\xx-\oxx)$ of double
Stanley symmetric functions in differences of consecutive alphabets is
nonzero.
\end{lemma}
\begin{proof}
One may use the double version of~\eqref{eq:FtoS} to write each term
$F_\ww(\xx-\oxx)$ as a nonnegative sum of products of Schur functions
in differences of consecutive alphabets.  Such products are linearly
independent by the remark in \cite[Section~2.2]{BF}.%
\end{proof}

Now we can finally identify the components of quiver degenerations.

\begin{thm} \label{t:components}
The components of the quiver degeneration $\Omega_\rr(0)$ are exactly
the orbit closures $\OO(\ww)$ for $\ww \in W(\rr)$.  Furthermore,
$\Omega_\rr(0)$ is generically reduced along each
component~$\OO(\ww)$, so its multiplicity there equals\/~$1$.
\end{thm}
\begin{proof}
By Corollary \ref{c:quivSchub} we find that
\begin{eqnarray} \label{KMS}
  \QQ_{m+\rr}(\xx-\oxx) &=& \sum_\ww c_{m+\ww}(\rr)
  \SS_{m+\ww}(\xx-\oxx)
\end{eqnarray}
for all nonnegative integers~$m$.  By~Proposition~\ref{p:stability},
for any ranks~$\rr$ and any lacing diagram~$\ww$, we know that
$c_\ww(\rr) = c_{m+\ww}(m+\rr)$, and also that \mbox{$c_\ww(m+\rr) =
0$} unless every partial permutation in the list $\ww$
fixes~$1,\ldots,m$.  Hence we can rewrite~\eqref{KMS} as
\begin{eqnarray*}
  \QQ_{m+\rr}(\xx-\oxx) &=& \sum_\ww c_{m+\ww}(m + \rr)
  \SS_{m+\ww}(\xx-\oxx)
\end{eqnarray*}
and note again that it holds for all integers $m \geq 0$.  The
alphabets in the previous two displayed equations are all finite, but
grow with~$m$.  Taking limits as $m \to \infty$ yields
\begin{eqnarray} \label{=}
  \lim_{m\to\infty} \QQ_{m+\rr}(\xx-\oxx)
  &=&
  \sum_\ww c_\ww(\rr) F_\ww(\xx-\oxx),
\end{eqnarray}
the sum of products of double Stanley symmetric functions
corresponding to components of the quiver
degeneration~$\Omega_\rr(0)$, each counted according to its
multiplicity.  This equation takes place in the ring of symmetric
functions, where the alphabets $\xx = \xx^0,\ldots,\xx^n$ are
infinite.

On the other hand, specializing Theorem~\ref{t:stanley} to $\yy =
\oxx$ implies that the limit on the left side of~\eqref{=} is a
multiplicity-free sum of Stanley symmetric functions
\mbox{$F_\ww(\xx-\oxx)$} indexed by $\ww \in W_\rp(\rr)$.  By
Lemma~\ref{l:nonzero}, it must therefore be that \mbox{$c_\ww(\rr) =
0$} unless $\ww \in W_\rp(\rr)$, in which case $c_\ww(\rr) = 1$.  But
we know by Proposition~\ref{p:lace} that $c_\ww(\rr) \geq 1$ for $\ww
\in W(\rr)$, so it must be that $W(\rr) \subseteq W_\rp(\rr)$.  We
conclude that $W_\rp(\rr) = W(\rr)$, and that all the multiplicities
$c_\ww(\rr)$ equal~$1$ for $\ww\in W(\rr)$.%
\end{proof}

\begin{cor}[Component formula for quiver polynomials---Schubert version]
\label{c:component}\mbox{}\\
The quiver polynomial $\QQ_\rr(\xx-\oxx)$ equals the sum
\begin{eqnarray*}
  \QQ_\rr(\xx-\oxx) &=& \sum_{\ww \in W(\rr)} \SS_\ww(\xx-\oxx)
\end{eqnarray*}
of products of double Schubert polynomials (in consecutive alphabets
$\xx = \xx^0,\ldots,\xx^n\!$) indexed by minimum length lacing
diagrams with rank array\/ $\rr$.%
\end{cor}

We have chosen to present this consequence of
Theorem~\ref{t:components} and Corollary~\ref{c:quivSchub} now, even
though it is a specialization of the stable double component formula
in the next section, to emphasize that it has a direct geometric
interpretation: the right-hand side is the sum of the equivariant
cohomology classes of components in the quiver degeneration.  The
double formula to come in Theorem~\ref{t:lace} currently lacks such a
geometric interpretation, and moreover the reason why it implies
Corollary~\ref{c:component} is somewhat subtle; see
Corollary~\ref{c:specialize} and Remark~\ref{rk:SSww}.  The subtlety
in this argument is related to the fact, discussed in
Remark~\ref{rk:SSww}, that replacing~$\oxx$ with~$\yy$ in
Corollary~\ref{c:component} {\em always}\/ yields a false statement
whenever there is more than one term on the right hand side.

In the course of proving Theorem~\ref{t:components}, we reached a
notable combinatorial result.

\begin{cor} \label{l:W}
Every minimal length lacing diagram for a rank array\/~$\rr$ occurs as
the lacing diagram derived from a reduced pipe dream for the
Zelevinsky permutation\/~$v(\rr)$:
\begin{eqnarray*}
  W(\rr) &=& W_\rp(\rr).
\end{eqnarray*}
\end{cor}

\begin{remark}
Theorem~\ref{t:lacing} and Corollary~\ref{l:W} give another
explanation (and proof) for Theorem~\ref{t:mindiag}.  That laces of
$\ww \in W(\rr)$ can cross at most once is exactly the statement that
pipes in reduced pipe dreams for~$v(\rr)$ cross at most once.  That
laces starting or ending in the same column do not cross at all
follows from the fact that Zelevinsky permutations and their inverses
have no descents within block~rows.%
\end{remark}

\section{Stable double component formula}

\label{sec:pos}

\begin{thm} \label{t:lace}
The double quiver function can be expressed as the sum
\begin{eqnarray*}
  \FF_\rr(\xx-\oyy) &=& \sum_{\ww \in W(\rr)} F_\ww(\xx-\oyy)
\end{eqnarray*}
of products of double Stanley symmetric functions.
\end{thm}
\begin{proof}
Using Corollary~\ref{l:W},
we can replace the sum in Theorem~\ref{t:stanley} over $W_\rp(\rr)$
by one over $W(\rr)$.
\end{proof}

\begin{remark} \label{rk:yong}
Definition~\ref{d:w(D)} reads pipes as if they flow northeast to
southwest (compare Definition~\ref{d:pipe}) because our proof of
Theorem~\ref{t:stanley} requires that we read antidiagonal blocks in
Zelevinsky pipe dreams as miniature pipe dreams with their southeast
corners as their origins.  This method suggests a direct combinatorial
proof of Corollary~\ref{l:W}, without appealing to the geometric lower
bound on multiplicities of components in quiver degenerations afforded
by Corollary~\ref{c:quivSchub}.  This idea has since been carried out
in \cite{Yong03}.  Thus the component formula in Theorem~\ref{t:lace}
can be proved in a purely combinatorial way using
Theorem~\ref{t:stanley}, if one assumes the ratio formula (whose proof
requires geometry).
\end{remark}


\begin{example}
In the case of {\em Fulton rank conditions}, the ranks
\mbox{$r_0,r_1,r_2,r_3,\ldots$} equal
$1,2,3,\ldots,n-1,n,n,n-1,\ldots,3,2,1$ and
$\rr = \rr_w$ is specified by a \mbox{permutation}~$w$ in~$S_{n+1}$.
Theorem~\ref{t:lace} contains
a combinatorial formula found independently in \cite{BKTY02}.
This is because the Zelevinsky permutation $v(\rr_w)$ has as many
diagonal $\times$~entries as will fit in superantidiagonal blocks, and
the southeast corner (last $n+1$ block rows and block columns) is a
block version of the permutation~$w$ itself (rotated by $180^\circ$).
See \cite{Yong03} for details.\qed%
\end{example}

\begin{cor}[Component formula for quiver polynomials---Stanley version]
\label{c:specialize}
The quiver polynomial can be expressed as the specialization
$$
\begin{array}{@{}rcccl@{}}
  \dis \QQ_\rr(\xx-\oxx) &=&\dis \FF_\rr(\xx_\rr-\oxx_\rr) &=& \dis
  \sum_{\ww \in W(\rr)} F_\ww(\xx_\rr-\oxx_\rr)
\end{array}
$$
of the \bem{quiver function}~$\FF_\rr(\xx-\oxx)$ to a sequence $\xx_\rr$
of finite alphabets of size $r_0,\ldots,r_{\!n}$.
\end{cor}
\begin{proof}
Since Theorem~\ref{t:BFtoKMS} holds for any ranks~$\rr$, it holds in
particular for $m+\rr$ when $m \gg 0$.  Moreover, it is shown
in~\cite{BF} that
\begin{eqnarray} \label{eq:clastab}
  c_\bla(m+\rr) &=& c_\bla(\rr)
\end{eqnarray}
for all $\bla$ and $m \geq 0$.  Specializing Theorem~\ref{t:stanley}
to $\yy = \xx$, we may take the limit there with $\sum c_\bla(\rr)
s_\bla(\xx_{m+\rr}-\oxx_{m+\rr})$ in place of~$\QQ_{m+\rr}(\xx-\oxx)$,
by Theorem~\ref{c:ratio} and Theorem~\ref{t:BFtoKMS}.  We conclude
that
\begin{eqnarray} \label{eq:Fx}
  \FF_\rr(\xx-\oxx) &=& \sum_\bla c_\bla(\rr) s_\bla(\xx-\oxx)
\end{eqnarray}
as power series in the ring of functions that are symmetric in each of
the infinite sets $\xx^0,\dotsc,\xx^n$ of variables.  Specializing
each~$\xx^j$ in~\eqref{eq:Fx} to have $r_{\!j}$ nonzero variables
yields $Q_\rr(\xx-\oxx)$ by Theorem~\ref{t:BFtoKMS}, and it yields
$\sum_{\ww \in W(\rr)} F_\ww(\xx-\oxx)$ by Theorem~\ref{t:lace}.%
\end{proof}

\begin{remark} \label{rk:stable}
As a consequence of Corollary~\ref{c:specialize}, there is no need to
define `stable quiver polynomials', just as there is no need to define
`stable Schur polynomials'.  The fact that Conjecture~\ref{m>0} is
actually true, as proved in \cite{BuchAltSign}, means that we are
correct to avoid the
term `stable double quiver polynomial' in favor of the more apt
`double quiver function'.
\end{remark}

\begin{remark} \label{rk:double}
We chose our definition of the double quiver polynomial
$\QQ_\rr(\xx-\oyy)$ via the ratio formula from among four
possibilities that we considered, three of which turned out to be
pairwise distinct.
In hindsight we could equally well have chosen the candidate obtained
by replacing~$\oxx$ with~$\oyy$ in the expression $\sum
c_\bla(\rr)s_\bla(\xx_\rr-\oxx_\rr)$ from Theorem~\ref{t:BFtoKMS},
% ...
or equivalently by~(\ref{eq:Fx}), in the expression $\sum
F_\ww(\xx_\rr-\oxx_\rr)$ from Corollary~\ref{c:specialize}.  However,
it was the algebraic and combinatorial properties derived from the
ratio formula that made $\QQ_\rr(\xx-\oyy)$ useful to us.  Moreover,
our methods did not prove the fact that this other candidate is equal
to our double quiver polynomials; this is the content of
Conjecture~\ref{m>0}, which we now know to hold only because of its
subsequent proof by Buch \cite{BuchAltSign}.

One of our other two candidates for double quiver polynomial was
obtained by replacing~$\oxx$ with~$\oyy$ on the right hand side of
Corollary~\ref{c:component}.  As we explain in Remark~\ref{rk:SSww},
below, this turned out to be less natural than we had expected from
seeing the geometric degeneration of quiver loci to unions of products
of matrix Schubert varieties.

The remaining candidate for double quiver polynomial was actually an
ordinary quiver polynomial, but for a quiver of twice the length as
the original.  Given a lacing diagram~$\ww$ with rank array~$\rr$, the
``doubled'' ranks $\rr^2$ were defined as the rank array of a
``doubled'' lacing diagram~$\ww^2$.  To get~$\ww^2$ from~$\ww$,
elongate each column of dots in~$\ww$ to a ``ladder'' with two
adjacent columns (of the same height) and horizontal rungs connecting
them.  Replacing the even-indexed alphabets in $\QQ_{\rr^2}(\xx-\oxx)$
with $\yy$ alphabets yielded our final candidate.

It is readily verified that the three candidates for double quiver
polynomial differ when $n=3$, the dimension vector is $(1,2,1)$, and
the Zelevinsky permutation is $2143 \in S_4$.%
\end{remark}

\begin{remark} \label{rk:SSww}
We now have formulae for quiver polynomials $\QQ_\rr(\xx-\oxx)$ in
terms of
\begin{itemize}
\item
Stanley symmetric functions $F_\ww(\xx_\rr-\oxx_\rr)$ in finite
alphabets, and

\item
Schubert polynomials $\SS_\ww(\xx-\oxx)$, which a~priori involve
finite alphabets.
\end{itemize}
We also have a formula for double quiver functions $\FF_\rr(\xx-\oyy)$
in terms of double Stanley symmetric functions $F_\ww(\xx-\oyy)$, in
infinite alphabets.  However, the double quiver
polynomial~$\QQ_\rr(\xx-\oyy)$ {\em never}\/ equals the sum $\sum
\SS_\ww(\xx-\oyy)$ of double Schubert polynomials for minimum length
lacing diagrams $\ww \in W(\rr)$, unless the sum only has one term,
even though
\begin{itemize}
\item
setting $\yy = \xx$ in this sum yields Corollary~\ref{c:component},
and

\item
taking limits for uniformly growing ranks yields Theorem~\ref{t:lace}.
\end{itemize}
This failure does not disappear by restricting to ranks $m+\rr$ for~$m
\gg 0$, either: indeed, $\SS_\ww(\xx-\oyy)$ equals the double Schubert
polynomial~$\SS_\ol\ww(\xx-\oyy)$ for a certain permutation $\ol\ww
\in S_d$ constructed from~$\ww$, so linear independence of double
Schubert polynomials prevents a direct double generalization of
Corollary~\ref{c:component}.  We view this as evidence for the
naturality of double quiver polynomials as we defined them (rather
than as $\sum{}_{\!\ww \in W(\rr)} \SS_\ww(\xx-\oyy)$, for example),
and even stronger evidence for the naturality of double quiver
functions.  The recent proof of our
Conjecture~\ref{m>0} in \cite{BuchAltSign} cements our belief in this
naturality.

The reason that $\SS_\ww(\xx-\oyy)$ equals an honest double Schubert
polynomial $\SS_\ol\ww(\xx-\oyy)$
%
\begin{excise}{%
  We will also need to consider some permutations that are not
  Zelevinsky permutations, but still the diagrams of these will contain
  all cells strictly above the block superantidiagonal.  More
  specifically, we shall require permutations constructed from partial
  permutation lists as follows.

  \begin{lemma} \label{l:ww}
  Suppose $\ww = (w_1,\ldots,w_n)$ is a list of partial permutation
  matrices that fit inside the superantidiagonal blocks in block rows
  \mbox{$0,\ldots,n-1$}.  Also denote by\/~$\ww$ the {$d \times d$}
  partial permutation matrix obtained by placing each~$w_i$ in the
  northwest corner of the superantidiagonal block in
  row~\mbox{$i-1$}. Then $\ww$ can be completed uniquely to a
  permutation $\ol \ww \in S_d$ by adding $\times$ entries in or below
  the block antidiagonal in such a way that all cells in the diagram
  of\/~$\ol \ww$ lie in or above the block superantidiagonal.
  \end{lemma}
  \begin{proof}
  Find a $\times$ entry in the top row (not block row) of~$\ww$, or
  else place a new $\times$ entry in the northwest corner of the top
  antidiagonal block.  Then delete the top row along with the column
  containg the top row~$\times$ entry, and continue by induction.%
  \end{proof}
}\end{excise}%
comes from a direct geometric connection between the orbit closure
$\OO(\ww)$ inside $\homv$ and the matrix Schubert variety $\ol
X_\ol\ww$ inside~$\endv$.  In fact, it can be shown that there is a
flat (but not Gr\"obner) degeneration of the matrix Schubert variety
$\ol X_{v(\rr)}$ to a generically reduced union $\ol X_{v(\rr)}(0)$ of
matrix Schubert varieties for permutations~$\ol\ww$ associated to
lacing diagrams $\ww \in W(\rr)$.  However, this flat family is only
equivariant for the multigrading by~$\ZZ^d$,
not~$\ZZ^{2d}$.
\begin{excise}{%
  Concretely, it is constructed by ``smearing'' the quiver
  degeneration $\til\Omega_\rr$ into~$\endv$ as follows.

  Let~$U$ be the block lower-triangular matrices in~$\gld$, and recall
  that~$B$ is the lower-triangular Borel group inside the Levi
  factor~$L$.  Define $N_+^\circ = {w_0}\ww_0 N_+ \ww_0{w_0}$ to be
  the unipotent upper-triangular matrices in the block-column reversed
  Levi factor ${w_0}\ww_0 L \ww_0{w_0}$.  The family $\til\Omega_\rr$
  sits inside $\homv \times \AA^1$, so we can write~$\til{Y}_\rr$ for
  its Zelevinsky image inside $\endv \times \AA^1$.  Then the family
  degenerating to~$\ol X_{v(\rr)}(0)$ is the closure inside $\endv
  \times \AA^1$ of the family $\UB \til{Y}_\rr N_+^\circ$.  The key
  point now rests in checking that $\UB \cZ(\OO(w)) N_+^\circ$ is a
  dense subvariety of a matrix Schubert variety~$\ol X_\ol\ww$, and
  similarly that $\UB Y_\rr N_+^\circ$ is a dense subvariety of~$\ol
  X_{v(\rr)}$.%
}\end{excise}%
\end{remark}

\part{Quiver constants}

\label{part:coeffs}
\setcounter{section}{0}
\setcounter{thm}{0}
\setcounter{equation}{0}

Using the
component formula for double quiver functions, we show that the quiver
constants~$c_\bla(\rr)$ arise as coefficients in the expansion of the
Zelevinsky Schubert polynomial~$\SS_{v(\rr)}$ into Demazure
characters, and consequently the quiver constants~$c_\bla(\rr)$ equal
the corresponding Stanley coefficients~$\alpha^\la_w$.  Using the
interpretation of Stanley coefficients as enumerating peelable
tableaux \cite{RSplact95,RSpeelable98}, we conclude a combinatorial
\mbox{formula for quiver constants}.

\section{Demazure characters}

\label{sec:demazure}

We recall here the type $A_\infty$ Demazure characters \cite{Dem74}
following \cite{Mac}.  Using the ordinary divided
difference~$\partial_i$ from~\eqref{eq:partial}, the \bem{Demazure
operator}~$\pi_i$, also called the \bem{isobaric divided difference},
is the linear operator on $\ZZ[x_1,x_2,\dotsc]$ defined by
\begin{eqnarray*}
  \pi_i f &=& \partial_i (x_i f).
\end{eqnarray*}
Let $\sigma_i$ be the transposition that exchanges $x_i$
and~$x_{i+1}$.  The operator~$\pi_i$ is idempotent and, like
$\partial_i$, commutes with multiplication by any
$\sigma_i$-symmetric polynomial:
\begin{eqnarray} \label{eq:pisymm}
\begin{array}{cccc}
  \partial_i (fg) = f \partial_i(g) &\text{and}& \pi_i (fg) = f
  \pi_i(g) \qquad\text{if } \sigma_i f = f.
\end{array}
\end{eqnarray}
Since $\pi_i 1 = 1$ it follows that
\begin{eqnarray} \label{eq:piinvariant}
  \pi_i f &=& f \qquad\text{if } \sigma_i f = f.
\end{eqnarray}

Given a permutation $w \in \bigcup_n S_n$, let
$w=\sigma_{i_1}\dotsm \sigma_{i_\ell}$ be any reduced
decomposition of~$w$, that is, a factorization of~$w$ into a
minimum number of simple reflections~$\sigma_i$.  Define the
operators $\pi_w = \pi_{i_1}\dotsm \pi_{i_\ell}$.  It is
independent of the factorization.

Let $\beta=(\beta_1,\dotsc,\beta_k)\in\NN^k$ be a sequence of
nonnegative integers, and $\beta_+$ the sequence obtained by
sorting~$\beta$ into weakly decreasing order.  If~$w\in S_k$ is the
shortest permutation satisfying $\beta=w\beta_+$, then define the
\bem{Demazure character} $\kp_\beta\in\ZZ[x_1,\dotsc,x_k]$~by
\begin{eqnarray*}
  \kp_\beta &=& \pi_w \xx^{\beta_+}.
\end{eqnarray*}
For a partition $\la\in\NN^k$ and $w_0^{(k)}\in S_k$ the longest
permutation, we have
\begin{eqnarray}
  \kp_\la &=& \xx^\la \\
\label{eq:DemSchur}
  \kp_{w_0^{(k)}\la} &=& s_\la(x_1,x_2,\dotsc,x_k).
\end{eqnarray}
Alternatively one may give the recursive definition
\begin{eqnarray} \label{eq:demrec}
\kp_\beta &=& \begin{cases}
  \xx^\beta & \text{if $\beta_1\ge\beta_2\ge\dotsm \ge \beta_k$} \\
  \pi_i \kp_{\sigma_i \beta} & \text{if $\beta_i < \beta_{i+1}$.}
\end{cases}
\end{eqnarray}
It can be shown that
\begin{eqnarray} \label{eq:dempi}
  \pi_i \kp_\beta &=& \begin{cases}
  \kp_\beta & \text{if } \beta_i \le \beta_{i+1} \\
  \kp_{\sigma_i \beta} & \text{if } \beta_i\ge\beta_{i+1}.
  \end{cases}
\end{eqnarray}

\begin{remark} \label{rem:KeyBasis}
The Demazure characters $\{\kp_\beta\mid \beta\in\NN^k\}$ form a
$\ZZ$-basis of the polynomial ring $\ZZ[x_1,x_2,\dotsc,x_k]$ since
\begin{eqnarray*}
  \kp_\beta &=& \xx^\beta + \text{terms lower in reverse lexicographic
  order}.
\end{eqnarray*}
\end{remark}

Demazure characters interpolate between dominant monomials and
Schur polynomials.  They have a tableau realization, due to
Lascoux and Sch\"utzenberger \cite{LS90}.  Before describing it,
we review some conventions regarding tableaux; for more background
and unexplained terminology, see \cite{Ful97}.

Our (semistandard) Young tableaux will be drawn with origin at the
northwest corner, weakly increasing along each row, and strictly
increasing down each column.  A \bem{word} is a sequence of positive
integers called ``letters''.  A \bem{column word} is one that strictly
decreases.  In what follows, each tableau~$t$ is identified with its
\bem{column reading word} $t_1t_2\dotsm$, where $t_i$ is the column
word obtained by reading the letters upward in the $i^\th$ column
of~$t$.  Given a word~$u$, denote by~$[u]$ the unique tableau
equivalent to~$u$ under the {\em Knuth-equivalence relations}\/ on
words: $acb \sim cab$ for $a \leq b < c$ and $bac \sim bca$ for $a < b
\leq c$ \cite[Section~2.1]{Ful97}.

The \bem{weight} of a tableau~$t$ is the sequence $\wgt(t) =
(\beta_1,\beta_2,\dotsc)$, where $\beta_j$ is the number of entries
equal to~$j$ in~$t$.  For each \bem{composition}
$\beta=(\beta_1,\beta_2,\dotsc)$, by which we mean a sequence of
nonnegative integers that is eventually zero, the \bem{key tableau}
$\key(\beta)$ of weight~$\beta$ is the unique tableau of
shape~$\beta_+$ and weight~$\beta$.

Suppose $t$ has shape~$\la$, so that its $j^\th$ column~$t_j$ has
length~$\la'_j$ for $1\le j\le r=\la_1$, where $\la'$ is the partition
\bem{conjugate} to~$\la$.  Given any reordering
$\gamma=(\gamma_1,\dotsc,\gamma_r)$ of these column lengths, there is
a unique word $u=u(t,\gamma)$ with a factorization $u=u_1u_2\dotsm
u_r$ into column words $u_j$ of length $\gamma_j$, such that $[u]=t$.
Such a word is called \bem{frank} in \cite{LS90}.

The leftmost column factor~$u_1$ in the word $u(t,\gamma)$ actually
depends only on~$t$ and~$\gamma_1$.  Define $\lc_j(t) = u_1$ to be the
$j^\th$ \bem{left column word} of~$t$ when $\gamma_1$ equals the
length~$\la'_j$ of the $j^\th$ column of~$t$.  Gathering together all
left column words produces the \bem{left key} of~$t$, which is the
tableau $\leftkey(t)$ whose $j^\th$ column is given by~$\lc_j(t)$ for
$j = 1,\ldots,r$.

Similarly, the rightmost column factor~$u_r$ in $u(t,\gamma)$ depends
only on~$t$ and~$\gamma_r$.  Define $\rc_j(t) = u_r$ to be the $j^\th$
\bem{right column word} of~$t$ when $\gamma_r = \la'_j$, and construct
the tableau $\rightkey(t)$, called the the \bem{right key} of~$t$,
whose $j^\th$ column is given by~$\rc_j(t)$ for $j = 1,\ldots,r$.  The
left and right keys of $t$ are key tableaux of the same shape as~$t$.

\begin{excise}{%
\begin{remark} \label{rk:swap}
It can be shown that if $\gamma$ and $\gamma'$ differ by exchanging
the $j^\th$ and $(j+1)^\st$ parts, then the frank words
$u(t,\gamma)=u_1u_2\dotsm u_r$ and $u(t,\gamma')=u'_1u'_2\dotsm u'_r$
satisfy $u_i=u'_i$ for $i\not\in\{j,j+1\}$, and the two-column words
$u_j u_{j+1}$ and $u'_j u'_{j+1}$ are easily computed from each other
using a jeu de taquin inside a two-column rectangle.  Indeed, one of
the two-column words is a tableau word and the other is the column
reading word of a two-column antitableau (skew tableau with unique
southeast corner).
\end{remark}
We refer to the process of passing between $u$ and $u'$ as
``column-swapping''.
}\end{excise}%

For tableaux $s$ and~$t$ of the same shape, write $s\le t$ to mean
that every entry of~$s$ is less than or equal to the corresponding
entry of~$t$.  In \cite{LS90}, the equation
\begin{eqnarray} \label{eq:keytableau}
  \kp_\beta &=& \sum_{\substack{\rightkey(t) \le \key(\beta)}}
  \xx^{\wgt(t)}
\end{eqnarray}
is given as a formula for the Demazure character.  This implies two
positivity properties:
\begin{eqnarray} \label{eq:dempos}
  \kp_\beta \in\NN[\xx] && \text{for all compositions } \beta,
\text{ and}\\    \label{eq:demadd}
  \kp_{w \la} - \kp_{v \la} \in \NN[\xx] &&\text{if } v\le w \text{
  and } \la \text{ is a partition}.
\end{eqnarray}

It is well-known that
\begin{eqnarray} \label{eq:schur}
  s_\la(\xx) &=& \sum_{\substack{\text{tableaux $t$} \\[.8mm]
  \shape(t)=\la}}
  \xx^{\wgt(t)}
\end{eqnarray}
is a formula for the Schur function~$s_\la$ in an infinite
alphabet~$\xx$.
{}From~\eqref{eq:keytableau} and~\eqref{eq:schur} it follows that
\begin{eqnarray} \label{eq:LimKey}
  s_{\beta_+} &=& \lim_{m\to\infty} \kp_{(0^m,\beta)},
\end{eqnarray}
where $(0^m,\beta)$ is the composition obtained from~$\beta$ by
prepending $m$~zeros.
\begin{excise}{%
  Applying Remark~\ref{rk:swap} yields the following easy result.
  
  \begin{lemma} \label{lem:LeftKeyWt}
  Let $i$ be a letter appearing in each of the first $k$ columns of
  the tableau~$t$.  If $\beta=\wgt(\leftkey(t))$ is the weight of the
  left key tableau of~$t$, then $\beta_i \geq k$.
  \end{lemma}
  
  \begin{proof}
  Equivalently it must be shown that $\lc_j(t)$ contains $i$ for $1\le
  j\le k$.  $\lc_j(t)$ can be computed by swapping the $j^\th$ column
  of $t$ all the way to the left.  During each column swap, both
  columns contain the letter $i$ before the swap and therefore must
  both contain $i$ afterwards.  In particular the leftmost column will
  contain $i$ throughout the entire process.
  \end{proof}
}\end{excise}%
%

\section{Schubert polynomials as sums of Demazure characters}
							\label{sec:sums}

We now recall the expansion of a Schubert polynomial as a positive
sum of Demazure characters and deduce some special properties of
this expansion for Zelevinsky permutations.

\begin{thm}[{\cite{LS89},\cite[Theorem~29]{RSplact95}}]\label{t:SchubKeyWeak}
For any (partial) permutation~$w$, there is a multiset $M(w)$ of
compositions such that
\begin{eqnarray*}
  \Schub_w(\xx) &=& \sum_{\beta \in M(w)} \kp_\beta.
\end{eqnarray*}
\end{thm}

The Schubert polynomial $\SS_w(\xx)$ here is obtained by setting $\oyy
= \0$ in the double Schubert polynomial $\SS_w(\xx-\oyy)$ from
Section~\ref{sec:schub}.  The multiset $M(w)$ has been described
explicitly using reduced word tableaux in \cite{LS89} and peelable
tableaux in~\cite{RSplact95}.  (We shall define the latter in
Section~\ref{sec:peelable}.)  The expansion in
Theorem~\ref{t:SchubKeyWeak} is a refinement of~\eqref{eq:FtoS}, in
the sense that taking limits in
Theorem~\ref{t:SchubKeyWeak}, expresses the Stanley
coefficient~$\alpha^\la_w$ in terms of~$M(w)$:

\begin{cor} \label{c:stanley2schur}
The coefficient $\alpha^\la_w$ in the expansion of a Stanley symmetric
function in Schur functions is the number of compositions $\beta\in
M(w)$ such that $\beta_+ = \la$.
\end{cor}
\begin{proof}
By the definition of the Stanley symmetric function
and~\eqref{eq:LimKey}, it is enough to show that there is a bijection
$M(w)\mapsto M(1+w)$ given by $\beta\mapsto (0,\beta)$.
%
Supposing that $w$ lies in~$S_k$, the definition of $\SS_w$ in
Section~\ref{sec:schub} says that $\SS_w(\xx) = \partial_{w'}
\SS_{w_0^{(k)}}(\xx)$, where $w' = w^{-1}w_0^{(k)}$ and $w_0^{(k)} \in
S_k$ is the longest permutation.
Note that
\begin{eqnarray*}
 1+w&=&\sigma_1\sigma_2\cdots\sigma_k w\sigma_k\cdots\sigma_2\sigma_1,
\end{eqnarray*}
and that $w_0^{(k+1)}=\sigma_1\sigma_2\cdots\sigma_i w_0^{(k)}$.
Also, in the expression
$(1+w)^{-1}w_0^{(k+1)}=\sigma_1\sigma_2\cdots\sigma_k w'$, multiplying
by each reflection~$\sigma_i$ lengthens the permutation.  It follows
that
\begin{eqnarray*}
\SS_{1+w}
  &=&\partial_{(1+w)^{-1} w_0^{(k+1)}}\big(\SS_{w_0^{(k+1)}}(\xx)\big)
\\&=&\partial_{\sigma_1\sigma_2\dotsm\sigma_k w'}
	\big(\SS_{w_0^{(k+1)}}(\xx)\big)
\\&=&\partial_1\partial_2 \dotsm \partial_k \partial_{w'}
	\big(x_1\cdots x_k \cdot \SS_{w_0^{(k)}}(\xx)\big)
\\&=&\partial_1 x_1 \partial_2 x_2 \dotsm\partial_k x_k
	\partial_{w'}\big(\SS_{w_0^{(k)}}(\xx)\big)
\\&=&\pi_1\dotsm\pi_k \SS_w
\end{eqnarray*}
using the definitions, \eqref{eq:pisymm}, and \eqref{eq:partial}.

Let $\beta$ be a composition in~$M(w)$.  Since $w$ lies in the
symmetric group~$S_k$ on~$k$ letters, the Schubert polynomial
$\Schub_w(\xx)$ involves at most the variables $x_1,\dotsc,x_{k-1}$.
Hence only the first $k-1$ parts of~$\beta$ can be nonzero.
By~\eqref{eq:dempi}, $\pi_1\dotsm\pi_k \kp_\beta = \kp_{(0,\beta)}$,
as a zero in the $(k+1)^\st$ position is swapped to the first
position.  Since the Demazure characters form a basis
(Remark~\ref{rem:KeyBasis}), the map $\beta\mapsto (0,\beta)$ is a
bijection $M(w)\mapsto M(1+w)$.%
\end{proof}


We now recall the dominance bounds on the shapes that occur in the
Schur function expansion of Stanley symmetric functions.  For any
diagram~$D$, by which we mean a set of pairs of positive integers, let
$D\north$ be the diagram obtained from $D$ by top-justifying the cells
in each column.  Similarly, let $D\west$ be obtained by
left-justifying each row of~$D$.  The \bem{dominance partial order} on
partitions of the same size, written as $\la \rdom \mu$, is defined by
the condition that $\la_1+\dotsm+\la_i\ge \mu_1+\dotsm+\mu_i$ for
all~$i$.

\begin{prop}[\cite{St}] \label{pp:StanDom}
If the Stanley coefficient $\alpha^\la_w$ is nonzero, then
\begin{equation*}
   D(w) \west\north\;\,\ldom\,\la\,\ldom\,D(w) \north\west.
\end{equation*}
\end{prop}

We require two more lemmata.

\begin{lemma} \label{lem:SchubDemAscent}
If $w(i)<w(i+1)$ then $\beta_i \le \beta_{i+1}$ for every composition
$\beta\in M(w)$.
\end{lemma}
\begin{proof}
Suppose $w(i)<w(i+1)$.  Since $\sigma_i \circ \partial_i =
\partial_i$, it follows from the definition of Schubert polynomial
that $\Schub_w(\xx)$ is $\sigma_i$-invariant.  Applying $\pi_i$ to
Theorem~\ref{t:SchubKeyWeak} and using \eqref{eq:piinvariant} along
with~\eqref{eq:dempi}, it follows that $\sum_{\beta\in M(w)} \kp_\beta
= \sum_{\beta\in M(w)} \kp_{\beta'}$, where $\beta'=\beta$ if $\beta_i
\le \beta_{i+1}$ and $\beta'=\sigma_i \beta$ otherwise.  But
$\kp_{\beta'}-\kp_\beta$ is a polynomial with nonnegative integer
coefficients by~\eqref{eq:demadd}.  It follows that $\beta'=\beta$ and
hence that $\beta_i\le \beta_{i+1}$ for all $\beta\in M(w)$.%
\end{proof}

\begin{lemma} \label{lem:SchubDemSub}
Suppose $\gamma$ is a composition such that $\xx^\gamma$
divides~$\SS_w(\xx)$.  Then $\gamma \leq \beta$ for all compositions
$\beta\in M(w)$.
\end{lemma}
\begin{proof}
Suppose $\xx^\gamma$ divides~$\SS_w(\xx)$.  Then $\xx^\gamma$ divides
the sum of Demazure characters on the right hand side of
Theorem~\ref{t:SchubKeyWeak}.  By~\eqref{eq:dempos} it follows that
$\xx^\gamma$ divides every monomial of every Demazure character
$\kp_\beta$ for $\beta\in M(w)$.  By Remark \ref{rem:KeyBasis} the
reverse lexicographic leading monomial in~$\kp_\beta$ is~$\xx^\beta$,
so $\gamma \leq \beta$.%
\end{proof}

We now apply the above results to Zelevinsky permutations.

\begin{prop} \label{pp:ZelDem}
For a fixed dimension vector $(r_0,r_1,\dotsc,r_n)$, identify the
diagrams $D_\homv$ and~$D(\Omega_0)$ from Definition~\ref{d:diagram}
with partitions of those shapes.  Let $\rr$ be a rank array with the
above dimension vector.  Then for all compositions $\beta\in
M(v(\rr))$,
\begin{eqnarray} \label{eq:DomBound}
  &\makebox[0ex]{$\dis\quad\ D_\homv\subset\beta_+\subset D(\Omega_0)$
	\makebox[0ex][l]{\quad and}}&
\\\label{eq:blocklong}
  \beta &=& \ww_0 \beta_+,
\end{eqnarray}
where $\ww_0$ is the block long permutation, reversing the row indices
within each block row.
\end{prop}
\begin{proof}
Let $\beta\in M(v(\rr))$.  Equation~\eqref{eq:blocklong} is a
consequence of Lemma~\ref{lem:SchubDemAscent} and
Proposition~\ref{p:zel}, as Zelevinsky permutations have no descents
in each block row.

Both of the partition diagrams \mbox{$D_\rr \north\west$} and
\mbox{$D_\rr\west\north$} are contained in~$D(\Omega_0)$, by
Definition~\ref{d:diagram} along with~\eqref{eq:heightsums},
\eqref{eq:widthsums}, \eqref{eq:heightincr}, and \eqref{eq:widthincr}.
Proposition \ref{pp:StanDom} implies that $\beta_+\subset
D(\Omega_0)$.

Since $D_\homv \subset D_\rr$, the monomial $\xx^{D_\homv}$ divides
$\Schub_{v(\rr)}$.  By Lemma~\ref{lem:SchubDemSub}, $D_\homv \le
\beta$. But the partition $D_\homv$ is constant in each block row.
With \eqref{eq:blocklong}, it follows that $D_\homv \subset \beta_+$.
\end{proof}

Finally, we observe that each Demazure character appearing in
$\Schub_{v(\rr)}$ is the monomial $\xx^{D_\homv}$ times a product
of Schur polynomials.

\begin{defn} \label{d:delete}
Assume that $\la=\beta_+$ for some composition $\beta\in M(v(\rr))$.
The skew shape $\la/D_\homv$ afforded by~\eqref{eq:DomBound} consists
of a union of partition diagrams in which the partition~$\la_i$ is
contained in the $r_{i-1}\times r_i$ rectangle comprising the $i^\th$
block row of $D(\Omega_0)/D_\homv$.  In this situation, we say that
the partition list $\bla=(\la_1,\la_2,\dotsc,\la_n)$ is obtained by
\bem{deleting~$D_\homv$ from~$\la$}, and we write $\bla =
\la-D_\homv$.
\end{defn}

\begin{prop} \label{pp:ZSchubDem}
Suppose that $\beta$ is a composition in~$M(v(\rr))$, that
$\la=\beta_+$, and that $\bla$ is obtained from~$\la$ by
deleting~$D_\homv$.  Then
\begin{eqnarray*}
  \kp_\beta &=& \xx^{D_\homv} \prod_{i=1}^n s_{\la_i}(\xx^{i-1}).
\end{eqnarray*}
\end{prop}
\begin{noqed}
Break $D_\homv$ into a sequence $(\mu_0,\dotsc,\mu_n)$ of shapes,
where $\mu_i\in\NN^{r_i}$ is rectangular (all of its parts are equal)
and $\mu_n=(0^{r_n})$.  Let $w_0^i$ be the longest element of the
symmetric group~$S_{r_i}$ acting on the $i^\th$ block of row
indices.  By~\eqref{eq:pisymm} and~\eqref{eq:DemSchur} we have
$$
\begin{array}{@{}rcl@{}}
\multicolumn{3}{@{}c@{}}{\mbox{}\hspace{\textwidth}\mbox{}}\\[-3ex]
\displaystyle\kp_\beta
  &=&\displaystyle
        \pi_{\ww_0} \xx^\la \cr
  &=&\displaystyle
        \prod_{i=0}^{n-1} \pi_{w_0^i} (\xx^i)^{\la_{i+1}+\mu_i} \cr
  &=&\displaystyle
        \prod_{i=0}^{n-1} (\xx^i)^{\mu_i} s_{\la_{i+1}}(\xx^i) \cr
\mbox{}\hspace{.422\textwidth}\mbox{}
  &=&\displaystyle
        \xx^{D_\homv} \prod_{i=1}^n s_{\la_i}(\xx^{i-1}).\hfill\square
\end{array}
$$
\end{noqed}

\section{Quiver constants are Stanley coefficients}

\label{sec:count}

Having an expansion of the double quiver function as a sum of products
of Stanley symmetric functions automatically produces an expansion as
a sum of products of Schur functions.  Next, in Theorem~\ref{t:bla},
we shall see that the coefficients in this expansion are the
Buch--Fulton quiver constants. Using this result, we show in
Theorem~\ref{t:qc=ss} that quiver constants are special cases of
Stanley coefficients.

For a list $\bla=(\la_1,\la_2,\dotsc,\la_n)$ of partitions and two
sequences $\xx=(\xx^0,\dotsc,\xx^n)$ and $\oyy=(\yy^n,\dotsc,\yy^0)$
of infinite alphabets, let
\begin{eqnarray*}
  s_\bla(\xx-\oyy) &=& \prod_{i=1}^n s_{\la_i}(\xx^{i-1}-\yy^i)
\end{eqnarray*}
be the product of Schur functions in differences of alphabets.
This notation parallels that with $\yy = \xx$ in~\eqref{eq:slax},
and for products of Stanley symmetric functions in
Section~\ref{sec:stanley}. We follow the conventions
after~\eqref{F} for finite alphabets.

\begin{thm} \label{t:bla}
If~$c_\bla(\rr)$ is the quiver constant from
Theorem~\ref{t:BFtoKMS}, then
\begin{eqnarray*}
  \FF_\rr(\xx-\oyy) &=& \sum_\bla c_\bla(\rr) s_\bla(\xx-\oyy).
\end{eqnarray*}
Moreover, using the Stanley coefficients $\alpha^\la_w$
from~\eqref{eq:FtoS} and writing $\alpha^\bla_\ww = \prod_{i=1}^n
\alpha^{\la_i}_{w_i}$, we get
\begin{eqnarray} \label{eq:QuiverCoefPerm}
  c_\bla(\rr) &=& \sum_{\ww\in W(\rr)} \alpha^\bla_\ww.
\end{eqnarray}
\end{thm}
\begin{proof}
Expanding the right hand side of Theorem~\ref{t:lace} into Schur
functions yields
\begin{eqnarray*}
  \FF_\rr(\xx-\oyy) &=& \sum_{\ww\in W(\rr)} \sum_{\bla}
  \alpha^\bla_\ww s_\bla(\xx-\oyy).
\end{eqnarray*}
Specializing $\oyy$ to~$\oxx$ yields an expression for
$\FF_\rr(\xx-\oxx)$ as a sum of Schur function products
$s_\bla(\xx-\oxx)$.  Uniqueness in the Main Theorem of~\cite{BF}
implies that the coefficient $\sum_{\ww\in W(\rr)} \sum_\bla
\alpha^\bla_\ww$ on $s_\bla(\xx-\oxx)$ in this expression
is~$c_\bla(\rr)$.%
\end{proof}

\begin{remark}
Theorem~\ref{t:bla} and Corollary~\ref{c:specialize} are the only
results in this paper that logically depend on the Main Theorem
of~\cite{BF}, other than Theorem~\ref{t:BFtoKMS} (which needs the
Main Theorem of~\cite{BF} for its statement).  In other words,
starting from the multidegree characterization of quiver
polynomials in Definition~\ref{d:quivpoly}, the statements and
proofs of all of our other combinatorial formulae---including all
double and stable versions, as well as our combinatorial formula
for the quiver constants to come in
Theorem~\ref{t:PeelQuiver}---are independent from \cite{BF}. It is
only to identify the constants on the right side
of~\eqref{eq:QuiverCoefPerm} as the quiver constants appearing in
the Conjecture from \cite{BF} that we apply
Theorem~\ref{t:BFtoKMS}, and hence the Main Theorem of~\cite{BF}.
\end{remark}

\begin{remark}
The quiver constants $c_\bla(\rr)$ may be computed fairly efficiently
using~\eqref{eq:QuiverCoefPerm}.
\end{remark}

We need a proposition, in which our choice of notation (using
$\xx-\oxx$ in arguments of quiver polynomials) should finally become
clear: we use $\QQ_\rr(\xx)$ and $\FF_\rr(\xx)$ to denote the $\yy =
\0$ specializations (as opposed to $\yy = \xx$ specializations) of the
double quiver polynomial~\mbox{$\QQ_\rr(\xx-\oyy)$} and the double
quiver function~$\FF_\rr(\xx-\oyy)$.  The ``doubleness'' of these
expressions is crucial here, because in the ordinary single case it is
impossible in the summation formulae to set just the second set of
variables in each factor of every summand~to~zero.

\begin{prop} \label{p:QQrr}
If\/~$\rr$ is any rank array, then $\QQ_\rr(\xx) =
\FF_\rr(\xx_\rr)$ is the finite $\xx$~alphabet specialization of
the double quiver function at~$\yy=\0$.
\end{prop}
\begin{proof}
Calculate $\QQ_{m+\rr}(\xx)$ and $\QQ_\rr(\xx)$ using
Theorem~\ref{t:QQrr}.
\begin{excise}{
  Reduced pipe dreams for~$v(m+\rr)$ that happen to lie entirely in
  rows indexed by~$\xx_\rr$ (as opposed to~$\xx_{m+\rr}$) are in
  canonical bijection with~$\rp(v(\rr))$ by
  Proposition~\ref{p:top}. Since the bijection is by horizontal shift,
  corresponding pipe dreams obviously yield the same monomial
  in~$\xx$.  Therefore setting all variables in the set $\xx_{m+\rr}
  \minus \xx_\rr$ to zero in $\QQ_{m+\rr}(\xx)$ yields
  $\QQ_\rr(\xx)$.  Now use Corollary~\ref{c:stanley}.%
}\end{excise}%
%
Consider a reduced pipe dream $D \in \rp(v(m+\rr))$ that happens to
lie entirely in rows indexed by~$\xx_\rr$ (as opposed
to~$\xx_{m+\rr}$). Restricting attention to the $(j-1)^\st$ block row
of~$D$ for the moment, let $w_j$ satisfy $w_j(D) = m+w_j$ using
Corollary~\ref{c:W}. No cross~$\textcross$ can be farther to the left
of the vertical dividing line in that block row (from after
Lemma~\ref{l:enlarge}) than the leftmost~$\textcross$ in the
upside-down (rotated by~$180^\circ$) top pipe dream for~$w_j$.
Therefore, the miniature pipe dream in block row $j-1$ of~$D$ has at
least $m$ blank columns on its left.  By Proposition~\ref{p:top},
horizontally shifting all crosses $m$ cells to the left
therefore induces a canonical bijection onto $\rp(v(\rr)$ from the
set of pipe dreams in $\rp(v(m+\rr))$ that happen to lie entirely
in rows indexed by~$\xx_\rr$.  Corresponding pipe dreams obviously
yield the same monomial in~$\xx$.  Therefore setting all variables
in the set $\xx_{m+\rr} \minus \xx_\rr$ to zero in
$\QQ_{m+\rr}(\xx)$ yields
$\QQ_\rr(\xx)$.  Now use Corollary~\ref{c:stanley}.%
\end{proof}

\begin{thm} \label{t:qc=ss}
The quiver constant $c_\bla(\rr)$ is the Stanley coefficient
$\alpha^\la_{v(\rr)}$, where $\lambda$ is the unique partition such
that the partition list~$\bla$ is obtained from~$\la$ by
deleting~$D_\homv$.
\end{thm}
\begin{proof}
Multiplying through by the implicit denominator in
Proposition~\ref{p:QQrr} yields
\begin{eqnarray} \label{eq:SchubquotBF}
  \SS_{v(\rr)}(\xx) &=& \SS_{v(\homv)}(\xx) \sum_\bla c_\bla(\rr)
  s_\bla(\xx_\rr),
\end{eqnarray}
where $c_\bla(\rr)$ is the quiver constant from
Theorem~\ref{t:BFtoKMS} by Theorem~\ref{t:bla}.
Proposition~\ref{pp:ZSchubDem} shows that the Demazure character
expansion of $\SS_{v(\rr)}(\xx)$ has the same form as the right
side of~\eqref{eq:SchubquotBF}.  Linear independence of the
polynomials $s_\bla(\xx_\rr)$ therefore implies that
\eqref{eq:SchubquotBF} is the Demazure character expansion. Equate
coefficients on each Demazure character in
Theorem~\ref{t:SchubKeyWeak} and~\eqref{eq:SchubquotBF} and apply
Corollary \ref{c:stanley2schur}.
\end{proof}

\begin{remark} \label{rk:y=0}
The comparison of~\eqref{eq:SchubquotBF} with
Theorem~\ref{t:SchubKeyWeak} via Proposition~\ref{pp:ZSchubDem} occurs
at the level of ordinary Schubert polynomials (in~a single set~$\xx$
of variables) rather than double Schubert polynomials (in~$\xx$
and~$\yy$).  This restriction is forced upon us by the fact that the
Demazure character theory behind
Proposition~\ref{pp:ZSchubDem} has not been sufficiently developed in
the double setting.  (See \cite{LasDoubKey} for the beginnings of such
a `double' theory.)  Since the comparison of~\eqref{eq:SchubquotBF}
with Theorem~\ref{t:SchubKeyWeak} is the key point
in the proof of Theorem~\ref{t:qc=ss}, it is our single most
important motivation for developing double versions of quiver
polynomials: we must be able to set $\yy = \0$ for this
comparison.
\end{remark}

\section{Peelable tableaux}

\label{sec:peelable}

Peelable tableaux \cite{RSplact95,RSpeelable98} arose from Stanley's
problem of counting reduced decompositions of permutations \cite{St},
the expansion by Lascoux and Sch\"utzenberger of the Schubert
polynomial as a sum of Demazure characters \cite{LS89}, and from
Magyar's character formula for the global sections of line bundles
over Bott--Samelson varieties~\cite{Mag}.  Here we use them to give
our first combinatorial formula for quiver constants.  Although other
combinatorial interpretations of Stanley coefficients would also
produce formulae, it is the peelable tableaux that will connect to
factor sequences in Section~\ref{part:fs}.

Let $D$ be a \bem{diagram}, meaning a finite set of pairs of
positive integers.  The diagram~$D$ is \bem{northwest} if
$(i_1,j_2)\in D$ and $(i_2,j_1)\in D$ implies $(i_1,j_1)\in D$ for
$i_1<i_2$ and $j_1<j_2$; thus, if two cells lie in~$D$, then so
does the northwest corner of the smallest rectangle containing
both.  Permutation diagrams are always northwest.  The diagram
$D(\la)$ of a partition $\la$ consists of $\la_i$ left-justified
cells in the $i^\th$ row for each~$i$.

For tableaux $Q$ and~$P$, write $Q\supset P$ if $Q$ contains~$P$
in its northwest corner.  Write $Q-P$ for the skew tableau
obtained by removing the subtableau~$P$ from~$Q$.

A column of a diagram can be identified with the set of indices of
rows in which it has a cell.  This subset can also be identified
with the decreasing word consisting of these row indices.  As in
Section~\ref{sec:demazure}, denote by $[u]$ the unique tableau
Knuth-equivalent to a word~$u$.

Now consider a northwest diagram~$D$.  A~tableau $Q$ is
\bem{$D$-peelable} provided that:
\begin{numbered}
\item $Q$ is the empty tableau when $D$ is the empty diagram; or
else \item $Q\supset C$ and $[Q-C]$ is $(D-C)$-peelable when $C$
is the first nonempty column~of~$D$.
\end{numbered}
We refer to the map $Q\mapsto [Q-C]$ as \bem{peeling}, and denote
by~$\Peel(D)$ the set of $D$-peelable tableaux of partition shape.

\begin{remark} \label{rem:PeelWeight}
The weight of any $D$-peelable tableau is the \bem{code} of~$D$,
meaning the sequence $\code(D)=(c_1,c_2,\dotsc)$ in which $c_i$ is
the number of cells in row~$i$ of~$D$.
\end{remark}

\begin{example} \label{ex:peel}
The tableau in Fig.~\ref{f:peel} is $D_\rr$-peelable as exhibited,
where $\rr$ is the rank array in Example~\ref{ex:rank}, and the
diagram~$D_\rr$ is in Example~\ref{ex:zel}.  The sequences of
removed columns~$C$ (appearing in boldface below) are given by the
columns of~$D_\rr$, after each cell has been labeled by its row
index (this labeling of cells in~$D_\rr$ is depicted explicitly in
\begin{figure}
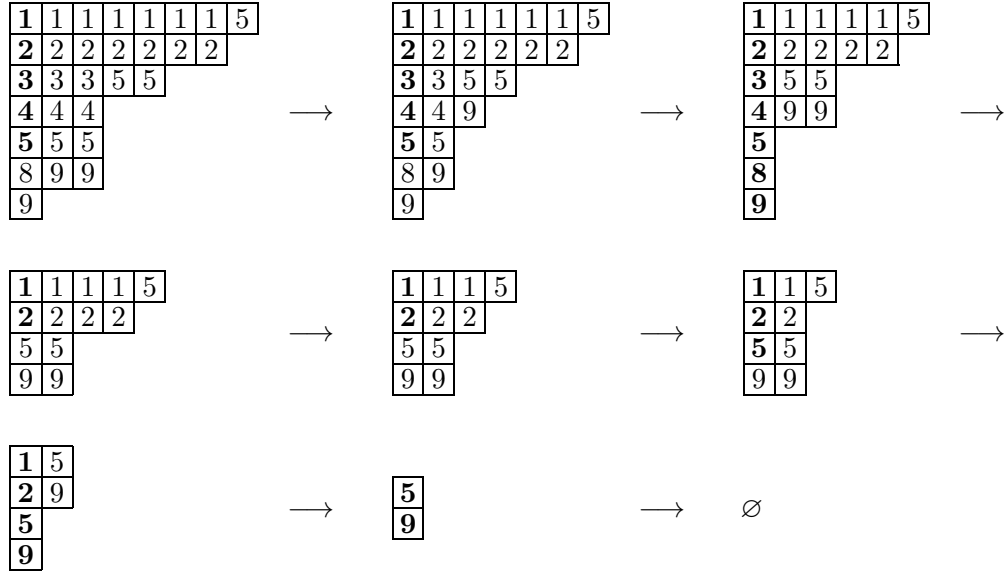

\newcommand{\pone}{\mathbf{1}}%
\newcommand{\ptwo}{\mathbf{2}}%
\newcommand{\pthree}{\mathbf{3}}%
\newcommand{\pfour}{\mathbf{4}}%
\newcommand{\pfive}{\mathbf{5}}%
\newcommand{\peight}{\mathbf{8}}%
\newcommand{\pnine}{\mathbf{9}}%
$$
\begin{array}{@{}lclclc@{}}
  \cyoung{\pone 1111115,\ptwo 222222,\pthree 3355,\pfour 44,\pfive 55,899,9}
    &\too\quad&
  \cyoung{\pone 111115,\ptwo 22222,\pthree 355,\pfour 49,\pfive 5,89,9}
    &\too\quad&
  \cyoung{\pone 11115,\ptwo 2222,\pthree 55,\pfour 99,\pfive,\peight,\pnine}
    &\too
\\\\
  \cyoung{\pone 1115,\ptwo 222,55,99}
    &\too\quad&
  \cyoung{\pone 115,\ptwo 22,55,99}
    &\too\quad&
  \cyoung{\pone 15,\ptwo 2,\pfive 5,99}
    &\too
\\\\
  \cyoung{\pone 5,\ptwo 9,\pfive,\pnine}
    &\too\quad&
  \cyoung{\pfive,\pnine}
    &\too\quad&
  \vn
\end{array}
$$
\caption{Peeling a peelable tableau}\label{f:peel}
\end{figure}%
Example~\ref{ex:tabarray},~below).
\end{example}

\begin{remark} \label{rem:all peelables}
The set of $D$-peelable tableaux may be constructed by
``unpeeling'' as follows.  Let $C$ be as in the definition of
$D$-peelability, and $p$ the number of cells in~$C$.  Suppose all
of the $(D-C)$-peelable tableaux have been constructed.  For each
pair $(T,V)$ where $T$ is a $(D-C)$-peelable tableau and $V$ is a
vertical $p$-strip (a skew shape $p$ cells, no two of which lie in
the same row) whose union with the shape of $T$ is a partition,
use the jeu de taquin \cite[Section~1.2]{Ful97} to slide~$T$ into
the topmost cell of~$V$, then into the next topmost cell of $V$,
and so~on.  This vacates $p$ cells at the top of the first column.
Place the single-column tableau~$C$ into the vacated cells.  If
the result is a (semistandard) tableau then it is $D$-peelable by
construction.  All $D$-peelable tableaux can be obtained in this
manner.
\end{remark}

\begin{thm}\cite[Theorem~29]{RSplact95} \label{t:SchubKey}
For any (partial) permutation~$w$,
\begin{eqnarray*}
  \Schub_w(\xx) &=& \sum_{Q\in\Peel(D(w))} \kp_{\wgt(\leftkey(Q))}.
\end{eqnarray*}
\end{thm}

\begin{remark} \label{rem:PeelMultiset}
Theorem \ref{t:SchubKey} gives a combinatorial description of the
multiset $M(w)$ in Theorem~\ref{t:SchubKeyWeak}:
it is the multiset of compositions given by $\wgt(\leftkey(Q))$ for
$Q\in\Peel(D(w))$.
\end{remark}

\begin{thm} \label{t:PeelQuiver}
The quiver constant $c_\bla(\rr)$ is the number of $D_\rr$-peelable
tableaux of shape~$\la$, where $\bla$ is obtained from~$\la$ by
deleting~$D_\homv$.
\end{thm}
\begin{proof}
This follows from Theorems~\ref{t:qc=ss} and~\ref{t:SchubKey} along
with Corollary~\ref{c:stanley2schur}, since for any tableau~$Q$, the
partition $\wgt(\leftkey(Q))_+$ is just the shape of~$Q$.%
\end{proof}



\part{Factor sequences from peelable tableaux}

\label{part:fs}
\setcounter{section}{0}
\setcounter{thm}{0}
\setcounter{equation}{0}

In \cite{BF} it was conjectured that the quiver constant~$c_\bla(\rr)$
counts certain lists of tableaux called `$\rr$-factor sequences of
shape~$\bla$'.  We already know by Theorem~\ref{t:PeelQuiver} that the
quiver constants
count peelable tableaux.  We shall verify the Buch--Fulton conjecture
by establishing a bijection from peelable tableaux to factor
sequences.

\section{Factor sequences}

\label{sec:factor}

The \bem{tableau array} $T(\rr)=(T_{ij}(\rr))$ for the rank
array~$\rr$ is defined as follows.  View the diagram $D_\rr =
D(v(\rr))$ of the Zelevinsky permutation as being filled with
integers, where every cell (either~$*$ or~$\sq$) in row~$k$ is filled
with $k \in \ZZ$.  Define $T_{ij}(\rr)$ to be the rectangular tableau
given by the $i^\th$ block column (from the right) and the $j^\th$
block row (from the top), defined for blocks on or below the
superantidiagonal.  The tableau $T_{ij}$ has shape $R_{i-1,j+1}$ by
Lemma~\ref{l:diagram}.

For later use, we define tableaux $Y_i$ and $K_i$ as follows.  Let $Y$
be the tableau given by restricting the filling of~$D_\rr$ to the
blocks above the superantidiagonal, and set~$Y_i$ equal to the part
of~$Y$ in the $i^\th$ block column.  Let $K_i$ be the tableau obtained
by stacking (left-justified) the tableaux in the $i^\th$ column of the
tableau array, with $T_{i,i-1}$ on top, $T_{i,i}$ below it, and so on.
Equivalently,
\begin{eqnarray} \label{eq:KT}
  K_i &=& [T_{i,n-1} T_{i,n-2} \dotsm T_{i,i} T_{i,i-1}].
\end{eqnarray}

Let $A_i$ be the interval of row indices occurring in the $i^\th$
block row, and then set $B_i=A_i\cup A_{i+1}\cup\dotsm$.  Given a set
$A$, let $A^\spot$ denote the set of words in the \mbox{alphabet~$A$.
Then}
\begin{equation} \label{eq:Yword}
  Y_i \in (A_0\cup\dotsm \cup A_{i-2})^\spot.
\end{equation}

\begin{example} \label{ex:tabarray}
For the ranks~$\rr$ in Example~\ref{ex:rank}, with diagram~$D_\rr$ in
Example~\ref{ex:zel}, Fig.~\ref{f:tabarray} depicts the filling
of~$D_\rr$, the tableau array $T=T(\rr)$, and the tableaux~$Y_i$~%
\begin{figure}
\begin{equation*}
\text{filled~}D_\rr\ \:=\:\
        \begin{array}{|%
        @{}  c@{\:}@{\:}c@{\:}@{\:}c@{\:}|%
        @{\:}c@{\:}@{\:}c@{\:}@{\:}c@{\:}@{\:}c@{\:}|%
        @{\:}c@{\:}@{\:}c@{\:}@{\:}c@{\:}|%
        @{\:}c@{\:}@{\:}c@{\:}|@{}}
\cline{1-12}
 \: 1 & 1 & 1 & 1 & 1 & 1 & 1 &\ti&\cd&\cd&\cd&\cd
\\
 \: 2 & 2 & 2 & 2 & 2 & 2 & 2 &\cd&\ti&\cd&\cd&\cd
\\\cline{1-12}
 \: 3 & 3 & 3 &\ti&\cd&\cd&\cd&\cd&\cd&\cd&\cd&\cd
\\
 \: 4 & 4 & 4 &\cd&\ti&\cd&\cd&\cd&\cd&\cd&\cd&\cd
\\
 \: 5 & 5 & 5 &\cd&\cd&5&5&\cd&\cd&5&\ti&\cd
\\\cline{1-12}
 \:\ti&\cd&\cd&\cd&\cd&\cd&\cd&\cd&\cd&\cd&\cd&\cd
\\
 \:\cd&\ti&\cd&\cd&\cd&\cd&\cd&\cd&\cd&\cd&\cd&\cd
\\
 \:\cd&\cd&8&\cd&\cd&\ti&\cd&\cd&\cd&\cd&\cd&\cd
\\
 \:\cd&\cd&9&\cd&\cd&\cd&9&\cd&\cd&9&\cd&\ti
\\\cline{1-12}
 \:\cd&\cd&\ti&\cd&\cd&\cd&\cd&\cd&\cd&\cd&\cd&\cd
\\
 \:\cd&\cd&\cd&\cd&\cd&\cd&\ti&\cd&\cd&\cd&\cd&\cd
\\
 \:\cd&\cd&\cd&\cd&\cd&\cd&\cd&\cd&\cd&\ti&\cd&\cd
\\\cline{1-12}
\end{array}\qquad
\qquad T \ \:=\:\
\begin{array}{ccc|c}
     3     &     2   &    1    & i \diagup j \\ \hline
           &         & {\varnothing \vphantom{\young(1)}} &     0
           \\[12pt]
           & \young(55) & \young(5) &     1       \\
 \young(8,9) & \young(9) & \young(9) &     2
\end{array}
\end{equation*}
$$
\begin{array}{rcl@{\qquad}rcl@{\qquad}rcl}
A_2 &=& \{6,7,8,9\} & A_1 &=& \{3,4,5\} & A_0 &=& \{1,2\}
\\[1ex]
Y_3 &=& \cyoung{111,222,333,444,555}
&
Y_2 &=& \cyoung{1111,2222}
&
Y_1 &=& \varnothing
\\\\[-1ex]
K_3 &=& \cyoung{8,9}
&
K_2 &=& \cyoung{55,9}
&
K_1 &=& \cyoung{5,9}
\end{array}
$$
\caption{Tableau array and other data from Example~\ref{ex:tabarray}}
\label{f:tabarray}
\end{figure}
and~$K_i$.\qed%
\end{example}

We now recall the recursive structure underlying the definition of an
$\rr$-factor sequence.  Let $\rrh$ be the rank array obtained by
removing the entries $r_{ii}$ for $0\le i\le n$.  Using notation from
Definition~\ref{d:diagram}, observe that $D_\rrh^*$ is obtained from
$D_\rr^*$ by removing the superantidiagonal cells and some empty rows
and columns.  Hence we may identify $D_\rrh^*$ with a subdiagram
of~$D_\rr^*$ by reindexing the nonempty rows.  Under this
identification the tableau array $T(\rrh)$ is obtained from $T(\rr)$
by removing the superantidiagonal tableaux.

\begin{example}
Continuing with Example~\ref{ex:tabarray}, we get $\rrh$, $\Rh$,
$D_\rrh$, and $T(\rrh)$ as
\begin{figure}
\begin{eqnarray*}
  \rrh \ \:=\
  \begin{matrix}
  & & & \cd \\
  & &\cd & 2\\
  & \cd&2&1\\
  \cd& 2&1&0
  \end{matrix}\quad
&\qquad&
  \Rh \ \:=\:
  \begin{matrix}
  & & \cd \\[1ex]
  & \cd & \yng(1) \\[1ex]
  \cd  & \yng(1) & \yng(1)
  \end{matrix}
\\[.5ex]\cr
  D_\rrh \ \:=\ \;
  \begin{array}{|cc|cc|cc|} \hline
   *  &  *  & \cd & \cd & \cd & \cd \\
   *  &  *  & \cd &  5  & \cd & \cd \\ \hline
  \cd & \cd & \cd & \cd & \cd & \cd \\
  \cd &  9  & \cd &  9  & \cd & \cd \\ \hline
  \cd & \cd & \cd & \cd & \cd & \cd \\
  \cd & \cd & \cd & \cd & \cd & \cd \\
  \hline
  \end{array}
&\qquad&
  T(\rrh) \ \:=
  \begin{array}{ccc}
  & & \cd \\[1ex]
  & \cd & \cyoung{5} \\[1ex]
  \cd & \cyoung{9} & \cyoung{9}
  \end{array}
\end{eqnarray*}
\caption{Data obtained by deleting the antidiagonal of~\(\rr\) from
Fig.~\ref{f:tabarray}}
\label{f:rrh}
\end{figure}
in Fig.~\ref{f:rrh}.\qed%
\end{example}

Returning to the general case, let $T_i=T_{i,i-1}$ for $1\le i\le n$.

\begin{defn} \label{d:factor}
The notion of \bem{$\rr$-factor sequence} is defined recursively, as
follows.
\begin{numbered}
\item
If $n=1$ then there is a unique $\rr$-factor sequence, namely~$(T_1)$.
\item
For $n\ge2$, an $\rr$-factor sequence is a tableau list
$(W_1,\dotsc,W_n)$ such that there exists an $\rrh$-factor sequence
$(U_1,\dotsc,U_{n-1})$ (for the tableau array $T(\rrh)$ as defined
above) and factorizations
\begin{alignat}{2}
\label{eq:U=PQ}
  U_i &\ \:=\:\ [P_i Q_i] &\qquad&\text{for $1\le i\le n-1$} \\
\intertext{such that}
\label{eq:W=QTP}
  W_i &\ \:=\:\ [Q_i T_i P_{i-1}]&\qquad&\text{for $1\le i\le n$,} \\
\intertext{where by convention}
\label{eq:QnP0}
Q_n&\ \:=\:\ P_0\ \:=\:\ \vn.
\end{alignat}
\end{numbered}
\end{defn}

The global computation of a factor sequence proceeds as follows.
Index the block antidiagonals by their distance (going northwest) to
the main superantidiagonal.
\begin{numbered}
\item
Initialize $W_{ij}:=T_{ij}$ for all $n \ge j+1\ge i \ge 1$.

\item
Set $k:=n-1$.

\item
At this point the $k^\th$ antidiagonal of~$W$ is an $\rr^{(k)}$-factor
sequence, where $\rr^{(k)}$ is obtained from $\rr$ by removing the
$0^\th$ through $(k-1)^\st$ antidiagonals.  Factor each tableau
$W_{ij}=[P_{ij}Q_{ij}]$ on the $k^\th$ antidiagonal.  Now move the
left factor $P_{ij}$ to the west and the right factor $Q_{ij}$ to the
north.

\item
At this point $W_{ij}=T_{ij}$ on the $(k-1)^\st$ antidiagonal.  For
every position on the $(k-1)^\st$ antidiagonal, multiply the current
entry $T_{ij}$ on the left with the tableau $Q_{i,j+1}$ coming from
the south and on the right with the tableau $P_{i-1,j}$ coming from
the east.  That is, set $W_{ij} := [Q_{i,j+1} T_{ij} P_{i-1,j}]$.

\item
Set $k:=k-1$.  If $k=0$, stop.  Otherwise go to step~3.
\end{numbered}

\begin{example}
In Fig.~\ref{f:fs} we compute a factor sequence as above, except that
after a tableau is factored and its factors moved, its position is
vacated.
\begin{figure}
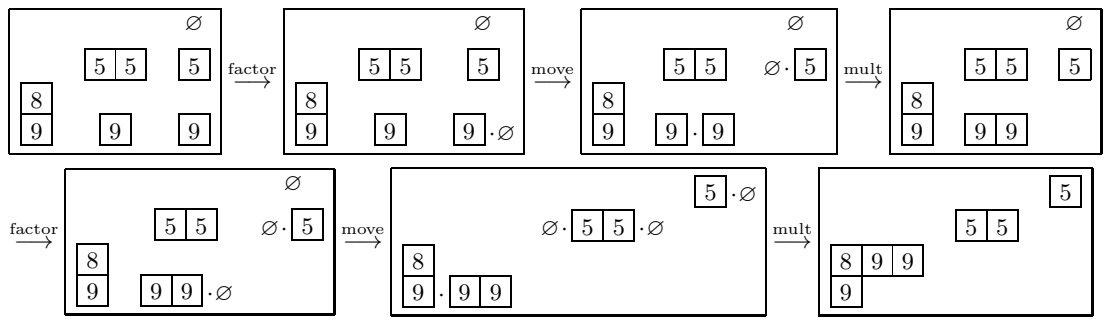

\begin{rcgraph}
\Yboxdim{12pt} \Yinterspace{1pt}
\begin{array}{@{}l@{}}
        \begin{array}{|@{\ }ccc@{\ }|}\hline
                     &             &    \vn     \\[1ex]
                     & \lyoung{55} & \lyoung{5} \\
        \lyoung{8,9} & \lyoung{9}  & \lyoung{9} \\[.5ex]\hline
        \end{array}
\overset{\mathrm{factor}}\longrightarrow
        \begin{array}{|@{\ }ccc@{\ }|}\hline
                     &             &     \vn    \\[1ex]
                     & \lyoung{55} & \lyoung{5} \\
        \lyoung{8,9} & \lyoung{9}  & \lyoung{9}\!\cdot\!\vn\\[.5ex]\hline
        \end{array}
\overset{\mathrm{move}}\longrightarrow
        \begin{array}{|@{\ }ccc@{\ }|}\hline
                     &                        & \vn \\[1ex]
                     & \lyoung{55}            & \vn\!\cdot\!\lyoung{5}\\
        \lyoung{8,9} & \lyoung{9}\!\cdot\!\lyoung{9}&\\[.5ex]\hline
        \end{array}
\overset{\mathrm{mult}}\longrightarrow
        \begin{array}{|@{\ }ccc@{\ }|}\hline
                     &             &    \vn     \\[1ex]
                     & \lyoung{55} & \lyoung{5} \\
        \lyoung{8,9} & \lyoung{99} &            \\[.5ex]\hline
        \end{array}
\\[7ex]
\overset{\mathrm{factor}}\longrightarrow
        \begin{array}{|@{\ }ccc@{\ }|}\hline
                     &                        & \vn \\[1ex]
                     & \lyoung{55}            &\vn\!\cdot\!\lyoung{5}\\
        \lyoung{8,9} & \lyoung{99}\!\cdot\!\vn&\\[.5ex]\hline
        \end{array}
\overset{\mathrm{move}}\longrightarrow
        \begin{array}{|@{\ }ccc@{\ }|}\hline&&\\[-2.2ex]
                     &                         &\lyoung{5}\!\cdot\!\vn\\
                     &\vn\!\cdot\!\lyoung{55}\!\cdot\!\vn&\\
        \lyoung{8,9}\!\cdot\!\lyoung{99}&      &          \\[.5ex]\hline
        \end{array}
\overset{\mathrm{mult}}\longrightarrow
        \begin{array}{|@{\ }ccc@{\ }|}\hline&&\\[-2.2ex]
                       &             & \lyoung{5} \\
                       & \lyoung{55} &            \\
        \lyoung{899,9} &             &            \\[.5ex]\hline
        \end{array}
\end{array}
\end{rcgraph}
\caption{Computation of a factor sequence}
\label{f:fs}
\end{figure}
Compare Example~\ref{ex:peelable}.\qed%
\end{example}

\begin{remark} \label{rem:fsNW}
It is clear from the global computation that for any $\rr$-factor
sequence $(W_1,\dotsc,W_n)$, the letters of~$W_i$ are a submultiset of
the letters from the tableaux~$T_{k\ell}$ to the south and east of the
$i^\th$ block column and $(i-1)^\st$ block row (that is, for
$\ell+1\ge i\ge k$).
\end{remark}

\section{Zelevinsky peelable tableaux}

\label{sec:zeltab}

This section contains some scattered but essential results regarding
peelable tableaux.  Our discussion culminates in
Definition~\ref{d:ptofs}, which says how to construct sequences of
tableaux from peelable tableaux associated to the diagram of a
Zelevinsky permutation.

The following lemma generalizes the dominance bounds for
$D(w)$-peelable tableaux implied by Proposition \ref{pp:StanDom} to
the case of peelables with respect to any northwest diagram.  This
added generality will be used in the proof of
Proposition~\ref{pp:ZPeelChar}.  Recall the definitions of~$\ldom$,
$D\north$, and~$D\west$ given just before
Proposition~\ref{pp:StanDom}.

\begin{lemma} \label{l:PeelDom}
If a $D$-peelable tableau for a northwest diagram~$D$ has
shape~$\la$,~then
\begin{equation*}
   D \west\north\;\,\ldom\,\la\,\ldom\,D\north\west.
\end{equation*}
\end{lemma}
\begin{noqed}
If a tableau has shape~$\la$ and weight~$\gamma$, then $\gamma_+
\ldom \la$ \cite[Section~2.2]{Ful97}.  Noting that $D
\west\north\: = \code(D)_+$, the first dominance condition holds
by Remark~\ref{rem:PeelWeight}.  The second condition is a
consequence of the first, using the following facts.
\begin{numbered}
\item The transposing of shapes gives an anti-isomorphism of the
poset of partitions under the dominance partial order.

\item There is a shape-transposing bijection $\Peel(D)\rightarrow
\Peel(D^t)$, where $D^t$ is the transpose diagram of $D$
\cite[Definition--Proposition~31]{RSpeelable98}.\qed
\end{numbered}
\end{noqed}

There are a number of operations $D\rightarrow D'$ one may perform
on diagrams that induce maps $\Peel(D)\rightarrow\Peel(D')$ on the
corresponding sets of peelable tableaux.  For notation, let $s_r
D$ be the diagram obtained by exchanging the $r^\th$
and~\mbox{$(r+1)^\st$} rows of~$D$, and similarly for~$D s_r$
exchanging columns.  Given a diagram~$D$ and a set~$I$ of row
indices, let~$D|_I$ be the diagram obtained by taking only the
rows of~$D$ indexed by~$I$.  If~$u$ is a word or a tableau, let
$u|_I$ be the word obtained from~$u$ by erasing all letters not
in~$I$, and set $[u]_I = [u|_I]$.

\begin{prop} \label{pp:peelprop} Let $D$ be a northwest diagram.
\begin{numbered}
\item
\label{it:rowswap}\cite[Definition--Proposition~13]{RSpeelable98}
If $s_r D$ is also northwest, then there is a shape-preserving
bijection $s_r:\Peel(D)\rightarrow \Peel(s_r D)$.  (It is given by
the action of the ``$r^\th$ type $A$ crystal reflection
operator''; see\/ \cite[Section~2]{RSpeelable98}\/ for the
definition of~$s_r$.)

\item \label{it:colswap} \cite[Lemma~51]{RSpeelable98} If both $D$
and $D s_r$ are northwest then $\Peel(D s_r)=\Peel(D)$.

\item \label{it:int} \cite[Lemmata~54 and~56]{RSpeelable98} There
is a surjective map $\Peel(D) \to \Peel(D|_I)$ given by $Q \mapsto
[Q]_I$.
\end{numbered}
\end{prop}

\begin{lemma} \label{lem:innershape}
Let $D$ be a northwest diagram such that $D\supset D(\mu)$ for a
partition~$\mu$.  Then $Q\supset \key(\mu)$ for all tableaux $Q\in
\Peel(D)$.
\end{lemma}
\begin{proof}
Suppose $\mu$ has $i$ nonempty rows.  By
Proposition~\ref{pp:peelprop}.\ref{it:int} with
$I=\{1,2,\dotsc,i\}$, we may assume that $D$ has $i$ rows.  Since
$D\supset D(\mu)$, the first column $C$ of $D$ is $C=i \dotsm 3 2
1$.  Let $\nu=(\mu_1-1,\mu_2-1,\dotsc,\mu_i-1)$.  The first column
of $Q$ contains $C$ and contains only letters in $I$, so it must
be equal to $C$.  Then $[Q-C]$ is merely $Q$ with its first column
erased.  By definition $[Q-C]$ is $(D-C)$-peelable.  By induction
$Q-C\supset \key(\nu)$.  But then $Q\supset C
\,\key(\nu)=\key(\mu)$.
\end{proof}

We say that a diagram $D$ is \bem{$\la$-partitionlike} for the
partition~$\la$ if $D$ is northwest and $D(\la)$ can be obtained
from~$D$ by repeatedly exchanging adjacent rows or columns,
staying within the family of northwest shapes.

\begin{lemma} \label{lem:partlike}
$\Peel(D)$ consists of one tableau, namely\/ $\key(\code(D))$, if
$D$ is partitionlike.
\end{lemma}
\begin{proof}
Using parts~\ref{it:rowswap} and~\ref{it:colswap} of
Proposition~\ref{pp:peelprop} to pass from $D$ to $D(\la)$, there
is a shape-preserving bijection from $\Peel(D)$
to~$\Peel(D(\la))$. Lemma~\ref{lem:innershape} and
Remark~\ref{rem:PeelWeight} imply that $\key(\la)$ is the unique
$D(\la)$-peelable tableau.  Thus $\Peel(D)$ consists of a single
tableau of shape~$\la$ (since $\key(\la)$ has shape $\la$) and
weight $\code(D)$ (by Remark~\ref{rem:PeelWeight}).
But $\key(\code(D))$ is the unique such tableau.%
\end{proof}

\begin{lemma} \label{lem:peelpairs}
Suppose $D$ is a northwest diagram such that its $r^\th$ row
(viewed as a subset of integers given by the column indices of its
cells) either contains or is contained in its $(r+1)^\st$ row.
Suppose that rows $r$ and~\mbox{$r+1$} have $a$ and~$b$ cells,
respectively. If~$Q$ is a $D$-peelable tableau, then
$[Q]_{\{r,r+1\}} = [(r+1)^b r^a]$.
\end{lemma}
\begin{proof}
Proposition~\ref{pp:peelprop}.\ref{it:int} with $I=\{r,r+1\}$ implies
$[Q]_I$ is $D|_I$-peelable.  By Lemma~\ref{lem:partlike} there is a
unique $D|_I$-peelable tableau, namely~\mbox{$[(r+1)^b r^a]$}.%
\end{proof}

In light of Remark~\ref{rem:PeelMultiset}, the following two results
are restatements of Lemma~\ref{lem:SchubDemAscent} and
Proposition~\ref{pp:ZelDem}.

\begin{lemma} \label{lem:peelkeypairs}
Let $w$ be a permutation such that $w(i)<w(i+1)$.  If
$Q\in\Peel(D(w))$ and $\beta=\wgt(\leftkey(Q))$, then $\beta_i \le
\beta_{i+1}$.
\end{lemma}

\begin{prop} \label{pp:ZPeelShape}
For a fixed dimension vector $(r_0,r_1,\dotsc,r_n)$, identify
$D_\homv$ and $D(\Omega_0)$ with partitions of those shapes.  Let
$\rr$ be a rank array with the above dimension vector, $Q\in
\Peel(D_\rr)$, and $\la=\shape(Q)$.  Then $D_\homv \subset \la
\subset D(\Omega_0)$.
\end{prop}

Let $Q\in\Peel(D_\rr)$ and $\la=\shape(Q)$.  By
Lemma~\ref{lem:innershape} $Q\supset\key(D_\homv)$.  The skew tableau
$Q-\key(\homv)$ has shape $\la/D_\homv\subset
D(\Omega_0)/D_\homv$. The skew shape $D(\Omega_0)/D_\homv$ is a
disjoint union of $r_{i-1}\times r_i$ rectangles for $1\le i\le n$ in
which the $r_{i-1}\times r_i$ rectangle resides in the $(i-1)^\st$
block row and $i^\th$ block column, where the indexing of block rows
and columns is the same as in Section~\ref{sec:zelperm}.
Proposition~\ref{pp:ZPeelShape} and Lemma~\ref{lem:innershape} imply
that the following map~$\ptofs_\rr$ from $\Peel(D_\rr)$ to length~$n$
sequences of tableaux is well-defined.  It is the peelable tableau
analogue of Definition~\ref{d:delete}.

\begin{defn} \label{d:ptofs}
Let~$Q$ be a $D_\rr$-peelable tableau.  The sequence $\ptofs_\rr(Q) =
(W_1,W_2,\dotsc,W_n)$ of tableaux is obtained by
\bem{deleting~$D_\homv$ from~$Q$} if $W_i$ is the subtableau of the
skew tableau $Q-\key(D_\homv)$ sitting inside the $r_{i-1}\times
r_i$~rectangle.
\end{defn}

\begin{example} \label{ex:peelable}
The image under $\ptofs_\rr$ of the $D_\rr$-peelable tableau in
Example~\ref{ex:peel} is the sequence $(W_1,W_2,W_3)$ depicted in
\begin{figure}
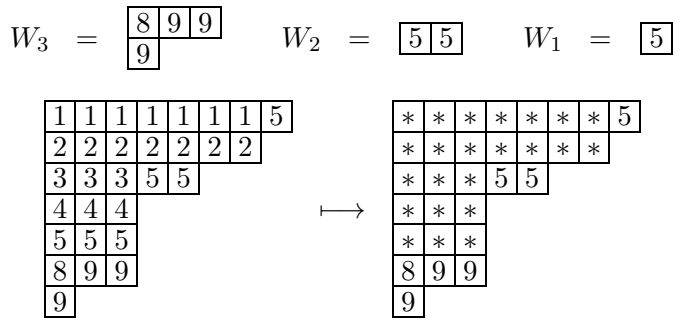

$$
\Yboxdim{12pt} \Yinterspace{1pt}
\begin{array}{@{}rcl@{\qquad}rcl@{\qquad}rcl@{}}
  W_3 &=& \cyoung{899,9}
& W_2 &=& \cyoung{55} & W_1 &=& \cyoung{5}
\end{array}
$$
$$
\Yboxdim{12pt} \Yinterspace{1pt}
\begin{array}{@{}c@{}}
\cyoung{11111115,2222222,33355,444,555,899,9}
\end{array}
\ \:\longmapsto\ \:
\begin{array}{@{}c@{}}
\cyoung{*******5,*******,***55,***,***,899,9}
\end{array}
$$
\caption{Peelable tableau to factor sequence under \(\ptofs_{\!\rr}\)}
\label{f:W}
\end{figure}
Fig.~\ref{f:W}.\qed%
\end{example}

\section{Bijection to factor sequences}

\label{sec:bijection}

Next we give a condition for a tableau list to be obtained by deleting
$D_\homv$ from a $D_\rr$-peelable tableau.  This will be indispensable
for proving the bijection from peelable tableaux to factor sequences,
the aim of this section.  The idea is to perform several peelings
together, so as to remove one block column at a time.  The `code
of~$D$' is defined in~Remark~\ref{rem:PeelWeight}.

\begin{lemma} \label{lem:peelseveral}
Let $D$ be a northwest diagram, divided by a vertical line into
diagrams~$D^w$ (``$D$-west'') and~$D^e$ (``$D$-east''), where $D^w$ is
the first several columns of~$D$.  Assume $D^w$ is
$\la$-partitionlike, and set $Y = \key(\code(D^w))$.
\begin{numbered}
\item
If $Q\in\Peel(D)$ then $Q\supset Y$.

\item
If $Q\supset Y$, then $Q\in\Peel(D)$ if and only if\/
$[Q-Y]\in\Peel(D^e)$.
\end{numbered}
\end{lemma}
\begin{proof}
By Proposition~\ref{pp:peelprop}.\ref{it:colswap} we assume that the
columns of~$D^w$ are decreasing with respect to containment.  Suppose
the rows of~$D^w$ are not decreasing with respect to containment.
Then some row~$r$ of~$D^w$ is a proper subset of row~\mbox{$r+1$}.
Since $D$ is northwest, row~$r$ of~$D^e$ is empty.  This given, apply
Proposition~\ref{pp:peelprop}.\ref{it:rowswap}.  It is straightforward
to check that the above conditions for~$D$ and~$Q$ are equivalent to
those for~$s_r D$ and~$s_r Q$.  Hence we may assume that $D^w=D(\la)$
for a partition~$\la$.

Now item~1 follows from Lemma~\ref{lem:innershape}.  For item~2 we use
the jeu de taquin \cite[Section~1.2]{Ful97}.  Let $C_j$ be the column
word that gives the (decreasing) set of row indices for the elements
in column~$j$ of~$D$.  The result $[Q-C_1]$ of one peeling step can be
computed by sliding the skew tableau $Q-C_1$ into the cells occupied
by~$C_1$.  Repeating this process for the columns of~$D^w$, the
tableau resulting from applying the corresponding peeling steps to~$Q$
is also obtained by removing the subtableau~$Y$ (which must be present
in~$Q$ by item~1) and sliding the skew tableau $Q-Y$ into the cells
that were occupied by~$Y$.%
\end{proof}

\begin{lemma} \label{lem:WIm}
If $(W_1,\dotsc,W_n)=\ptofs_\rr(Q)$ for $Q\in\Peel(D_\rr)$, then
\begin{eqnarray}
\label{eq:Wword}
  W_i &\in& B_{i-1}^\spot  \\
\label{eq:Wht}
  \hgt(W_i) &\leq& r_{i-1} \\
\label{eq:Wwd}
  \wid(W_i) &\leq& r_i,
\end{eqnarray}
and $Q$ has the block column factorization
\begin{eqnarray} \label{eq:QW}
  Q &=& (W_n Y_n)\dotsm (W_2 Y_2) (W_1 Y_1).
\end{eqnarray}
\end{lemma}
\begin{proof}
Statement~\eqref{eq:Wword} follows from \eqref{eq:Yword} and the fact
that $W_i$ sits below $Y_i$ in a block column of the tableau~$Q$.  The
rest follows by definition, using Proposition~\ref{pp:ZPeelShape}.%
\end{proof}

\begin{prop} \label{pp:ZPeelChar}
The list $(W_1,\dotsc,W_n)$ of tableaux is obtained from some
$D_\rr$-peelable tableau by deleting~$D_\homv$ if and only if there
exist tableaux $X_1,\ldots,X_n$ such that
\begin{align}
\label{eq:Xn}
  X_n &\ \:=\ \: W_n,
\\\label{eq:XK}
  X_i &\ \:\supset\ \: K_i \quad\rlap{for all\/ $i$,}
\\\label{eq:Xi}
  X_i &\ \:=\ \:[(X_{i+1}-K_{i+1})W_i]
        \quad\rlap{for\/ $i = 1,\ldots,n-1$,}\quad
\\\label{eq:X1}
  \llap{and\quad}X_1&\ \:=\ \: K_1;
\end{align}
and if
\begin{eqnarray} \label{eq:Zdef}
  Z_i &=& [X_i W_{i-1}\dotsm W_2 W_1]_{B_{i-1}}
\end{eqnarray}
then
\begin{eqnarray}
\label{eq:Zht} \hgt(Z_i) &\le& r_{i-1}, \\[-3ex]\cr
\label{eq:Zwd} \wid(Z_i) &\le& r_i, \\[-3ex]\cr
\label{eq:Zal} \llap{and\quad}\wgt(Z_i) &=& \wgt(T_{\le i,\,\ge i-1}),
\end{eqnarray}
where
$\displaystyle\wgt(T_{\le i,\,\ge j}) = \sum_{\substack{k\ge i \\ \ell
\le j}} \wgt(T_{k,\ell})$.
\end{prop}

The proof will come shortly.

\begin{lemma} \label{l:XZ}
The conditions in Proposition~\ref{pp:ZPeelChar} imply the following:
\begin{eqnarray}
\label{eq:Xword} X_i &\in& B_{i-1}^\spot \\
\label{eq:Xht} \hgt(X_i) &\le& r_{i-1} \\
\label{eq:Xwd} \wid(X_i) &\le& r_i \\
\label{eq:Xal} \wgt(X_i) &\le& \wgt(T_{\le i,\ge i-1})
\end{eqnarray}
\end{lemma}
\begin{proof}
\eqref{eq:Xword} can be proved by induction using
\eqref{eq:Wword}.  $X_i$ is a factor of $Z_i$ by \eqref{eq:Xword}
and \eqref{eq:Zdef}.  The last three assertions then follow from
\eqref{eq:Zht}, \eqref{eq:Zwd}, and \eqref{eq:Zal}.
\end{proof}

\begin{lemma} \label{l:XZ'}
If \eqref{eq:Xn} through~\eqref{eq:X1} hold, then the conditions
\eqref{eq:Xword} through \eqref{eq:Xwd} imply \eqref{eq:Wword} through
\eqref{eq:QW}, plus
\begin{eqnarray*}
  \wgt(W_i) &\le& \wgt(T_{\le i,\,\ge i-1}).
\end{eqnarray*}
\end{lemma}
This last inequality is the same condition on $(W_1,\dotsc,W_n)$
appearing in Remark~\ref{rem:fsNW}.

\begin{example}
The tableaux $X_i$ and $Z_i$ associated with
the above tableaux $W_i$ are given by
\begin{alignat}{3}
X_3 &= \cyoung{899,9} &\qquad  X_2 &= \cyoung{55,99} &\qquad  X_1 &=
\cyoung{5,9} \\
Z_3 &= \cyoung{899,9} & Z_2 &= \cyoung{555,99} & Z_1 &= \cyoung{5,9}
\end{alignat}
Note that the $X$ tableaux appear at the bottom of the leftmost block
columns in the first, fourth, and eighth tableaux in
Example~\ref{ex:peel}.\qed%
\end{example}

\begin{proof}[Proof of Proposition~\ref{pp:ZPeelChar}]
Set $D=D_\rr$, and let $D_i$ be the diagram obtained from~$D$ by
removing the $n^\th$ through $(i+1)^\st$ block columns (or
equivalently the leftmost $r_n+\dotsm+r_{i+1}$ columns).  Let
$D_{i,j}$ be obtained from $D_i$ by removing the $0^\th$ through
$(j-1)^\st$ block rows.  In other words, $D_{i,j}$ consists of the
block at the intersection of the $i^\th$ block column and $j^\th$
block row, along with all blocks to its south and east.

Suppose $(W_1,\dotsc,W_n) = \ptofs_\rr(Q)$ for a tableau $Q \in
\Peel(D)$.  Let $Q_i\in\Peel(D_i)$ be the tableau obtained from $Q$ by
$r_n+\dotsm+r_{i+1}$ peelings.  We shall proceed by induction, with
hypothesis (for descending~$i$) being that $Q_j$ has the block
column~factorization
\begin{eqnarray} \label{eq:Qfact}
  Q_j &=& (X_j Y_j) (W_{j-1} Y_{j-1}) \dotsm (W_2Y_2)(W_1Y_1)
\end{eqnarray}
and that $X_j,Z_j$ have been defined and have the desired properties
for $n\ge j\ge i$.

Suppose $0\le i\le n$ and the induction hypothesis holds for indices
greater than $i$.  It must be shown that the induction hypothesis
holds for~$i$.  The tableau $X_i$ is defined by~\eqref{eq:Xn} when
$i=n$, or by~\eqref{eq:Xi} when $i<n$ (note that \eqref{eq:XK} holds
for~\mbox{$i+1$}).

Apply Lemma~\ref{lem:peelseveral} to $Q_{i+1}\in\Peel(D_{i+1})$, with
$D_{i+1}^w$ the west block column (the leftmost $r_{i+1}$ columns) of
$D_{i+1}$.  By \eqref{eq:Qfact} for $i+1$, the first block column of
$Q_{i+1}$ is $X_{i+1} Y_{i+1}$.  Also
$\key(\code(D_{i+1}^w))=K_{i+1}Y_{i+1}$.  Recalling that $Q_i$ is
defined by iterated peelings from $Q$, we have
\begin{equation} \label{eq:Qprod}
\begin{array}{rcl}
  Q_i &=& [(X_{i+1}Y_{i+1}-K_{i+1}Y_{i+1})(W_iY_i)\dotsm(W_1Y_1)] \cr
      &=& [(X_{i+1}-K_{i+1})W_i Y_i \dotsm (W_1Y_1)] \cr
      &=& [ (X_i Y_i) (W_{i-1}Y_{i-1})\dotsm(W_1Y_1)],
\end{array}
\end{equation}
the last equality by definition of~$X_i$.  Apply
Proposition~\ref{pp:peelprop}.\ref{it:int} to $Q_i\in \Peel(D_i)$ and
the interval~$B_{i-1}$.  Note that $(D_i)|_{B_j}=D_{i,j}$, so in
particular $[Q_i]_{B_{i-1}}\in\Peel(D_{i,i-1})$.~~Thus
$$
  [Q_i]_{B_{i-1}}\ =\ [X_i W_{i-1} \dotsm W_1]_{B_{i-1}}\ =\ Z_i,
$$
by \eqref{eq:Qprod}, \eqref{eq:Yword}, and \eqref{eq:Zdef}.  We
have proved that
\begin{eqnarray} \label{eq:ZiPeel}
  Q\in \Peel(D) &\implies& Z_i \in \Peel(D_{i,i-1}).
\end{eqnarray}

The sum of the heights of the rectangles in the $i^\th$ block column
of $D_{i,i-1}\subset D$ is at most $r_{i-1}$ by \eqref{eq:heightsums}.
The sum of the widths of the rectangles in the $(i-1)^\st$ block row
of $D_{i,i-1}\subset D$ is at most $r_i$ by \eqref{eq:widthsums}.
$Z_i$ fits inside a $r_{i-1}\times r_i$ rectangle, by
Lemma~\ref{l:PeelDom}, \eqref{eq:heightincr}, and
\eqref{eq:widthincr}.  This proves \eqref{eq:Zht} and \eqref{eq:Zwd}.
Equation~\eqref{eq:Zal} follows from Remark~\ref{rem:PeelWeight}.  By
Lemma~\ref{l:XZ}, the corresponding properties \eqref{eq:Xword}
through~\eqref{eq:Xal} hold for the $X$ tableaux.  The block column
factorization \eqref{eq:Qfact} for $j=i$ then follows from
\eqref{eq:Qprod}, \eqref{eq:Xht}, \eqref{eq:Xwd}, and
\eqref{eq:Xword}.

To prove \eqref{eq:XK}, apply Lemma~\ref{lem:peelseveral} with
$Q_i\in \Peel(D_i)$ and $D_i^w$ given by the westmost block column
of $D_i$.  This yields that $Q_i \supset \key(\code(D_i^w))=K_i
Y_i$.  By \eqref{eq:Qfact} the first block column of $Q_i$ is $X_i
Y_i$, from which \eqref{eq:XK} follows.

For the converse (i.e.\ that \eqref{eq:Xn}--\eqref{eq:Zal} imply
$(W_1,\ldots,W_n)$ lies in the image of~$\ptofs_\rr$), suppose the
$X_i$ and $Z_i$ exist and have the desired properties.  By
Lemmas~\ref{l:XZ} and~\ref{l:XZ'}, conditions
\eqref{eq:Xword}--\eqref{eq:Xal} and \eqref{eq:Wword}--\eqref{eq:Wwd}
all hold.  Let $Q_i$ be defined by the block column factorization
\eqref{eq:Qfact}.  The aforementioned properties guarantee that $Q_i$
is a tableau.  It will be shown by induction on increasing $i$ that
$Q_i\in \Peel(D_i)$.

For $i=0$, $Q_0$ is empty by~\eqref{eq:Qfact}, and $D_0$ is the empty
diagram so $Q_0\in\Peel(D_0)$.  Suppose $i>0$.  Apply Lemma
\ref{lem:peelseveral} to $D_i$ and $Q_i$ with $D_i^w$ the west block
column of $D_i$.  Then $Q_i\in\Peel(D_i)$ if and only if $Q_i\supset
\key(\code(D_i^w))=K_iY_i$ and $[Q_i-(K_iY_i)]\in\Peel(D_{i-1})$.  But
$Q_i$ has first block column $X_iY_i$ so by \eqref{eq:XK} it contains
$K_iY_i$.  Also
\begin{eqnarray*}
  [Q_i-(K_iY_i)]&=&[(X_i-K_i)Y_i (W_{i-1}Y_{i-1})\dotsm(W_1Y_1)] \cr
  &=& [X_{i-1}Y_{i-1}(W_{i-2}Y_{i-2}\dotsm(W_1Y_1)] \cr
  &=&Q_{i-1}
\end{eqnarray*}
by \eqref{eq:Qfact} for $i$ and $i-1$, and \eqref{eq:Xi} for $i-1$.
Since by induction $Q_{i-1}\in \Peel(D_{i-1})$, we conclude that
$Q_i\in\Peel(D_i)$.  Taking $i=n$ we have $Q_n\in\Peel(D_n)$.

The proof is complete because $Q_n=Q$ and $D_n=D$.%
\end{proof}

For the record, the next result is ``rank stability for peelable
tableaux'', which is analogous to the rank stability statements we
proved for components of quiver degenerations
(Proposition~\ref{p:stability}) and Zelevinsky pipe dreams
(Proposition~\ref{p:top}).

\begin{cor} \label{c:PeelLim}
There is a bijection $\Peel(D_\rr)\to \Peel(D_{m+\rr})$ inducing a
shape-preserving bijection $\Image(\ptofs_\rr)\rightarrow
\Image(\ptofs_{m+\rr})$.
\end{cor}
\begin{proof}
With notation as in Definition~\ref{d:diagram}, the diagrams $D_\rr^*$
and~$D_{m+\rr}^*$ are the same up to removing some empty columns and
rows.  By Proposition~\ref{pp:ZPeelChar} the images of $\ptofs_\rr$
and $\ptofs_{m+\rr}$ differ by a trivial relabeling of their entries.%
\end{proof}

Corollary~\ref{c:PeelLim} could be used instead of
Proposition~\ref{p:QQrr} as the conduit through which
Corollary~\ref{c:stanley} enters into the proof of
Theorem~\ref{t:PeelQuiver}.


\section{The Buch--Fulton conjecture}

\label{sec:bf}

\begin{thm} \label{t:ptofs}
The map $\ptofs_\rr$ taking each $D_\rr$-peelable tableau $Q \in
\Peel(D_\rr)$ to the list $\ptofs_\rr(Q)$ of tableaux obtained by
deleting~$D_\homv$ from~$Q$ is a bijection from the set of
$D_\rr$-peelable tableaux to $\rr$-factor sequences.
\end{thm}
\begin{proof}
It suffices to show that the conditions on a tableau list
$(W_1,\dots,W_n)$ in Proposition~\ref{pp:ZPeelChar} for membership in
the image of $\ptofs_\rr$ are equivalent to the conditions defining an
$\rr$-factor sequence.  The case $n=1$ is trivial, so suppose $n\ge2$.
Write $\Kh_i$ to denote the analogue of the tableau~$K_i$ for~$\rrh$,
for $1\le i\le n-1$.  Identifying the tableau array of~$\rrh$ with
that of~$\rr$ but with $T_i=T_{i,i-1}$ removed for all~$i$, we have
\begin{alignat}{2}
 \label{eq:KKh}
  K_i&\ \:=\:\ [\Kh_i T_i] &\quad&\text{for $1\le i\le n$,} \\
\intertext{where by convention}
\label{eq:Khn}
  \Kh_n &\ \:=\:\ \vn. &&
\end{alignat}

Suppose first that $(W_1,\dotsc,W_n)$ is an $\rr$-factor sequence.
Let $(U_1,\dotsc,U_{n-1})$ be an $\rrh$-factor sequence, with
factorizations given by \eqref{eq:U=PQ} and such that
\eqref{eq:W=QTP} and \eqref{eq:QnP0} hold.  By induction the
$\rrh$-factor sequence $(U_1,\dotsc,U_{n-1})$ satisfies the
conditions of Proposition~\ref{pp:ZPeelChar} with associated
tableaux $\Xh_i$ for $1\le i\le n-1$.  For convenience define
\begin{eqnarray} \label{eq:Xhn}
  \Xh_n&=&\vn
\end{eqnarray}
The $\Xh_i$ satisfy the following analogues of \eqref{eq:Xi}
through~\eqref{eq:Xal}:
\begin{eqnarray}
\label{eq:Xhi}
  \Xh_i &=& [(\Xh_{i+1}-\Kh_{i+1}) U_i]
\\[-3ex]\cr
\label{eq:XKh}
  \Xh_i &\supset& \Kh_i  \qquad\text{for } 1\le i\le n-1
\\[-3ex]\cr
\label{eq:Xh1}
  \Xh_1 &=& \Kh_1.
\end{eqnarray}
Defining the tableau $\Zh_i$ by
\begin{eqnarray} \label{eq:Zhdef}
  \Zh_i &=& [\Xh_i U_{i-1}\dotsm U_1]_{B_i}
\end{eqnarray}
we have
\begin{eqnarray}
\label{eq:Zhht} \hgt(\Zh_i) &\le& r_{i-1,i} \\
\label{eq:Zhwd} \wid(\Zh_i) &\le& r_{i,i+1} \\
\label{eq:Zhal} \wgt(\Zh_i) &=& \wgt(T_{\le i,\,\ge i}).
\end{eqnarray}
Note that the analogue $\Xh_{n-1}=U_{n-1}$ of \eqref{eq:Xn}
follows from \eqref{eq:Khn}, \eqref{eq:Xhn} and \eqref{eq:Xhi}.

Define $X_i$ by
\begin{eqnarray} \label{eq:XXh}
  X_i &=& [\Xh_i \,T_i \, P_{i-1}].
\end{eqnarray}
By \eqref{eq:XKh} and \eqref{eq:KKh} we have \eqref{eq:XK} and
\begin{eqnarray} \label{eq:XKrec}
  [X_i - K_i] &=& [(\Xh_i - \Kh_i) P_{i-1}].
\end{eqnarray}
We shall show by descending induction on $i$ that $X_i$ and $Z_i$ (the
latter being defined by \eqref{eq:Zdef}) have the properties specified
by Proposition~\ref{pp:ZPeelChar}.

We have $X_n=[\vn T_n P_{n-1}]=[Q_n T_n P_{n-1}]=W_n$ by
\eqref{eq:Xhn}, \eqref{eq:QnP0}, and \eqref{eq:W=QTP},
proving \eqref{eq:Xn}.  Let $n>i\ge 1$.  We have
\begin{eqnarray*}
  X_i &=& [\Xh_i T_i P_{i-1}] \\
  &=& [(\Xh_{i+1}-\Kh_{i+1}) U_i T_i P_{i-1}] \\
  &=& [(\Xh_{i+1}-\Kh_{i+1}) P_i Q_i T_i P_{i-1}] \\
  &=& [(\Xh_{i+1}-\Kh_{i+1}) P_i W_i] \\
  &=& [(X_{i+1}-K_{i+1})W_i]
\end{eqnarray*}
by \eqref{eq:XXh}, \eqref{eq:Xhi}, \eqref{eq:U=PQ}, \eqref{eq:W=QTP},
and \eqref{eq:XKrec}, proving \eqref{eq:Xi}.  We have $X_1=[\Xh_1 T_1
P_0] = [\Kh_1 T_1] = K_1$ by \eqref{eq:XXh}, \eqref{eq:Xh1},
\eqref{eq:QnP0}, and \eqref{eq:KKh}, proving \eqref{eq:X1}.

By construction $Z_i$ satisfies \eqref{eq:Zal}.  To see this,
observe first that $\sum_{i=1}^n \wgt(W_i)$ is the
weight of the entire tableau array.  Then by \eqref{eq:Xn} and
\eqref{eq:Xi}, $\wgt(X_i W_{i-1}\dotsm W_1)$ is the weight of the entire
tableau array minus $\sum_{j=i+1}^m \wgt(K_j)$, which is the part
of the tableau array strictly west of the $i^\th$ block column.
Finally, the restriction of $X_i W_{i-1}\dotsm W_1$ to $B_{i-1}$
restricts to the part of the tableau array weakly below the
$(i-1)^\st$ block row.

In particular, $Z_i$ has at most $r_i$ distinct letters appearing
in it since the corresponding part of the tableau array has that
property (see the proof of Proposition~\ref{pp:ZPeelChar}).  Since
$Z_i$ is a tableau, \eqref{eq:Zht} follows.

For \eqref{eq:Zwd}, we have
\begin{equation} \label{eq:Zi}
\begin{array}{rcl}
Z_i
  &=& [(\Xh_i T_i P_{i-1}) (Q_{i-1} T_{i-1}
      P_{i-2})\dotsm(Q_1T_1P_0)]_{B_{i-1}}
\\&=& [\Xh_i T_i U_{i-1} U_{i-2}\dotsm U_1 ]_{B_{i-1}}
\end{array}
\end{equation}
by \eqref{eq:Zdef}, \eqref{eq:XXh}, \eqref{eq:W=QTP}, \eqref{eq:U=PQ},
\eqref{eq:QnP0}, and the fact that $T_j\in A_{j-1}^\spot$.  We get
\begin{equation} \label{eq:ZZh}
\begin{array}{rcl}
  \Zh_{i-1} &=& [\Xh_{i-1} U_{i-2}\dotsm U_1]_{B_{i-1}} \\
     &=& [ (\Xh_i-\Kh_i)U_{i-1}\dotsm U_1]_{B_{i-1}} \\
     &=& [ (\Xh_i T_i - \Kh_i T_i) U_{i-1}\dotsm U_1]_{B_{i-1}} \\
     &=& [ (\Xh_i T_i - K_i) U_{i-1}\dotsm U_1]_{B_{i-1}} \\
     &=& [Z_i-K_i]
\end{array}
\end{equation}
by \eqref{eq:Zhdef}, \eqref{eq:Xhi}, \eqref{eq:KKh}, and
\eqref{eq:Zi}.  Here we have used $\Xh_i \supset \Kh_i$ so that
$[\Xh_i T_i] \supset [\Kh_i T_i] = K_i$.  By \eqref{eq:ZZh} we have
\begin{eqnarray*}
  \wid(\Zh_{i-1}) &\ge& \wid(Z_i) -\wid(K_i).
\end{eqnarray*}
But $\wid(K_i)=\wid(T_i)=\wid(R_{i-1,i})=r_i-r_{i-1,i}$; together with
\eqref{eq:Zhwd} (applied to $i-1$ instead of~$i$), this
implies~\eqref{eq:Zwd}.

Conversely, suppose $(W_1,\dots,W_n)$ satisfies Proposition
\ref{pp:ZPeelChar} with tableaux $X_i$.  We must show that it is an
$\rr$-factor sequence.  By Lemma~\ref{l:XZ}, conditions
\eqref{eq:Xword}--\eqref{eq:Xal} hold.  We claim that
\begin{alignat}{2}
\label{eq:WsupT}
W_i&\supset T_i &\qquad&\text{for all $1\le i\le n$.} \\
\intertext{We have}
\label{eq:XsupT}
X_i&\supset T_i &\qquad&\text{for all $1\le i\le n$}
\end{alignat}
by \eqref{eq:XK} and \eqref{eq:KT}.  For $i=n$, \eqref{eq:Xn} and
\eqref{eq:XsupT} gives \eqref{eq:WsupT}.  For $n>i\ge 1$,
\eqref{eq:Xword} for $i+1$, \eqref{eq:Xi}, and \eqref{eq:XsupT}
imply \eqref{eq:WsupT}.

Consider the factorization of $W_i$ pictured below.
\begin{eqnarray*}
W_i &=&
\begin{array}{@{}c@{}}
\pspicture[.15](-.5,-.5)(8.5,5.5)
\psline(0,5)(8,5)
\psline(0,0)(0,5)
\psline(0,3)(4,3)
\psline(0,0)(3,0)(3,1)(7,1)(7,3)(8,3)(8,5)
\psline(4,1)(4,5)
\rput(2,4){$T_i$}
\rput(2,1.5){$Q_i$}
\rput(5.5,3){$P_{i-1}$}
\endpspicture
\end{array}
\end{eqnarray*}
By construction \eqref{eq:W=QTP} holds.

Since $T_n$ has the largest letters, $Q_n=\vn$.  To obtain
$P_0=\vn$ (and finish up \eqref{eq:QnP0}), note that $[\Kh_1 T_1 ]
= K_1 = X_1 = [(X_2-K_2)W_1]$ by \eqref{eq:KKh}, \eqref{eq:X1},
and \eqref{eq:Xi}.  But $K_1$ has the same width as $T_1$.  It
follows that $W_1$ does also, and hence that $P_0=\vn$.

Define $U_i$ by \eqref{eq:U=PQ}.  By induction it suffices
to prove that there exist tableaux $\Xh_i$ such that
$(U_1,\dotsc,U_{n-1})$ satisfies
the conditions of Proposition~\ref{pp:ZPeelChar}.
Define $\Xh_n$ by \eqref{eq:Xhn} and
\begin{eqnarray} \label{eq:Xhdef}
  \Xh_i &=& [(X_{i+1}-K_{i+1})Q_i] \qquad\text{for } 1\le i\le n-1.
\end{eqnarray}
We have
$$
  [\Xh_i T_i P_{i-1}]\ =\ [(X_{i+1}-K_{i+1})Q_i T_i P_{i-1}]
  \ =\ [(X_{i+1}-K_{i+1}) W_i]\ =\ X_i
$$
by \eqref{eq:Xhdef}, \eqref{eq:W=QTP}, and \eqref{eq:Xi}, proving
\eqref{eq:XXh}.  This implies \eqref{eq:XKh}.  We get
\begin{eqnarray*}
  \Xh_i &=& [(X_{i+1}-K_{i+1})Q_i] \cr
  &=& [(\Xh_{i+1}T_{i+1}P_i-\Kh_{i+1}T_{i+1}) Q_i] \cr
  &=& [(\Xh_{i+1}-\Kh_{i+1})P_i Q_i ] \cr
  &=& [(\Xh_{i+1}-\Kh_{i+1})U_i]
\end{eqnarray*}
by \eqref{eq:Xhdef}, by \eqref{eq:XXh}, \eqref{eq:KKh}, and
\eqref{eq:XKh} for $i+1$ (which hold by induction), and
\eqref{eq:U=PQ}, proving \eqref{eq:Xhi}.  We have $[\Kh_1 T_1] =
K_1 = X_1 = [\Xh_1 T_1 P_0] = [\Xh_1 T_1]$ by \eqref{eq:KKh},
\eqref{eq:X1}, \eqref{eq:XXh}, and \eqref{eq:QnP0}, from which
\eqref{eq:Xh1} follows.

Define $Z_i$ by \eqref{eq:Zdef} and $\Zh_i$ by \eqref{eq:Zhdef}.  Note
that \eqref{eq:Zhal} holds by construction.  This implies
\eqref{eq:Zhht} because $\Zh_i$ is a tableau in an alphabet with at
most $r_{i,i+1}$ distinct letters.  To show that \eqref{eq:Zhwd}
holds, note that \eqref{eq:ZZh} still holds.  Let $D_{i,j}$ be as in
the proof of Proposition~\ref{pp:ZPeelChar}.  Recall that also that
$Z_i\in\Peel(D_{i,i-1})$ by \eqref{eq:ZiPeel}.  Apply Lemma
\ref{lem:peelseveral} with $D_{i,i-1}^w$ given by the west block
column of $D_{i,i-1}$.  Note that $K_i=\key(\code(D_{i,i-1}^w))$.  It
follows that the $r_i$-fold peeling of $Z_i$ (which is none other than
$\Zh_i$ by \eqref{eq:ZZh}) is $D_{i,i}$-peelable.  The inequality
\eqref{eq:Zhwd} follows.  Since the conditions of
Proposition~\ref{pp:ZPeelChar} have been satisfied,
$(U_1,\dotsc,U_{n-1})$ is in the image of~$\ptofs_\rrh$.  By
induction $(U_1,\dotsc,U_{n-1})$ is an $\rrh$-factor sequence.  Thus
$(W_1,\dotsc,W_n)$ is an $\rr$-factor sequence.
\end{proof}

Our final result is an alternate combinatorial version of the
formula in Theorem~\ref{t:PeelQuiver}.  The list of shapes in a
tableau list is called the \bem{shape} of the tableau~list.

\begin{cor}[Buch--Fulton conjecture \cite{BF}] \label{c:bla}
The quiver constant\/ $c_\bla(\rr)$ in the expression\/
$\QQ_\rr(\xx-\oxx) = \sum_\bla c_\bla(\rr)s_\bla(\xx_\rr-\oxx_\rr)$ of
the quiver polynomial in terms of Schur functions (or by
Theorem~\ref{t:bla}, in the expression of double quiver functions in
terms of double Schur functions) equals the number of\/ $\rr$-factor
sequences~of~shape\/~$\bla$.
\end{cor}
\begin{proof}
Theorem~\ref{t:ptofs} and Theorem~\ref{t:PeelQuiver}.
\end{proof}



\begin{thebibliography}{ADK81}

\bibitem[AD80]{AD}
S.~Abeasis and A.~Del\thinspace{}\thinspace{}Fra, \emph{Degenerations
  for the representations of a quiver of type $A_m$},
  Boll. Un. Mat. Ital. Suppl. (1980) no. \textbf{2}, 157--171.

\bibitem[ADK81]{ADK}
S.~Abeasis, A.~Del\thinspace{}\thinspace{}Fra, and H.~Kraft, \emph{The
  geometry of representations of ${A}\sb{m}$}, Math. Ann. \textbf{256}
  (1981), no.~3, 401--418.

\bibitem[BB93]{BB}
Nantel Bergeron and Sara Billey, \emph{{RC}-graphs and {Schubert}
  polynomials}, Experimental Math. \textbf{2} (1993), no.~4,
  257--269.

\bibitem[BJS93]{BJS}
Sara~C. Billey, William Jockusch, and Richard~P. Stanley, \emph{Some
  combinatorial properties of {S}chubert polynomials}, J. Algebraic
  Combin. \textbf{2} (1993), no.~4, 345--374.

\bibitem[Bri97]{Bri97}
Michel Brion, \emph{Equivariant {Chow} groups for torus actions},
  Transform. Groups \textbf{2} (1997), no.~3, 1--43.

\bibitem[BF99]{BF}
Anders~Skovsted Buch and William Fulton, \emph{Chern class formulas
  for quiver varieties}, Invent. Math. \textbf{135} (1999), no.~3,
  665--687.

\bibitem[BFR03]{BFR03}
Anders S.\thinspace{}Buch, L{\'a}szl{\'o} M.\thinspace{}Feh{\'e}r and
  Rich{\'a}rd Rim{\'a}nyi, \emph{Positivity of quiver coefficients
  through Thom polynomials}, in preparation, 2003.

\bibitem[BKTY02]{BKTY02}
Anders~S. Buch, Andrew Kresch, Harry Tamvakis, and Alexander Yong,
  \emph{Schubert polynomials and quiver formulas}, 2002,
  \textsf{arXiv:math.AG/0211300}.


\bibitem[Buc01a]{Bu01}
Anders~Skovsted Buch, \emph{Stanley symmetric functions and quiver
  varieties}, J. Algebra \textbf{235} (2001), no. 1, 243--260.

\bibitem[Buc01b]{Bu01a}
Anders~Skovsted Buch, \emph{On a conjectured formula for quiver
  varieties}, J. Algebraic Combin. \textbf{13} (2001), no. 2,
  151--172.

\bibitem[Buc02]{Buch02}
Anders~Skovsted Buch, \emph{Grothendieck classes of quiver varieties},
  Duke Math. J. \textbf{115} (2002), no.~1, 75--103.

\bibitem[Buc03]{BuchAltSign}
Anders~Skovsted Buch, \emph{Alternating signs of quiver coefficients},
  preprint, 2003.

\bibitem[Dem74]{Dem74}
M. Demazure, \emph{Une nouvelle formule des caract\`eres},
  Bull. Sci. Math. (2) \textbf{98} (1974) 163--172.

\bibitem[EG87]{EGn87}
Paul Edelman and Curtis Greene, \emph{Balanced tableaux}, Adv. in
  Math. \textbf{63} (1987), no.~1, 42--99.

\bibitem[EG98]{EG98}
Dan Edidin and William Graham, \emph{Equivariant intersection theory},
  Invent. Math. \textbf{131} (1998), no.~3, 595--634.

\bibitem[Eis95]{Eis}
David Eisenbud, \emph{Commutative algebra, with a view toward
  algebraic geometry}, Graduate Texts in Mathematics, vol. 150,
  Springer-Verlag, New York, 1995.

\bibitem[FR02]{FRthomPoly}
L{\'a}szl{\'o} Feh{\'e}r and Rich{\'a}rd Rim{\'a}nyi, \emph{Classes of
  degeneracy loci for quivers: the {T}hom polynomial point of view},
  Duke Math. J. \textbf{114} (2002), no.~2, 193--213.

\bibitem[FK94]{FK94}
Sergey Fomin and Anatol~N. Kirillov, \emph{Grothendieck polynomials
  and the {Y}ang--{B}axter equation}, 1994, Proceedings of the Sixth
  Conference in Formal Power Series and Algebraic Combinatorics,
  DIMACS, pp.~183--190.

\bibitem[FK96]{FK96}
Sergey Fomin and Anatol~N. Kirillov, \emph{The {Y}ang-{B}axter
  equation, symmetric functions, and {S}chubert polynomials}, Discrete
  Math. \textbf{153} (1996), no.~1-3, 123--143, Proceedings of the 5th
  Conference on Formal Power Series and Algebraic Combinatorics
  (Florence, 1993).

\bibitem[FS94]{FS}
Sergey Fomin and Richard~P. Stanley, \emph{Schubert polynomials and
  the nil-{C}oxeter algebra}, Adv. Math. \textbf{103} (1994), no.~2,
  196--207.

\bibitem[Ful92]{Ful92}
William Fulton, \emph{Flags, {S}chubert polynomials, degeneracy loci,
  and determinantal formulas}, Duke Math. J. \textbf{65} (1992),
  no.~3, 381--420.

\bibitem[Ful97]{Ful97}
William Fulton, \emph{Young tableaux, with applications to
  representation theory and geometry}, London Mathematical Society
  Student Texts \textbf{35}, Cambridge University Press, Cambridge,
  1997.

\bibitem[Ful98]{FulIT}
William Fulton, \emph{Intersection theory}, second ed.,
  Springer-Verlag, Berlin, 1998.

\bibitem[FP98]{FP}
William Fulton and Piotr Pragacz, \emph{Schubert varieties and
  degeneracy loci}, Springer-Verlag, Berlin, 1998, Appendix J by the
  authors in collaboration with I. Ciocan-Fontanine.

\bibitem[Ful99]{Ful99}
William Fulton, \emph{Universal {S}chubert polynomials}, Duke
  Math. J. \textbf{96} (1999), no.~3, 575--594.

\bibitem[Gia04]{giambelli}
G. Z. Giambelli, \emph{Ordine di una variet\`a pi\`u ampia di quella
  rappresentata coll'annullare tutti i minori di dato ordine estratti
  da una data matrice generica di forme}, Mem. R. Ist. Lombardo
  \textbf{3} (1904), no. 11, 101--135.

\bibitem[Hai92]{Ha92}
Mark Haiman, \emph{Dual equivalence with applications, including a
  conjecture of Proctor}, Discrete Math. \textbf{99} (1992), no. 1--3,
  79--113.

\bibitem[Jos84]{joseph}
Anthony Joseph, \emph{On the variety of a highest weight module},
  J. Algebra \textbf{88} (1984), no.~1, 238--278.

\bibitem[Kaz97]{Kaz97}
M.~{\'E}. Kazarian, \emph{Characteristic classes of singularity
  theory}, The Arnold-Gelfand mathematical seminars, Birkh\"auser
  Boston, Boston, MA, 1997, pp.~325--340.

\bibitem[KM03a]{grobGeom}
Allen Knutson and Ezra Miller, \emph{{Gr\"{o}bner} geometry of
  {Schubert} polynomials}, to appear in Ann.\ of Math~(2), 2003.

\bibitem[KM03b]{subword}
Allen Knutson and Ezra Miller, \emph{Subword complexes in {Coxeter}
  groups}, to appear in Adv.\ in Math., 2003.

\bibitem[Kog00]{KoganThesis}
Mikhail Kogan, \emph{Schubert geometry of flag varieties and
  Gel\/$'\!$fand--Cetlin theory}, Ph.D. thesis, Massachusetts
  Institute of Technology, 2000.

\bibitem[Las90]{Las90}
Alain Lascoux, \emph{Anneau de {G}rothendieck de la vari\'et\'e de
  drapeaux}, The Grothendieck Festschrift, Vol.\ III, Birkh\"auser
  Boston, Boston, MA, 1990, pp.~1--34.

\bibitem[Las01]{Las01}
Alain Lascoux, \emph{Transition on Grothendieck polynomials}, Physics
  and combinatorics, 2000 (Na\-goya), 164--179, World Sci. Publishing,
  River Edge, NJ, 2001.

\bibitem[Las03]{LasDoubKey}
Alain Lascoux, \emph{Double crystal graphs}, Studies in Memory of
  {Issai} {Schur} (Anthony Joseph, Anna Melnikov, and Rudolf
  Rentschler, eds.), Birkh{\"a}user, Boston, 2003, pp.~95--114.

\bibitem[Len02]{Len02}
Cristian Lenart, \emph{A unified approach to combinatorial formulas
  for Schubert polynomials}, Preprint at
  \textsf{http:/$\!$/www.math.albany.edu:8000/math/pers/lenart/articles/schubert.html},
  2002.

\bibitem[LM98]{LM}
V.~Lakshmibai and Peter Magyar, \emph{Degeneracy schemes, quiver
  schemes, and {S}chubert varieties}, Internat. Math. Res. Notices
  (1998), no.~12, 627--640.

\bibitem[LS82]{LS82}
Alain Lascoux and Marcel-Paul Sch{\"u}tzenberger, \emph{Structure de
  {H}opf de l'anneau de cohomologie et de l'anneau de {G}rothendieck
  d'une vari\'et\'e de drapeaux}, C. R. Acad. Sci. Paris S\'er. I
  Math. \textbf{295} (1982), no.~11, 629--633.

\bibitem[LS85]{LS85}
Alain Lascoux and Marcel-Paul Sch{\"u}tzenberger, \emph{Schubert
  polynomials and the Littlewood--Richardson rule},
  Lett. Math. Phys. \textbf{10} (1985), no. 2-3, 111--124.

\bibitem[LS89]{LS89}
Alain Lascoux and Marcel-Paul Sch{\"u}tzenberger, \emph{Tableaux and
  noncommutative Schubert polynomials}, Func. Anal. Appl. \textbf{23}
  (1989) 63--64.

\bibitem[LS90]{LS90}
Alain Lascoux and Marcel-Paul Sch{\"u}tzenberger, \emph{Keys and
  standard bases}, Invariant theory and tableaux (Minneapolis, MN,
  1988), 125--144, IMA Vol. Math. Appl., 19, Springer, New York,
  1990.

\bibitem[Lus90]{Lu90}
G. Lusztig, \emph{Canonical bases arising from quantized enveloping
  algebras}, J. Amer. Math. Soc.\ \textbf{3} (1990), no. 2, 447--498.

\bibitem[Mac95]{Mac}
I. G. Macdonald, \emph{Notes on Schubert polynomials}, Publ. LACIM
  \textbf{6}, UQAM, Montr\'eal, 1991.

\bibitem[Mag98]{Mag}
Peter Magyar, \emph{Schubert polynomials and Bott-Samelson varieties},
  Commentarii Mathematici Helvetici \textbf{73} (1998), 603--636.

\bibitem[Mil03]{altSign}
Ezra Miller, \emph{Alternating formulae for \K-theoretic quiver
  polynomials}, preprint, 2003.

\bibitem[MS04]{CCA}
Ezra Miller and Bernd Sturmfels, \emph{Combinatorial commutative
  algebra}, Springer Verlag, in preparation, to appear in 2004.

\bibitem[Ram85]{Ram85}
A.~Ramanathan, \emph{Schubert varieties are arithmetically
  {C}ohen-{M}acaulay}, Invent. Math. \textbf{80} (1985), no.~2,
  283--294.

\bibitem[RR85]{RR85}
S.~Ramanan and A.~Ramanathan, \emph{Projective normality of flag
  varieties and {S}chubert varieties}, Invent. Math. \textbf{79}
  (1985), no.~2, 217--224.

\bibitem[Ros89]{rossmann}
W.~Rossmann, \emph{Equivariant multiplicities on complex varieties},
  Ast\'erisque (1989), no.~173--174, 11, 313--330, Orbites unipotentes
  et repr\'esentations, III.

\bibitem[RS95a]{RSplact95}
Victor Reiner and Mark Shimozono, \emph{Plactification}, J. Algebraic
  Combin. \textbf{4} (1995), no.~4, 331--351.

\bibitem[RS95b]{RS95K}
\emph{Key polynomials and a flagged Littlewood-Richardson rule},
  J. Combin. Theory Ser. A \textbf{70} (1995), no.~1, 107--143.

\bibitem[RS98]{RSpeelable98}
\bysame, \emph{Percentage-avoiding, northwest shapes and peelable
  tableaux}, J. Combin. Theory Ser. A \textbf{82} (1998), no.~1, 1--73.

\bibitem[Sta84]{St}
R. P. Stanley, \emph{On the number of reduced decompositions of
  elements of Coxeter groups}, Eur. J. Math. \textbf{5} (1984)
  359--372.

\bibitem[Tho55]{thom55}
R.~Thom, \emph{Les singularit\'es des applications diff\'erentiables},
  Ann. Inst. Fourier, Grenoble \textbf{6} (1955--1956), 43--87.

\bibitem[Yon03]{Yong03}
Alexander Yong, \emph{On combinatorics of quiver component formulas},
  preprint, 2003.

\bibitem[Zel85]{Zel85}
A.~V. Zelevinski\u{\i}, \emph{Two remarks on graded nilpotent
  classes}, Uspekhi Mat. Nauk \textbf{40} (1985), no.~1(241), 199--200.

\end{thebibliography}
\def\cprime{$'$}
\providecommand{\bysame}{\leavevmode\hbox to3em{\hrulefill}\thinspace}


\end{document}